\documentclass[12pt,a4paper]{amsart}
\usepackage[cmtip,arrow]{xy}
\usepackage{pb-diagram,pb-xy}
\usepackage{amsxtra}

\calclayout

\DeclareMathSymbol{\subsetneq}{\mathrel}{AMSb}{"28}
\DeclareMathSymbol{\rightrightarrows}{\mathrel}{AMSa}{"13}

\newcommand{\+}{\nobreakdash-}
\renewcommand{\:}{\colon}
\renewcommand{\.}{\mskip.5\thinmuskip\relax}
\renewcommand{\;}{,\medspace}

\newcommand{\rarrow}{\longrightarrow}
\newcommand{\ot}{\otimes}
\newcommand{\oc}{\mathbin{\text{\smaller$\square$}}}

\newcommand{\lrarrow}{\.\relbar\joinrel\relbar\joinrel\rightarrow\.}
\newcommand{\bu}{{\text{\smaller\smaller$\scriptstyle\bullet$}}}

\DeclareMathOperator{\cone}{cone}
\DeclareMathOperator{\id}{id}
\DeclareMathOperator{\Id}{Id}
\DeclareMathOperator{\coker}{coker}
\DeclareMathOperator{\Hom}{Hom}
\DeclareMathOperator{\Supp}{Supp}
\DeclareMathOperator{\supp}{supp}
\DeclareMathOperator{\cHom}{\mathcal{H}\mathit{om}}
\DeclareMathOperator{\Tot}{Tot}
\DeclareMathOperator{\Ext}{Ext}
\DeclareMathOperator{\Tor}{Tor}
\DeclareMathOperator{\Spec}{Spec}

\newcommand{\op}{{\mathrm{op}}}

\newcommand{\Z}{{\mathbb Z}}
\newcommand{\boC}{{\mathbb C}}

\renewcommand{\O}{{\mathcal O}}
\newcommand{\A}{{\mathcal A}}
\newcommand{\B}{{\mathcal B}}
\newcommand{\M}{{\mathcal M}}
\newcommand{\N}{{\mathcal N}}
\renewcommand{\L}{{\mathcal L}}
\newcommand{\K}{{\mathcal K}}
\renewcommand{\P}{{\mathcal P}}
\newcommand{\Q}{{\mathcal Q}}
\newcommand{\R}{{\mathcal R}}
\newcommand{\E}{{\mathcal E}}
\newcommand{\F}{{\mathcal F}}
\newcommand{\G}{{\mathcal G}}
\renewcommand{\H}{{\mathcal H}}
\newcommand{\U}{{\mathcal U}}
\newcommand{\V}{{\mathcal V}}
\newcommand{\W}{{\mathcal W}}

\newcommand{\D}{{\mathcal D}}
\newcommand{\J}{{\mathcal J}}
\newcommand{\I}{{\mathcal I}}

\newcommand{\sE}{{\mathsf E}}
\newcommand{\sF}{{\mathsf F}}
\newcommand{\sD}{{\mathsf D}}
\newcommand{\sC}{{\mathsf C}}
\newcommand{\sT}{{\mathsf T}}

\newcommand{\sDG}{{\mathsf{DG}}}
\newcommand{\MF}{{\mathit{MF}}}

\newcommand{\Perf}{\mathit{Perf}}
\newcommand{\Sing}{{\mathit{Sing}}}

\newcommand{\boL}{{\mathbb L}}
\newcommand{\boR}{{\mathbb R}}

\newcommand{\Ups}{\Upsilon}

\newcommand{\vect}{{\operatorname{\mathsf{--vect}}}}
\newcommand{\modl}{{\operatorname{\mathsf{--mod}}}}
\newcommand{\modr}{{\operatorname{\mathsf{mod--}}}}
\newcommand{\qcoh}{{\operatorname{\mathsf{--qcoh}}}}
\newcommand{\qcohr}{{\operatorname{\mathsf{qcoh--}}}}
\newcommand{\coh}{{\operatorname{\mathsf{--coh}}}}
\newcommand{\wfl}{{w\!\operatorname{\mathsf{--fl}}}}
\newcommand{\fwfl}{{f^*w\!\operatorname{\mathsf{--fl}}}}
\newcommand{\fadj}{{f\!\operatorname{\mathsf{--adj}}}}
\newcommand{\qc}{{\operatorname{\mathsf{--qc}}}}
\newcommand{\flc}{{\operatorname{\mathsf{fl-c}}}}
\newcommand{\zdim}{{\operatorname{\mathsf{0-dim}}}}

\newcommand{\fl}{{\mathsf{fl}}}
\newcommand{\ffd}{{\mathsf{ffd}}}
\newcommand{\lfd}{{\mathsf{lfd}}}
\newcommand{\fid}{{\mathsf{fid}}}
\newcommand{\lf}{{\mathsf{lf}}}
\newcommand{\lp}{{\mathsf{lp}}}
\newcommand{\inj}{{\mathsf{inj}}}
\newcommand{\co}{{\mathsf{co}}}
\newcommand{\ctr}{{\mathsf{ctr}}}
\newcommand{\abs}{{\mathsf{abs}}}
\newcommand{\st}{{\mathsf{st}}}
\renewcommand{\b}{{\mathsf{b}}}
\newcommand{\sop}{{\mathsf{op}}}
\newcommand{\ch}{{\mathsf{coh}}}

\newcommand{\rmodl}{\text{-}\mathrm{mod}}
\newcommand{\modld}{\text{-}\mathrm{mod}^{\mathrm{dg}}}
\newcommand{\modlc}{\text{-}\mathrm{mod}^{\mathrm{cdg}}}

\newcommand{\modrd}{\mathrm{mod}^{\mathrm{dg}}\text{-}}
\newcommand{\modrc}{\mathrm{mod}^{\mathrm{cdg}}\text{-}}

\newcommand{\modrcfgp}{\mathrm{mod}^{\mathrm{cdg}}_
                     {\mathrm{fgp}}\text{-}}
\newcommand{\modrdfgp}{\mathrm{mod}^{\mathrm{dg}}_
                     {\mathrm{fgp}}\text{-}}

\newcommand{\modrdfp}{\mathrm{mod}^{\mathrm{dg}}_
                     {\mathrm{fp}}\text{-}}
\newcommand{\modrcfp}{\mathrm{mod}^{\mathrm{cdg}}_
                     {\mathrm{fp}}\text{-}}

\newcommand{\rinj}{{\mathrm{inj}}}
\newcommand{\rfl}{{\mathrm{fl}}}
\newcommand{\rprj}{{\mathrm{prj}}}
\newcommand{\rfg}{{\mathrm{fg}}}
\newcommand{\rabs}{{\mathrm{abs}}}
\newcommand{\rco}{{\mathrm{co}}}
\newcommand{\rctr}{{\mathrm{ctr}}}
\newcommand{\rfgp}{{\mathrm{fgp}}}

\newcommand{\rD}{{\mathrm D}}

\newcommand{\Section}[1]{\bigskip
\addtocontents{toc}{\medskip}\section{#1}\medskip}
\theoremstyle{plain}
\newtheorem*{thm}{Theorem}
\newtheorem*{lem}{Lemma}
\newtheorem*{lema}{Lemma A}
\newtheorem*{lemb}{Lemma B}
\newtheorem*{lemc}{Lemma C}
\newtheorem*{lemd}{Lemma D}
\newtheorem*{leme}{Lemma E}
\newtheorem*{lemf}{Lemma F}
\newtheorem*{lemg}{Lemma G}
\newtheorem*{lemh}{Lemma H}
\newtheorem*{cor}{Corollary}
\newtheorem*{prop}{Proposition}
\newtheorem*{propa}{Proposition A}
\newtheorem*{propb}{Proposition B}
\newtheorem*{propc}{Proposition C}
\newtheorem*{propd}{Proposition D}

\theoremstyle{definition}
\newtheorem*{rem}{Remark}

\begin{document}

\title{Coherent analogues of matrix factorizations \\
and relative singularity categories}
\author{Alexander I.~Efimov \ and \ Leonid Positselski}

\address{Steklov Mathematical Institute of RAS, Gubkin str.~8,
GSP-1, Moscow 119991, Russia \newline
\indent AG Laboratory, HSE, 7~Vavilova str., Moscow, Russia, 117312}
\email{efimov@mccme.ru}

\address{Sector of Algebra and Number Theory, Institute for
Information Transmission Problems, Moscow 127994; and \newline
\indent Laboratory of Algebraic Geometry, National Research University
Higher School of Economics, Moscow 117312, Russia; and \newline
\indent Mathematics Department, Technion --- Israel Institute of
Technology, Haifa 32000, Israel}
\email{posic@mccme.ru}

\begin{abstract}
 We define the triangulated category of relative singularities
of a closed subscheme in a scheme.
 When the closed subscheme is a Cartier divisor, we consider
matrix factorizations of the related section of a line bundle,
and their analogues with locally free sheaves replaced by
coherent ones.
 The appropriate exotic derived category of coherent matrix
factorizations is then identified with the triangulated category
of relative singularities, while the similar exotic derived
category of locally free matrix factorizations is its full
subcategory.
 The latter category is identified with the kernel of the direct image
functor corresponding to the closed embedding of the zero locus and
acting between the conventional (absolute) triangulated categories
of singularities.
 Similar results are obtained for matrix factorizations of infinite
rank; and two different ``large'' versions of the triangulated
category of relative singularities, corresponding to the approaches
of Orlov and Krause, are identified in the case of a Cartier divisor.
 A version of the Thomason--Trobaugh--Neeman localization theory is
proven for coherent matrix factorizations and disproven for
locally free matrix factorizations of finite rank.
 Contravariant (coherent) and covariant (quasi-coherent) versions
of the Serre--Grothendieck duality theorems for matrix factorizations
are established, and pull-backs and push-forwards of matrix
factorizations are discussed at length.
 A number of general results about derived categories of the second
kind for CDG\+modules over quasi-coherent CDG\+algebras are
proven on the way.
 Hochschild (co)homology of matrix factorization categories are
discussed in an appendix.
\end{abstract}

\maketitle

\tableofcontents

\section*{Introduction}
\medskip

 A \emph{matrix factorization} of an element $w$ in a commutative
ring~$R$ is a pair of square matrices $(\Phi,\Psi)$ of the same
size, with entries from $R$, such that both the products $\Phi\Psi$
and $\Psi\Phi$ are equal to $w$ times the identity matrix.
 In the coordinate-free language, a matrix factorization is a pair
of finitely generated free $R$\+modules $M^0$ and $M^1$ together with
$R$\+module homomorphisms $M^0\rarrow M^1$ and $M^1\rarrow M^0$ such
that both the compositions $M^0\rarrow M^1\rarrow M^0$ and
$M^1\rarrow M^0\rarrow M^1$ are equal to the multiplication with~$w$.
 Matrix factorizations were introduced by Eisenbud~\cite{Ei} and used
by Buchweitz~\cite{Bu} for the purposes of study of the maximal
Cohen--Macaulay modules over hypersurface local rings.

 Another name for this notion is ``D\+branes in the Landau--Ginzburg
B model'' (as suggested by Kontsevich) \cite{KL}; in this context,
the element~$w$ is called the \emph{potential}.
 One generalizes the above definition, replacing free modules with
projective modules~\cite{KL,Or1}, with locally free sheaves~\cite{Or3},
and finally with coherent sheaves~\cite{PL}.
 The importance of the latter generalization is emphasized in
the present paper.

 Being particular cases of curved DG\+modules over a curved
DG\+ring~\cite{KL,Pkoszul}, matrix factorizations form a DG\+category.
 So one can consider the corresponding category of closed degree-zero
morphisms up to chain homotopy, which is a triangulated category.
 Generally speaking, however, the homotopy category is ``too big''
for most purposes, and one would like to pass from it to
an appropriately defined derived category.
 One can use the homotopy category in lieu of the derived one when
dealing with projective modules~\cite{KL,Or1}; for locally free
matrix factorizations over a nonaffine scheme, there is an option of
working with the quotient category of the homotopy category by
the locally contractible objects~\cite[Definition~3.13]{PV}.
 When dealing with coherent (analogues of) matrix factorizations,
having some kind of derived category construction
is apparently unavoidable.

 The relevant concept of a derived category is that of the derived
category of the second kind, as developed in~\cite{Pkoszul,Psemi}.
 There are several versions of this notion; the appropriate one
for quasi-coherent sheaves is called the \emph{coderived category}
and for coherent sheaves it is the \emph{absolute derived category}.
 The absolute derived category of locally free matrix factorizations
was studied in~\cite{Or3}; for coherent matrix factorizations over
a smooth variety, it was considered in~\cite{PL}.
 These two absolute derived categories are equivalent for regular
schemes, but can be different otherwise (as we show with an explicit
counterexample).

 The \emph{triangulated category of singularities} of a Noetherian
scheme was defined by D.~Orlov in~\cite{Or1} as the quotient category
of the bounded derived category of coherent sheaves by its full
triangulated subcategory of perfect complexes, i.~e., the objects
locally presentable as finite complexes of locally free sheaves.
 This triangulated category vanishes if and only if the Noetherian
scheme is regular.
 It was shown in~\cite[Theorem~3.9]{Or1}, under mild assumptions on
an affine regular Noetherian scheme $X$ and a potential
(regular function) $w$ on it, that the homotopy category of
locally free matrix factorizations of $w$ over $X$ is equivalent to
the triangulated category of singularities of the zero locus
$X_0$ of~$w$ in~$X$.

 In his recent paper~\cite{Or3}, Orlov shows that the affineness
assumption on $X$ can be dropped in this result if one replaces
the homotopy category of locally free matrix factorizations with
their absolute derived category.
 He also considers the general case of a nonaffine singular scheme $X$,
for which he obtains a fully faithful functor from the absolute derived
category of locally free matrix factorizations over $X$ to
the triangulated category of singularities of~$X_0$.
 The problem of studying the difference between these two triangulated
categories was posed in the introduction to~\cite{Or3}.

 The first aim of the present paper is to provide an alternative proof
of these results of Orlov for regular schemes, an alternative
generalization of them to singular schemes, and a more precise version
of Orlov's original generalization.
 We replace the triangulated category at the source of Orlov's fully
faithful functor by a ``larger'' category (containing the original one)
and the triangulated category at the target by a ``smaller'' category
(a quotient of the original one), thereby transforming this
functor into an equivalence of triangulated categories.
 We also describe the image of Orlov's fully faithful functor as
the kernel of a certain other triangulated functor.

 More precisely, we show that the absolute derived category of coherent
matrix factorizations of~$w$ over $X$ is equivalent to what we call
the \emph{triangulated category of singularities of $X_0$ relative
to~$X$}.
 The latter category is a certain quotient category of
the triangulated category of singularities of $X_0$; it measures,
roughly speaking, how much worse are the singularities of $X_0$
compared to those of~$X$.
 As to the image of Orlov's fully faithful embedding, it consists
precisely of those objects of the conventional (absolute)
triangulated category of singularities of $X_0$ whose direct images
vanish in the triangulated category of singularities of~$X$.

 The paper consists of three sections and two appendices.
 In Section~1, we prove three rather general technical assertions about
derived categories of the second kind for CDG\+modules over
a quasi-coherent CDG\+algebra with a restriction on the homological
dimension.
 One of them, claiming that certain embeddings of DG\+categories of
CDG\+modules induce equivalences of the derived categories of
the second kind, is a generalization of~\cite[Theorem~3.2]{PP} based
on a modification of the same argument, originally introduced for
the proof of~\cite[Theorem~7.2.2]{Psemi}.

 The idea of the proof of the other assertion, according to which
certain natural functors between derived categories of the second kind
are fully faithful, is new.
 The third technical assertion explains when the coderived category
coincides with the absolute derived category of the same class of
CDG\+modules: e.~g., for the locally projective CDG\+modules this
is true.

 A version of (the former two of) these results is used in Section~2
in order to extend Orlov's cokernel functor from the absolute derived
category of locally free matrix factorizations to the absolute derived
category of coherent ones.
 This extension of the cokernel functor admits a simple construction
of a functor in the opposite direction, suggested in~\cite{PL}.
 We use these constructions to obtain a new proof of Orlov's theorem,
and our own generalization of it to the singular case.

 When $X$ is regular, Orlov's and our results amount to the same
assertion, since the absolute derived categories of locally free
and coherent matrix factorizations are equivalent by our
Theorem~\ref{flat-dimension-thm}.
 When $X$ is singular, the natural functor between these two
absolute derived categories is fully faithful by our
Proposition~\ref{embedding-prop}, and Orlov's full-and-faithfulness
theorem follows from ours by virtue of an appropriate
semiorthogonality property.

 We also compare a ``large'' version of the triangulated category
of relative singularities with the coderived category of
quasi-coherent matrix factorizations, strengthening some results
of Polishchuk--Vaintrob~\cite{PV}.
 A ``large'' version of the absolute triangulated category of
singularities, defined by Orlov in~\cite{Or1}, is identified with
H.~Krause's stable derived category~\cite{Kr} in the case of
a divisor in a regular scheme.
 A similar result is proven in the case of a Cartier divisor in
a singular scheme, where we extend Krause's theory by defining
the \emph{relative stable derived category}.
 For any closed subscheme of finite flat dimension in a separated
Noetherian scheme, the relative stable derived category is
compactly generated by its full triangulated subcategory equivalent
to the triangulated category of relative singularities.

 The homotopy categories of unbounded complexes of projective modules
over a ring and injective quasi-coherent sheaves over a scheme were
studied in the papers by J\o rgensen~\cite{Jo} and Krause~\cite{Kr};
subsequently, Iyengar and Krause have constructed an equivalence
between these two categories for rings with dualizing
complexes~\cite{IK}.
 These results were extended to quasi-coherent sheaves over schemes
by Neeman~\cite{N-f} and Murfet~\cite{M-th}, who found a way to define
a replacement of the homotopy category of (nonexistent) projective
sheaves in terms of the flat ones.
 The equivalence between these two categories is a covariant version
of the Serre--Grothendieck duality~\cite{Har}.
 It is also very similar to the derived comodule-contramodule
correspondence theory, developed by the second author
in~\cite{Pkoszul,Psemi}.

 The Serre--Grothendieck duality for matrix factorizations in
the situation of a smooth variety $X$ (and an isolated singularity
of $X_0$) was studied in~\cite{M-d}.
 In this paper we extend the duality to matrix factorizations over
much more general schemes $X$, constructing an equivalence between
two ``large'' exotic derived categories, namely, the coderived category
of flat (or locally free) matrix factorizations of possibly infinite
rank and the coderived category of quasi-coherent matrix factorizations.
 Unless $X$ is Gorenstein, this equivalence is \emph{not} provided by
the natural functor induced by the embedding of DG\+categories, but
rather differs from it in that the tensor product with
the dualizing complex has to be taken along the way.
 A contravariant Serre duality in the form of an auto-anti-equivalence
of the absolute derived category of coherent matrix factorizations
is also obtained.

 There was some attention paid to pull-backs and push-forwards of
matrix factorizations recently~\cite{PV,DM,PV2}.
 In Section~3, we approach this topic with our techniques,
constructing the push-forwards of locally free matrix factorizations
of infinite rank for any morphism of finite flat dimension between
schemes of finite Krull dimension, and the push-forwards of
locally free matrix factorizations of finite rank for any such
morphism for which the induced morphism of the zero loci
of~$w$ is proper.
 At the price of having to adjoin the images of idempotent
endomorphisms, the preservation of finite rank under push-forwards
is proven assuming only the support of the matrix
factorization~\cite{PV} to be proper over the base.

 Push-forwards of quasi-coherent matrix factorizations are well-defined
for any morphism of Noetherian schemes, and push-forwards of coherent
matrix factorizations exist under properness assumptions similar to
the above.
 A general study of category-theoretic and set-theoretic supports
of quasi-coherent and coherent CDG\+modules is undertaken in this
paper in order to obtain an independent proof of the preservation of
coherence under the push-forwards not based on the passage to
the triangulated categories of singularities.

 The compatibility with pull-backs and push-forwards is an organic
part of the Serre--Grothendieck duality theory.
 The contravariant duality agrees with push-forwards of coherent
sheaves (or matrix factorizations) with respect to proper
morphisms~\cite{Har}, while the covariant duality transforms
the conventional inverse image of flat matrix factorizations into
the extraordinary inverse image of quasi-coherent ones~\cite{Pcosh}.
 We use the latter result in order to construct the Hartshorne--Deligne
extraordinary inverse image functor, which is denoted by~$f^!$
in~\cite{Har} and which we denote by~$f^+$, in the case of
quasi-coherent matrix factorizations.

 Appendix~A contains proofs of some basic facts about flat,
locally projective, and injective quasi-coherent graded modules
which are occasionally used in the main body of the paper.
 Appendix~B can be viewed as a complement to the paper~\cite{PP}.
 While Section~\ref{loc-free-hochschild} contains some variations of
and improvements on the results about Hochschild (co)homology of
(C)DG\+categories and (locally free) matrix factorizations in~\cite{PP},
Section~\ref{coherent-hochschild} presents an alternative approach to
the Hochschild (co)homology of coherent matrix factorizations based on
the techniques developed in the main body of this paper.

\medskip
\textbf{Acknowledgment.}
 The second author is grateful to Daniel Pomerleano and Kevin Lin for
sending him a preliminary version of their manuscript~\cite{PL}, which
inspired the present work.
 We also wish to thank Roman Bezrukavnikov, Alexander Kuznetsov, and
Amnon Neeman for helpful conversations, Alexander Polishchuk and
Henning Krause for stimulating questions, and Dmitri Orlov for
his insightful remarks on an early version of this manuscript.
 Finally, the second author would like to thank an anonymous referee
of an early version for detailed suggestions on the improvement of
the exposition.
 It was initially my attempt to implement those suggestions that
eventually resulted in the growth of the size of this paper by
a factor of almost seven.

 The first author was partially supported by RFBR grant 2998.2014.1
and RFBR research project 15-51-50045.
 The second author was supported in part by a Simons Foundation grant
and RFBR grants in Moscow and by a fellowship from the Lady Davis
Foundation at the Technion.
 The paper was prepared within the framework of a subsidy granted to
the HSE by the Government of the Russian Federation for
the implementation of the Global Competitiveness Program.

\Section{Exotic Derived Categories of Quasi-Coherent CDG-Modules}

\subsection{CDG-rings and CDG-modules}
 A \emph{CDG\+ring} (curved differential graded ring) $B=(B,d,h)$ is
defined as a graded ring $B=\bigoplus_{i\in\Z}B^i$ endowed with
an odd derivation $d\:B\rarrow B$ of degree~$1$ and an element
$h\in B^2$ such that $d^2(b)=[h,b]$ for all $b\in B$ and $d(h)=0$.
 So one should have $d\:B^i\rarrow B^{i+1}$ and
$d(ab) = d(a)b + (-1)^{|a|}ad(b)$; the brackets $[-,-]$ denote
the supercommutator $[a,b]=ab-(-1)^{|a||b|}ba$.
 The element~$h$ is called the \emph{curvature element}.

 A morphism of CDG\+rings $B\rarrow A$ is a pair $(f,a)$, with
a morphism of graded rings $f\:B\rarrow A$ and an element
$a\in A^1$, such that $f(d_Bb) = d_Af(b) + [a,f(b)]$ for all
$b\in B$ and $f(h_B) = h_A + d_Aa + a^2$.
 The composition of morphisms of CDG\+rings is defined by the obvious
rule $(f,a)\circ(g,b)=(f\circ g\;a+f(b))$.
 The element~$a$ is called the \emph{change-of-connection element}.
 A discussion of the origins of these definitions can be found
in the paper~\cite{Pcurv}, where the above terminology first
appeared (see also an earlier paper~\cite{GJ}, where the motivation
was entirely different).

 A \emph{left CDG\+module} $M=(M,d_M)$ over a CDG\+ring $B$ is
a graded $B$\+module endowed with an odd derivation $d_M\:M\rarrow M$
compatible with the derivation~$d$ on $B$ such that $d^2_M(m)=hm$
for all $m\in M$.
 Given a morphism of CDG\+rings $(f,a)\:B\rarrow A$ and
a CDG\+module $(M,d)$ over $A$, the CDG\+module $(M,d')$
over $B$ is defined by the rule $d'(m) = d(m) + am$.

 Given graded left $B$\+modules $M$ and $N$, homogeneous $B$\+module
morphisms $f\:M\rarrow N$ of degree~$n$ are defined as
homogeneous maps supercommuting with the action of $B$, i.~e.,
$f(bm)=(-1)^{n|b|}bf(m)$.
 When $M$ and $N$ are CDG\+modules, the homogeneous $B$\+module
morphisms $M\rarrow N$ form a complex of abelian groups with
the differential $d(f)(m)=d(f(m))-(-1)^{|f|}f(d(m))$.
 The curvature-related terms cancel out in the computation of
the square of this differential, so one has $d^2(f)=0$.
 Therefore, left CDG\+modules over $B$ form a DG\+category.

 Two aspects of the above definitions are worth to be pointed out.
 First, the CDG\+rings or modules have \emph{no} cohomology modules,
as their differentials do not square to zero.
 Second, given a CDG\+ring $B$, there is \emph{no} natural way to
define a CDG\+module structure on the free graded $B$\+module $B$
(though $B$ is naturally a CDG\+bimodule over itself, in
the appropriate sense).

 We refer the reader to~\cite[Section~3.1]{Pkoszul} or
\cite[Sections~0.4.3\+-0.4.5]{Psemi} for more detailed
discussions of the above notions.
 We will not need to consider any gradings different from
$\Z$\+gradings in this paper, though all the general results
will be equally applicable in the $\Gamma$\+graded situation
in the sense of~\cite[Section~1.1]{PP}.

\subsection{Quasi-coherent CDG-algebras}  \label{qcoh-cdg}
 Throughout this paper, unless specified otherwise, $X$ is
a separated Noetherian scheme with enough vector bundles;
in other words, it is assumed that every coherent sheaf on $X$ is
the quotient sheaf of a locally free sheaf of finite rank.
 Note that the class of all schemes satisfying these conditions
is closed under the passages to open and closed
subschemes~\cite[Section~1.2]{Or1} and contains all regular
separated Noetherian schemes~\cite[Exercise~III.6.8]{Har2}.

 Recall the definition of a \emph{quasi-coherent CDG\+algebra}
from~\cite[Appendix B]{Pkoszul}.
 A quasi-coherent CDG\+algebra $\B$ over $X$ is a graded
quasi-coherent $\O_X$\+algebra such that for each affine open
subscheme $U\subset X$ the graded ring $\B(U)$ is endowed with
a structure of CDG\+ring, i.~e., a (not necessarily $\O_X$\+linear)
odd derivation $d\:\B(U)\rarrow\B(U)$ of degree~$1$ and an element
$h\in\B^2(U)$.
 For each pair of embedded affine open subschemes $U\subset
V\subset X$, an element $a_{UV}\in\B^1(U)$ is fixed such that
the restriction morphism $\B(V)\rarrow\B(U)$ together with
the element $a_{UV}$ form a morphism of CDG\+rings.
 The obvious compatibility condition is imposed for triples of
embedded affine open subschemes $U\subset V\subset W\subset X$.

 A \emph{quasi-coherent left CDG\+module} $\M$ over $\B$ is
an $\O_X$\+quasi-coherent (or, equivalently, $\B$\+quasi-coherent)
sheaf of graded left modules over $\B$ together with
a family of differentials $d\:\M(U)\rarrow\M(U)$
defined for all affine open subschemes $U\subset X$ such that
$\M(U)$ is a CDG\+module over $\B(U)$ and the appropriate
compatibility condition holds with respect to the restriction
morphisms of CDG\+rings $\B(V)\rarrow\B(U)$.
 Specifically, for a quasi-coherent left CDG\+module $\M$ one
should have $d(s)|_U = d(s|_U) + a_{UV}s|_U$ for any $s\in\M(V)$.

 Quasi-coherent left CDG\+modules over a quasi-coherent CDG\+algebra
$\B$ form a DG\+category~\cite{Pkoszul}.
 The complex of morphisms between CDG\+modules $\N$ and $\M$ is
the graded abelian group of homogeneous $\B$\+module morphisms
$f\:\N\rarrow\M$ with the differential $d(f)$ defined locally as
the supercommutator of~$f$ with the differentials in
$\N(U)$ and $\M(U)$.
 We denote this DG\+category by $\B\qcoh$.

 We will call a quasi-coherent graded algebra $\B$ over $X$
\emph{Noetherian} if the graded ring $\B(U)$ is left Noetherian
for any affine open subscheme $U\subset X$.
 Equivalently, $\B$ is Noetherian if the abelian category of
quasi-coherent graded left $\B$\+modules is a locally Noetherian
Grothendieck category.
 In this case, the full DG\+subcategory in $\B\qcoh$ formed by
CDG\+modules whose underlying graded $\B$\+modules are coherent
(i.~e., finitely generated over~$\B$) is denoted by $\B\coh$. 

 Given a quasi-coherent graded left $\B$\+module $\M$ and 
a quasi-coherent graded right $\B$\+module $\N$, one can define
their tensor product $\N\ot_\B\M$, which is a quasi-coherent
graded $\O_X$\+module.
 A quasi-coherent graded left $\B$\+module $\M$ is called \emph{flat}
if the functor ${-}\ot_\B\M$ is exact on the abelian category of
quasi-coherent graded right $\B$\+modules.
 Equivalently, $\M$ is flat if the graded left $\B(U)$\+module
$\M(U)$ is flat for any affine open subscheme $U\subset X$.
 The \emph{flat dimension} of a quasi-coherent graded module $\M$
is the minimal length of its flat left resolution.

 The full DG\+subcategory in $\B\qcoh$ formed by CDG\+modules
whose underlying graded $\B$\+modules are flat is denoted by
$\B\qcoh_\fl$, and the full subcategory formed by CDG\+modules
whose underlying graded $\B$\+modules have finite flat dimension
is denoted by $\B\qcoh_\ffd$.
 The similarly defined DG\+categories of coherent CDG\+modules
are denoted by $\B\coh_\fl$ and $\B\coh_\ffd$.

 All the above DG\+categories of quasi-coherent CDG\+modules
(and the similar ones defined below in this paper) admit shifts
and twists, and, in particular, cones.
 It follows that their homotopy categories $H^0(\B\qcoh)$,
\ $H^0(\B\qcoh_\fl)$, \ $H^0(\B\coh)$, etc.\ are triangulated.
 Besides, to any finite complex (of objects and closed
morphisms) in one of these DG\+categories one can assign its
total object, which is an object of (i.~e., a CDG\+module belonging
to) the same DG\+category~\cite[Section~1.2]{Pkoszul}.

 The DG\+categories $\B\qcoh$ and $\B\qcoh_\fl$ also admit
infinite direct sums.
 Hence in these two DG\+categories one can totalize even
an unbounded complex by taking infinite direct sums along
the diagonals.

 The DG\+category $\B\qcoh$ also admits infinite products (which
one can obtain using the coherator construction
from~\cite[Section~B.14]{TT}), but these are not well-behaved
(neither exact nor local), so we will not use them.

\subsection{Derived categories of the second kind}
\label{second-kind}
 The nonexistence of the cohomology groups for curved structures
stands in the way of the conventional definition of the derived
category of CDG\+modules, which therefore does not seem to
make sense.
 The suitable class of constructions of derived categories for
CDG\+modules is that of the \emph{derived categories of
the second kind}~\cite{Psemi,Pkoszul}.

 Let $\B$ be a quasi-coherent CDG\+algebra over $X$; assume that
the quasi-coherent graded algebra $\B$ is Noetherian.
 Then a coherent CDG\+module over $\B$ is called \emph{absolutely
acyclic} if it belongs to the minimal thick subcategory of
the homotopy category of coherent CDG\+modules $H^0(\B\coh)$
containing the total CDG\+modules of all the short exact sequences
of coherent CDG\+modules over $\B$ (with closed morphisms between
them).
 The quotient category of $H^0(\B\coh)$ by the thick subcategory
of absolutely acyclic CDG\+modules is called the \emph{absolute
derived category} of coherent CDG\+modules over $\B$ and denoted
by $\sD^\abs(\B\coh)$ \cite{Pkoszul}.

 For any quasi-coherent CDG\+algebra $\B$ over $X$, a quasi-coherent
CDG\+module over $\B$ is called \emph{coacyclic} if it belongs to
the minimal triangulated subcategory of the homotopy category of
quasi-coherent CDG\+modules $H^0(\B\qcoh)$ containing the total
CDG\+modules of all the short exact sequences of quasi-coherent
CDG\+modules over $\B$ and closed under infinite direct sums.
 The quotient category of $H^0(\B\coh)$ by the thick subcategory
of coacyclic CDG\+modules is called the \emph{coderived category} of
quasi-coherent CDG\+modules over $\B$ and denoted by
$\sD^\co(\B\qcoh)$ \cite{Psemi,Pkoszul}.

 Given an exact subcategory $\sE$ in the abelian category of
quasi-coherent graded left $\B$\+modules, one can define
the \emph{absolute derived category of left CDG\+modules over $\B$
with the underlying graded $\B$\+modules belonging to\/ $\sE$} as
the quotient category of the corresponding homotopy category by
its minimal thick subcategory containing the total CDG\+modules of
all the exact triples of CDG\+modules with the underlying graded
$\B$\+modules belonging to~$\sE$.
 The objects of the latter subcategory are called \emph{absolutely
acyclic with respect to\/ $\sE$} (or with respect to the DG\+category
of CDG\+modules with the underlying graded modules belonging to~$\sE$)
\cite{PP}.

 So one defines the absolute derived categories $\sD^\abs(\B\coh_\ffd)$
and $\sD^\abs(\B\coh_\fl)$ as the quotient categories of the homotopy
categories $H^0(\B\coh_\ffd)$ and $H^0(\B\coh_\fl)$ by the thick
subcategories of CDG\+modules absolutely acyclic with respect to 
$\B\coh_\ffd$ and $\B\coh_\fl$, respectively.

 When the exact subcategory $\sE$ is closed under infinite direct sums,
the thick subcategory of CDG\+modules \emph{coacyclic with respect
to\/ $\sE$} is the minimal triangulated subcategory of the homotopy
category CDG\+modules with the underlying graded modules belonging
to $\sE$, containing the total CDG\+modules of all the exact triples
of CDG\+modules with the underlying graded modules belonging to $\sE$
and closed under infinite direct sums.
 The quotient category by this thick subcategory is called
the \emph{coderived category of left CDG\+modules over $\B$ with
the underlying graded modules belonging to\/ $\sE$}
\cite{Psemi,PP}.

 Thus one defines the coderived category $\sD^\co(\B\qcoh_\fl)$ as
the quotient categories of the homotopy category $H^0(\B\qcoh_\fl)$
by the thick subcategory of CDG\+modules coacyclic with respect to
$\B\qcoh_\fl$.
 The definition of the coderived category $\sD^\co(\B\qcoh_\ffd)$
requires a little more care, since the class of graded modules of
finite flat dimension is not in general closed under infinite
direct sums.
 An object $\M\in H^0(\B\qcoh_\ffd)$ is said to be \emph{coacyclic
with respect to $\B\qcoh_\ffd$} if there exists an integer $d\ge0$ 
such that $\M$ is coacyclic with respect to the exact category of
quasi-coherent CDG\+modules of flat dimension~$\le d$.
 The coderived category of quasi-coherent CDG\+modules of finite
flat dimension is, by the definition, the quotient category of
$H^0(\B\qcoh_\ffd)$ by the above-defined thick subcategory of
coacyclic CDG\+modules~\cite[Section~3.2]{PP}.

\begin{rem}
 One may wonder whether coacyclicity (absolute acyclicity) of
quasi-coherent CDG\+modules (of a certain class) is a local notion.
 One general approach to this kind of problems is to consider
the Mayer--Vietoris/\v Cech exact sequence
$$\textstyle
 0\lrarrow\M\lrarrow \bigoplus_\alpha j_{U_\alpha *}j_{U_\alpha}^*
 \M\lrarrow\bigoplus_{\alpha<\beta} j_{U_\alpha\cap\U_\beta *}
 j_{U_\alpha\cap U_\beta}^*\M\lrarrow\dotsb\lrarrow0
$$
for a finite affine open covering $U_\alpha$ of~$X$.
 Since the inverse and direct images with respect to affine open
embeddings are exact and compatible with direct sums, they preserve
coacyclicity (absolute acyclicity).
 Hence if the restrictions of $\M$ to all $U_\alpha$ are coacyclic
(absolutely acyclic), then so is $\M$~itself.

 Alternatively, one can base this kind of argument on
the implications of the Noetherianness assumption, rather than
the separatedness assumption.
 For this purpose, one replaces a quasi-coherent CDG\+module $\M$
with its injective resolution (see Lemma~\ref{gorenstein-case}(b))
before writing down its \v Cech resolution.
 In this approach, the covering need not be affine, as injective
coacyclic objects are contractible, and direct images preserve
contractibility; but it is important that the restrictions to open
subschemes should preserve injectivity of quasi-coherent graded
$\B$\+modules (see~\cite[Theorem~II.7.18]{Har} and
Theorem~\ref{injective-sheaves}; cf.~\cite[Appendix~B]{TT}).

 When one is working with coherent CDG\+modules, the \v Cech
sequence argument is to be used in conjuction with
Proposition~\ref{embedding-prop} below.
 (Cf.\ Sections~\ref{cdg-supports} and~\ref{locality}.)
\end{rem}

\subsection{Finite flat dimension theorem}  \label{flat-dimension-thm}
 The next theorem is our main technical result on which the proofs
in Section~2 are based.

 Though we generally prefer the coderived categories
of (various classes of) infinitely generated CDG\+modules over
their absolute derived categories, technical considerations
sometimes force us to deal with the latter (see
Remark~\ref{embedding-prop}).
 Therefore, let $\sD^\abs(\B\qcoh_\fl)$, \ $\sD^\abs(\B\qcoh_\ffd)$,
and $\sD^\abs(\B\qcoh)$ denote the absolute derived categories of
(flat, of finite flat dimension, or arbitrary) quasi-coherent
CDG\+modules over a quasi-coherent CDG\+algebra~$\B$.

\begin{thm}
\textup{(a)} For any quasi-coherent CDG\+algebra $\B$ over $X$,
the functor\/ $\sD^\co(\B\qcoh_\fl)\rarrow\sD^\co(\B\qcoh_\ffd)$
induced by the embedding of DG\+categories $\B\qcoh_\fl\rarrow
\B\qcoh_\ffd$ is an equivalence of triangulated categories.
{\hbadness=2000\par}
\textup{(b)} For any quasi-coherent CDG\+algebra $\B$ over $X$,
the functor\/ $\sD^\abs(\B\qcoh_\fl)\rarrow\sD^\abs(\B\qcoh_\ffd)$
induced by the embedding of DG\+categories $\B\qcoh_\fl\rarrow
\B\qcoh_\ffd$ is an equivalence of triangulated categories. \par
\textup{(c)} For any quasi-coherent CDG\+algebra $\B$ over $X$
such that the underlying quasi-coherent graded algebra $\B$ is
Noetherian, the functor\/ $\sD^\abs(\B\coh_\fl)\rarrow\sD^\abs
(\B\coh_\ffd)$ induced by the embedding of DG\+categories
$\B\coh_\fl\rarrow\B\coh_\ffd$ is an equivalence of triangulated
categories.
\end{thm}

\begin{proof}
 The proof follows that of~\cite[Theorem~3.2]{PP} (see
also~\cite[Theorem~7.2.2]{Psemi}) with some modifications.
 We will prove part~(a); the proofs of parts~(b\+c) are completely
similar.
 (Alternatively, parts~(b\+c) can be deduced from
Proposition~\ref{embedding-prop}(a\+b) below.)

 Given an affine open subscheme $U\subset X$ and a graded module
$P$ over the graded ring $\B(U)$, one can construct
the freely generated CDG\+module $G^+(P)$ over the CDG\+ring
$\B(U)$ in the way explained in~\cite[proof of Theorem~3.6]{Pkoszul}.
 The elements of $G^+(P)$ are formal expressions of the form $p+dq$,
where $p$, $q\in P$.
 Given a quasi-coherent graded module $\P$ over $\B$,
the CDG\+modules $G^+(\P(U))$ glue together to form a quasi-coherent
CDG\+module $G^+(\P)$ over~$\B$.
 For any quasi-coherent CDG\+module $\M$ over $\B$, there is
a bijective correspondence between morphisms of graded $\B$\+modules
$\P\rarrow\M$ and closed morphisms of CDG\+modules
$G^+(\P)\rarrow\M$ over~$\B$.
 There is a natural short exact sequence of quasi-coherent graded
$\B$\+modules $\P\rarrow G^+(\P)\rarrow \P[-1]$.
 The quasi-coherent CDG\+module $G^+(\P)$ is naturally contractible
with the contracting homotopy~$t_\P$ given by the composition
$G^+(\P)\rarrow\P[-1]\rarrow G^+(\P)[-1]$.

 Due to our assumption on $X$, for any quasi-coherent $\O_X$\+module
$\K$ over $X$ there exists a surjective morphism $\E\rarrow\K$
onto $\K$ from a direct sum $\E$ of locally free sheaves of finite
rank on~$X$.
 Hence for any quasi-coherent graded $\B$\+module $\M$ there is
a surjective morphism onto $\M$ from a flat quasi-coherent graded
$\B$\+module $\P=\bigoplus_n\B\ot_{\O_X}\E_n[n]$, and for any
quasi-coherent CDG\+module $\M$ over $\B$ there is a surjective
closed morphism onto $\M$ from the CDG\+module $G^+(\P)\in\B\qcoh_\fl$.
 (In fact, parts~(a\+b) of Theorem can be proven without
the assumption of enough vector bundles on~$X$, since there are always
enough flat sheaves; see Remark~\ref{w-flat-cor} and
Lemma~\ref{flat-sheaves}.)

 Now the construction from~\cite[proof of Theorem~3.6]{Pkoszul}
provides for any object $\M$ of $\B\qcoh_\ffd$ a closed morphism
onto $\M$ from an object of $\B\qcoh_\fl$ with the cone absolutely
acyclic with respect to $\B\qcoh_\ffd$.
 To obtain this morphism, one picks a finite left resolution of
$\M$ consisting of objects from $\B\qcoh_\fl$ with closed morphisms
between them, and takes the total CDG\+module of this resolution.
 By~\cite[Lemma~1.6]{Pkoszul}, it follows that the triangulated
category $\sD^\co(\B\qcoh_\ffd)$ is equivalent to the quotient category
of $H^0(\B\qcoh_\fl)$ by its intersection in $H^0(\B\qcoh_\ffd)$
with the thick subcategory of CDG\+modules coacyclic with respect
to $\B\qcoh_\ffd$.
 It only remains to show that any object of $H^0(\B\qcoh_\fl)$ that is
coacyclic with respect to $\B\qcoh_\ffd$ is coacyclic with respect
to $\B\qcoh_\fl$.

 Let us call a quasi-coherent CDG\+module $\M$ over $\B$
\emph{$d$\+flat} if its underlying quasi-coherent graded
$\B$\+module $\M$ has flat dimension not exceeding~$d$.
 A $d$\+flat quasi-coherent CDG\+module is said to be
\emph{$d$\+coacyclic} if it is homotopy equivalent to
a CDG\+module obtained from the total CDG\+modules of exact
triples of $d$\+flat CDG\+modules using the operations of
cone and infinite direct sum.
 Our goal is to show that any $0$\+flat $d$\+coacyclic CDG\+module
is $0$\+coacyclic.
 For this purpose, we will prove that any $(d-1)$-flat
$d$\+coacyclic CDG\+module is $(d-1)$-coacyclic; the desired
assertion will then follow by induction.

 It suffices to construct for any $d$\+coacyclic CDG\+module $\M$
a $(d-1)$-coacyclic CDG\+module $\L$ with a $(d-1)$-coacyclic
CDG\+submodule $\K$ such that the quotient CDG\+module $\L/\K$
is isomorphic to~$\M$.
 Then if $\M$ is $(d-1)$-flat, it would follow that both the cone
of the morphism $\K\rarrow\L$ and the total CDG\+module of
the exact triple $\K\rarrow\L\rarrow\M$ are $(d-1)$-coacyclic, so
$\M$ also is.
 The construction is based on four lemmas similar to those
in~\cite[Section~3.2]{PP}.

\begin{lema}
 Let $\M$ be the total CDG\+module of an exact triple
of $d$\+flat quasi-coherent CDG\+modules $\M'\rarrow\M''\rarrow
\M'''$ over~$\B$.
 Then there exists a surjective closed morphism onto $\M$ from
a contractible $0$\+flat CDG\+module $\P$ with a $(d-1)$-coacyclic
kernel~$\K$.
\end{lema}

\begin{proof}
 Choose $0$\+flat quasi-coherent CDG\+modules $\P'$ and $\P'''$
such that there exist surjective closed morphisms $\P'\rarrow\M'$
and $\P'''\rarrow\M''$.
 Then there exists a surjective morphism from the exact
triple of CDG\+modules $\P'\rarrow\P'\oplus\P'''\rarrow\P'''$
onto the exact triple $\M'\rarrow\M''\rarrow\M'''$.
 The rest of the proof is similar to that in~\cite{PP}.
\end{proof}

\begin{lemb}
\textup{(a)} Let $\K'\subset\L'$ and $\K''\subset\L''$ be
$(d-1)$-coacyclic CDG\+submodules in $(d-1)$-coacyclic CDG\+modules,
and let $\L'/\K'\rarrow\L''/\K''$ be a closed morphism of CDG\+modules.
 Then there exists a $(d-1)$-coacyclic CDG\+module $\L$ with
a $(d-1)$-coacyclic CDG\+submodule $\K$ such that $\L/\K\simeq
\cone(\L'/\K'\to\L''/\K'')$. \par
\textup{(b)} In the situation of \textup{(a)}, assume that the morphism
$\L'/\K'\rarrow\L''/\K''$ is injective with a $d$\+flat cokernel $\M_0$.
 Then there exists a $(d-1)$-coacyclic CDG\+module $\L_0$ with
a $(d-1)$-coacyclic CDG\+submodule $\K_0$ such that $\L_0/\K_0
\simeq\M_0$.
\end{lemb}

\begin{proof}
 The proof is similar to that in~\cite{PP}.
\end{proof}

\begin{lemc}
 For any contractible $d$\+flat CDG\+module $\M$ there exists
a surjective closed morphism onto $\M$ from a contractible $0$\+flat
CDG\+module $\L$ with a $(d-1)$-coacyclic kernel~$\K$.
\end{lemc}

\begin{proof}
 Let $p\:\P\rarrow\M$ be a surjective morphism onto the quasi-coherent
graded $\B$\+module $\M$ from a flat quasi-coherent graded
$\B$\+module $\P$, and $\tilde p\:G^+(\P)\rarrow\M$ be the induced
surjective closed morphism of quasi-coherent CDG\+modules.
 Let $t\:\M\rarrow\M$ be a contracting homotopy for $\M$ and
$t_\P\:G^+(\P)\rarrow G^+(\P)$ be the natural contracting homotopy
for $G^+(\P)$.
 Then $\tilde u = \tilde p t_\P - t\tilde p\:G^+(\P)\rarrow\M$ is
a closed morphism of quasi-coherent CDG\+modules of degree~$-1$.
 Denote by~$u$ the restriction of~$\tilde u$ to $\P\subset G^+(\P)$.
 There exists a surjective morphism from a flat quasi-coherent
graded $\B$\+module $\Q$ onto the fibered product of
the morphisms $p\:\P\rarrow\M$ and $u\:\P\rarrow \M$.
 Hence we obtain a surjective morphism of quasi-coherent graded
$\B$\+modules $q\:\Q\rarrow\P$ and a morphism of quasi-coherent
graded $\B$\+modules $v\:\Q\rarrow\P$ of degree~$-1$
such that $uq=pv$.

 The morphism~$q$ induces a surjective closed morphism of
quasi-coherent CDG\+mod\-ules $\tilde q\: G^+(\Q)\rarrow G^+(\P)$.
 The morphism $\tilde q$ is homotopic to zero with the natural
contracting homotopy $\tilde q t_\Q = t_\P\tilde q$.
 The morphism~$v$ induces a closed morphism of CDG\+modules
$\tilde v\:G^+(\Q)\rarrow G^+(\P)$ of degree~$-1$.
 The morphism $t_\P\tilde q - \tilde v$ is another contracting
homotopy for~$\tilde q$.
 The latter homotopy forms a commutative square with the morphisms
$\tilde p$, \ $\tilde p\tilde q$, and the contracting homotopy~$t$
for the CDG\+module $\M$.

 Let $\N$ be the kernel of the morphism $\tilde p\tilde q\:
G^+(\Q)\rarrow\M$ and $\K$ be the kernel of the morphism
$\tilde p\:G^+(\P)\rarrow\M$.
 Then the natural surjective closed morphism $r\:\N\rarrow\K$ is
homotopic to zero; the restriction of the map $t_\P\tilde q -
\tilde v$ provides the contracting homotopy that we need.
 In addition, the kernel $G^+(\ker q)$ of the morphism~$r$
is contractible.
 So the cone of the morphism~$r$ is isomorphic to $\K\oplus\N[1]$,
and on the other hand there is an exact triple
$G^+(\ker q)[1]\rarrow\cone(r)\rarrow\cone(\id_\K)$.
 Since $\K$ is $(d-1)$-flat and $\ker q$ is flat, this proves that
$\K$ is $(d-1)$-coacyclic.
 It remains to take $\L=G^+(\P)$.
\end{proof}

\begin{lemd}
 Let $\M\rarrow\M'$ be a homotopy equivalence of $d$\+flat CDG\+modules
such that $\M'$ is the quotient CDG\+module of a $(d-1)$-coacyclic
CDG\+module by a $(d-1)$-coacyclic CDG\+submodule.
 Then $\M$ is also such a quotient.
\end{lemd}

\begin{proof}
 The proof is similar to that in~\cite{PP}.
\end{proof}

 It is clear that the property of a CDG\+module to be presentable
as the cokernel of an injective closed morphism of $(d-1)$-coacyclic
CDG\+modules is stable under infinite direct sums.
 This finishes our construction and the proof of Theorem.
\end{proof}

\begin{rem}
 The assertion of part~(c) of Theorem~\ref{flat-dimension-thm} can be
equivalently rephrased with flat modules replaced by locally projective
ones.
 Indeed, a finitely presented module over a ring is flat if and only if
it is projective.

 In the infinitely generated situation of parts~(a\+b), flatness of
quasi-coherent sheaves is different from their local projectivity
(which is a stronger condition), but the assertions remain true
after one replaces the former with the latter.
 The same applies to Proposition~\ref{embedding-prop}(a) below.
 Indeed, by Theorem~\ref{locally-projective}, for any quasi-coherent
graded algebra $\B$ over an affine scheme $U$, projectivity of
a graded module over the graded ring $\B(U)$ is a local notion.
 Taking this fact into account, our proof goes through for locally
projective quasi-coherent graded modules in place of flat ones and
the locally projective dimension (defined as the minimal length of
a locally projective resolution) in place of the flat dimension.

 When $\B=\O_X$, local projectivity of quasi-coherent modules is
equivalent to local freeness~\cite[Corollary~4.5]{Ba}.
 Furthermore, in this case, assuming additionally that $X$ has
finite Krull dimension, the classes of quasi-coherent sheaves of
finite flat dimension and of finite locally projective dimension
coincide~\cite[Corollaire~II.3.3.2]{RG}.
\end{rem}

\subsection{Fully faithful embedding} \label{embedding-prop}
 The next proposition is stronger than Theorem~\ref{flat-dimension-thm}
in some respects, and is proven by an entirely different technique.

\begin{prop}
\textup{(a)} For any quasi-coherent CDG\+algebra $\B$ over $X$,
the functor\/ $\sD^\abs(\B\qcoh_\fl)\rarrow \sD^\abs(\B\qcoh)$
induced by the embedding of DG\+categories $\B\qcoh_\fl\rarrow
\B\qcoh$ is fully faithful. {\hbadness=1300\par}
 Furthermore, let $\B$ be a quasi-coherent CDG\+algebra over $X$
such that the underlying quasi-coherent graded algebra
$\B$ is Noetherian.
 Then \par
\textup{(b)} the functor $\sD^\abs(\B\coh_\fl)\rarrow\sD^\abs(\B\coh)$
induced by the embedding of DG\+categories $\B\coh_\fl\rarrow
\B\coh$ is fully faithful; \par
\textup{(c)} the functor\/ $\sD^\abs(\B\coh)\rarrow\sD^\abs(\B\qcoh)$
induced by the embedding of DG\+categories $\B\coh\rarrow\B\qcoh$
is fully faithful; \par
\textup{(d)} the functor\/ $\sD^\abs(\B\coh)\rarrow\sD^\co(\B\qcoh)$
induced by the embedding of DG\+categories $\B\coh\rarrow\B\qcoh$
is fully faithful and its image forms a set of compact generators
for\/ $\sD^\co(\B\qcoh)$.
\end{prop}

\begin{proof}
 The proof of part~(d) in the case when $X$ is affine can be
found in~\cite[Section~3.11]{Pkoszul} (the part concerning
compact generation belongs to D.~Arinkin).
 The proof in the general case is similar; and part~(c) can be
also proven in the way similar to~\cite[Theorem~3.11.1]{Pkoszul}.
 Part~(b) in the affine case is easy and follows from
the semiorthogonality property of CDG\+modules with projective
underlying graded modules and absolutely acyclic/contraacyclic
CDG\+modules~\cite[Theorem~3.5(b)]{Pkoszul}, since finitely
generated flat modules over a Noetherian ring are projective.
 A detailed proof of part~(b) in the general case is given below;
and the proof of part~(a) (which does not automatically simplify
in the affine case) is similar.

 We will show that any morphism $\E\rarrow\L$ from a CDG\+module
$\E\in H^0(\B\coh_\fl)$ to a CDG\+module $\L\in H^0(\B\coh)$
absolutely acyclic with respect to $\B\coh$ can be annihilated
by a morphism $\P\rarrow\E$ from a CDG\+module
$\P\in H^0(\B\coh_\fl)$ with a cone of the morphism $\P\rarrow\E$
being absolutely acyclic with respect to $\B\coh_\fl$.
 By the definition, the CDG\+module $\L$ is a direct summand of
a CDG\+module homotopy equivalent to a CDG\+module obtained from
the totalizations of exact triples of CDG\+modules in $\B\coh$
using the operation of passage to the cone of a closed morphism
repeatedly.
 It suffices to consider the case when $\L$ itself is obtained
from totalizations of exact triples using cones.
 We proceed by induction in the number of operations of passage
to the cone in such a construction of~$\L$.

 So we assume that there is a distinguished triangle
$\K\rarrow\L\rarrow\M\rarrow\K[1]$ in $H^0(\B\coh)$ such that
$\M$ is the total CDG\+module of an exact triple of
CDG\+modules in $\B\coh$, while the CDG\+module $\K$ has
the desired property with respect to morphisms into it from
all CDG\+modules $\F\in H^0(\B\coh_\fl)$.
 If we knew that the object $\M$ also has the same property,
it would follow that the composition $\E\rarrow\L\rarrow\M$ can be
annihilated by a morphism $\F\rarrow\E$ with $\F\in H^0(\B\coh_\fl)$
and a cone absolutely acyclic with respect to $\B\coh_\fl$.
 The composition $\F\rarrow\E\rarrow\L$ then factorizes through
$\K$, and the morphism $\F\rarrow\K$ can be annihilated by
a morphism $\P\rarrow\F$ with $\P\in H^0(\B\coh_\fl)$ and a cone
absolutely acyclic with respect to $\B\coh_\fl$.
 The composition $\P\rarrow\F\rarrow\E$ provides the desired
morphism $\P\rarrow\E$.

 Thus it remains to construct a morphism $\F\rarrow\E$ with
the required properties annihilating a morphism $\E\rarrow\M$,
where $\M$ is the total CDG\+module of an exact triple of
CDG\+modules $\U\rarrow\V\rarrow\W$.
 For any graded module $\N$ over $\B$, morphisms of graded
$\B$\+modules $\N\rarrow\M$ of degree~$n$ are represented by triples
$(f,g,h)$, where $f\:\N\rarrow\U$ is a morphism of degree~$n+1$,
\ $g\:\N\rarrow\V$ is a morphism of degree~$n$, and
$h\:\N\rarrow\W$ is a morphism of degree~$n-1$.
 Denote the closed morphisms in the exact triple $\U\rarrow\V
\rarrow\W$ by $j\:\U\rarrow\V$ and $k\:\V\rarrow\W$.

\begin{leme}
 Let $\N$ be a CDG\+module over $\B$ and $\M$ be the total CDG\+module
of an exact triple of CDG\+modules $\U\rarrow\V\rarrow\W$ as above.
Then \par
\textup{(a)} the differential of a morphism of graded
$\B$\+modules $\N\rarrow\M$ of degree~$n$ rep\-resented by
a triple $(f,g,h)$ is given by the rule $d(f,g,h) =
(-df\;-jf+dg\;kg-dh)$; \par
\textup{(b)} when $(f,g,h)$ is a closed morphism of CDG\+modules of
degree~$n$ and the morphism of graded $\B$\+modules $h\:\N\rarrow\W$
can be lifted to a morphism of graded $\B$\+modules $t\:\N\rarrow\V$
of degree~$n-1$, the morphism $(f,g,h)$ is homotopic to zero.
\end{leme}

\begin{proof}
 The complex of morphisms in the DG\+category of CDG\+modules
$\Hom_\B(\N,\M)$ is the total complex of the bicomplex
of abelian groups $\Hom_\B(\N,\U)\rarrow\Hom_\B(\N,\V)\rarrow
\Hom_\B(\N,\W)$.
 The formula in~(a) is the formula for the differential
of a total complex. {\hbadness=1500\par}

 Furthermore, the sequence $0\rarrow\Hom_\B(\N,\U)\rarrow
\Hom_\B(\N,\V)\rarrow\Hom_\B(\N,\W)$ is exact.
 Let $\Hom'_\B(\N,\W)$ denote the cokernel of the morphism of
complexes $\Hom_\B(\N,\U)\rarrow\Hom_\B(\N,\V)$; then
$\Hom'_\B(\N,\W)$ is a subcomplex of $\Hom_\B(\N,\W)$ and
the total complex of the bicomplex $\Hom_\B(\N,\U)\rarrow
\Hom_\B(\N,\V)\rarrow\Hom'_\B(\N,\W)$ is an acyclic subcomplex
of $\Hom_\B(\N,\M)$.
 Hence any cocycle in $\Hom_\B(\N,\M)$ that belongs to this
subcomplex is a coboundary.

 To present the same argument using our letter notation for
morphisms, assume that $kt=h$.
 Then $k(dt-g) = dh-kg = 0$, so there exists a morphism of graded
$\B$\+modules $s\:\N\rarrow\U$ of degree~$n$ such that $dt-g = js$.
 Then $jds = -dg = -jf$, hence $ds=-f$ and $d(s,t,0)=(f,g,h)$.
\end{proof}

 Recall the notation $G^+(\Q)$ for the CDG\+module freely generated
by a graded $\B$\+module $\Q$ (see the beginning of the proof of
Theorem~\ref{flat-dimension-thm}).

\begin{lemf}
 Let $\M$ be the total CDG\+module of an exact triple of
CDG\+modules $\U\rarrow\V\rarrow\W$ as above, and let $\Q$
be a graded $\B$\+module.
 Assume that a morphism of graded $\B$\+modules $p\:\Q\rarrow\M$
of degree~$n$ with the components $(f,g,h)$ is given such that
the component $h\:\Q\rarrow\W$ can be lifted to a morphism of
graded $\B$\+modules $t\:\Q\rarrow\V$ of degree~$n-1$.
 Let $\tilde p\:G^+(\Q)\rarrow\M$ be the induced closed morphism
of CDG\+modules of degree~$n$ and $(\tilde f,\tilde g,\tilde h)$
be its three components.
 Then the morphism of graded $\B$\+modules $\tilde h\:G^+(\Q)
\rarrow\W$ can be lifted to a morphism of graded $\B$\+modules
$\tilde t\:G^+(\Q)\rarrow\V$ of degree~$n-1$.
\end{lemf}

\begin{proof}
 Notice that any closed morphism of CDG\+modules $G^+(\Q)\rarrow\M$
is homotopic to zero, since the CDG\+module $G^+(\Q)$
is contractible.
 The conclusion of the lemma is stronger, and we will need its full
strength.
 The argument consists in a computation in the letter notation for
morphisms.

 For any CDG\+module $\N$ over $\B$, morphisms of graded $\B$\+modules
$\tilde r\: G^+(\Q)\rarrow\N$ of degree~$n-1$ are uniquely determined
by their restriction to $\Q$ and the restriction to $\Q$ of their
differential $d\tilde r$, which can be arbitrary morphisms of graded
$\B$\+modules $\Q\rarrow\N$ of the degrees $n-1$ and~$n$, respectively.
 Extend our morphism $t\:\Q\rarrow\V$ to a morphism of graded
$\B$\+modules $\tilde t\:G^+(\Q)\rarrow\V$ of degree~$n-1$ such that
$(d\tilde t)|_\Q = g$.
 Then $k\tilde t|_\Q = kt = h = \tilde h|_\Q$ and 
$(d(k\tilde t))|_\Q = k(d\tilde t)|_\Q = kg = k\tilde g|_\Q =
(d\tilde h)|_\Q$ by Lemma~E(a), hence $k\tilde t = \tilde h$.
\end{proof}

 Now represent a closed morphism $\E\rarrow\M$ by a triple $(f,g,h)$
of morphisms of the degrees~$1$, $0$, and~$-1$, respectively.
 Let $\Q$ be a flat coherent graded $\B$\+module mapping surjectively
onto the fibered product of the morphisms $k\:\V\rarrow\W$ and
$h\:\E\rarrow\W$ (see the beginning of the proof of
Theorem~\ref{flat-dimension-thm} again).
 Then there is a surjective morphism of graded $\B$\+modules
$q\:\Q\rarrow\E$ and its composition with the morphism
$h\:\E\rarrow\W$ can be lifted to a morphism of graded $\B$\+modules
$t\:\Q\rarrow\V$ of degree~$-1$.
 Consider the induced morphism of CDG\+modules $\tilde q\:
G^+(\Q)\rarrow\E$.
 By Lemma~F, the composition $h\tilde q\:G^+(\Q)\rarrow\W$
can be lifted to a morphism of graded $\B$\+modules
$\tilde t\:G^+(\Q)\rarrow\V$ of degree~$-1$.

 Let $\R$ denote the kernel of the closed morphism~$\tilde q$.
 Then the cone $\F$ of the embedding $\R\rarrow G^+(\Q)$ maps
naturally onto $\E$ with the cone absolutely acyclic with respect to
$\B\coh_\fl$.
 As a graded $\B$\+module, the CDG\+module $\F$ is isomorphic
to $G^+(\Q)\oplus\R[1]$; the composition $\F\rarrow\E\rarrow\M$
factorizes through the direct summand $G^+(\Q)$, where it is
defined by the triple $(f\tilde q\;g\tilde q\;h\tilde q)$.
 Since the morphism $h\tilde q$ can be lifted to $\V$, so can
the corresponding component $\F\rarrow\W$ of the morphism
$\F\rarrow\M$.
 Thus the latter morphism is homotopic to zero by Lemma~E(b).
\end{proof}

 In some cases the use of Lemma~F in the above proof of part~(b)
can be avoided.
 Assume that $X$ is a projective scheme over a Noetherian ring and
the category of coherent graded $\B$\+modules is equivalent to
the category of coherent modules over some coherent (graded)
$\O_X$\+algebra $\A$.
 In this situation, one takes $\Q$ to be the graded $\B$\+module
corresponding to the (graded) $\A$\+module induced from a large enough
finite direct sum of (shifts of) copies of a sufficiently negative
invertible $\O_X$\+module; then there is a surjective morphism of
graded $\B$\+modules $\Q\rarrow\E$ and any morphism of graded
$\B$\+modules $G^+(\Q)\rarrow\W$ lifts to~$\V$.

\begin{rem}
 We do \emph{not} know how to extend the proof of 
Proposition~\ref{embedding-prop}(a\+b) to the coderived categories
of quasi-coherent CDG\+modules.
 Instead, this argument appears to be well-suited for use with
the \emph{contraderived} categories (see~\cite[Section~3.3]{Pkoszul}
for the definition).
 In particular, it allows to show that the contraderived category
of left CDG\+modules over a CDG\+ring $B$ with a right coherent
underlying graded ring is equivalent to the contraderived category
of CDG\+modules whose underlying graded $B$\+modules are flat
(cf.~\cite[paragraph after the proof of Theorem~3.8]{Pkoszul}).

 This is the main reason why we sometimes find it easier to deal
with the absolute derived rather than the coderived categories of
infinitely generated CDG\+modules (cf.\ Remark~\ref{infinite-matrix}).
 On the other hand, for the coderived category of quasi-coherent
CDG\+modules we have the compact generation result (part~(d) of
Proposition), the results and arguments of
Sections~\ref{gorenstein-case}, \ref{cdg-supports},
\ref{serre-duality}, \ref{stable-derived}, etc.
 The conditions under which these two versions of the construction
of the derived category of the second kind for a given class of
CDG\+modules lead to the same triangulated category are
discussed below in Section~\ref{homol-dimension-thm}.
\end{rem}

\subsection{Finite homological dimension theorem}
\label{homol-dimension-thm}
 Let $\B\qcoh_\lp$ denote the DG\+cate\-gory of quasi-coherent
CDG\+modules over $\B$ whose underlying graded $\B$\+modules are
locally projective (see Remark~\ref{flat-dimension-thm} and
Theorem~\ref{locally-projective}).
 Denote by $\sD^\co(\B\qcoh_\lp)$ and $\sD^\abs(\B\qcoh_\lp)$
the corresponding coderived and absolute derived categories.

\begin{thm}
 The triangulated categories\/ $\sD^\co(\B\qcoh_\lp)$ and\/
$\sD^\abs(\B\qcoh_\lp)$ coincide, i.~e., every CDG\+module
over $\B$ that is coacyclic with respect to $\B\qcoh_\lp$
is also absolutely acyclic with respect to $\B\qcoh_\lp$.
\end{thm}

\begin{proof}
 The reason for this assertion to be true is that the exact
category of locally projective graded $\B$\+modules has finite
homological dimension~\cite[Lemma~1.12]{Or1} and exact functors
of infinite direct sums.
 If this exact category also had enough injectives, the simple
argument from~\cite[Theorem~3.6(a) and Remark~3.6]{Pkoszul}
would suffice to establish the desired $\sD^\co=\sD^\abs$
isomorphism for it (see also~\cite[Remark~2.1]{Psemi}).
 The lengthy argument below is designed to provide a way around
the injective objects issue in this kind of proof.

 Our aim is to show that for any closed morphism $\P\rarrow\L$
from a CDG\+module $\P\in\B\qcoh_\lp$ to a CDG\+module $\L$ absolutely
acyclic with respect to $\B\qcoh_\lp$ there exists an exact
sequence $0\rarrow\Q_d\rarrow\Q_{d-1}\rarrow\dotsb\rarrow\Q_0
\rarrow\P\rarrow 0$ of CDG\+modules and closed morphisms in
$\B\qcoh_\lp$ such that the induced morphism from the total
CDG\+module of $\Q_d\rarrow\dotsb\rarrow\Q_0$ to $\L$ is
homotopic to zero.
 Here $d$~is a fixed integer equal to the homological dimension
of the exact category of locally projective graded $\B$\+modules,
which does not exceed the number of open subsets in an affine
covering of~$X$ minus one.

 Taking $\P=\L$ and the morphism $\P\rarrow\L$ to be the identity,
we will then conclude that $\P$ is isomorphic to a direct summand
of the total CDG\+module of 
$\Q_d\rarrow\dotsb\rarrow\Q_0\rarrow\P$ in $H^0(\B\qcoh_\lp)$.
 Hence an object of $H^0(\B\qcoh_\lp)$ is absolutely acyclic with
respect to $\B\qcoh_\lp$ if and only if it is isomorphic to
a direct summand of the total CDG\+module of a $(d+2)$-term
exact sequence of CDG\+modules from $\B\qcoh_\lp$ with closed
morphisms between them.
 It will immediately follow that the class of CDG\+modules 
absolutely acyclic with respect to $\B\qcoh_\lp$ is closed under
infinite direct sums, so it coincides with the class of
coacyclic CDG\+modules.

 We can suppose that there exists a sequence of distinguished
triangles $\K_{i-1}\rarrow\K_i\rarrow\M_i\rarrow\K_{i-1}[1]$ 
in $H^0(\B\qcoh_\lp)$ such that $\K_0=0$, \ $\K_n=\L$, and
$\M_i$ is the total CDG\+module of an exact triple
$\U_i\rarrow\V_i\rarrow\W_i$ of CDG\+modules from $\B\qcoh_\lp$
for all $1\le i\le n$.
 We will start with constructing an exact sequence
$0\rarrow\Q'_n\rarrow\dotsb\rarrow\Q'_0\rarrow\P \rarrow0$ with
the above properties, but of the length~$n$ rather than~$d$.
 Then we will use the finite homological dimension property of
locally projective graded $\B$\+modules in order to obtain 
the desired resolution $\Q_\bu$ of a fixed length~$d$ from
a resolution~$\Q'_\bu$.

\begin{lemg}
 Let $\M$ be the total CDG\+module of an exact triple
$\U\rarrow\V\rarrow\W$ of CDG\+modules from $\B\qcoh_\lp$
and $\K\rarrow\L\rarrow\M\rarrow\K[1]$ be a distinguished
triangle in $H^0(\B\qcoh_\lp)$.
 Then for any CDG\+module $\P\in\B\qcoh_\lp$ and a morphism
$\P\rarrow\L$ in $H^0(\B\qcoh_\lp)$ there exists an exact
triple $\R\rarrow\Q\rarrow\P$ of CDG\+modules from $\B\qcoh_\lp$
and a morphism $\R[1]\rarrow\K$ in $H^0(\B\qcoh_\lp)$
such that the composition $\F\rarrow\P\rarrow\L$, where
$\F$ is the cone of the closed morphism $\R\rarrow\Q$,
is equal to the composition $\F\rarrow\R[1]\rarrow\K\rarrow\L$
in $H^0(\B\qcoh_\lp)$.
\end{lemg}

\begin{proof}
 The argument is based on Lemmas~E\+-F from
Section~\ref{embedding-prop}.
 We can assume that $\L$ is the cone of a closed morphism
$\M[-1]\rarrow\K$ and fix a closed morphism $\P\rarrow\L$
representing the given morphism in the homotopy category.
 Arguing as in the proof of Proposition~\ref{embedding-prop},
we can construct a surjective closed morphism $\Q'\rarrow\P$
onto $\P$ from a CDG\+module $\Q'\in\B\qcoh_\lp$ such that
the composition $\Q'\rarrow\P\rarrow\L\rarrow\M\rarrow\W[-1]$
lifts to a morphism of graded $\B$\+modules $\Q'\rarrow\V[-1]$.
 Here it suffices to apply the functor $G^+$ to the fibered
product of the morphisms of graded $\B$\+modules
$\P\rarrow\W[-1]$ and $\V[-1]\rarrow\W[-1]$, and use Lemma~F.

 Then the morphism $\Q'\rarrow\M$ is homotopic to zero with
a natural contracting homotopy (provided by the proof of Lemma~E),
so the morphism $\Q'\rarrow\L$ factorizes, up to a homotopy,
as the composition of a naturally defined closed morphism
$\Q'\rarrow\K$ and the closed morphism $\K\rarrow\L$.
 Set $\Q$ to be the cocone of the closed morphism $\Q'\rarrow\K$;
then we have a surjective closed morphism $\Q\rarrow\Q'$ such
that the composition $\Q\rarrow\Q'\rarrow\K$ is homotopic
to zero.

 Let $\R$ be kernel of the morphism $\Q\rarrow\P$ and $\F$
be the cone of the morphism $\R\rarrow\Q$; then there is
a natural closed morphism $\F\rarrow\P$.
 Using Lemma~E and arguing as in the end of the proof of
Proposition~\ref{embedding-prop} again, we can conclude that
the composition $\F\rarrow\P\rarrow\L\rarrow\M$ is homotopic
to zero.
 Indeed, the composition $\F\rarrow\M\rarrow\W[-1]$ lifts to
a graded $\B$\+module morphism $\F\rarrow\V[-1]$, since
$\F\simeq\Q\oplus\R[-1]$ as a graded $\B$\+module, the morphism
$\F\rarrow\M$ factorizes through the projection of $\F$ onto
$\Q$, and the morphism $\Q\rarrow\Q'\rarrow\W[-1]$ lifts to
a graded $\B$\+module morphism $\Q\rarrow\Q'\rarrow\V[-1]$ by our
construction.

 Notice that the contracting homotopy that we have obtained for
the closed morphism $\F\rarrow\M$ forms a commutative diagram
with the closed morphisms $\Q\rarrow\F$, \ $\Q\rarrow\Q'$, and
the contracting homotopy that we have previously had for
the closed morphism $\Q'\rarrow\M$ (since so do the liftings
$\F\rarrow\V[-1]$ and $\Q'\rarrow\V[-1]$).
 This allows to factorize, up to a homotopy, the closed morphism
$\F\rarrow\L$ as the composition of a closed morphism
$\F\rarrow\K$ and the closed morphism $\K\rarrow\L$ in
such a way that the morphism $\F\rarrow\K$ forms
a commutative diagram with the closed morphisms $\Q\rarrow\F$,
\ $\Q\rarrow\Q'$, and the closed morphism $\Q'\rarrow\K$
that we have previously constructed.
 The composition $\Q\rarrow\F\rarrow\K$, being equal to
the composition $\Q\rarrow\Q'\rarrow\K$, is homotopic
to zero; hence the morphism $\F\rarrow\K$ factorizes through
the closed morphism $\F\rarrow\R[1]$ in $H^0(\B\qcoh_\lp)$.
\end{proof}

 Applying Lemma~G to the morphism $\P\rarrow\L$ and the distinguished
triangle $\K_{n-1}\rarrow\L\rarrow\M_n\rarrow\K_{n-1}$, we obtain
an exact triple $\R'_0\rarrow\Q'_0\rarrow\P$ and a morphism
$\R'_0[1]\rarrow\K_{n-1}$ in $H^0(\B\qcoh_\lp)$.
 Applying the same lemma again to the morphism $\R'_0[1]
\rarrow\K_{n-1}$ and the distinguished triangle $\K_{n-2}\rarrow
\K_{n-1}\rarrow\M_{n-1}\rarrow\K_{n-2}[1]$, we construct an exact
triple $\R'_1\rarrow\Q'_1\rarrow\R'_0$ and a morphism
$\R'_1[2]\rarrow\K_{n-2}$, etc.
 Finally we obtain an exact triple $\R'_{n-1}\rarrow\Q'_{n-1}
\rarrow\R'_{n-2}$ and a morphism $\R'_{n-1}[n]\rarrow\K_0=0$.

 Let us check that the natural morphism from the total CDG\+module
of the complex $0\rarrow\R'_{n-1}\rarrow\Q'_{n-1}\rarrow\dotsb
\rarrow\Q'_0$ to the CDG\+module $\L$ is homotopic to zero.
 Denote this morphism by~$f_n$.
 It factorizes naturally through the cone $\F_0$ of the closed
morphism $\R'_0\rarrow\Q'_0$, and the morphism $\F_0\rarrow\L$
is homotopic to the composition $\F_0\rarrow\R'_0[1]\rarrow\K_{n-1}
\rarrow\L$.
 Hence, up to the homotopy, the morphism~$f_n$ factorizes through
the morphism~$f_{n-1}$ from the total CDG\+module of the complex
$0\rarrow\R'_{n-1}\rarrow\Q'_{n-1}\rarrow\dotsb\rarrow\Q'_1$
to $\K_{n-1}$ induced by the morphism $\R'_0[1]\rarrow\K_{n-1}$.
 Continuing to argue in this way, we conclude that the morphism~$f$
factorizes, up to a homotopy, through the morphism $f_0\:\R'_{n-1}[n]
\rarrow\K_0=0$.

 It remains to ``cut'' our exact sequence of an unknown length~$n$
to a fixed size~$d$.
 For this purpose, we will assume that $n>d$ and construct from
our exact sequence of length~$n$ another exact sequence with
the same properties, but of the length $n-1$.
 This part of the argument is based on the following lemma.

\begin{lemh}
 For any CDG\+module $\M\in\B\qcoh_\lp$, locally projective graded
$\B$\+module $\E$, and a homogeneous surjective morphism of
locally projective graded $\B$\+modules $\E\rarrow\M$, there exist
a CDG\+module $\Q\in\B\qcoh_\lp$, a surjective closed morphism
of CDG\+modules $\Q\rarrow\M$, and a homogeneous surjective morphism
of locally projective graded $\B$\+modules $\Q\rarrow\E$, such that
the triangle $\Q\rarrow\E\rarrow\M$ commutes.
\end{lemh}

\begin{proof}
 For any open subscheme $U\subset X$, one can simply define
$\Q^i(U)$ as the abelian group of all pairs $(e'\in\E^{i+1}(U)\;
e\in\E^i(U))$ such that $df(e)=f(e')$, where $f$~denotes
the morphism of graded $\B$\+modules $\E\rarrow\M$ and $d$~is
the differential in~$\M$.
 The action of $\B$ in $\Q$ is defined by the formula $b(e',e) =
((-1)^{|b|}be'+d(b)e\;be)$; the differential in $\Q$ is given
by the obvious rule $d(e',e)=(he,e')$.
 The morphism $\Q\rarrow\E$ is defined as $(e',e)
\longmapsto e$; the morphism $\Q\rarrow\M$, given by
$(e',e)\longmapsto f(e)$, obviously commutes with
the differentials.

 It remains to check that the graded $\B$\+module $\Q$ is locally
projective.
 This can be done by comparing the above construction with
the constructions of the freely (co)generated CDG\+modules
$G^+(\E)$ and $G^-(\E)$ from~\cite[proof of Theorem~3.6]{Pkoszul}
(see the beginning of the proof of Theorem~\ref{flat-dimension-thm}).
 One can simply define $G^-(\E)$ as being isomorphic to $G^+(\E)[1]$.
 Since $\M$ is a CDG\+module, there is a natural closed morphism
of CDG\+modules $\M\rarrow G^-(\M)$.
 The CDG\+module $\Q$ is the fibered product of the surjective
closed morphism of CDG\+modules $G^-(\E)\rarrow G^-(\M)$ and
the closed morphism $\M\rarrow G^-(\M)$; hence the graded
$\B$\+module $\Q$ is locally projective.
 The morphism $\Q\rarrow\E$ is induced by the natural morphism of
graded $\B$\+modules $G^-(\E)\rarrow\E$.
 It forms a commutative diagram with the morphism $\E\rarrow\M$,
since the composition $\M\rarrow G^-(\M)\rarrow\M$ is the identity
morphism.
\end{proof}

 The exact sequence of CDG\+modules $0\rarrow\R'_{n-1}\rarrow
\Q'_{n-1}\rarrow\dotsb\rarrow\Q'_0\rarrow\P\rarrow 0$ represents
a certain Yoneda Ext class of degree~$n$ between the locally
projective graded $\B$\+modules $\P$ and $\R'_{n-1}$.
 Since the homological dimension of the exact category of such
$\B$\+modules is equal to~$d$ and we assume that $n>d$, this
Ext class has to vanish.
 This means that there exists an exact sequence of locally
projective graded $\B$\+modules $0\rarrow\R'_{n-1}\rarrow\E_{n-1}
\rarrow\dotsb\rarrow\E_0\rarrow\P\rarrow0$ mapping to our original
exact sequence, with the maps on the rightmost and leftmost terms
being the identity maps, such that the embedding of $\B$\+modules
$\R'_{n-1}\rarrow\E_{n-1}$ splits.

 As explained in~\cite[proof of Lemma~4.4]{Partin}, one can
assume the morphisms $\E_i\rarrow\Q'_i$ to be surjective.
 Applying Lemma~H, we obtain a surjective closed morphism of
CDG\+modules $\Q_0\rarrow\Q'_0$ and a morphism of graded
$\B$\+modules $\Q_0\rarrow\E_0$ forming a commutative triangle
with the morphism $\E_0\rarrow\Q'_0$.
 Applying Lemma~H to the surjective morphism of fibered products
$\Q_0\times_{\E_0}\E_1\rarrow\Q_0\times_{\Q'_0}\Q'_1$, we obtain
a surjective closed morphism $\Q_1\rarrow\Q'_1$ and a closed
morphism $\Q_1\rarrow\Q_0$ forming a commutative square with
the closed morphisms $\Q_0\rarrow\Q'_0$ and $\Q'_1\rarrow\Q'_0$.
 Besides, the sequence $\Q_1\rarrow\Q_0\rarrow\P$ is exact
at~$\Q_0$.
 We also obtain a morphism of graded $\B$\+modules $\Q_1\rarrow\E_1$
forming a commutative triangle with the morphisms to $\Q'_1$ and
a commutative square with the morphisms to~$\E_0$.

 Proceeding in this way, we construct a sequence $\Q_{n-2}\rarrow
\dotsb\rarrow\Q_0\rarrow\P\rarrow0$, which is exact at all
the middle terms, maps onto the sequence $\Q'_{n-2}\rarrow
\dotsb\rarrow\Q'_0\rarrow\P$ by closed morphisms, and maps into
the sequence $\E_{n-2}\rarrow\dotsb\rarrow\E_0\rarrow\P$ so that
the triangle of the maps of sequences commutes.
 Finally, notice that $\E_{n-1}\simeq\E_{n-2}\times_{\Q'_{n-2}}
\Q'_{n-1}$, and set $\Q_{n-1}=\Q_{n-2}\times_{\Q'_{n-2}}\Q'_{n-1}$.
 Then the exact sequence of CDG\+modules $0\rarrow\R'_{n-1}\rarrow
\Q_{n-1}\rarrow \dotsb\rarrow\Q_0\rarrow\P\rarrow0$ maps onto
the exact sequence $0\rarrow\R'_{n-1}\rarrow\Q'_{n-1}\rarrow\dotsb
\rarrow\Q'_0 \rarrow\P\rarrow0$ by closed morphisms, and this map
of exact sequences factorizes through the exact sequence of graded
$\B$\+modules $0\rarrow\R'_{n-1}\rarrow\E_{n-1}\rarrow\dotsb\rarrow
\E_0\rarrow\P\rarrow0$.
 The composition of the morphism $\Q_{n-1}\rarrow\E_{n-1}$ with
the splitting $\E_{n-1}\rarrow\R'_{n-1}$ of the embedding
$\R'_{n-1}\rarrow\E_{n-1}$ provides a graded $\B$\+module splitting
$\Q_{n-1}\rarrow\R'_{n-1}$ of the embedding of CDG\+modules
$\R'_{n-1}\rarrow\Q_{n-1}$.

 Denote by $\R_{n-2}$ the image of the morphism of CDG\+modules
$\Q_{n-1}\rarrow\Q_{n-2}$.
 The morphism from the total CDG\+module of the complex
$\R'_{n-1}\rarrow\Q'_{n-1}\rarrow\dotsb\rarrow\Q'_0$ to
the CDG\+module $\L$ is homotopic to zero, hence so is
the morphism to $\L$ from the total CDG\+module of the complex
$\R'_{n-1}\rarrow\Q_{n-1}\rarrow\dotsb\rarrow\Q_0$.
 The latter morphism factorizes naturally through the total
CDG\+module of the complex $\R_{n-2}\rarrow\Q_{n-2}\rarrow\dotsb
\rarrow\Q_0$.
 The cone of this closed morphism between two total CDG\+modules
is homotopy equivalent to the total CDG\+module of the exact
triple $\R'_{n-1}\rarrow\Q_{n-1}\rarrow\R_{n-2}$.
 Since this exact triple splits as an exact triple of graded
$\B$\+modules, its total CDG\+module is contractible.
 Consequently, the morphism between the total CDG\+modules of
$\R'_{n-1}\rarrow\Q_{n-1}\rarrow\dotsb\rarrow\Q_0$ and
$\R_{n-2}\rarrow\Q_{n-2}\rarrow\dotsb\rarrow\Q_0$ is a homotopy
equivalence.

 It follows that the natural morphism from the total CDG\+module of
the resolution $\R_{n-2}\rarrow\Q_{n-2}\rarrow\dotsb\rarrow\Q_0$
of the CDG\+module $\P$ to the CDG\+module $\L$ is homotopic
to zero, and we are done.
\end{proof}

 So far we have only considered flat coherent CDG\+modules over
quasi-coherent CDG\+algebras $\B$ whose underlying quasi-coherent
graded algebras are Noetherian.
 But the latter restriction is not necessary, as flat and locally
finitely presented (or, which is equivalent, locally projective
and finitely generated) quasi-coherent graded $\B$\+modules always
form an exact subcategory of flat (or locally projective)
graded $\B$\+modules.
 The notation $\B\coh_\lp$ (understood in the obvious sense as
the DG\+category of CDG\+modules over $\B$ with coherent and
locally projective underlying graded $\B$\+modules) is synonymous
to $\B\coh_\fl$ (see Remark~\ref{flat-dimension-thm}).

\begin{cor}
 The functor\/ $\sD^\abs(\B\coh_\lp)\rarrow\sD^\co(\B\qcoh_\lp)$
induced by the embedding of DG\+categories $\B\coh_\lp\rarrow
\B\qcoh_\lp$ is fully faithful.
\end{cor}

\begin{proof}
 When $\B$ is Noetherian, one can show that the functor
$\sD^\abs(\B\coh_\lp)\rarrow \sD^\abs(\B\qcoh_\lp)$ is fully faithful
by comparing parts~(a\+c) of Proposition~\ref{embedding-prop} (with
the flatness condition replaced by the local projectivity).
 In the general case, one proves this assertion directly, using
an argument similar to the proof of
Proposition~\ref{embedding-prop}(a\+b).
 Then it remains to use the above Theorem.
\end{proof}

 When every flat quasi-coherent graded module over $\B$ has
finite locally projective dimension
(see Remark~\ref{flat-dimension-thm}), one has
$\sD^\co(\B\qcoh_\lp)\simeq\sD^\co(\B\qcoh_\fl)\simeq
\sD^\co(\B\qcoh_\ffd)$ and $\sD^\abs(\B\qcoh_\lp)
\simeq\sD^\abs(\B\qcoh_\fl)\simeq\sD^\abs(\B\qcoh_\ffd)$
by appropriate versions of Theorem~\ref{flat-dimension-thm}.
 Consequently, it follows from Theorem above that
$\sD^\abs(\B\qcoh_\fl) = \sD^\co(\B\qcoh_\fl)$ and
$\sD^\abs(\B\qcoh_\ffd) = \sD^\co(\B\qcoh_\ffd)$ in this case.
 Thus the functor $\sD^\abs(\B\coh_\fl)\rarrow\sD^\co(\B\qcoh_\fl)$
is fully faithful; when $\B$ is Noetherian, so is the functor
$\sD^\abs(\B\coh_\ffd)\rarrow\sD^\co(\B\qcoh_\ffd)$.

\subsection{Gorenstein case}  \label{gorenstein-case}
 Here we establish a sufficient condition for the functor
$\sD^\co(\B\qcoh_\fl)\rarrow\sD^\co(\B\qcoh)$ to be
an equivalence of triangulated categories.

 Let $\B\qcoh_\inj$ denote the full DG\+subcategory in $\B\qcoh$
consisting of the CDG\+mod\-ules whose underlying quasi-coherent
graded $\B$\+modules are injective.
 Furthermore, let $\B\qcoh_\fid$ be the full DG\+subcategory in
$\B\qcoh$ consisting of the CDG\+modules whose underlying
quasi-coherent graded $\B$\+modules have finite injective
dimension (i.~e., admit a finite right resolution by
injective quasi-coherent graded $\B$\+modules).
 Let $\sD^\abs(\B\qcoh_\fid)$ and $\sD^\co(\B\qcoh_\fid)$ denote
the corresponding derived categories of the second kind.
 (The difficulty in the definition of the latter category, similar
to the difficulty in the definition of $\sD^\co(\B\qcoh_\ffd)$
discussed in Section~\ref{second-kind}, does not actually arise,
as it is clear from part~(a) of the next lemma.) 
{\hbadness=1100\par}

\begin{lem}
\textup{(a)} \hbadness=2100
 For any quasi-coherent CDG\+algebra $\B$ over $X$,
the natural functors $H^0(\B\qcoh_\inj)\rarrow
\sD^\abs(\B\qcoh_\fid)\rarrow\sD^\co(\B\qcoh_\fid)$
are equivalences of triangulated categories. \par
\textup{(b)}
 Let $\B$ be a quasi-coherent CDG\+algebra over $X$ whose
underlying quasi-coherent graded algebra $\B$ is Noetherian.
 Then the functor $H^0(\B\qcoh_\inj)\rarrow\sD^\co(\B\qcoh)$
induced by the embedding $\B\qcoh_\inj\rarrow\B\qcoh$ is
an equivalence of triangulated categories.
\end{lem}

\begin{proof}
 Part~(a) is provided by~\cite[Theorem and
Remark in Section~3.6]{Pkoszul}.
 Part~(b) is a particular case of~\cite[Theorem and Remark
in Section~3.7]{Pkoszul}, since the class of injective
quasi-coherent graded $\B$\+modules is closed under infinite
direct sums in its assumptions. 
 (Cf.~\cite[Proposition~2.4]{PL}.)
\end{proof}

\begin{prop}
 Let $\B$ be a quasi-coherent CDG\+algebra over $X$ such that
the quasi-coherent graded algebra $\B$ is Noetherian and
the classes of quasi-coherent graded $\B$\+modules of finite
flat dimension and of finite injective dimension coincide.
 Then the functors\/ $\sD^\abs(\B\qcoh_\fl)\rarrow
\sD^\co(\B\qcoh_\fl)\rarrow\sD^\co(\B\qcoh)$ induced by
the embedding\/ $\B\qcoh_\fl\rarrow\B\qcoh$ are equivalences
of triangulated categories.
\end{prop}

\begin{proof}
 Since $\B\qcoh_\ffd=\B\qcoh_\fid$, the isomorphism of categories
$\sD^\abs(\B\qcoh_\ffd)=\sD^\co(\B\qcoh_\ffd)$ follows from
part~(a) of Lemma.
 Applying Theorem~\ref{flat-dimension-thm}, we obtain
the isomorphism of categories $\sD^\abs(\B\qcoh_\fl)\rarrow
\sD^\co(\B\qcoh_\fl)$.
 Similarly, it suffices to compare parts~(a) and (b) of Lemma
in order to conclude that the functor $\sD^\co(\B\qcoh_\fid)
\rarrow\sD^\co(\B\qcoh)$ is an equivalence of categories,
hence so are the functors $\sD^\co(\B\qcoh_\fl)\rarrow
\sD^\co(\B\qcoh_\ffd)\rarrow\sD^\co(\B\qcoh)$.
 (Cf.~\cite[Section~3.9]{Pkoszul}.)
\end{proof}

\subsection{Pull-backs and push-forwards}  \label{pull-push}
 Let $f\:Y\rarrow X$ be a morphism of separated Noetherian schemes,
$\B_X$ be a quasi-coherent CDG\+algebra over $X$, and $\B_Y$
a quasi-coherent CDG\+algebra over~$Y$.
 A \emph{morphism of quasi-coherent CDG\+algebras} $\B_X\rarrow
\B_Y$ \emph{compatible with the morphism $Y\rarrow X$} is
the data of a CDG\+ring morphism $\B_X(U)\rarrow\B_Y(V)$ for
each pair of affine open subschemes $U\subset X$ and $V\subset Y$
such that $f(V)\subset U$.
 This data should satisfy the obvious compatibility condition:
for any affine open subschemes $U'\subset U$ and $V'\subset V$
such that $f(V')\subset U'$, the square diagram of CDG\+ring
morphisms between the CDG\+rings $\B_X(U)$, \ $\B_X(U')$, \
$\B_Y(V)$, and $\B_Y(V')$ must be commutative.

 Let $\B_X\rarrow\B_Y$ be a morphism of quasi-coherent CDG\+algebras
compatible with a morphism of schemes $Y\rarrow X$.
 Then for any quasi-coherent left CDG\+module $\M$ over $\B_X$
the quasi-coherent graded left module
$f^*\M=\B_Y\ot_{f^{-1}\B_X}f^{-1}\M$ over $\B_Y$ has a natural
structure of quasi-coherent CDG\+module over~$\B_Y$.
 Similarly, for any quasi-coherent left CDG\+module $\N$ over $\B_Y$
the quasi-coherent graded left module 
$f_*\N$ over $\B_X$ has a natural structure of quasi-coherent
CDG\+module over~$\B_X$.
 These CDG\+module structures are defined in terms of the CDG\+ring
morphisms $\B_X(U)\rarrow\B_Y(V)$.
 The above constructions provide the underived direct and inverse
image functors, which can be viewed as triangulated functors
$f^*\:H^0(\B_X\qcoh)\rarrow H^0(\B_Y\qcoh)$ and $f_*\:
H^0(\B_Y\qcoh)\rarrow H^0(\B_X\qcoh)$.
 The functor~$f_*$ is right adjoint to the functor~$f^*$.

 The derived inverse image functor $\boL f^*$ is in general only
defined on CDG\+modules satisfying certain finite flat dimension
conditions.
 Restricting the functor~$f^*$ to flat CDG\+modules, we obtain
a triangulated functor $H^0(\B_X\qcoh_\fl)\rarrow
H^0(\B_Y\qcoh_\fl)$, which takes objects coacyclic with respect
to $\B_X\qcoh_\fl$ to objects coacyclic with respect to
$\B_Y\qcoh_\fl$, since the inverse image preserves infinite direct
sums and short exact sequences of flat quasi-coherent graded modules.
 Hence there is the induced triangulated functor
$\sD^\co(\B_X\qcoh_\fl)\rarrow\sD^\co(\B_Y\qcoh_\fl)$. 
 Applying Theorem~\ref{flat-dimension-thm}(a), we construct
the derived inverse image functor
$$
 \boL f^*\:\sD^\co(\B_X\qcoh_\ffd)\lrarrow\sD^\co(\B_Y\qcoh_\ffd).
$$

 Assuming that there are enough vector bundles on $X$ and $Y$, and
restricting the functor~$f^*$ to flat coherent CDG\+modules, we
obtain a triangulated functor $H^0(\B_X\coh_\fl)\rarrow
H^0(\B_Y\coh_\fl)$, which induces a triangulated functor
$\sD^\abs(\B_X\coh_\fl)\rarrow\sD^\abs(\B_Y\coh_\fl)$.
 Assuming additionally that the quasi-coherent graded algebras $\B_X$
and $\B_Y$ are Noetherian and applying
Theorem~\ref{flat-dimension-thm}(c), we construct the derived
inverse image functor
$$
 \boL f^*\:\sD^\abs(\B_X\coh_\ffd)\lrarrow\sD^\abs(\B_Y\coh_\ffd).
$$

 When $f$~is an affine morphism, the direct image of quasi-coherent
sheaves is an exact functor (preserving also infinite direct sums),
so the functor~$f_*\:H^0(\B_Y\qcoh)\rarrow H^0(\B_X\qcoh)$
induces a triangulated functor $\sD^\co(\B_Y\qcoh)\rarrow\sD^\co
(\B_X\qcoh)$.
 To construct the derived direct image functor between
the coderived categories in the general case, we need to use
injective resolutions.

 From now on we assume that $\B_X$ and $\B_Y$ are Noetherian;
so Lemma~\ref{gorenstein-case}(b) is applicable to~$\B_Y$.
 Restricting the functor $f_*$ to the full subcategory
$H^0(\B_Y\qcoh_\inj)\subset H^0(\B_Y\qcoh)$ and composing it with
the localization functor $H^0(\B_X\qcoh)\rarrow\sD^\co(\B_X\qcoh)$,
we obtain the derived direct image functor
$$
 \boR f_*\:\sD^\co(\B_Y\qcoh)\lrarrow\sD^\co(\B_X\qcoh).
$$

\begin{prop}
 Assume that there are enough vector bundles on $X$ and~$Y$.
 Then the functors\/ $\boL f^*\:\sD^\abs(\B_X\coh_\ffd)\rarrow
\sD^\abs(\B_Y\coh_\ffd)$ and\/ $\boR f_*\:\sD^\co\allowbreak
(\B_Y\qcoh)\rarrow\sD^\co(\B_X\qcoh)$ are ``partially adjoint''
to each other in the following sense: for any objects
$\M\in\sD^\abs(\B_X\coh_\ffd)$ and $\N\in\sD^\co(\B_Y\qcoh)$
there is a natural isomorphism of abelian groups
$$
 \Hom_{\sD^\co(\B_X\qcoh)}(\iota_X\M\;\boR f_*\N)\simeq
 \Hom_{\sD^\co(\B_Y\qcoh)}(\iota_Y\boL f^*\M\;\N),
$$
where\/ $\iota_X\:\sD^\abs(\B_X\coh_\ffd)\rarrow\sD^\co(\B_X\qcoh)$
and\/ $\iota_Y\:\sD^\abs(\B_X\coh_\ffd)\rarrow
\sD^\co\allowbreak(\B_Y\qcoh)$ are the natural fully faithful
triangulated functors.
\end{prop}

\begin{proof}
 The functors $\iota_X$ and $\iota_Y$ are fully faithful by
Theorem~\ref{flat-dimension-thm}(c) and
Proposition~\ref{embedding-prop}(b,\,d).
 Using Theorem~\ref{flat-dimension-thm}(c), let us assume that
$\M\in\sD^\abs(\B_X\coh_\fl)$.
 We can also assume that $\N\in H^0(\B_Y\qcoh_\inj)$.

 Then the left-hand side is the (filtered) inductive limit of
$\Hom_{H^0(\B_X\qcoh)}(\M'',f_*\N)$ over all morphisms $\M''\rarrow\M$
in $H^0(\B_X\qcoh)$ with a cone coacyclic with respect to
$\B_X\qcoh$.
 According to the proofs of Proposition~\ref{embedding-prop}(b)
and~\cite[Theorem~3.11.1]{Pkoszul}, any morphism from $\M$
to an object coacyclic with respect to $\B_X\qcoh$ factorizes
through an object absolutely acyclic with respect to $\B_X\coh_\fl$.
 Thus the above inductive limit coincides with the similar limit
taken over all morphisms $\M'\rarrow \M$ in $H^0(\B_X\coh_\fl)$
with a cone absolutely acyclic with respect to $\B_X\coh_\fl$.

 By~\cite[Theorem~3.5(a), Remark~3.5, and Lemma~1.3]{Pkoszul},
the right-hand side is isomorphic to $\Hom_{H^0(\B_Y\qcoh)}(f^*\M,\N)$
and to $\Hom_{H^0(\B_Y\qcoh)}(f^*\M',\N)$, since the objects of
$H^0(\B_Y\qcoh_\inj)$ are right orthogonal to any coacyclic objects
in $H^0(\B_Y\qcoh)$.
 So the assertion follows from the adjointness of the functors
$f^*$ and~$f_*$ on the level of the homotopy categories of
quasi-coherent CDG\+modules.
\end{proof}

\begin{rem}
 It is not immediately obvious from the above construction that
the derived functor $\boR f_*$ is compatible with the compositions,
i.~e., for $g\:Z\rarrow Y$ and $f\:Y\rarrow X$ one has
$\boR (fg)_*\simeq \boR f_*\circ\boR g_*$.
 The problem is that the direct image functor $f_*$ does not
preserve injectivity of quasi-coherent graded modules in general.
 When the derived direct image functors are adjoint to appropriately
defined derived inverse images (see Section~\ref{finite-dim-morphisms}
below for some results of this kind), the problem reduces to checking
that the derived inverse images are compatible with the compositions,
which may be easier to see from our definitions.

 One general approach to this problem is to replace injective
quasi-coherent graded $\B$\+modules with quasi-coherent graded
$\B$\+modules that are flabby as sheaves of graded abelian groups
in our construction of the derived direct images.
 The class of flabby sheaves of abelian groups is closed under
infinite direct sums, since the underlying topological space of
the scheme is Noetherian; it is also always closed under extensions
and cokernels of injective morphisms.
 Whenever the quasi-coherent graded algebra $\B$ is Noetherian,
all injective quasi-coherent graded $\B$\+modules are
flabby by Theorem~\ref{injective-sheaves}.
 Therefore, the coderived category of flabby quasi-coherent
CDG\+modules over $\B$ is equivalent to the homotopy category
$H^0(\B\qcoh_\inj)$ by a version of Lemma~\ref{gorenstein-case}(b),
hence it is also equivalent to the coderived category of all
quasi-coherent CDG\+modules $\sD^\co(\B\qcoh)$ (cf.\ the proof of
Proposition~\ref{gorenstein-case}).

 The direct images preserve exact triples of flabby sheaves, so
derived direct images can be defined using flabby resolutions.
 The direct images also take flabby sheaves to flabby sheaves, hence
the desired compatibility of their derived functors with
the compositions of scheme morphisms follows.

 Moreover, assuming additionally that the scheme has finite Krull
dimension, the absolute derived category of flabby quasi-coherent
CDG\+modules is equivalent to $\sD^\abs(\B\qcoh)$ by a dual version
of Theorem~\ref{flat-dimension-thm}(b), as the ``flabby dimension''
of any quasi-coherent graded $\B$\+module is finite.
 This allows to define the derived direct images on the absolute
derived categories of quasi-coherent CDG\+modules (another approach
to this question is to use the construction from the proof of
Proposition~\ref{finite-dim-morphisms} below).
 Notice that all our constructions of derived inverse images are
also applicable to the categories $\sD^\abs(\B\qcoh)$.
\end{rem}

 Finally, let us point out that for any morphism of quasi-coherent
CDG\+algebras $\B_X\rarrow\B_Y$ with Noetherian underlying
quasi-coherent graded algebras $\B_X$ and $\B_Y$ compatible with
a morphism of separated Noetherian schemes $f\:Y\rarrow X$
the functor $\boR f_*$ has a right adjoint functor
$$
 f^!\:\sD^\co(\B_X\qcoh)\lrarrow\sD^\co(\B_Y\qcoh).
$$
 Indeed, the triangulated category $\sD^\co(\B_Y\qcoh)$ is compactly
generated by Proposition~\ref{embedding-prop}(d), and the functor
$\boR f_*$ preserves infinite direct sums, since the class of
injective quasi-coherent graded $\B_Y$\+modules is closed under
infinite direct sums, due to Noetherianness of~$\B_Y$.
 So it remains to apply~\cite[Theorem~4.1]{Neem}.

 There is a special situation when one can construct the above
functor~$f^!$ explicitly.
 Assume that $f\:Y\rarrow X$ is an affine morphism.
 Let us say that the quasi-coherent graded algebra $\B_Y$ is
\emph{finite} over $\B_X$ if for any affine open subscheme
$U\subset X$ the graded $\B_X(U)$\+module $\B_Y(f^{-1}(U))$ is
finitely generated, or in other words, if the quasi-coherent
graded $\B_X$\+module $f_*\B_Y$ is coherent.

 Let $\B_X\rarrow\B_Y$ be a morphism of Noetherian quasi-coherent
CDG\+algebras compatible with an affine morphism of separated
Noetherian schemes $f\:Y\rarrow X$ such that the quasi-coherent
graded algebra $\B_Y$ is finite over~$\B_X$.
 Given a quasi-coherent graded left module $\M$ over $\B_X$, we
set $(f^!\M)(f^{-1}(U))$ to be the graded left module
of homogeneous morphisms (of various degrees)
$\Hom_{\B_X(U)}(\B_Y(f^{-1}(U)),\M(U))$ over the graded ring
$\B_Y(f^{-1}(U))$ for any affine open subscheme $U\subset X$.
 Due to the finiteness condition on $\B_Y$ over $\B_X$, for any
affine open subscheme $V\subset U$ there are natural isomorphisms
$(f^!\M)(f^{-1}(V))\simeq\O_X(V)\ot_{\O_X(U)}(f^!\M)(f^{-1}(U))
\simeq\O_Y(f^{-1}(V))\ot_{\O_Y(f^{-1}(U))}(f^!\M)(f^{-1}(U))$,
which allow to extend the assignment $f^{-1}(U)\longmapsto
(f^!\M)(f^{-1}(U))$ to a quasi-coherent graded module $f^!(\M)$
over the quasi-coherent graded algebra~$\B_Y$.

 Given a quasi-coherent CDG\+module $\M$ over $\B_X$,
the conventional rule $d(g)(m)=d(g(m))-(-1)^{|g|}g(d(m))$
(with the usual change-of-connection modifications) defines
the structure of a quasi-coherent CDG\+module over $\B_Y$
on the quasi-coherent graded module $f^!(\M)$.
 This construction provides a triangulated functor $f^!\:H^0(\B_X\qcoh)
\allowbreak\rarrow H^0(\B_Y\qcoh)$ right adjoint to the triangulated
functor $f_*\:H^0(\B_Y\qcoh)\rarrow H^0(\B_X\qcoh)$.
 Restricting the functor $f^!\:H^0(\B_X\qcoh)\rarrow H^0(\B_Y\qcoh)$
to the full subcategory of injective quasi-coherent CDG\+modules in
$H^0(\B_X\qcoh)$ and taking into account
Lemma~\ref{gorenstein-case}(b), we obtain the right derived functor
$$
 \boR f^!\:\sD^\co(\B_X\qcoh)\lrarrow\sD^\co(\B_Y\qcoh),
$$
which is right adjoint to the (underived, as the morphism~$f$ is
affine) direct image functor $f_*\:\sD^\co(\B_Y\qcoh)\rarrow
\sD^\co(\B_X\qcoh)$.
 In other words, the functor $\boR f^!$ coincides with the above
adjoint functor $f^!\:\sD^\co(\B_X\qcoh)\rarrow\sD^\co(\B_Y\qcoh)$
in our special case.

\subsection{Morphisms of finite flat dimension} 
\label{finite-dim-morphisms}
 Let $f\:Y\rarrow X$ be a morphism of separated Noetherian schemes,
and let $\B_X\rarrow\B_Y$ be a compatible morphism of quasi-coherent
CDG\+algebras.
 We will say that the quasi-coherent graded algebra $\B_Y$
\emph{has finite flat dimension over\/ $\B_X$} if (the left derived
functor of) the functor of inverse image~$f^*$ acting between
the abelian categories of quasi-coherent graded modules over $\B_X$
and $\B_Y$ has finite homological dimension.
 Equivalently, for any affine open subschemes $U\subset X$ and
$V\subset Y$ such that $f(V)\subset U$ the graded right
$\B_X(U)$\+module $\B_Y(V)$ should have finite flat dimension.

 A quasi-coherent graded $\B_X$\+module is said to be
\emph{adjusted to\/~$f^*$} if its derived inverse image under~$f$,
as an object of the derived category of the abelian category of
quasi-coherent graded $\B_Y$\+modules, coincides with
the underived inverse image.
 Denote the DG\+category of quasi-coherent CDG\+modules over
$\B_X$ whose underlying graded $\B_X$\+modules are adjusted
to~$f^*$ by $\B_X\qcoh_\fadj$.
 When $\B_X$ is Noetherian, let $\B_X\coh_\fadj$ denote the similarly
defined DG\+category of coherent CDG\+modules.
 We will use our usual notation for the absolute derived and
coderived categories of these DG\+categories of CDG\+modules.

\begin{lem}
 Assume that the quasi-coherent graded algebra $\B_Y$ has finite
flat dimension over~$\B_X$.  Then \par
\textup{(a)} the functor\/ $\sD^\co(\B_X\qcoh_\fadj)\rarrow
\sD^\co(\B_X\qcoh)$ induced by the embedding of DG\+categories
$\B_X\qcoh_\fadj\rarrow\B_X\qcoh$ is an equivalence of
triangulated categories; \par
\textup{(b)} the functor\/ $\sD^\abs(\B_X\qcoh_\fadj)\rarrow
\sD^\abs(\B_X\qcoh)$ induced by the embedding of DG\+categories
$\B_X\qcoh_\fadj\rarrow\B_X\qcoh$ is an equivalence of
triangulated categories; \par
\textup{(c)} if there are enough vector bundles on $X$ and $\B_X$
is Noetherian, the functor\/ $\sD^\abs(\B_X\coh_\fadj)\rarrow
\sD^\abs(\B_X\coh)$ induced by the embedding of DG\+categories
$\B_X\coh_\fadj\rarrow\B_X\coh$ is an equivalence of
triangulated categories.
\end{lem}

\begin{proof}
 This is a version of Theorem~\ref{flat-dimension-thm}, provable
in the same way (cf.\ Corollary~\ref{w-flat-cor} below).
 The assertions hold, because any quasi-coherent graded
$\B_X$\+module has a finite left resolution consisting of
quasi-coherent CDG\+modules adjusted to~$f^*$, and similarly
for coherent CDG\+modules.
\end{proof}

 The functor of inverse image $f^*\:H^0(\B_X\qcoh)\rarrow
H^0(\B_Y\qcoh)$ takes CDG\+mod\-ules coacyclic with respect to
$\B_X\qcoh_\fadj$ to CDG\+modules coacyclic with respect to
$\B_Y\qcoh$, and hence induces a triangulated functor
$\sD^\co(\B_X\qcoh_\fadj)\rarrow\sD^\co(\B_Y\qcoh)$.
 Taking Lemma into account, we construct the derived inverse
image functor
$$
 \boL f^*\:\sD^\co(\B_X\qcoh)\lrarrow\sD^\co(\B_Y\qcoh).
$$
 One shows that this functor is left adjoint to the functor
$\boR f_*$ constructed in~\ref{pull-push} in the way analogous
to (but simpler than) the proof of Proposition~\ref{pull-push}.

 When there are enough vector bundles on $X$, and $\B_X$ and $\B_Y$
are Noetherian, we construct the derived inverse image functor
$$
 \boL f^*\:\sD^\abs(\B_X\coh)\lrarrow\sD^\abs(\B_Y\coh)
$$
in the similar way.

 Let $\B_X^\op$ and $\B_Y^\op$ denote the quasi-coherent graded
algebras with the opposite multiplication to $\B_X$ and~$\B_Y$.

\begin{prop}
 When $\B_Y^\op$ has finite flat dimension over $\B_X^\op$,
the derived inverse image functor\/ $\boL f^*\:
\sD^\co(\B_X\qcoh_\ffd)\rarrow\sD^\co(\B_Y\qcoh_\ffd)$ constructed
in\/~\textup{\ref{pull-push}} has a right adjoint functor
$$
 \boR f_*\:\sD^\co(\B_Y\qcoh_\ffd)\lrarrow\sD^\co(\B_X\qcoh_\ffd).
$$
\end{prop}

\begin{proof}
 Let $\{U_\alpha\}$ be a finite affine covering of~$Y$.
 To any object $\N\in\B_Y\qcoh_\ffd$, assign the total CDG\+module
$\boR_{\{U_\alpha\}} f_*\N$ of the finite \v Cech complex
$$\textstyle
 \bigoplus_\alpha f|_{U_\alpha}{}_*(\N|_{U_\alpha})\lrarrow
 \bigoplus_{\alpha<\beta} f|_{U_\alpha\cap\U_\beta}{}_*
 (\N|_{U_\alpha\cap U_\beta})\lrarrow\dotsb
$$
of CDG\+modules over~$\B_X$.

 The terms of this complex belong to $\B_X\qcoh_\ffd$, since
the morphism $f|_V\:V\rarrow X$ is affine for any intersection $V$
of a nonempty subset of affine open subschemes $U_\alpha\subset Y$
and the quasi-coherent graded algebra $\B_Y^\op$ has
finite flat dimension over $\B_X^\op$.
 Hence one has $\boR_{\{U_\alpha\}} f_*\N\in\B_X\qcoh_\ffd$; it is
clear that $\boR_{\{U_\alpha\}} f_*$ is a DG\+functor
$\B_Y\qcoh_\ffd\rarrow\B_X\qcoh_\ffd$ taking coacyclic objects
to coacyclic objects.
 So we have the induced functor $\boR f_*$ between the
coderived categories.

 It remains to obtain the adjunction isomorphism
$$
 \Hom_{\sD^\co(\B_X\qcoh_\ffd)}(\M,\boR f_*\N)\simeq
 \Hom_{\sD^\co(\B_Y\qcoh_\ffd)}(\boL f^*\M,\N)
$$
for $\M\in\sD^\co(\B_X\qcoh_\ffd)$.
 Denote by $\N_+$ the total CDG\+module of the finite complex
$$\textstyle
 C_{\{U_\alpha\}}^\bu\N \.=\.
 \big(\bigoplus_\alpha j_{U_\alpha *}j_{U_\alpha}^*
 \N\rarrow\bigoplus_{\alpha<\beta} j_{U_\alpha\cap\U_\beta *}
 j_{U_\alpha\cap U_\beta}^*\N\rarrow\dotsb)
$$
of CDG\+modules over $\B_Y$ (where $j_V\:V\rarrow Y$ denotes
the embedding of an affine open subscheme).
 Then we have $\boR_{\{U_\alpha\}} f_*\N\simeq f_*\N_+$.
 There is a natural closed morphism $\N\rarrow\N_+$ of CDG\+modules
over $\B_Y$ with the cone coacyclic (and even absolutely acyclic)
with respect to $\B_Y\qcoh_\ffd$.

 For any CDG\+module $\Q\in\B_Y\qcoh_\ffd$ such that
$f_*\Q\in\B_X\qcoh_\ffd$, there is a natural map
$$
 \psi\:\Hom_{\sD^\co(\B_X\qcoh_\ffd)}(\M,f_*\Q)\lrarrow
 \Hom_{\sD^\co(\B_Y\qcoh_\ffd)}(\boL f^*\M,\Q).
$$
 Indeed, by (the proof of) Theorem~\ref{flat-dimension-thm}(a), any
morphism $\M\rarrow f_*\Q$ in $\sD^\co(\B_X\qcoh_\ffd)$ can be represented
as a fraction formed by a morphism $\M'\rarrow\M$ in
$H^0(\B_X\qcoh_\ffd)$ with $\M'\in\B_X\qcoh_\fl$ and a cone coacyclic
with respect to $\B_X\qcoh_\ffd$, and a morphism $\M'\rarrow f_*\Q$
in $H^0(\B_X\qcoh_\ffd)$.
 To such a fraction, the map~$\psi$ assigns the related morphism
$\boL f^*\M=f^*\M'\rarrow\Q$. {\hbadness=1200\par}

 For a fixed $\M$, the map~$\psi$ is a morphism of cohomological
functors of the argument $\Q\in H^0(\B_Y\qcoh_\ffd)$ with
$f_*\Q\in H^0(\B_X\qcoh_\ffd)$.
 Thus in order to show that it is an isomorphism for $\Q=\N_+$,
it suffices to check that it is an isomorphism for $\Q=j_{V*}\P$
for every affine $V\subset Y$ and $\P\in\B_Y|_V\qcoh_\ffd$.
 This follows from the adjunction isomorphism
$$
 \Hom_{\sD^\co(\B_X\qcoh_\ffd)}(\M,f|_V{}_*\P)\simeq
 \Hom_{\sD^\co(\B_Y|_V\qcoh_\ffd)}(\boL f|_V^*\M,\P)
$$
and the similar isomorphism for the embedding~$j_V$, which hold
because the functors $f|_V{}_*$ and $j_{V*}$ are exact,
the morphisms $f|_V$ and $j_V$ being affine.
\end{proof}

\begin{rem}
 One can also use the above \v Cech complex approach in order to
construct a version of the derived functor $\boR f_*\:\sD^\co
(\B_Y\qcoh)\rarrow\sD^\co(\B_X\qcoh)$.
 One can check that this construction agrees with the injective
resolution construction from Section~\ref{pull-push}, using
the fact that the restrictions of injective quasi-coherent graded
$\B_Y$\+modules to open subschemes are injective
(Theorem~\ref{injective-sheaves}).
 Alternatively, in the assumption of finite flat dimension of $\B_Y$
over $\B_X$, one checks that both constructions provide functors
right adjoint to $\boL f^*$, hence they are isomorphic.

 This allows to conclude that the derived functors $\boR f_*$
acting on arbitrary quasi-coherent CDG\+modules and quasi-coherent
CDG\+modules of finite flat dimension form a commutative diagram
with the natural functors from the coderived categories of
the latter to the coderived categories of the former.
\end{rem}

\subsection{Supports of CDG-modules} \label{cdg-supports}
 Let $X$ be a Noetherian scheme.
 The \emph{set-theoretic support} of a quasi-coherent sheaf $\M$
on $X$ is the minimal closed subset $T\subset X$ such that
the restriction of $\M$ to the open subscheme $X\setminus T$
vanishes.
 Given a Noetherian quasi-coherent graded algebra $\B$ over $X$
and a quasi-coherent graded $\B$\+module $\M$, the set-theoretic
support $T=\Supp\M$ of $\M$ is defined similarly.
 It only depends on the underlying quasi-coherent $\O_X$\+module
of~$\M$.

 Let $\B$ be a quasi-coherent CDG\+algebra over $X$ whose
underlying quasi-coherent graded algebra $\B$ is Noetherian.
 Fix a closed subset $T\subset X$.
 Denote by $\B\qcoh_T$ the full DG\+subcategory in $\B\qcoh$
consisting of all the quasi-coherent CDG\+modules whose underlying
quasi-coherent graded $\B$\+modules have their set-theoretic
supports contained in~$T$.
 The DG\+category $\B\coh_T$ of coherent CDG\+modules with
the set-theoretic support in $T$ is defined similarly.

 Let $\sD^\co(\B\qcoh_T)$ and $\sD^\abs(\B\coh_T)$ denote
the coderived and the absolute derived category of these
DG\+categories of CDG\+modules.
 Finally, let $\B\qcoh_{T,\inj}$ denote the DG\+category of
quasi-coherent CDG\+modules over $\B$ whose underlying
quasi-coherent graded modules are injective objects of
the abelian category of quasi-coherent graded $\B$\+modules
with the set-theoretic support contained in~$T$.

\begin{prop}
\textup{(a)} The functor $H^0(\B\qcoh_{T,\inj})\rarrow\sD^\co
(\B\qcoh_T)$ induced by the embedding of DG\+categories
$\B\qcoh_{T,\inj}\rarrow\B\qcoh_T$ is an equivalence of
triangulated categories. \par
\textup{(b)} The functor\/ $\sD^\abs(\B\coh_T)\rarrow\sD^\co
(\B\qcoh_T)$ induced by the embedding of DG\+categories
$\B\coh_T\rarrow\B\qcoh_T$ is fully faithful and its image
is a set of compact generators of the target category. \par
\textup{(c)} The functor\/ $\sD^\co(\B\qcoh_T)\rarrow\sD^\co
(\B\qcoh)$ induced by the embedding of DG\+categories
$\B\qcoh_T\rarrow\B\qcoh$ is fully faithful. \par
\textup{(d)} The functor\/ $\sD^\abs(\B\coh_T)\rarrow\sD^\abs
(\B\coh)$ induced by the embedding of DG\+categories
$\B\coh_T\rarrow\B\coh$ is fully faithful.
\end{prop}

\begin{proof}
 Part~(a) is essentially a particular case of~\cite[Theorem
and Remark in Section~3.7]{Pkoszul}.
 It is only important here that there are enough injective
objects in the abelian category of quasi-coherent graded
$\B$\+modules supported set-theoretically in~$T$ and the class
of such injective objects is closed under infinite direct sums.
 This is so because the abelian category in question is
a locally Noetherian Grothendieck category (since $X$ and
$\B$ are Noetherian).
 Part~(b) can be proven in the same way as the results
of~\cite[Section~3.11]{Pkoszul}.
 Part~(d) follows from parts~(b\+c) and
Proposition~\ref{embedding-prop}(d).

 Finally, part~(c) follows from part~(a),
Lemma~\ref{gorenstein-case}(b), and the fact that any injective
object $\J$ in the category of quasi-coherent graded $\B$\+modules
supported set-theoretically in~$T$ is also an injective object in
the category of arbitrary quasi-coherent graded $\B$\+modules.
 The latter is essentially a reformulation of the Artin--Rees
lemma.

 Indeed, it suffices to check that for any coherent graded
$\B$\+module $\M$ and its coherent graded $\B$\+submodule $\N$,
any morphism of quasi-coherent graded $\B$\+modules $\phi\:
\N\rarrow\J$ can be extended to~$\M$.
 Let $Z$ be a closed subscheme structure on the closed
subset $T\subset X$.
 Then there is an integer $n\ge0$ such that the morphism~$\phi$
annihilates $\I_Z^n\N$ (where $\I_Z$ is the sheaf of ideals of
the closed subscheme~$Z$).
 By Lemma~\ref{injective-sheaves}, there exists $m\ge0$ such that
$\I_Z^m\M\cap\N\subset\I_Z^n\N$.
 Then there exists a morphism $\M/\I_Z^m\M\rarrow\J$ of quasi-coherent
graded $\B$\+modules supported set-theoretically in~$T$ which extends
the given morphism into $\J$ from the quasi-coherent graded
$\B$\+submodule $\N/(\I_Z^m\M\cap\N)\subset\M/\I_Z^m\M$.
\end{proof}

 Let $U\subset X$ denote the open subscheme $X\setminus T$.

\begin{thm}\hfuzz=6.4pt
\textup{(a)} The functor of restriction to the open subscheme\/
$\sD^\co(\B\qcoh)\rarrow\sD^\co(\B|_U\qcoh)$ is the Verdier
localization functor by the thick subcategory\/ $\sD^\co(\B\qcoh_T)
\allowbreak\subset\sD^\co(\B\qcoh)$.
 In particular, the kernel of the restriction functor coincides
with the subcategory\/ $\sD^\co(\B\qcoh_T)$. \par
\textup{(b)} The functor of restriction to the open subscheme\/
$\sD^\abs(\B\coh)\rarrow\sD^\abs(\B|_U\coh)$ is the Verdier
localization functor by the triangulated subcategory\/
$\sD^\abs(\B\coh_T)\subset\sD^\abs(\B\coh)$.
 In particular, the kernel of the restriction functor coincides
with the thick envelope of (i.~e., the minimal thick subcategory
containing) $\sD^\abs(\B\coh_T)$ in $\sD^\abs(\B\coh)$.
\end{thm}

\begin{proof}
 Let $j\:U\rarrow X$ denote the natural open embedding.
 To prove part~(a), consider the functor $\boR j_*\:
\sD^\co(\B|_U\qcoh)\rarrow\sD^\co(\B\qcoh)$ as constructed
in Section~\ref{pull-push}.
 The quasi-coherent graded algebra $\B|_U$ being flat over~$\B$,
the functor $\boR j_*$ is right adjoint to the restriction
functor~$j^*\:\sD^\co(\B\qcoh)\rarrow\sD^\co(\B|_U\qcoh)$.
 Obviously, the composition $j^*\boR j_*$ is the identity
functor.
 It follows that the functor $j^*$~is a Verdier localization functor
by its kernel, which is the full subcategory consisting of
all the cones of the adjunction morphisms
$\M\rarrow \boR j_*j^*\M$, where $\M\in\sD^\co(\B\qcoh)$.

 Represent the object $\M$ by a CDG\+module with an injective
underlying quasi-coherent graded $\B$\+module.
 By Theorem~\ref{injective-sheaves}, the quasi-coherent graded
$\B|_U$\+module $j^*\M$ is then also injective, so we have
$\boR j_*j^*\M = j_*j^*\M$.
 Obviously, both the kernel and the cokernel of the closed
morphism of CDG\+modules $\M\rarrow j_*j^*\M$ belong to
$\B\qcoh_T$, and it follows, in view of part~(c) of Proposition,
that the cone also belongs to $\sD^\co(\B\qcoh_T)$.

 To prove part~(b), notice first that any coherent
CDG\+module over $\B|_U$ can be extended to a coherent
CDG\+module over~$\B$ (because a coherent sheaf $\K$ on $U$
can be extended to a coherent subsheaf of $j_*\K$), so
the restriction functor is essentially surjective.
 Taking this observation into account, part~(b) follows from
part~(a), part~(b) of the above Proposition,
Proposition~\ref{embedding-prop}(d), and the standard results
about localization of compactly generated triangulated
categories~\cite[Lemma~2.5 to Theorem~2.1]{Neem0}.
\end{proof}

 Define the \emph{category-theoretic support} $\supp\M$ of
a quasi-coherent CDG\+module $\M$ over $\B$ as the minimal closed
subset $T\subset X$ such that the restriction $\M|_U$ of $\M$ to
the open subscheme $U=X\setminus T$ is a coacyclic CDG\+module
over $\B|_U$.
 In other words, $X\setminus\supp\M$ is the union of all open
subschemes $V\subset X$ such that $\M|_V$ is a coacyclic
CDG\+module over $\B|_V$ (see Remark~\ref{second-kind}).
 Obviously, one has $\supp\M\subset\Supp\M$.

 The category-theoretic support of a coherent CDG\+module
$\M$ over $\B$ can be equivalently defined as the minimal
closed subset $T\subset X$ such that the restriction $\M|_U$
of $\M$ to the open subscheme $U=X\setminus T$ is absolutely
acyclic.
 Indeed, any CDG\+module from $\B|_U\coh$ that is coacyclic with
respect to $\B|_U\qcoh$ is also absolutely acyclic with respect
to $\B|_U\coh$ by Proposition~\ref{embedding-prop}(d).

\begin{cor}
\textup{(a)} For any quasi-coherent CDG\+module $\M$ over $\B$
with the category-theoretic support\/ $\supp\M$ contained in~$T$,
there exists a quasi-coherent CDG\+module $\M'$ over $\B$
such that $\M$ is isomorphic to $\M'$ in\/ $\sD^\co(\B\qcoh)$
and the set-theoretic support\/ $\Supp\M'$ is contained in~$T$. \par
\textup{(b)} For any coherent CDG\+module $\M$ over $\B$ with
the category-theoretic support\/ $\supp\M$ contained in~$T$,
there exists a coherent CDG\+module $\M'$ over $\B$ such that
$\M$ is isomorphic to a direct summand of $\M'$ in\/
$\sD^\abs(\B\coh)$ and the set-theoretic support\/ $\Supp\M'$
is contained in~$T$. 
\end{cor}

\begin{proof}
 Follows immediately from Theorem.
\end{proof}

\begin{rem}
 One can prove that the restriction functor in part~(a) of Theorem
is a Verdier localization functor without assuming the quasi-coherent
graded algebra $\B$ to be Noetherian.
 Indeed, one can construct a right adjoint functor $\boR j_*$ to
the restriction functor~$j^*$ in the way similar to that of
Proposition~\ref{finite-dim-morphisms}; then it is easy to see that
$j^*\boR j_*$ is the identity functor. 

 When $\B$ is Noetherian, the above Theorem can be generalized as
follows.
 Let $S$ and $T$ be closed subsets in $X$; set $U=X\setminus T$.
Then the restriction functor $\sD^\co(\B\qcoh_S)\rarrow\sD^\co
(\B|_U\qcoh_{U\cap S})$ is the Verdier localization functor by
the thick subcategory $\sD^\co(\B\qcoh_{T\cap S})$, and
the restriction functor $\sD^\abs(\B\coh_S)\rarrow
\sD^\abs(\B|_U\coh_{U\cap S})$ is the Verdier localization functor
by the triangulated subcategory $\sD^\abs(\B\coh_{T\cap S})$.
 The proof is similar to the above.

 It is not difficult to deduce from the latter assertions, using
the result of~\cite[Theorem~2.1(5)]{Neem}, that the property of
an object of $\sD^\co(\B\qcoh)$ to belong to the thick envelope
of $\sD^\abs(\B\coh)$ is local in~$X$.
 Using the \v Cech exact sequence as in Remark~\ref{second-kind},
one can easily see that the property of an object of
$\sD^\abs(\B\qcoh)$ to belong to $\sD^\abs(\B\qcoh_\fl)$ is
also local.

 We do \emph{not} know whether the property of an object
of $\sD^\abs(\B\coh)$ or $\sD^\abs(\B\qcoh_\fl)$ to belong to
$\sD^\abs(\B\coh_\fl)$ is local in general.
 In the particular case of matrix factorizations, such results
will be proven in Section~\ref{locality} using the connection with
singularity categories (cf.\ Remark~\ref{finite-dim-matrix-push}).
\end{rem}

 Notice that the theory of localization for compactly generated
triangulated categories, on which the proof of part~(b) of Theorem
is based, was originally applied in algebraic geometry for
the purposes of the Thomason--Trobaugh localization theory of
perfect complexes.
 In this section we apply it to coherent CDG\+modules.
 In fact, we will see in Section~\ref{nonlocalization} that
the localization theory \emph{fails} for locally free matrix
factorizations of finite rank.

\Section{Triangulated Categories of Relative Singularities}

\subsection{Relative singularity category}  \label{relative-sing}
 Recall that $X$ denotes a separated Noetherian scheme with enough
vector bundles.
 The \emph{triangulated category of singularities} $\sD^\b_\Sing(X)$
of the scheme $X$ is defined~\cite[Section~1.2]{Or1} as the quotient
category of the bounded derived category $\sD^\b(X\coh)$ of
coherent sheaves on $X$ by its thick subcategory $\Perf(X)$ of
perfect complexes on~$X$.

 The perfect complexes, in our assumptions, can be simply defined as
bounded complexes of locally free sheaves of finite rank, so $\Perf(X)
= \sD^\b(X\coh_\lf)$ is the bounded derived category of the exact
category $X\coh_\lf$ of locally free sheaves of finite rank
on~$X$.
 Equivalently, the perfect complexes are the compact objects of
the unbounded derived category of quasi-coherent sheaves
$\sD(X\qcoh)$ on the scheme~$X$
\cite[Examples~1.10\+-1.11 and Corollary~2.3]{Neem}.

 Let $Z\subset X$ be a closed subscheme such that $\O_Z$ has
finite flat dimension as an $\O_X$\+module.
 In this case the derived inverse image functor $\boL i^*$
for the closed embedding $i\:Z\rarrow X$ acts on the bounded
derived categories of coherent sheaves, $\sD^\b(X\coh)\rarrow
\sD^\b(Z\coh)$.
 We call the quotient category of $\sD^\b(Z\coh)$ by the thick
subcategory generated by the objects in the image of this functor
the \emph{triangulated category of singularities of\/ $Z$ relative
to $X$} and denote it by $\sD^\b_\Sing(Z/X)$.

 Note that the triangulated category of relative singularities
$\sD^\b_\Sing(Z/X)$ is a quotient category of the conventional
(absolute) triangulated category of singularities $\sD^\b_\Sing(Z)$
of the scheme~$Z$.
 Indeed, the thick subcategory $\Perf(Z)\subset\sD^\b(Z\coh)$ is
generated by any ample family of vector bundles on $Z$, since any
such family is a set of compact generators of the unbounded derived
category of quasi-coherent sheaves $\sD(Z\qcoh)$ on $Z$ \cite{Neem};
in particular, it is generated by the restrictions to $Z$ of
vector bundles from~$X$ (see also Lemma~\ref{infinite-matrix}).
 
 The functor $\boL i^*\:\sD^\b(X\coh)\rarrow\sD^\b(Z\coh)$ induces
a triangulated functor $i^\circ\:\sD^\b_\Sing(X)\rarrow
\sD^\b_\Sing(Z)$.
 Furthermore, since the sheaf $i_*\O_Z$ belongs to $\Perf(X)$,
the functor $i_*\:\sD^\b(Z\coh)\rarrow\sD^\b(X\coh)$ takes $\Perf(Z)$
to $\Perf(X)$ (cf.~\cite[paragraphs before Proposition~1.14]{Or1}).
 Hence the functor~$i_*$ induces a triangulated functor
$i_\circ\:\sD^\b_\Sing(Z)\rarrow\sD^\b_\Sing(X)$ right adjoint
to~$i^\circ$.
 The triangulated category $\sD^\b_\Sing(Z/X)$ is the quotient
category of $\sD^\b_\Sing(Z)$ by the thick subcategory generated
by the image of the functor~$i^\circ$.

 When $X$ is regular, any coherent sheaf on $X$ has a finite resolution
by locally free sheaves of finite rank.
 So $\sD^\b_\Sing(X)=0$, hence the triangulated categories
$\sD^\b_\Sing(Z)$ and $\sD^\b_\Sing(Z/X)$ coincide. 
 The converse is also true: the structure sheaf of the reduced scheme
structure on the closure of any singular point of $X$ is not a perfect
complex on $X$, so $\sD^\b_\Sing(X)\ne0$ when $X$ is not regular.

\begin{rem}
 Roughly speaking, the triangulated category of relative
singularities $\sD^\b_\Sing(Z/X)$ measures how much worse are
the singularities of $Z$ compared to the singularities of $X$
in a neighborhood of~$Z$.

 The basic formal properties of $\sD^\b_\Sing(Z/X)$ are similar to
those of $\sD^\b_\Sing(Z)$.
 When the $\O_X$\+module $\O_Z$ has finite flat dimension,
the derived category $\sD^\b(X\coh)$ is generated by coherent
sheaves adjusted to~$i^*$.
 Let $\sE_{Z/X}$ denote the minimal full subcategory of the abelian
category of coherent sheaves on~$Z$ containing the restrictions
of such coherent sheaves from $X$ and closed under extensions
and the kernels of epimorphisms of sheaves. 
 Then $\sE_{Z/X}$ is naturally an exact category and its bounded
derived category $\sD^\b(\sE_{Z/X})$ is equivalent to the thick
subcategory of $\sD^\b(Z\coh)$ generated by the derived restrictions
of coherent sheaves from~$X$, so $\sD^\b_\Sing(Z/X) = 
\sD^\b(Z\coh)/\sD^\b(\sE_{Z/X})$.
 One can define the $\sE$\+homological dimension of a coherent
sheaf (or bounded complex) on $Z$ as the minimal length of
a left resolution consisting of objects from~$\sE_{Z/X}$.
 This dimension does not depend on the choice of a resolution
(in the same sense as the conventional flat dimension doesn't).
 The thick subcategory $\sD^\b(\sE_{Z/X})$ consists of those
objects of $\sD^\b(Z\coh)$ that have finite $\sE$\+homological
dimensions.

 Unlike in the case of perfect complexes, we do not know whether
the property to belong to $\sE_{Z/X}$ or $\sD^\b(\sE_{Z/X})$
is local, though.
 In the case when $Z$ is a Cartier divisor, locality can be
established using Theorem~\ref{main-theorem} below and
Remark~\ref{second-kind}.
\end{rem}

\subsection{Matrix factorizations}  \label{matrix-subsect}
 Following~\cite{PV}, we will consider matrix factorizations of
a global section of a line bundle.
 So let $\L$ be a line bundle (invertible sheaf) on $X$ and
$w\in \L(X)$ be a fixed section, called the \emph{potential}.

 Let $\B=(X,\L,w)$ denote the following $\Z$\+graded quasi-coherent
CDG\+algebra over~$X$.
 The component $\B^n$ is isomorphic to $\L^{\ot n/2}$ for $n\in 2\Z$
and vanishes for $n\in 2\Z+1$, the multiplication in $\B$ being
given by the natural isomorphisms $\L^{\ot n/2}\ot_{\O_X}\L^{\ot m/2}
\rarrow \L^{\ot (n+m)/2}$.
 For any affine open subscheme $U\subset X$, the differential
on $\B(U)$ is zero, and the curvature element is $w|_U\in\B^2(U)=
\L(U)$.
 The elements $a_{UV}$ defining the restriction morphisms of
CDG\+rings $\B(V)\rarrow\B(U)$ all vanish.

 The category of quasi-coherent $\Z$\+graded $\B$\+modules is
equivalent to the category of quasi-coherent $\Z/2$\+graded
$\O_X$\+modules, the equivalence assigning to a graded $\B$\+module
$\M$ the pair of $\O_X$\+modules which we denote symbolically
by $\U^0=\M^0$ and $\U^1\ot\L^{\ot 1/2}=\M^1$.
 Conversely, $\M^n \simeq \U^{n\bmod 2}\ot_{\O_X}\L^{\ot n/2}$ for
all $n\in\Z$ (the meaning of the notation in the right-hand side
being the obvious one).
 This equivalence of abelian categories preserves all the properties
of coherence, flatness, flat dimension, local projectivity/local
freeness, etc.\ that we have been interested in in Section~1.

 Following~\cite{PL}, we will consider CDG\+modules over $\B=(X,\L,w)$
whose underlying graded $\B$\+modules correspond to coherent or
quasi-coherent $\O_X$\+modules, rather than just locally free sheaves
(as in the conventional matrix factorizations).
 A quasi-coherent CDG\+module over $(X,\L,w)$ is the same thing as
a pair of quasi-coherent $\O_X$\+modules $\U^0$ and
$\U^1\ot\L^{\ot 1/2}$ endowed with $\O_X$\+linear morphisms
$\U^0\rarrow\U^1\ot\L^{\ot 1/2}$ and $\U^1\ot\L^{\ot 1/2}\rarrow
\U^0\ot_{\O_X}\L$ such that both compositions $\U^0\rarrow\U^1\ot
\L^{\ot 1/2}\rarrow\U^0\ot_{\O_X}\L$ and $\U^1\ot\L^{\ot 1/2}\rarrow
\U^0\ot_{\O_X}\L\rarrow\U^1\ot_{\O_X}\L^{\ot 3/2}$ are equal to
the multiplications with~$w$.

\subsection{Exotic derived categories of matrix factorizations}
\label{exotic-derived-matrix-cor}

 The following corollary is a restatement of the results
of Section~1 in the application to the quasi-coherent CDG\+algebra
$\B=(X,\L,w)$.
 We will use the notation $(X,\L,w)\coh_\lf$ (instead of
the previously introduced $\B\coh_\fl$) for the DG\+category of
locally free matrix factorizations of finite rank, and
the notation $(X,\L,w)\qcoh_\lf$ (instead of the previously
introduced $\B\qcoh_\lp$) for the DG\+category of locally free
matrix factorizations of possibly infinite rank (see
Remark~\ref{flat-dimension-thm}).
 The rest of our notation system for various classes of
quasi-coherent CDG\+modules over $\B=(X,\L,w)$ remains in use.

 In addition, we also denote by $(X,\L,w)\qcoh_\lfd$
the DG\+category of quasi-coherent CDG\+modules of finite
locally free/locally projective dimension over $(X,\L,w)$
(see Remark~\ref{flat-dimension-thm} again).
 Let $\sD^\co((X,\L,w)\qcoh_\lfd)$ and $\sD^\abs((X,\L,w)
\qcoh_\lfd)$ be the corresponding derived categories of
the second kind.

\begin{cor}
\textup{(a)} The functor\/ $\sD^\co((X,\L,w)\qcoh_\fl)\rarrow
\sD^\co((X,\L,w)\qcoh_\ffd)$ induced by the embedding of
DG\+categories $(X,\L,w)\qcoh_\fl\rarrow(X,\L,w)\qcoh_\ffd$
is an equivalence of triangulated categories. \par
\textup{(b)} The functor\/ $\sD^\abs((X,\L,w)\qcoh_\fl)\rarrow
\sD^\abs((X,\L,w)\qcoh_\ffd)$ induced by the embedding of
DG\+categories $(X,\L,w)\qcoh_\fl\rarrow(X,\L,w)\qcoh_\ffd$
is an equivalence of triangulated categories. \par
\textup{(c)} The functors\/ $\sD^\co((X,\L,w)\qcoh_\lf)\rarrow\
\sD^\co((X,\L,w)\qcoh_\lfd)$ and\/ $\sD^\abs((X,\L,\allowbreak
w)\qcoh_\lf)\rarrow\sD^\abs((X,\L,w)\qcoh_\lfd)$ induced by
the embedding of DG\+categories $(X,\L,w)\qcoh_\lf\rarrow(X,\L,w)
\qcoh_\lfd$ are equivalences of triangulated categories. \par
\textup{(d)} The triangulated categories\/ $\sD^\co((X,\L,w)
\qcoh_\lf)$ and\/ $\sD^\abs((X,\L,w)\qcoh_\lf)$ coincide, as do
the categories\/ $\sD^\co((X,\L,w)\qcoh_\lfd)$ and\/
$\sD^\abs((X,\L,w)\qcoh_\lfd)$.  
 The natural functors between these four categories form
a commutative square of equivalences of triangulated categories. \par
\textup{(e)} When the scheme $X$ has finite Krull dimension,
the functors\/ $\sD^\co((X,\L,w)\allowbreak\qcoh_\lf)\rarrow
\sD^\co((X,\L,w)\qcoh_\fl)$ and\/ $\sD^\abs((X,\L,w)\qcoh_\lf)\rarrow
\sD^\abs((X,\L,w)\allowbreak\qcoh_\fl)$ induced by the embedding
of DG\+categories $(X,\L,w)\qcoh_\lf\rarrow(X,\L,w)\allowbreak
\qcoh_\fl$ are equivalences of triangulated categories.
 The natural functors between these four categories form
a commutative square of equivalences. \par
\textup{(f)} When the scheme $X$ has finite Krull dimension,
the triangulated category\/ $\sD^\co((X,\L,w)\qcoh_\fl)$
coincides with\/ $\sD^\abs((X,\L,w)\qcoh_\fl)$ and the triangulated
category\/ $\sD^\co((X,\L,w)\qcoh_\ffd)$ coincides with\/
$\sD^\abs((X,\L,w)\qcoh_\ffd)$.
 The natural functors between these four categories form
a commutative square of equivalences. {\hfuzz=1.1pt\par}
\textup{(g)} The functor\/ $\sD^\abs((X,\L,w)\coh_\lf)\rarrow
\sD^\abs((X,\L,w)\coh_\ffd)$ induced by the embedding of
DG\+categories $(X,\L,w)\coh_\lf\rarrow(X,\L,w)\coh_\ffd$
is an equivalence of triangulated categories. \par
\textup{(h)} The triangulated functors\/ $\sD^\abs((X,\L,w)\qcoh_\lf)
\rarrow\sD^\abs((X,\L,w)\qcoh_\fl)\allowbreak\rarrow
\sD^\abs((X,\L,w)\qcoh)$ induced by the embeddings of DG\+categories
$(X,\L,w)\allowbreak\qcoh_\lf\rarrow(X,\L,w)\qcoh_\fl\rarrow
(X,\L,w)\qcoh$ are fully faithful. \par
\textup{(i)} The triangulated functor\/ $\sD^\abs((X,\L,w)\coh_\lf)
\rarrow\sD^\abs((X,\L,w)\coh)$ induced by the embedding of
DG\+categories $(X,\L,w)\coh_\lf\rarrow(X,\L,w)\coh$
is fully faithful. \par
\textup{(j)} The triangulated functor\/ $\sD^\abs((X,\L,w)\coh_\lf)
\rarrow\sD^\co((X,\L,w)\qcoh_\lf)$ induced by the embedding
of DG\+categories $(X,\L,w)\coh_\lf\rarrow(X,\L,w)\qcoh_\lf$
is fully faithful. \par
\textup{(k)} The triangulated functor\/ $\sD^\abs((X,\L,w)\coh)
\rarrow\sD^\abs((X,\L,w)\qcoh)$ induced by the embedding of
DG\+categories $(X,\L,w)\coh\rarrow(X,\L,w)\qcoh$ is fully faithful.
\par
\textup{(l)} The triangulated functor\/ $\sD^\abs((X,\L,w)\coh)
\rarrow\sD^\co((X,\L,w)\qcoh)$ induced by the embedding of
DG\+categories $(X,\L,w)\coh\rarrow(X,\L,w)\qcoh$ is fully faithful
and its image forms a set of compact generators for\/
$\sD^\co((X,\L,w)\qcoh)$.
\end{cor}

\begin{proof}
 Parts (a\+b) and~(g) are particular cases of
Theorem~\ref{flat-dimension-thm}, and the proof of part~(c) is
similar (see Remark~\ref{flat-dimension-thm}).
 Part~(g) also essentially follows from
Proposition~\ref{embedding-prop}(b) (and part~(b) can be
proven similarly).
 Parts (h\+i) and~(k\+l) are particular cases of
Proposition~\ref{embedding-prop} (except for ``locally free half''
of part~(h), which is similar to the ``flat half'').
 Part~(d) is Theorem~\ref{homol-dimension-thm} together
with part~(c).
 Part~(j) is Corollary~\ref{homol-dimension-thm}.
 Part~(e) follows from parts~(a\+c) and
Remark~\ref{flat-dimension-thm} (cf.\ the discussion in
the end of Section~\ref{homol-dimension-thm}).
 Part~(f) follows from parts (a\+b) and~(d\+e); alternatively,
it can be proven directly in the way similar to part~(d),
using the fact that the exact category of flat 
quasi-coherent sheaves on $X$ has finite homological dimension
when the Krull dimension of $X$ is finite.
\end{proof}

\subsection{Regular and Gorenstein scheme cases}
\label{regular-cor}
 When the scheme $X$ is regular or Gorenstein, the assertions of
Corollary~\ref{exotic-derived-matrix-cor} simplify as follows.

\begin{cor}
\textup{(a)} When the scheme $X$ is Gorenstein of finite Krull
dimension, the functors\/ $\sD^\abs((X,\L,w)\qcoh_\fl)\rarrow
\sD^\co((X,\L,w)\qcoh_\fl)\rarrow\sD^\co((X,\L,w)\qcoh)$
induced by the embedding of DG\+categories $(X,\L,w)\qcoh_\fl\rarrow
(X,\L,w)\qcoh$ are equivalences of triangulated categories. \par
\textup{(b)} When the scheme $X$ is regular of finite Krull
dimension, the natural functors between the categories\/
$\sD^\abs((X,\L,w)\qcoh_\fl)$, \ $\sD^\co((X,\L,w)\qcoh_\fl)$,
\ $\sD^\abs((X,\L,w)\allowbreak\qcoh)$, and\/
$\sD^\co((X,\L,w)\qcoh)$ form a commutative square of equivalences of
triangulated categories. \par
\textup{(c)} When the scheme $X$ is regular, the natural functor\/
$\sD^\abs((X,\L,w)\coh_\lf)\rarrow\sD^\abs((X,\L,w)\coh)$
is an equivalence of triangulated categories.
\end{cor}

\begin{proof}
 Part~(a) is a particular case of Proposition~\ref{gorenstein-case}.
 Part~(c) follows from Corollary~\ref{exotic-derived-matrix-cor}(g),
since any coherent sheaf on a regular scheme has finite flat
dimension.
 In the assumptions of part~(b), the functor $\sD^\abs((X,\L,w)\qcoh)
\rarrow\sD^\co((X,\L,w)\qcoh)$ is an isomorphism of triangulated
categories by~\cite[Theorem~3.6(a) and Remark~3.6]{Pkoszul},
since the abelian category of quasi-coherent sheaves on a regular
scheme of finite Krull dimension has finite homological dimension
and enough injectives (cf.\ Theorem~\ref{homol-dimension-thm}).
 The remaining assertions of part~(b) follow from 
Corollary~\ref{exotic-derived-matrix-cor}(a\+b), or alternatively
from part~(a).
\end{proof}

 Assuming that $X$ has finite Krull dimension, the assertions
of Corollaries~\ref{exotic-derived-matrix-cor}\+-\ref{regular-cor}
may be summarized by the following commutative diagram of
triangulated functors.
 Here, as above, $\B$ denotes the quasi-coherent CDG\+algebra
$(X,\L,w)$:
\small
$$\dgARROWLENGTH=-10.75em %1.7em
\begin{diagram}
\node{\sD^\abs(\B\coh_\lf)}\arrow[6]{s,V}\arrow[2]{e,=}
\node[2]{\sD^\abs(\B\coh_\ffd)}\arrow[6]{s,V}\arrow[2]{see,t,V}
{\text{$=$ when $X$ regular}} \\ \\
\node[7]{\sD^\abs(\B\coh)}\arrow[4]{s,r,V} %{\text{c.g.}}
{\text{\smaller\smaller\smaller$\begin{matrix}
\text{comp.}\\ \text{gener.}\end{matrix}$}}
\arrow[2]{sssw,V} \\ \\
\node{}\node{\!\!\sD^{\co=\abs}(\B\qcoh_\lfd)\!\!}\arrow{ssw,=}
\arrow{sse,=}
\\ \\
\node{\sD^{\co=\abs}(\B\qcoh_{\lf})}\arrow[2]{e,=}
\node[2]{\sD^{\co=\abs}(\B\qcoh_\ffd)}\arrow[4]{e,t}
{\text{$=$ when $X$ Gorenstein\qquad\qquad\ \ }}
\arrow[2]{se,b,V}{\text{\smaller\smaller\smaller$\begin{matrix}
\text{$=$ when $X$}\\ \text{regular}\end{matrix}$}}
\node[4]{\sD^\co(\B\qcoh)}\\ \\
\node{}
\node{\sD^{\co=\abs}(\B\qcoh_\fl)}\arrow{nnw,=}\arrow{nne,=}
\node{}\node[2]{\sD^\abs(\B\qcoh)}\arrow[2]{ne,b,A}
{\text{\smaller\smaller\smaller$\begin{matrix}
\text{$=$ when $X$}\\ \text{regular}\end{matrix}$}}
\end{diagram}
$$\normalsize

 The four categories in the left lower area are coderived
categories coinciding with absolute derived categories
(of the same classes of quasi-coherent CDG\+modules).
 The five double lines between these four categories are
equivalences, as is the upper left horizontal line. 
 All the arrows going down are fully faithful functors.
 The image of the rightmost vertical arrow is a set of compact
generators in the target category.
 The only arrow going up is a Verdier localization functor.

 Nothing is claimed about the long horizontal arrow in the right
lower area of the diagram in general; but when $X$ is Gorenstein,
this functor is an equivalence of categories.
 When $X$ is regular, all the arrows going right are equivalences
of categories (so the whole diagram reduces to one triangulated
category with infinite direct sums, containing a full triangulated
subcategory of compact generators).

 Recall also that, by Lemma~\ref{gorenstein-case}, for any $X$
we have a commutative diagram of triangulated functors
\small
$$\dgARROWLENGTH=2.25em
\begin{diagram}
\node{H^0(\B\qcoh_\inj)}\arrow{e,=}
\node{\sD^{\co=\abs}(\B\qcoh_\fid)}\arrow[2]{e,=}
\arrow{se,V}\node[2]{\sD^\co(\B\qcoh)} \\
\node{}\node{}\node{\sD^\abs(\B\qcoh)}\arrow{ne,A}
\end{diagram}
$$\normalsize
with equivalences of categories in the upper line.
 The fully faithful embedding $\sD^\abs(\B\qcoh_\fid)\rarrow
\sD^\abs(\B\qcoh)$, which in the Gorenstein case (of finite Krull
dimension) coincides with the embedding $\sD^\abs(\B\qcoh_\ffd)
\rarrow\sD^\abs(\B\qcoh)$, is always right adjoint to
the localization functor $\sD^\abs(\B\qcoh)\rarrow\sD^\co(\B\qcoh)$.

\begin{rem}
 When $X$ is an affine Noetherian scheme of finite Krull
dimension, the embeddings of DG\+categories $(X,\L,w)\qcoh_\lp
\rarrow(X,\L,w)\qcoh_\fl\rarrow(X,\L,w)\qcoh$ induce
equivalences $H^0(\B\qcoh_\lp)\simeq\sD^\abs(\B\qcoh_\fl)
\simeq\sD^\ctr(\B\qcoh)$ between the homotopy category of (locally)
projective matrix factorizations of infinite rank, the absolute
derived category of flat matrix factorizations, and
the contraderived category of arbitrary quasi-coherent matrix
factorizations (see~\cite[Section~3.8]{Pkoszul};
cf.\ Remark~\ref{embedding-prop}).
\end{rem}

\subsection{Serre--Grothendieck duality}  \label{serre-duality}
 The aim of this section is to show that the somewhat mysterious
long horizontal arrow in the above large diagram is actually
a functor between two equivalent triangulated categories, for
a rather wide class of schemes~$X$.
 The functor $\sD^\co((X,\L,w)\qcoh_\fl)\rarrow\sD^\co((X,\L,w)\qcoh)$
in the above diagram, which is induced by the embedding of
DG\+categories $(X,\L,w)\qcoh_\fl\rarrow(X,\L,w)\qcoh$, is \emph{not}
the equivalence that we have in mind, however (unless the scheme
is Gorenstein).
 Instead, the equivalence between the categories
$\sD^\co((X,\L,w)\qcoh_\fl)$ and $\sD^\co((X,\L,w)\qcoh)$ is
constructed using a dualizing complex on~$X$ \cite[Section~V.2]{Har}.

 Before recalling the definition of a dualizing complex, let us
discuss the notion of the \emph{quasi-coherent internal Hom}.
 Given quasi-coherent sheaves $\M$ and $\N$ over $X$,
the quasi-coherent sheaf $\cHom_{X\qc}(\M,\N)$ is defined by
the isomorphism $\Hom_{\O_X}({-}\ot_{\O_X}\M\;\N)
\simeq\Hom_{\O_X}({-},\cHom_{X\qc}(\M,\N))$ of functors from
the category of quasi-coherent sheaves to the category of abelian
groups.
 Equivalently, the quasi-coherent sheaf $\cHom_{X\qc}(\M,\N)$
can be obtained by applying the coherator
functor~\cite[Sections~B.12--B.14]{TT} to the sheaf of $\O_X$\+modules
$\cHom_{\O_X}(\M,\N)$.
 Whenever $\M$ is a coherent sheaf, the sheaf $\cHom_{\O_X}(\M,\N)$
of $\O_X$\+module internal Hom is quasi-coherent, and
$\cHom_{X\qc}(\M,\N)\simeq\cHom_{\O_X}(\M,\N)$.

 Notice that the construction of the sheaf $\cHom_{X\qc}(\M,\N)$
is \emph{not} local in general, i.~e., it does not commute
with the restrictions of quasi-coherent sheaves to open subschemes;
when the sheaf $\M$ is coherent, it does.

\begin{lem}
\textup{(a)} For any injective quasi-coherent sheaf $\J$ over
a separated Noetherian scheme $X$, the functor
$\M\longmapsto\cHom_{X\qc}(\M,\J)$ is exact. \par
\textup{(b)} For any flat quasi-coherent sheaf $\F$ and injective
quasi-coherent sheaf $\J$ over $X$, the quasi-coherent sheaves
$\F\ot_{\O_X}\J$ and\/ $\cHom_{X\qc}(\F,\J)$ are injective. \par
\textup{(c)} For any injective quasi-coherent sheaves $\J'$ and
$\J$ over $X$, the quasi-coherent sheaf\/ $\cHom_{X\qc}(\J',\J)$
is flat.
\end{lem}

\begin{proof}
 The second assertion of part~(b) is obvious from the universal
property defining $\cHom_{X\qc}$.
 To prove the first one, notice that injectivity of quasi-coherent
sheaves over a Noetherian scheme is a local
property (\cite[Lemma~II.7.16 and Theorem~II.7.18]{Har} or
Theorem~\ref{injective-sheaves}), a flat quasi-coherent sheaf over
an affine scheme is a filtered inductive limit of locally free
sheaves of finite rank~\cite[No.~1.5--6]{Bour}, and injectivity of
modules over a Noetherian ring is preserved by filtered
inductive limits.

 The proof of parts~(a) and~(c) follows the argument
in~\cite[Lemma~8.7]{M-th}.
 Choose a finite affine covering $U_\alpha$ of the scheme $X$ and
consider the morphism $\J\rarrow\bigoplus_\alpha j_{U_\alpha *}
j_{U_\alpha}^*\J$.
 Being an embedding of injective quasi-coherent sheaves,
it splits, so $\J$ is a direct summand of the direct sum of
$j_{U_\alpha *}j_{U_\alpha}^*\J$.
 Hence it suffices to prove both assertions in the case when
$\J=j_{V*}\J''$, where $\J''$ is an injective quasi-coherent
sheaf on an affine open subscheme $V\subset X$.

 Now we have $\cHom_{X\qc}(\M,j_{V*}\J'')\simeq
j_{V*}\cHom_{V\qc}(j_V^*\M,\J'')$.
 Since $V\rarrow X$ is a flat affine morphism, the functor $j_V{}_*$
is exact and preserves flatness of quasi-coherent sheaves.
 This proves part~(a), and reduces part~(c) to the case of an affine
scheme $X=V$.
 Then it remains to apply~\cite[Proposition~VI.5.3]{CE}.
\end{proof}

 For our purposes, a \emph{dualizing complex} $\D_X^\bu$ on $X$ is
a finite complex of injective quasi-coherent sheaves such that
the cohomology sheaves of $\D_X^\bu$ are coherent and for any
coherent sheaf $\M$ over $X$ the natural morphism of finite complexes
of quasi-coherent sheaves $\M\rarrow\cHom_{X\qc}(\cHom_{X\qc}
(\M,\D_X^\bu),\D_X^\bu)$ is a quasi-isomorphism.
 Note that it follows from the former two conditions on $\D_X^\bu$
that the complex $\cHom_{X\qc}(\M,\D_X^\bu)$ has coherent cohomology
sheaves.
 This makes the conditions imposed on $\D_X^\bu$ actually local
in~$X$, so the restriction $\D^\bu_U=\D_X^\bu|_U$ of the complex
of sheaves $\D_X^\bu$ to an open subscheme
$U\subset X$ is a dualizing complex on~$U$.

 Given a quasi-coherent CDG\+algebra $\B$ over $X$, a quasi-coherent
left CDG\+module $\M$ over $\B$, and a complex of quasi-coherent
sheaves $\F^\bu$ on $X$, one can consider the complexes of
quasi-coherent left CDG\+modules $\F^\bu\ot_{\O_X}\M$ and
$\cHom_{X\qc}(\F^\bu,\M)$ over~$\B$.
 Taking their totalizations (formed, if necessary, by taking
infinite direct sums along the diagonals), one constructs two
triangulated functors $H^0(\B\qcoh)\rarrow H^0(\B\qcoh)$
depending on a complex~$\F^\bu$.
 Given a right CDG\+module $\N$ over $\B$ (see~\cite[Sections~3.1
and~B.1]{Pkoszul}), one can similarly construct a complex of
quasi-coherent left CDG\+modules $\cHom_{X\qc}(\N,\F^\bu)$
over~$\B$, obtaining a triangulated functor from the homotopy
category of right CDG\+modules $H^0(\qcohr\B)$ to $H^0(\B\qcoh)$.

 In the particular case of matrix factorizations, we conclude
that the covariant functors $\F^\bu\ot_{\O_X}{-}$ and
$\cHom_{X\qc}(\F^\bu,{-})$ take quasi-coherent matrix
factorizations of a potential $w\in\L(X)$ to (complexes
of) quasi-coherent matrix factorizations of~$w$, while
the contravariant functor $\cHom_{X\qc}({-},\F^\bu)$ transforms
quasi-coherent matrix factorizations of the opposite potential
$-w\in\L(X)$ into (complexes of) quasi-coherent matrix
factorizations of~$w$.
 Of course, the quasi-coherent CDG\+algebras $(X,\L,w)$ and
$(X,\L,-w)$ over a scheme $X$ are naturally isomorphic,
but we prefer to keep the distinction between the two.

 The next proposition provides the matrix factorization version of
the conventional (contravariant) Serre--Grothendieck duality for
bounded complexes of coherent sheaves. 
 We assume that $X$ is a separated Noetherian scheme with
a dualizing complex $\D_X^\bu$.
 Recall that any such scheme has finite
Krull dimension~\cite[Corollary~V.7.2]{Har}.
 We denote by $\sD^\sop$ the opposite category to a category~$\sD$.

\begin{prop}
 The triangulated functor\/ $\cHom_{X\qc}({-},\D_X^\bu)\:
H^0((X,\L,-w)\qcoh)^\sop\allowbreak\rarrow H^0((X,\L,w)\qcoh)$
induces a well-defined triangulated functor between the absolute
derived categories\/ $\sD^\abs((X,\L,-w)\qcoh)^\sop\rarrow
\sD^\abs((X,\L,w)\qcoh)$ taking the full triangulated
subcategory\/ $\sD^\abs((X,\L,-w)\coh)^\sop\subset
\sD^\abs((X,\L,-w)\qcoh)^\sop$ into the full subcategory\/
$\sD^\abs((X,\L,w)\coh)\subset\sD^\abs((X,\L,w)\qcoh)$.
 The composition of the duality functors\/ $\sD^\abs((X,\L,w)
\coh)\rarrow\sD^\abs((X,\L,-w)\coh)^\sop\rarrow\sD^\abs((X,\L,w)
\coh)$ is the identity functor.
\end{prop}

\begin{proof}
 The functor $\cHom_{X\qc}({-},\D_X^\bu)$ preserves absolute
acyclicity, because $\D_X^\bu$ is a complex of injective
quasi-coherent sheaves, so part~(a) of Lemma applies.
 Given a coherent matrix factorization $\M$, the finite complex
of matrix factorizations $\cHom_{X\qc}({-},\D_X^\bu)$ has
coherent cohomology matrix factorizations, so one can use its
canonical truncations in order to prove by induction that its
totalization belongs to the triangulated subcategory
$\sD^\abs((X,\L,w)\coh)$.

 Finally, for any quasi-coherent matrix factorization $\M$
consider the bicomplex of matrix factorizations
$\cHom_{X\qc}(\cHom_{X\qc}(\M,\D_X^\bu),\D_X^\bu)$ and take
its totalization in the two directions where it is a complex,
obtaining a complex of matrix factorizations.
 Then there is a natural morphism of finite complexes of matrix
factorizations $\M\rarrow\cHom_{X\qc}(\cHom_{X\qc}(\M,\D_X^\bu),
\D_X^\bu)$, which is a quasi-isomorphism of complexes of matrix
factorizations when $\M$ is coherent.
 The induced closed morphism of the total matrix factorizations
is an isomorphism in $\sD^\abs((X,\L,w)\qcoh)$, since 
the totalization of a finite acyclic complex of matrix
factorizations is absolutely acyclic.
 It remains to use the fact that the functor $\sD^\abs((X,\L,w)
\coh)\rarrow\sD^\abs((X,\L,w)\qcoh)$ is fully faithful 
(see Corollary~\ref{exotic-derived-matrix-cor}(k)) again.
\end{proof}
 
 The next result is our covariant Serre--Grothendieck duality
theorem for matrix factorizations.
 It is the matrix factorization analogue of the similar results
for complexes of projective and injective
modules~\cite[Theorem~4.2]{IK} and sheaves~\cite[Theorem~8.4]{M-th}.
 It also strongly resembles the \emph{derived comodule-contramodule
correspondence} theory (see~\cite[Theorem~5.2]{Pkoszul},
\cite[Corollaries~5.4 and~6.3]{Psemi}; cf.\
Remark~\ref{regular-cor} above).
 Notice that our proof is more akin to the arguments
in~\cite{Pkoszul,Psemi} than those of~\cite{IK,M-th} in that we
give a direct proof of the covariant duality independent of both
the contravariant duality and any descriptions of the compact
objects in the categories to be compared.

\begin{thm}\hfuzz=1.7pt
 The functors\/ $\D_X^\bu\ot_{\O_X}{-}\:H^0((X,\L,w)\qcoh_\fl)
\rarrow H^0((X,\L,w)\qcoh_\inj)$ and\/ $\cHom_{X\qc}(\D_X^\bu,{-})
\:H^0((X,\L,w)\qcoh_\inj)\rarrow H^0((X,\L,w)\qcoh_\fl)$
induce mutually inverse equivalences between the coderived
categories\/ $\sD^\co((X,\L,w)\qcoh_\fl)$ and\/
$\sD^\co((X,\L,w)\qcoh)$.
\end{thm}

\begin{proof}
 Recall that $H^0((X,\L,w)\qcoh_\inj)\simeq\sD^\co((X,\L,w)\qcoh)$
by Lemma~\ref{gorenstein-case}(b) and
$\sD^\abs((X,\L,w)\qcoh_\fl)=\sD^\co((X,\L,w)\qcoh_\fl)$ by
Corollary~\ref{exotic-derived-matrix-cor}(f) (though we will
reprove the latter fact rather than use it in the following
argument; see also Remark~\ref{w-flat-cor} below and
Lemma~\ref{flat-sheaves}).
 The functor $\D_X^\bu\ot_{\O_X}{-}\:H^0((X,\L,w)\qcoh_\fl)
\rarrow H^0((X,\L,w)\qcoh_\inj)$ obviously takes matrix
factorizations coacyclic with respect to $(X,\L,w)\qcoh_\fl$
to matrix factorizations coacyclic with respect to
$(X,\L,w)\qcoh_\inj$, which are all contractible.
 It remains to check that the induced functors are mutually inverse.

 Let $\E$ be a matrix factorization from $(X,\L,w)\qcoh_\fl$.
 As in the previous proof, consider the bicomplex of matrix
factorizations $\cHom_{X\qc}(\D_X^\bu\;\D_X^\bu\ot_{\O_X}\E)$
and take its total complex of matrix factorizations.
 Then there is a natural morphism $\E\rarrow\cHom_{X\qc}(\D_X^\bu\;
\D_X^\bu\ot_{\O_X}\E)$ of finite complexes of matrix factorizations
from $(X,\L,w)\qcoh_\fl$.
 To prove that the induced morphism of the total matrix factorizations
is an isomorphism in $\sD^\co((X,\L,w)\qcoh_\fl)$,
we once again use the fact that the totalization of a finite
acyclic complex of matrix factorizations is absolutely acyclic.
 So it suffices to check that for any flat quasi-coherent sheaf $\F$
over $X$ the natural morphism $\F\rarrow\cHom_{X\qc}(\D_X^\bu\;
\D_X^\bu\ot_{\O_X}\F)$ is a quasi-isomorphism of complexes of flat
quasi-coherent sheaves.
 This will be done below.

 Similarly, let $\M$ be a matrix factorization from
$(X,\L,w)\qcoh_\inj$.
 Consider the morphism of finite complexes of injective matrix
factorizations $\D_X^\bu\ot_{\O_X}\cHom_{X\qc}\allowbreak
(\D_X^\bu,\M)\rarrow\M$.
 To prove that the cone of the induced morphism of the total
matrix factorizations is contractible, it suffices to check
that for any injective quasi-coherent sheaf $\J$ over $X$
the natural morphism of complexes of injective sheaves
$\D_X^\bu\ot_{\O_X}\cHom_{X\qc}(\D_X^\bu,\J)\rarrow\J$
is a quasi-isomorphism.

 Let ${}'\D_X^\bu$ denote a finite complex of coherent sheaves
over~$X$ endowed with a quasi-isomorphism
${}'\D_X^\bu\rarrow\D_X^\bu$.
 Then the morphism $\cHom_{X\qc}(\D_X^\bu\;\D_X^\bu\ot_{\O_X}\F)
\rarrow\cHom_{X\qc}({}'\D_X^\bu\;\D_X^\bu\ot_{\O_X}\F)$ is
a quasi-isomorphism for any flat quasi-coherent sheaf~$\F$.
 The construction of the composition $\F\rarrow\cHom_{X\qc}
(\D_X^\bu\;\D_X^\bu\ot_{\O_X}\F)\rarrow\cHom_{X\qc}({}'\D_X^\bu\;
\D_X^\bu\ot_{\O_X}\F)$ is local in $X$, so it suffices to check
that the composition is a quasi-isomorphism when $X$ is affine.
 Then, using the passage to the filtered inductive limit, we may
assume that $\F$ is locally free of finite rank, and further
that $\F=\O_X$.
 It remains to recall that the morphism $\O_X\rarrow\cHom_{X\qc}
({}'\D_X^\bu,\D_X^\bu)$ is a quasi-isomorphism by the definition
of~$\D_X^\bu$.

 Let ${}''\D_X^\bu$ be a bounded above complex of flat quasi-coherent
sheaves mapping quasi-isomorphically to ${}'\D_X^\bu$.
 Then for any injective quasi-coherent sheaf $\J$ over $X$ there are
quasi-isomorphisms ${}''\D_X^\bu\ot_{\O_X}\cHom_{X\qc}(\D_X^\bu,\J)
\rarrow\D_X^\bu\ot_{\O_X}\cHom_{X\qc}(\D_X^\bu,\J)$ and
${}''\D_X^\bu\ot_{\O_X}\cHom_{X\qc}(\D_X^\bu,\J)\rarrow
{}''\D_X^\bu\ot_{\O_X}\cHom_{X\qc}({}'\D_X^\bu,\J)$ forming
a commutative diagram with the evaluation morphisms into~$\J$.
 Hence it remains to check that the morphism ${}''\D_X^\bu\ot_{\O_X}
\cHom_{X\qc}({}'\D_X^\bu,\J)\rarrow\J$ is a quasi-isomorphism, which
is a local question.
 Assume further that ${}''\D_X^\bu$ is a bounded above complex of
locally free sheaves of finite rank.
 Then there is a natural isomorphism of complexes of sheaves
${}''\D_X^\bu\ot_{\O_X}\cHom_{X\qc}({}'\D_X^\bu,\J)\simeq
\cHom_{X\qc}(\cHom_{X\qc}({}''\D_X^\bu,{}'\D_X^\bu),\J)$.
 The related morphism $\cHom_{X\qc}(\cHom_{X\qc}({}''\D_X^\bu,
{}'\D_X^\bu),\J)\rarrow\J$ is induced by the natural morphism
of complexes $\O_X\rarrow\cHom_{X\qc}({}''\D_X^\bu,{}'\D_X^\bu)$.
 The latter is again a quasi-isomorphism essentially by
the definition of $\D_X^\bu$.
\end{proof}

 From this point on we resume assuming that $X$ has enough
vector bundles.

 Notice that the equivalence functor $\D_X^\bu\ot_{\O_X}{-}\:
\sD^\co((X,\L,w)\qcoh_\lf)\rarrow\sD^\co((X,\L,w)\qcoh)$
that we have constructed takes the full triangulated subcategory
$\sD^\abs((X,\L,w)\coh_\lf)\subset\sD^\co((X,\L,w)\qcoh_\lf)$
into the full triangulated subcategory $\sD^\abs((X,\L,w)\coh)
\subset\sD^\co((X,\L,w)\qcoh)$.
 This is so because the dualizing complex $\D_X^\bu$ has
bounded coherent cohomology sheaves. {\hbadness=1500\par}

 Now we will use the above Proposition and Theorem in order to
construct compact generators of the triangulated category
$\sD^\co((X,\L,w)\qcoh_\lf)$ (cf.~\cite{Jo,N-f}).

 Consider the abelian category $Z^0((X,\L,-w)\coh)$ of coherent
matrix factorizations of $-w$ and closed morphisms of degree~$0$
between them, and its exact subcategory of locally free matrix
factorizations of finite rank $Z^0((X,\L,-w)\coh_\lf)$.
 The natural functor between the bounded above derived
categories of our abelian category and its exact subcategory
$\sD^-(Z^0((X,\L,-w)\coh_\lf))\rarrow\sD^-(Z^0((X,\L,-w)\coh))$
is an equivalence of triangulated categories.
 The vector bundle duality functor $\cHom_{X\qc}({-},\O_X)\:
Z^0((X,\L,-w)\coh_\lf)^\sop\rarrow Z^0((X,\L,w)\coh_\lf)$
induces a triangulated functor $\sD^-(Z^0((X,\L,-w)\coh_\lf))^\sop
\rarrow\sD^+(Z^0((X,\L,w)\coh_\lf))$ taking bounded above
complexes to bounded below ones.

 Let $\sD^+(Z^0((X,\L,w)\qcoh_\lf))$ denote the bounded below derived
category of the exact category of locally free matrix factorizations
of possibly infinite rank.
 Since the bounded below acyclic complexes over any exact category
with infinite direct sums are coacyclic~\cite[Lemma~2.1]{Psemi},
there is a well-defined, triangulated direct sum totalization functor
$\sD^+(Z^0((X,\L,w)\qcoh_\lf))\rarrow\sD^\co((X,\L,w)\qcoh_\lf)$.
 Consider the composition
\begin{multline*}
 Z^0((X,\L,-w)\coh)^\sop\lrarrow\sD^-(Z^0((X,\L,-w)\coh))^\sop \\
 \.\simeq\.\sD^-(Z^0((X,\L,-w)\coh_\lf))^\sop
 \lrarrow \sD^+(Z^0((X,\L,w)\coh_\lf))\lrarrow \\
 \sD^+(Z^0((X,\L,w)\qcoh_\lf))\lrarrow\sD^\co((X,\L,w)\qcoh_\lf),
\end{multline*} 
where two of the functors are the duality and the totalization
discussed above, while the other two are the natural embedding
and the functor induced by such.

 One easily checks that this composition takes cones of closed
morphisms in $Z^0((X,\L,-w)\coh)$ to cocones in
$\sD^\co((X,\L,w)\qcoh_\lf)$, hence induces a triangulated
functor $H^0((X,\L,-w)\coh)^\sop\rarrow\sD^\co((X,\L,w)\qcoh_\lf)$.
 Similarly, the above composition takes the totalizations of
short exact sequences in $(X,\L,-w)\coh$ to objects corresponding
to the totalizations of short exact sequences in $(X,\L,w)\qcoh_\lf$;
one checks this by considering a left locally free resolution of
a short exact sequence of coherent matrix factorizations.
 Thus we obtain a triangulated functor
$$
 \Omega\:\sD^\abs((X,\L,-w)\coh)^\sop\lrarrow
 \sD^\co((X,\L,w)\qcoh_\lf).
$$

\begin{cor}
 The functor $\Omega$ is fully faithful, and its image forms
a set of compact generators in $\sD^\co((X,\L,w)\qcoh_\lf)$.
 The following diagram of triangulated functors is commutative:
$$ \dgARROWPARTS=6\dgARROWLENGTH=2.6em
\begin{diagram}
\node{\sD^\abs((X,\L,-w)\coh_\lf)^\sop}
\arrow[4]{e,t,V}{\upsilon^\sop} \arrow{s,l,=}{\cHom_{X\qc}({-},\O_X)}
\node[4]{\sD^\abs((X,\L,-w)\coh)^\sop} 
\arrow{s,r,=}{\cHom_{X\qc}({-},\D_X^\bu)}
\arrow[2]{sww,t,V,2}{\Omega} \\
\node{\sD^\abs((X,\L,w)\coh_\lf)}
\arrow[4]{e,t,V}{\D_X^\bu\ot_{\O_X}{-}\qquad\qquad}
\arrow{s,l,V}{\kappa}
\node[4]{\sD^\abs((X,\L,w)\coh)}\arrow{s,lr,V}{\gamma}
{\text{\rm\smaller\smaller\smaller$\begin{matrix}
\text{comp.}\\ \text{gener.}\end{matrix}$}} \\
\node{\sD^\co((X,\L,w)\qcoh_\lf)}
\arrow[4]{e,tb,=}{\D_X^\bu\ot_{\O_X}{-}}{\cHom_{X\qc}(\D_X^\bu,{-})}
\node[4]{\sD^\co((X,\L,w)\qcoh)}
\end{diagram}
$$
 Here $\upsilon$, $\kappa$, and~$\gamma$ denote the fully faithful
functors induced by the natural embeddings of DG\+categories of
CDG\+modules.
 The two upper vertical lines are the natural contravariant dualities
(anti-equivalences) on the (absolute derived) categories of locally
free matrix factorizations of finite rank and coherent matrix
factorizations.
 The lower horizontal line is the equivalence of categories from
Theorem, and the middle horizontal arrow is the fully faithful
functor discussed after the proof of Theorem.
\end{cor}

 The above diagram is to be compared with the following
subdiagram of the large diagram in the end of
Section~\ref{regular-cor}:
$$ \dgARROWLENGTH=2.6em
\begin{diagram}
\node{\sD^\abs((X,\L,w)\coh_\lf)}
\arrow[4]{e,t,V}{\upsilon} \arrow{s,l,V}{\kappa}
\node[4]{\sD^\abs((X,\L,w)\coh)} \arrow{s,lr,V}{\gamma}
{\text{\rm\smaller\smaller\smaller$\begin{matrix}
\text{comp.}\\ \text{gener.}\end{matrix}$}} \\
\node{\sD^\co((X,\L,w)\qcoh_\lf)}
\arrow[4]{e,t}{\lambda}\node[4]{\sD^\co((X,\L,w)\qcoh)}
\end{diagram}
$$
 Here $\lambda$ denotes the triangulated functor induced by
the embedding of DG\+categories of CDG\+modules
$(X,\L,w)\qcoh_\lf\rarrow(X,\L,w)\qcoh$.

 Notice that it is clear from these two diagrams that
the functor~$\lambda$ is an equivalence of triangulated
categories whenever the functor $\upsilon$ is.
 Indeed, if $\upsilon$ is an equivalence of categories,
then the image of $\kappa$ is a set of compact generators
in the target category, and $\lambda$ is an infinite direct
sum-preserving triangulated functor identifying triangulated
subcategories of compact generators, hence an equivalence.
 In this case, the functor $\D^\bu_X\ot_{\O_X}{-}$ becomes
an auto-equivalence of the triangulated category
$\sD^\co((X,\L,w)\qcoh)$ and restricts to an auto-equivalence
of its full subcategory of compact generators
$\sD^\abs((X,\L,w)\coh)$.

\begin{proof}[Proof of Corollary]
 The assertions in the first sentence follow from the second one,
as we know $\gamma$ to be fully faithful and its image to be
a set of compact generators by
Corollary~\ref{exotic-derived-matrix-cor}(l).
 The commutativity of both squares and the upper left triangle
is clear.
 To check commutativity of the lower right triangle, consider
a coherent matrix factorization $\M$ of the potential~$-w$;
let $\E_\bu$ be its left resolution in the abelian category
$Z^0((X,\L,-w)\coh)$ whose terms $\E_n$ belong to
$Z^0((X,\L,-w)\coh_\lf)$.
 Then the finite complex of matrix factorizations
$\cHom_{X\qc}(\M,\D_X^\bu)$ maps quasi-isomorphically
to the bounded below complex of injective matrix factorizations
$\cHom_{X\qc}(\E_\bu,\D_X^\bu)\simeq\D_X^\bu\ot_{\O_X}
\cHom_{X\qc}(\E_\bu,\O_X)$, so the cone of the corresponding
morphism of the total matrix factorizations is coacyclic. 
\end{proof}

\subsection{$w$-flat matrix factorizations}  \label{w-flat-cor}
 From now on we will assume that for any affine open subscheme
$U\subset X$ the element $w|_U$ is not a zero divisor in
the $\O(U)$\+module $\L(U)$; in other words, the morphism of
sheaves $w\:\O_X\rarrow\L$ is injective.

 The following results will be used in the proof of the main
theorem and its analogues below.
 Let us call a quasi-coherent $\O_X$\+module $\E$ \ \emph{$w$\+flat}
if the map $w\:\E\rarrow\E\ot_{\O_X}\L$ is injective.
 Notice that any submodule of a $w$\+flat module is $w$\+flat,
so the ``$w$\+flat dimension'' of a quasi-coherent sheaf
over $X$ never exceeds~$1$.

 Denote by $(X,\L,w)\coh_\wfl$ the DG\+category of coherent
CDG\+modules over $(X,\L,w)$ with $w$\+flat underlying graded
$\O_X$\+modules and by $(X,\L,w)\qcoh_\wfl$ the similar
DG\+category of quasi-coherent CDG\+modules.
 Let $\sD^\abs((X,\L,w)\coh_\wfl)$, \ $\sD^\abs((X,\L,w)\qcoh_\wfl$,
and $\sD^\co((X,\L,w)\qcoh_\wfl)$ denote the corresponding
derived categories of the second kind.

 Furthermore, denote by $(X,\L,w)\coh_{\wfl\cap\ffd}$
the DG\+category of coherent CDG-modules over $(X,\L,w)$ whose
underlying graded $\O_X$\+modules are both $w$\+flat
and of finite flat dimension, and by $(X,\L,w)\qcoh_{\wfl\cap\lfd}$
the DG\+category of $w$\+flat quasi-coherent CDG\+modules of
finite locally free dimension.
 Let $\sD^\abs((X,\L,w)\coh_{\wfl\cap\ffd})$, \
$\sD^\abs((X,\L,w)\qcoh_{\wfl\cap\lfd})$,
and $\sD^\co((X,\L,w)\qcoh_{\wfl\cap\lfd})$ denote
the corresponding exotic derived categories.
{\emergencystretch=0em\par} \emergencystretch=2em

\begin{cor}
\textup{(a)} The functor\/ $\sD^\co((X,\L,w)\qcoh_\wfl)\rarrow
\sD^\co((X,\L,w)\qcoh)$ induced by the embedding of
DG\+categories $(X,\L,w)\qcoh_\wfl\rarrow(X,\L,w)\qcoh$
is an equivalence of triangulated categories. \par
\textup{(b)} The functor\/ $\sD^\abs((X,\L,w)\qcoh_\wfl)\rarrow
\sD^\abs((X,\L,w)\qcoh)$ induced by the embedding of
DG\+categories $(X,\L,w)\qcoh_\wfl\rarrow(X,\L,w)\qcoh$
is an equivalence of triangulated categories. \par
\textup{(c)} The functor\/ $\sD^\abs((X,\L,w)\coh_\wfl)\rarrow
\sD^\abs((X,\L,w)\coh)$ induced by the embedding of
DG\+categories $(X,\L,w)\coh_\wfl\rarrow(X,\L,w)\coh$
is an equivalence of triangulated categories. \par
\textup{(d)} The functor\/ $\sD^\co((X,\L,w)\qcoh_{\wfl\cap\lfd})
\rarrow\sD^\co((X,\L,w)\qcoh_\lfd)$ induced by the embedding of
DG\+categories $(X,\L,w)\qcoh_{\wfl\cap\lfd}\rarrow(X,\L,w)
\qcoh_\lfd$ is an equivalence of triangulated categories. \par
\textup{(e)} The functor\/ $\sD^\abs((X,\L,w)\qcoh_{\wfl\cap\lfd})
\rarrow\sD^\abs((X,\L,w)\qcoh_\lfd)$ induced by the embedding of
DG\+categories $(X,\L,w)\qcoh_{\wfl\cap\lfd}\rarrow(X,\L,w)
\qcoh_\lfd$ is an equivalence of triangulated categories. \par
\textup{(f)} The functor\/ $\sD^\abs((X,\L,w)\coh_{\wfl\cap\ffd})
\rarrow\sD^\abs((X,\L,w)\coh_\ffd)$ induced by the embedding of
DG\+categories $(X,\L,w)\coh_{\wfl\cap\ffd}\rarrow(X,\L,w)
\coh_\ffd$ is an equivalence of triangulated categories.
\end{cor}

\begin{proof}
 The proofs are analogous to those of
Corollary~\ref{exotic-derived-matrix-cor}(a\+c) and~(g)
(except that no induction in~$d$ is needed, as it suffices
to consider the case $d=1$).
 Parts~(d), (e), (f) are analogous to parts~(a), (b), (c),
respectively.
 Parts~(b\+c) and~(e\+f) can be also proven in the way similar
to Corollary~\ref{exotic-derived-matrix-cor}(h-i).
\end{proof}

\begin{rem}
 The assertions of parts~(a\+b) hold under somewhat weaker assumptions
than above: namely, one does not need to assume the existence of
enough vector bundles on~$X$.
 And one can make parts~(d\+e) hold without vector bundles by
replacing the finite locally free dimension condition in their
formulation with the finite flat dimension condition.
 The reason is that there are enough flat sheaves on any reasonable
scheme (see Lemma~\ref{flat-sheaves}).

 In fact, even part~(c) does not depend on the existence of vector
bundles, since a surjective morphism onto a given coherent sheaf $\M$
from a $w$\+flat coherent sheaf can be easily constructed, e.~g.,
by starting from a surjective morphism onto $\M$ from a flat
quasi-coherent sheaf $\F$ and picking a large enough coherent
subsheaf in~$\F$.
 Accordingly, one does not need vector bundles to prove
the equivalence of categories in the lower horizontal line in
Theorem~\ref{main-theorem} below and the other two equivalences
in Theorem~\ref{infinite-matrix}.
 Replacing locally free sheaves with flat ones in the relevant
definitions and assuming the Krull dimension to be finite, one can
have the whole of Proposition~\ref{infinite-matrix} hold without
vector bundles as well.

 Another alternative is to use \emph{very flat} quasi-coherent
sheaves, which there are always enough of and which always form
a category of finite homological dimension on a quasi-compact
semi-separated scheme~\cite[Section~4.1]{Pcosh}.
 Similarly, the existence of vector bundles is not needed for
the validity of Theorem~\ref{flat-dimension-thm}(a\+b),
Proposition~\ref{embedding-prop}(a,\,c\+d), all
the assertions of Sections~\ref{gorenstein-case}
and~\ref{cdg-supports}, Corollaries
\ref{exotic-derived-matrix-cor}(a\+b,\,f,\,k\+l)
and~\ref{regular-cor}(a\+b), Proposition~\ref{serre-duality} and
Theorem~\ref{serre-duality}, and some other results.
\end{rem}

\subsection{Main theorem}  \label{main-theorem}
 Let $X_0\subset X$ be the closed subscheme defined locally by
the equation $w=0$, and $i\:X_0\rarrow X$ be the natural
closed embedding.
 The next theorem is the main result of this paper.

\begin{thm}
 There is a natural equivalence of triangulated categories
$$
 \sD^\abs((X,\L,w)\coh)\.\simeq\.\sD^\b_\Sing(X_0/X).
$$
 Together with the functor\/
$\Sigma\:\sD^\abs((X,\L,w)\coh_\lf)\rarrow\sD^\b_\Sing(X_0)$
constructed in~\cite{Or3}, this equivalence forms
the following diagram of triangulated functors
$$\dgARROWLENGTH=1.75em
\begin{diagram}
\node[10]{\sD^\b_\Sing(X)}\arrow[2]{sw,tb,<>}
{i_\bu,\,i_\circ\!}{i^\circ} \\ \\
\node{0}\arrow[3]{e}
\node[3]{\sD^\abs((X,\L,w)\coh_\lf)}\arrow[4]{e,t,V}{\Sigma}
\arrow[3]{s,V} \node[4]{\sD^\b_\Sing(X_0)} \arrow[3]{s,A} \\ \\ \\
\node[4]{\sD^\abs((X,\L,w)\coh)} \arrow[4]{e,tb,=}
{\vphantom{b_b}\boL\Xi}{\vphantom{\tilde b}\Ups}
\node[4]{\sD^\b_\Sing(X_0/X)} \arrow[2]{s} \\ \\ 
\node[8]{0}
\end{diagram}
$$
where the upper horizontal arrow\/ $\Sigma$ is fully faithful,
the left vertical arrow is fully faithful, the right vertical arrow
is a Verdier localization functor, and the lower horizontal line\/
$\boL\Xi=\Ups^{-1}$ is an equivalence of categories.
 The square is commutative;
the three diagonal arrows $i_\bu$, $i^\circ$, $i_\circ$ (the middle
one pointing down and the two other ones pointing up) are adjoint.

 Furthermore, the image of the functor\/~$\Sigma$ is precisely
the full subcategory of objects annihilated by the functor~$i_\circ$,
or equivalently, by the functor~$i_\bu$.
 In other words, the image of\/~$\Sigma$ is equal both to the left
and to the right orthogonal complements to the thick subcategory
generated by the image of the functor~$i^\circ$,
i.~e., an object $\F\in\sD^\b_\Sing(X_0)$ is isomorphic
to\/~$\Sigma(\M)$ for some $\M\in\sD^\abs((X,\L,w)\coh_\lf)$ if
and only if for every $\E\in\sD^\b_\Sing(X)$ one has\/
$\Hom_{\sD^\b_\Sing(X_0)}(i^\circ\E,\F)=0$, or equivalently,
for every~$\E\in\sD^\b_\Sing(X)$ one has\/
$\Hom_{\sD^\b_\Sing(X_0)}(\F,i^\circ\E)=0$.
\end{thm}

 The thick subcategory generated by the image of the functor~$i^\circ$
is the kernel of the right vertical arrow.
 So the upper horizontal arrow and the right vertical arrow are
included into ``exact sequences'' of triangulated categories
(as marked by the zeroes at the ends; there is no exactness at
the uppermost rightmost end).

 When $X$ is a regular scheme, the functor $\sD^\abs((X,\L,w)\coh_\lf)
\rarrow\sD^\abs((X,\L,w)\allowbreak\coh)$ is an equivalence of
categories by Corollary~\ref{regular-cor}(c), and so is
the functor $\sD^\b_\Sing(X_0)\rarrow\sD^\b_\Sing(X_0/X)$ (as
explained in Section~\ref{relative-sing}).
 Hence it follows that the functor~$\Sigma$ is an equivalence of
categories, too.
 Thus we recover the result of Orlov~\cite[Theorem~3.5]{Or3}
claiming the equivalence of triangulated categories
$\sD^\abs((X,\L,w)\coh_\lf)\simeq\sD^\b_\Sing(X_0)$ for
a regular~$X$.

 A counterexample in Section~\ref{nonlocalization} will show that
when $X$ is not regular, the functor $\sD^\abs((X,\L,w)\coh_\lf)
\rarrow\sD^\abs((X,\L,w)\coh)$ does not have to be an equivalence,
and indeed, the thick subcategory generated by
$\sD^\abs((X,\L,w)\coh_\lf)$ can be a proper strictly full
subcategory in $\sD^\abs((X,\L,w)\coh)$.

\begin{proof}[Proof of the lower horizontal equivalence]
 To obtain the equivalence of triangulated categories in
the lower horizontal line, we will construct triangulated functors
in both directions, and then check that they are mutually inverse.
 Given a bounded complex of coherent sheaves $\F^\bu$ over $X_0$,
consider the CDG\+module $\Ups(\F^\bu)$ over $(X,\L,w)$ with
the underlying coherent graded module given by the rule
$$
 \Ups^n(\F^\bu)\.=\.\textstyle\bigoplus_{m\in\Z}\.
 i_*\F^{n-2m}\ot_{\O_X}\L^{\ot m}
$$
and the differential induced by the differential on~$\F^\bu$.
 Since $d^2=0$ on $\F^\bu$ and $w$ acts by zero in $i_*\F^j$,
this is a CDG\+module.
 It is clear that $\Ups$ is a well-defined triangulated functor
$\sD^\b(X_0\coh)\rarrow\sD^\abs((X,\L,w)\coh)$, since the derived
category of bounded complexes over an abelian category coincides
with their absolute derived category.

 Let us check that $\Ups$ annihilates the image of
the functor~$\boL i^*$.
 It suffices to consider a $w$\+flat coherent sheaf $\E$ on
$X$ and check that $\Ups(\coker w)=0$, where $w\:\E\ot_{\O_X}\L^{\ot-1}
\rarrow\E$.
 Indeed, $\Ups(\coker w)$ is the cokernel of the injective morphism
of contractible coherent CDG\+modules $\N\rarrow\M$, where
$\N^{2n+1} = \M^{2n+1} = \E\ot_{\O_X}\L^{\ot n}$ and
$\N^{2n} = \E\ot_{\O_X}\L^{\ot n-1}$, while $\M^{2n} =
\E\ot_{\O_X}\L^{\ot n}$ for $n\in\Z$.

 This provides the desired triangulated functor
$$
 \Ups\:\.\sD^\b_\Sing(X_0/X)\lrarrow\sD^\abs((X,\L,w)\coh).
$$
 The functor in the opposite direction is a version of Orlov's
cokernel functor, but in our situation it has to be constructed as
a derived functor, since the functor of cokernel of an arbitrary
morphism is not exact.
 Recall the equivalence of triangulated categories
$\sD^\abs((X,\L,w)\coh_\wfl)\rarrow\sD^\abs((X,\L,w)\coh)$ from
Corollary~\ref{w-flat-cor}(c).

 Define the functor $\Xi\:Z^0((X,\L,w)\coh_\wfl)\rarrow
\sD^\b_\Sing(X_0/X)$ from the category of $w$\+flat coherent
CDG\+modules over $(X,\L,w)$ and closed morphisms of degree~$0$
between them to the triangulated category of relative singularities
by the rule
$$
 \Xi(\M)\.=\.\coker(d\:\M^{-1}\to\M^0)\.=\.
 \coker(i^*d\:\.i^*\M^{-1}\to i^*{\M^0}),
$$
where the former cokernel, which is by definition a coherent sheaf
on $X$ annihilated by~$w$, is considered as a coherent sheaf on~$X_0$.
 One can immediately see that the functor $\Xi$ transforms
morphisms homotopic to zero into morphisms factorizable through
the restrictions to $X_0$ of $w$\+flat coherent sheaves on~$X$.
 Hence the functor $\Xi$ factorizes through the homotopy category
$H^0((X,\L,w)\coh_\wfl)$.

 It is explained in~\cite[Lemma~3.12]{PV} (see also
Lemma~\ref{finite-dim-matrix-push} below) that the functor
$\Xi$ is triangulated and in~\cite[Proposition~3.2]{Or3} that
the functor $\Xi$ factorizes through $\sD^\abs((X,\L,w)\coh_\wfl)$.
 The latter assertion can be also deduced by considering
the complex~(1.3) from~\cite{PV}.
 Indeed, the complex $i^*\M$ corresponding to the total CDG\+module
$\M$ of an exact triple in $\B\coh_\wfl$ is the total complex of
an exact triple of complexes in the exact category $\sE_{X_0/X}$
from Remark~\ref{relative-sing}, hence the complex $i^*\M$ is exact
with respect to $\sE_{X_0/X}$ and the cokernels of its differentials
belong to this exact subcategory in the abelian category of
coherent sheaves over~$X_0$.
 So we obtain the triangulated functor
$$
 \Xi\:\.\sD^\abs((X,\L,w)\coh_\wfl)\lrarrow\sD^\b_\Sing(X_0/X),
$$
and consequently, the left derived functor
$$
 \boL\Xi\:\.\sD^\abs((X,\L,w)\coh)\lrarrow\sD^\b_\Sing(X_0/X).
$$

 Let us check that the two functors $\Ups$ and $\boL\Xi$ are
mutually inverse.
 For any $w$\+flat coherent CDG\+module $\M$ over $(X,\L,w)$,
there is a natural surjective closed morphism of CDG\+modules
$\phi\:\M\rarrow\Ups\Xi(\M)$ with a contractible kernel.
 Clearly, $\phi\:\Id\rarrow\Ups\boL\Xi$ is an (iso)morphism
of functors.

 Conversely, any object of $\sD^\b_\Sing(X_0/X)$ can be represented
by a coherent sheaf on~$X_0$, and any morphism in
$\sD^\b_\Sing(X_0/X)$ is isomorphic to a morphism coming from
the abelian category of such coherent sheaves.
 Indeed, the bounded above derived category $\sD^-(X_0\coh)$ of
coherent sheaves over $X_0$ is equivalent to the bounded above
derived category $\sD^-(X_0\coh_\lf)$ of locally free sheaves;
using a truncation far enough to the left, one can represent any
object or morphism in $\sD^\b_\Sing(X_0/X)$ by a long enough shift
of a coherent sheaf or a morphism of coherent sheaves.
 Now for any coherent sheaf $\F$ on $X_0$ there is a natural
distinguished triangle $\F\ot_{\O_{X_0}}i^*\L^{\ot-1}[1]\rarrow
\boL i^* i_*\F\rarrow\F\rarrow\F\ot_{\O_{X_0}}i^*\L^{\ot-1}[2]$ in
$\sD^\b(X_0\coh)$, which provides a natural isomorphism $\F\simeq
\F\ot_{\O_{X_0}}i^*\L^{\ot-1}[2]$ in $\sD^\b_\Sing(X_0/X)$.

 Let $\F$ be a coherent sheaf on $X_0$; pick a vector bundle
$\E$ on $X$ together with a surjective morphism $\E\rarrow i_*\F$
with the kernel~$\E'$.
 Then the CDG\+module $\M$ over $(X,\L,w)$ with the components
$\M^{2n}=\E\ot_{\O_X}\L^{\ot n}$ and $\M^{2n-1}=\E'\ot_{\O_X}
\L^{\ot n}$ maps surjectively onto $\Ups(\F)$ with a contractible
kernel, and $\boL\Xi\Ups(\F)=\Xi(\M)=\F$ (cf.~\cite[Lemma~2.18]{PL}).
 Denote the isomorphism we have constructed by
$\psi\:\boL\Xi\Ups(\F)\rarrow\F$.
 The composition $\Ups\psi\circ\phi\Ups\:\Ups(\F)\rarrow\Ups\boL\Xi
\Ups(\F)\rarrow\Ups(\F)$ is clearly the identity morphism.
 It is obvious that $\psi$~commutes with any morphisms of coherent
sheaves $\F$ on $X_0$, but checking that it commutes with all
morphisms, or all isomorphisms, in $\sD^\b_\Sing(X_0/X)$ is a little
delicate (cf.\ Remark below).

 Notice that $\Ups\psi$ is an (iso)morphism of functors since
$\phi\Ups$ is, and consequently $\boL\Xi\Ups\psi$ is
an (iso)morphism of functors.
 Thus it remains to check that the functor $\boL\Xi\Ups$ is faithful,
i.~e., does not annihilate any morphisms.
 Indeed, any morphism in $\sD^\b_\Sing(X_0/X)$ is isomorphic to
a morphism coming from the abelian category of coherent sheaves
on $X_0$, and the functor $\boL\Xi\Ups$ transforms such morphisms
into isomorphic ones.
 The construction of the equivalence of categories in the lower
horizontal line is finished.
 One still has to check that the isomorphisms~$\phi$ commute
with the isomorphisms $\Ups\Xi(\M[1])\simeq\Ups\Xi(\M)[1]$, but
this is straightforward.

 Alternatively, one can use $w$\+flat coherent sheaves on $X$
or objects of the exact category $\sE_{X_0/X}$ of coherent
sheaves on~$X_0$ (as applicable) instead of the locally free
sheaves everywhere in the above argument.
\end{proof}

\begin{proof}[Proof of ``exactness'' in the upper line]
 We start with a discussion of the three adjoint functors in
the right upper corner.
 The functor $i_\circ$ right adjoint to the functor
$i^\circ\:\sD^\b_\Sing(X)\rarrow\sD^\b_\Sing(X_0)$
was constructed in Section~\ref{relative-sing}.

 To construct the left adjoint functor to~$i^\circ$, notice that
the right derived functor of subsheaf with the scheme-theoretic
support in the closed subscheme $\boR i^!\:\sD^\b(X\coh)\rarrow
\sD^\b(X_0\coh)$ only differs from the functor $\boL i^*$ by
a shift and a twist, $\boR i^!\E^\bu\simeq\boL i^*\E^\bu
\ot_{\O_{X_0}}\L|_{X_0}[-1]$.
 One can check this first for $w$\+flat coherent sheaves $\E$,
when both objects to be identified are shifts of sheaves, so
it suffices to compare their direct images under~$i$, which are
both computed by the same two-term complex $\E\rarrow\E\ot_{\O_X}\L$;
then replace a complex $\E^\bu$ with a finite complex of
$w$\+flat coherent sheaves (for a general result of this kind,
see~\cite[Theorem~5.4]{Neem}). {\hbadness=1150\par}

 Hence the functor $\boR i^!$ takes $\Perf(X)$ to $\Perf(X_0)$ and
induces a triangulated functor $i^\bu\:\sD^\b_\Sing(X)\rarrow
\sD^\b_\Sing(X_0)$ right adjoint to~$i_\circ$.
 It follows that the functor $i_\bu(\F)=i_\circ(\F)\ot_{\O_X}
\L[-1]$ is left adjoint to the functor~$i^\circ$.

 To prove the vanishing of the composition of functors in the upper
line and the orthogonality assertions, notice that
$$
 \Hom_{\sD^\b_\Sing(X_0)}(i^\circ\E,\Sigma\M) \.\simeq\.
 \Hom_{\sD^\b_\Sing(X)}(\E,i_\circ\Sigma\M)
$$
and $i_*\Sigma(\M)=\coker(\M^{-1}\to\M^0)\in\Perf(X)$
for any $\M\in\sD^\abs((X,\L,w)\coh_\lf)$, since the morphism
$\M^{-1}\rarrow\M^0$ of locally free sheaves on $X$ is injective.
 Similarly,
$$
 \Hom_{\sD^\b_\Sing(X_0)}(\Sigma\M,i^\circ\E)\.\simeq\.
 \Hom_{\sD^\b_\Sing(X)}(i_\bu\Sigma\M,\E)
$$
and $i_\bu\Sigma(\M) = i_\circ\Sigma(\M)\ot_{\O_X}
\L[-1] = 0$ in $\sD^\b_\Sing(X)$. 

 Obviously, our derived cokernel functor $\boL\Xi$ makes
a commutative diagram with the cokernel functor $\Sigma$
from~\cite{Or3}.
 The left vertical arrow is fully faithful by
Corollary~\ref{exotic-derived-matrix-cor}(i).
 The assertion that the upper horizontal arrow is fully faithful
is due to Orlov~\cite[Theorem~3.4]{Or3}.
 We have just obtained a new proof of it with our methods.
 Indeed, it follows from the orthogonality that the functor
$\sD^\b_\Sing(X_0)\rarrow\sD^\b_\Sing(X_0/X)$ induces isomorphisms
on the groups of morphisms between any two objects one of which
comes from $\sD^\abs((X,\L,w)\coh_\lf)$.
 Conversely, Orlov's theorem together with the orthogonality argument
and the equivalence of categories in the lower horizontal line
imply that the left vertical arrow is fully faithful.

 Now assume that $i_\circ\F=0$ for some $\F\in\sD^\b_\Sing(X_0)$.
 Clearly, there exists $m\ge0$ and a coherent sheaf $\K$ on $X_0$
such that $\F\simeq\K[m]$ in $\sD^\b_\Sing(X_0)$.
 Then $i_*\K$ is a perfect complex, i.~e., a coherent sheaf of
finite flat dimension on~$X$.
 Let us view it as an object of $(X,\L,w)\coh_\ffd$, i.~e.,
consider the CDG\+module $\N$ over $(X,\L,w)$ with the components
$\N^{2n}=i_*\K\ot_{\O_X}\L^{\ot n}$ and $\N^{2n+1}=0$.

 The construction of the cokernel functor $\Sigma$ can be
straightforwardly extended to $w$\+flat coherent matrix factorizations
of finite flat dimension, providing a triangulated functor
$\widetilde\Sigma\:\sD^\abs((X,\L,w)\coh_{\wfl\cap\ffd})\rarrow
\sD^\b_\Sing(X_0)$.
 The functor $\widetilde\Sigma$ is well-defined, since one has
$i^*\M\in\Perf(X_0)$ for any $w$\+flat coherent sheaf $\M$ of finite
flat dimension on~$X$.
 Using the equivalence of triangulated categories
$\sD^\abs((X,\L,w)\coh_{\wfl\cap\ffd})\simeq\sD^\abs((X,\L,w)
\coh_\ffd)$ from Corollary~\ref{w-flat-cor}(f), one constructs
the derived functor $\boL\widetilde\Sigma\:\sD^\abs((X,\L,w)\coh_\ffd)
\rarrow\sD^\b_\Sing(X_0)$ in the same way as it was done above
for the derived functor~$\boL\Xi$.
 Since the functor $\sD^\abs((X,\L,w)\coh_\lf)\rarrow
\sD^\abs((X,\L,w)\coh_\ffd)$ is an equivalence of categories
by Corollary~\ref{exotic-derived-matrix-cor}(g), the (essential)
images of the functors $\Sigma$ and $\boL\widetilde\Sigma$ coincide.

 Let us check that $\boL\widetilde\Sigma(\N)\simeq\K$ as
an object of $\sD^\b_\Sing(X_0)$.
 We argue as above, picking a vector bundle $\E$ on $X$ together
with a surjective morphism $\E\rarrow i_*\K$ with the kernel~$\E'$.
 Then the CDG\+module $\M$ over $(X,\L,w)$ with the components
$\M^{2n}=\E\ot_{\O_X}\L^{\ot n}$ and $\M^{2n-1}=\E'\ot_{O_X}
\L^{\ot n}$ maps surjectively onto $\N$ with a contractible kernel.
 Hence the object $\M\in (X,\L,w)\coh_{\wfl\cap\ffd}$ is
isomorphic to $\N$ in $\sD^\abs((X,\L,w)\coh_\ffd)$, and we have
$\boL\widetilde\Sigma(\N)=\widetilde\Sigma(\M)=\K$.
 Therefore, the object $\K\in\sD^\b_\Sing(X_0)$ belongs to
the (essential) image of the functor $\Sigma$, and it follows that
so does the object $\F\simeq\K[m]$.

 One can strengthen the above argument so as to obtain
a construction of the (partial) inverse functor $\Delta$ to
the functor $\Sigma$ similar to the above construction of
the functor~$\Ups$ inverse to the functor~$\boL\Xi$.
 Consider the full subcategory $\sF_{X_0/X}\subset X_0\coh$ in
the abelian category of coherent sheaves on $X_0$ consisting
of all the sheaves $\F$ such that the sheaf $i_*\F$ has finite
flat dimension (i.~e., is a perfect complex) on~$X$.
 The category $\sF_{X_0/X}$ contains all the locally free
sheaves on $X_0$ and is closed under the kernels of surjections,
the cokernels of embeddings, and the extensions.

 Hence $\sF_{X_0/X}$ is an exact subcategory in $X_0\coh$.
 The natural functor $\sD^\b(\sF_{X_0/X})\allowbreak
\rarrow\sD^\b(X_0\coh)$ is fully faithful; its image
coincides with the kernel of the composition of the direct
image and Verdier localization functors
$\sD^\b(X_0\coh)\rarrow\sD^\b(X\coh)\rarrow\sD^\b_\Sing(X)$.
 Accordingly, the quotient category $\sD^\b(\sF_{X_0/X})/
\sD^\b(X_0\coh_\lf)$ is identified with the kernel of the direct
image functor $i_\circ\:\sD^\b_\Sing(X_0)\rarrow\sD^\b_\Sing(X)$.

 Now the functor
$$
 \Delta\:\.\sD^\b(\sF_{X_0/X})/\sD^\b(X_0\coh_\lf)\lrarrow
 \sD^\abs((X,\L,w)\coh_\ffd)
$$
is constructed in the way similar to the construction of
the functor~$\Ups$, by taking the direct image from $X_0$ to $X$
and applying the periodicity summation.
 That is
$$
 \Delta^n(\F^\bu)\.=\.\textstyle\bigoplus_{m\in\Z}\.
 i_*\F^{n-2m}\ot_{\O_X}\L^{\ot m}
$$
for any $\F^\bu\in\sD^\b(\sF_{X_0/X})$.
 One checks that the functor $\Delta$ is inverse to the functor
$\boL\widetilde\Sigma$, the latter being viewed as a functor
taking values in the triangulated subcategory $\sD^\b(\sF_{X_0/X})/
\sD^\b(X_0\coh_\lf)\subset\sD^\b_\Sing(X_0)$, in the same way as
it was done above for the functors $\Ups$ and~$\boL\Xi$.
 This provides yet another proof of the fact that the functor
$\Sigma$ is fully faithful, together with another proof of
our description of its image.
 It is also obvious from the constructions that the functor
$\Delta$ makes a commutative diagram with the functor~$\Ups$. 
\end{proof}

\begin{rem}
 The somewhat tricky technical argument in the first part of the above
proof can be clarified and generalized using the approach developed
by the first author in~\cite[Appendix~A]{Ef2}.

 Let $\sC$ be an abelian category, $L\:\sC\rarrow\sC$ be
its covariant autoequivalence, and $w\:\Id\rarrow L$ be a natural
transformation commuting with~$L$ (that is for any object $B\in\sC$
one has $w_{L(B)}=L(w_B)$).
 Let $\MF(\sC,L,w)$ denote the abelian category of ``matrix
factorizations of~$w$ in~$\sC$'', i.~e., pairs of objects $U^0$, \
$L^{1/2}(U^1)\in\sC$ endowed with pairs of morphisms
$U^0\rarrow L^{1/2}(U^1)$, \ $L^{1/2}(U^1)\rarrow L(U^0)$ such that
the compositions $U^0\rarrow L^{1/2}(U^1)\rarrow L(U^0)$ and
$L^{1/2}(U^1)\rarrow L(U^0)\rarrow L^{3/2}(U^1)$ are equal to
$w_{U^0}$ and $w_{L^{1/2}(U^1)}$, respectively.
 Given a matrix factorization $M=(U^0,U^1)$, one sets
$M^n=L^{n/2}(U^{n\bmod 2})$.
 Passing to the quotient category by the ideal of morphisms homotopic
to zero, one obtains the homotopy category of matrix factorizations
of~$w$ in $\sC$; and their absolute derived category, denoted by
$\sD^\abs(\sC,L,w)$, is produced by the Verdier localization procedure
similar to the one discussed in Section~\ref{second-kind}.
 (Cf.~\cite[Remark~4.3]{Partin}.)

 Let $\sC_0\subset\sC$ denote the full subcategory formed by all
the objects $A\in\sC$ for which $w_A=0$; so $\sC_0$ is an abelian
subcategory in $\sC$ closed under subobjects and quotient objects.
 An object $B\in\sC$ is said to have no $w$\+torsion if
the morphism~$w_B$ is injective; and one says that the potential
(natural transformation) $w$ does not divide zero in $\sC$ if every
object of $\sC$ is the quotient object of an object without $w$\+torsion.
 Let $i_*\:\sC_0\rarrow\sC$ denote the exact identity embedding functor,
and $i^*\:\sC\rarrow\sC_0$ be the functor left adjoint to~$i_*$, so
$i^*(B)=\coker(w_{L^{-1}(B)}\:L^{-1}(B)\to B)$.
 Assuming that $w$~does not divide zero (as we do in the sequel),
one can construct the left derived functor $\boL i^*\:\sD^\b(\sC)
\rarrow\sD^\b(\sC_0)$ with $\boL_s i^*(B)=0$ for all $s\ne0$,~$1$
and any object $B\in\sC$.
 The functor $\boL i^*$ is left adjoint to the triangulated functor
$i_*\:\sD^\b(\sC_0)\rarrow\sD^\b(\sC)$ induced by the identity embedding
$i_*\:\sC_0\rarrow\sC$.

 Similarly, let $\upsilon^n\:\sC_0\rarrow\MF(\sC,L,w)$ denote the exact
functor assigning to an object $A\in\sC_0$ the matrix factorization $M$
with $M^n=A$ and $M^{n+1}=0$; and let $\xi^n\:\MF(\sC,L,w)\rarrow\sC_0$
be the functor left adjoint to~$\upsilon^n$, assigning the object
$\coker(M^{n-1}\to M^n)\in\sC_0$ to a matrix factorization~$M$.
 Considering the bounded derived category $\sD^\b\MF(\sC,L,w)$ of
the abelian category $\MF(\sC,L,w)$, one can construct the left derived
functor $\boL\xi^n\:\sD^\b\MF(\sC,L,w)\rarrow\sD^\b(\sC_0)$; once
again, the functor $\boL\xi^n$ is left adjoint to $\upsilon^n\:
\sD^\b(\sC_0)\rarrow\sD^\b\MF(\sC,L,w)$ and one has $\boL_s\xi^n(M)=0$
for all $s\ne0$,~$1$ and any matrix factorization $M\in\MF(\sC,L,w)$.

 Then the composition of the functor $\upsilon^n\:
\sD^\b(\sC_0)\rarrow\sD^\b\MF(\sC,L,w)$ with the totalization functor
$\sD^\b\MF(\sC,L,w)\rarrow\sD^\abs(\sC,L,w)$ induces an equivalence of
triangulated categories $\Ups^n\:\sD^\b(\sC_0)/
\langle \boL i^*\sD^\b(\sC)\rangle\rarrow\sD^\abs(\sC,L,w)$
between the quotient category of the derived category $\sD^\b(\sC_0)$
by the thick subcategory generated by the image of the functor
$\boL i^*$ and the absolute derived category of matrix factorizations.
 The composition of the functor $\boL\xi^n\:\sD^\b\MF(\sC,L,w)\rarrow
\sD^\b(\sC_0)$ with the Verdier localization functor
$\sD^\b(\sC_0)\rarrow\sD^\b(\sC_0)/\langle \boL i^*\sD^\b(\sC)\rangle$
factorizes through the totalization functor $\sD^\b\MF(\sC,L,w)\rarrow
\sD^\abs(\sC,L,w)$, providing the triangulated functor $\boL\Xi^n\:
\sD^\abs(\sC,L,w)\rarrow\sD^\b(\sC_0)/\langle \boL i^*\sD^\b(\sC)
\rangle$ inverse to~$\Ups^n$.

 Indeed, let $F^n\:\MF(\sC,L,w)\rarrow\sC$ denote the forgetful functor
taking a matrix factorization $M$ to the object $M^n\in\sC$, and let
$G^{n-}\:\sC\rarrow\MF(\sC,L,w)$ denote the functor left adjoint to
$F^{n-1}$ (and right adjoint to $F^n$); so the functor $G^{n-}$ takes
an object $B\in\sC$ to a contractible matrix factorization $M$ with
$M^{n-1}=M^n=B$ (cf.\ the constructions of the functors $G^+$ and $G^-$
in the proofs in Sections~\ref{flat-dimension-thm}
and~\ref{homol-dimension-thm}).
 It is claimed that the induced triangulated functors $G^{n-}\:
\sD^\b(\sC)\rarrow\sD^\b\MF(\sC,L,w)$ and $\upsilon^n\:\sD^\b(\sC_0)
\rarrow\sD^\b\MF(\sC,L,w)$ are fully faithful and their images form
a semiorthogonal decomposition of the derived category
$\sD^\b\MF(\sC,L,w)$.

 To check the first assertion, it suffices to notice that
the triangulated functor $G^{n-}$ is left adjoint to the functor
$F^{n-1}\:\sD^\b\MF(\sC,L,w)\rarrow\sD^\b(\sC)$, and their composition
$F^{n-1}\circ G^{n-}$ is the identity endofunctor on $\sD^\b(\sC)$.
 Similarly, the composition of triangulated functors $\boL\xi^n\circ
\upsilon^n$ is the identity endofunctor on $\sD^\b(\sC_0)$, so
$\upsilon^n$~is fully faithful as a functor between the derived
categories.
 Furthermore, one has $F^{n-1}\circ\upsilon^n=0=\boL\xi^n\circ G^{n-}$,
implying the semiorthogonality.
 Finally, for any matrix factorization $M$ whose terms are objects
without $w$\+torsion there is a short exact sequence
$0\rarrow G^{n-}F(M)\rarrow M\rarrow \upsilon^n\xi^nM\rarrow0$ in
$\MF(\sC,L,w)$ and $\boL\xi^nM=\xi^nM$, proving the decomposition claim.

 Now we notice that for any object $B\in\sC$ having no $w$\+torsion
there is a short exact sequence $0\rarrow G^{(n+2)-}(B)\rarrow
G^{(n+1)-}(B)\rarrow\upsilon^nj^*B\rarrow0$ in $\MF(\sC,L,w)$.
 According to (the proof of) \cite[Proposition~A.3(1\+2)]{Ef2},
the totalization functor $\sD^\b\MF(\sC,L,w)\rarrow\sD^\abs(\sC,L,w)$
is the Verdier localization functor by the thick subcategory generated
by the objects of the form $G^{n-}(B)=G^{(n+2)-}L(B)$ and
$G^{(n+1)-}(B)\in\MF(\sC,L,w)\subset\sD^\b\MF(\sC,L,w)$.
 The assertions about the existence of triangulated functors
$\Ups^n$ and $\boL\Xi^n$ and their being mutually inverse equivalences
of categories follow from these observations.

 Returning to a separated Noetherian scheme $X$ with enough vector
bundles and the Cartier divisor $X_0\subset X$ of a global section $w$
of a line bundle $\L$ on $X$, the above approach based
on~\cite[Proposition~A.3]{Ef2} provides an elegant construction of
Orlov's triangulated cokernel functor
$\Sigma\:\sD^\abs((X,\L,w)\coh_\lf)\rarrow\sD^\b_\Sing(X_0)$ in
addition to a proof of our equivalence of categories
$\sD^\abs((X,\L,w)\coh)\simeq\sD^\b_\Sing(X_0/X)$.
\end{rem}

\subsection{Infinite matrix factorizations}  \label{infinite-matrix}
 Following~\cite[paragraphs after Remark~1.9]{Or1}, one can define
a ``large'' version of the triangulated category of singularities
$\sD'_\Sing(X)$ of a scheme $X$ as the quotient category of
the bounded derived category of quasi-coherent sheaves
$\sD^\b(X\qcoh)$ by the thick subcategory $\sD^\b(X\qcoh_\lf)$
of bounded complexes of locally free sheaves (of infinite rank).
 When $X$ has finite Krull dimension, the latter subcategory
coincides with the thick subcategory $\sD^\b(X\qcoh_\fl)$ of
bounded complexes of flat sheaves (see
Remark~\ref{flat-dimension-thm}).

 Similarly, let $Z\subset X$ be a closed subscheme such that
$\O_Z$ has finite flat dimension as an $\O_X$\+module.
 Let us define a ``large'' triangulated category of relative
singularities $\sD'_\Sing(Z/X)$ as the quotient category
of $\sD^\b(Z\qcoh)$ by the minimal thick subcategory containing
the image of the functor $\boL i^*\:\sD^\b(X\qcoh)\rarrow
\sD^\b(Z\qcoh)$ and closed under those infinite direct sums that
exist in $\sD^\b(Z\qcoh)$.
 The quotient category of $\sD^\b(Z\qcoh)$ by the minimal thick
subcategory containing $\boL i^*\sD^\b(X\qcoh)$ (without
the direct sum closure) will be also of interest to us; let us
denote it by $\sD''_\Sing(Z/X)$.

\begin{lem}
 The triangulated categories\/ $\sD'_\Sing(Z/X)$ and\/
$\sD''_\Sing(Z/X)$ are quotient categories of\/ $\sD'_\Sing(Z)$.
 When the scheme $X$ is regular of finite Krull dimension, these
three triangulated categories coincide.
\end{lem}

\begin{proof}
 To prove the first assertion, let us show that any locally free
sheaf on $Z$, considered as an object of $\sD^\b(Z\qcoh)$, is
a direct summand of a bounded complex whose terms are direct sums
of locally free sheaves of finite rank restricted from~$X$.
 Indeed, pick a finite left resolution of a given locally free sheaf
on $Z$ with the middle terms as above, long enough compared to
the number of open subsets in an affine covering of~$Z$.
 Then the corresponding Ext class between the cohomology sheaves at
the rightmost and leftmost terms has to vanish in view of
the Mayer--Vietoris sequence for Ext groups between quasi-coherent
sheaves~\cite[Lemma~1.12]{Or1}.
 Hence the rightmost term is a direct summand of the complex
formed by the middle terms.

 The second assertion holds for the categories $\sD''_\Sing(Z/X)$
and $\sD'_\Sing(Z)$, since any quasi-coherent sheaf on
a regular scheme of finite Krull dimension has a finite left
resolution consisting of locally free sheaves.
 To identify these two categories with $\sD'_\Sing(Z/X)$, one needs
to know that the subcategory of bounded complexes of locally free
sheaves on $Z$ is closed under those infinite direct sums that
exist in $\sD^\b(Z\qcoh)$.
 The latter is true for any Noetherian scheme $Z$ of finite Krull
dimension with enough vector bundles, since the finitistic
projective dimension of a commutative ring of finite Krull dimension
is finite~\cite[Th\'eor\`eme~II.3.2.6]{RG}.
\end{proof}

 Now let $\L$ be a line bundle on $X$, \ $w\in\L(X)$ be a global
section corresponding to an injective morphism of sheaves
$\O_X\rarrow\L$, and $X_0\subset X$ be the locus of $w=0$.

\begin{prop}
 There is a natural equivalence of triangulated categories
$$
 \sD^\abs((X,\L,w)\qcoh)\.\simeq\.\sD''_\Sing(X_0/X).
$$
 Together with the infinite-rank version\/ $\Sigma'\:
\sD^\abs((X,\L,w)\qcoh_\lf)\rarrow\sD'_\Sing(X_0)$ of Orlov's
cokernel functor\/~$\Sigma$ from~\cite{Or3}, this equivalence
forms the following diagram of triangulated functors
$$\dgARROWLENGTH=1.75em
\begin{diagram}
\node[10]{\sD'_\Sing(X)}\arrow[2]{sw,tb,<>}
{i_\bu,\,i_\circ\!}{i^\circ} \\ \\
\node{0}\arrow[3]{e}
\node[3]{\sD^\abs((X,\L,w)\qcoh_\lf)}\arrow[4]{e,t,V}{\Sigma'}
\arrow[3]{s,V} \node[4]{\sD'_\Sing(X_0)} \arrow[3]{s,A} \\ \\ \\
\node[4]{\sD^\abs((X,\L,w)\qcoh)} \arrow[4]{e,=}
\node[4]{\sD''_\Sing(X_0/X)} \arrow[2]{s} \\ \\ 
\node[8]{0}
\end{diagram}
$$
where the upper horizontal arrow\/ $\Sigma'$ is fully faithful,
the left vertical arrow is fully faithful, the right vertical arrow
is the Verdier localization functor by the thick subcategory
generated by the image of the diagonal down arrow~$i^\circ$, and
the lower horizontal line is an equivalence of categories.
 The square is commutative; the three diagonal arrows
$i_\bu$, $i^\circ$, $i_\circ$ are adjoint.

 Furthermore, the image of the functor\/~$\Sigma'$ is precisely
the full subcategory of objects annihilated by the functor~$i_\circ$,
or equivalently, by the functor~$i_\bu$.
 In other words, the image of\/~$\Sigma'$ is equal both to the left
and to the right orthogonal complements to (the thick subcategory
generated by) the image of the functor~$i^\circ$.
\end{prop}

\begin{proof}
 The proof is completely similar to that of
Theorem~\ref{main-theorem}.
 It uses Corollaries~\ref{w-flat-cor}(b), 
\ref{exotic-derived-matrix-cor}(h), \ref{w-flat-cor}(e),
and~\ref{exotic-derived-matrix-cor}(c).
 The first assertion can be also obtained as a particular case
of the result of Remark~\ref{main-theorem}.

 Alternatively, one can prove that the functor $\Sigma'$
is fully faithful in the same way as it was done for
the functor~$\Sigma$ in~\cite[Theorem~3.4]{Or3}, and deduce
the assertion that the left vertical arrow is fully faithful
from the orthogonality.

 Note that one can check in a straightforward way that the functor
$\Sigma'$ annihilates the objects coacyclic with respect to
$(X,\L,w)\qcoh_\lf$.
 This provides another proof of
Corollary~\ref{exotic-derived-matrix-cor}(d), working
in the assumption that $w$~is a local nonzero-divisor.
\end{proof}

 The functors $\Sigma$ and $\Sigma'$ together with
the direct image functors~$i_\circ$ form the commutative diagram
of an embedding of ``exact sequences'' of triangulated functors
$$\dgARROWLENGTH=2em
\begin{diagram}
\node{0}\arrow[2]{e}\node[2]{\sD^\abs((X,\L,w)\coh_\lf)}
\arrow[3]{e,t,V}{\Sigma}\arrow{s,V}\node[3]{\sD^\b_\Sing(X_0)}
\arrow[2]{e,t}{i_\circ}\arrow{s,V}\node[2]{\sD^\b_\Sing(X)}
\arrow{s,V} \\
\node{0}\arrow[2]{e}\node[2]{\sD^\co((X,\L,w)\qcoh_\lf)}
\arrow[3]{e,t,V}{\Sigma'}\node[3]{\sD'_\Sing(X_0)}
\arrow[2]{e,t}{i_\circ}\node[2]{\sD'_\Sing(X)}
\end{diagram}
$$

 The leftmost vertical arrow is fully faithful by
Corollary~\ref{exotic-derived-matrix-cor}(j).
 The other two vertical arrows are fully faithful by Orlov's
theorem~\cite[Proposition~1.13]{Or1} claiming that the functor
$\sD^\b_\Sing(X)\rarrow\sD'_\Sing(X)$ is fully faithful
for any separated Noetherian scheme~$X$ with enough vector bundles.
 The leftmost nontrivial terms in both lines are the kernels
of the rightmost arrows by Theorem~\ref{main-theorem}
and Proposition above.

\begin{thm}
 There is a natural equivalence of triangulated categories
$$
 \sD^\co((X,\L,w)\qcoh)\.\simeq\.\sD'_\Sing(X_0/X)
$$
forming a commutative diagram of triangulated functors
$$\dgARROWLENGTH=6.4em
\begin{diagram}
\node{}\arrow{s,l,V}{\text{\smaller\smaller\rm$\begin{matrix}
\text{comp.}\\ \text{gener.}\end{matrix}$}}\\ \node{}
\end{diagram}
\dgARROWLENGTH=2em
\begin{diagram}
\node{\sD^\abs((X,\L,w)\coh)} \arrow[3]{e,=} \arrow{s,V}
\node[3]{\sD^\b_\Sing(X_0/X)} \arrow{s,V} \\
\node{\sD^\abs((X,\L,w)\qcoh)} \arrow[3]{e,=} \arrow{s,A}
\node[3]{\sD''_\Sing(X_0/X)} \arrow{s,A} \\ 
\node{\sD^\co((X,\L,w)\qcoh)} \arrow[3]{e,=}
\node[3]{\sD'_\Sing(X_0/X)}
\end{diagram}
\dgARROWLENGTH=6.4em
\begin{diagram}
\node{}\arrow{s,r,V}{\text{\,\smaller\smaller\rm$\begin{matrix}
\text{comp.}\\ \text{gener.}\end{matrix}$}}\\ \node{}
\end{diagram}
$$
with the equivalences of categories from
Theorem~\textup{\ref{main-theorem}} and the above Proposition.
 The upper vertical arrows are fully faithful, the lower ones
are Verdier localization functors, and the vertical compositions
are fully faithful.
 The categories in the lower line admit arbitrary direct
sums, and the images of the vertical compositions are sets
of compact generators in the target categories.
\end{thm}

\begin{proof}
 The construction of the desired equivalence of categories is very
similar to the construction of the equivalence of categories in
Theorem~\ref{main-theorem} and the Proposition.
 Using Corollary~\ref{w-flat-cor}(a), one defines
the infinite-rank version of the functor $\boL\Xi$, then
shows that the obvious infinite-rank version of
the functor $\Ups$ is inverse to it.
 Notice that the functor $\Xi\:Z^0((X,\L,w)\qcoh_\wfl)\rarrow
\sD^\b(X_0\qcoh)$ preserves infinite direct sums and the functor
$\Ups\:\sD^\b(X_0\qcoh)\rarrow\sD^\co((X,\L,w)\qcoh)$ preserves
those infinite direct sums that exist in $\sD^\b(X_0\qcoh)$, so
the functors $\Xi\:\sD^\co((X,\L,w)\qcoh_\wfl)\allowbreak\rarrow
\sD'_\Sing(X_0/X)$ and $\Ups\:\sD'_\Sing(X_0/X)\rarrow
\sD^\co((X,\L,w)\qcoh)$ are well-defined.

 The upper left vertical arrow is fully faithful by
Corollary~\ref{exotic-derived-matrix-cor}(k); it follows that
the upper right vertical arrow is fully faithful, too.
 The assertions about the vertical compositions are proved similarly.
 The category $\sD'_\Sing(X_0/X)$ admits arbitrary direct sums,
since the category $\sD^\co((X,\L,w)\qcoh)$ does.
 By Corollary~\ref{exotic-derived-matrix-cor}(l), the left vertical
composition is fully faithful and its image is a set of compact
generators in the target, so the right vertical composition has
the same properties.
\end{proof}

 The following square diagram of triangulated functors is
commutative:
$$
\begin{diagram}
\node{\sD^\co((X,\L,w)\qcoh_\lf)}\arrow[3]{e,t,V}{\Sigma'}
\arrow[2]{s}\node[3]{\sD'_\Sing(X_0)}\arrow[2]{s,A} \\ \\
\node{\sD^\co((X,\L,w)\qcoh)}\arrow[3]{e,=}
\node[3]{\sD'_\Sing(X_0/X)}
\end{diagram}
$$
 The upper horizontal arrow $\Sigma'$ is fully faithful;
the right vertical arrow is a Verdier localization functor.
 The lower line is an equivalence of triangulated categories.
 Nothing is claimed about the left vertical arrow in general.

 When the scheme $X$ is Gorenstein of finite Krull dimension,
the left vertical arrow is an equivalence of categories by
Corollary~\ref{regular-cor}(a).
 When $X$ is also regular, the right vertical arrow is
an equivalence of categories by the above Lemma.
 So $\Sigma'$ is an equivalence of categories $\sD^\abs((X,\L,w)
\qcoh_\lf)\simeq\sD'_\Sing(X_0)$ and we have obtained
a strengthened version of~\cite[Theorem~4.2]{PV}
(in the scheme case).

\begin{rem}
 It is well-known that the Verdier localization functor of
a triangulated category with infinite direct sums by a thick
subcategory closed under infinite direct sums preserves
infinite direct sums~\cite[Lemma~3.2.10]{N-book}.
 This result is not applicable to the localization functors
$\sD^\b(X\qcoh)\rarrow\sD'_\Sing(X)$ and $\sD^\b(Z\qcoh)\rarrow
\sD'_\Sing(Z/X)$, as the category $\sD^\b(X\qcoh)$ does not
admit arbitrary infinite direct sums.

 Using the equivalence of categories from the above Theorem
and the observation that the functor $\Ups$ preserves infinite
direct sums, one can show that the localization functor
$\sD^\b(X_0\qcoh)\rarrow\sD'_\Sing(X_0/X)$ takes those infinite
direct sums that exist in $\sD^\b(X_0\qcoh)$ into direct sums
in the triangulated category of relative singularities
$\sD'_\Sing(X_0/X)$ of the zero locus of~$w$ in~$X$.
 However, there is \emph{no} obvious reason why the localization
functor $\sD^\b(X_0\qcoh)\rarrow\sD'_\Sing(X_0)$ should take those
infinite direct sums that exist in $\sD^\b(X_0\qcoh)$ into
direct sums in the absolute triangulated category of singularities
$\sD'_\Sing(X_0)$.

 That is the problem one encounters attempting to prove that
the kernel of the localization functor $\sD'_\Sing(X_0)\rarrow
\sD'_\Sing(X_0/X)$ is semiorthogonal to the image of
the functor~$\Sigma'$.
\end{rem}

\subsection{Stable derived category}  \label{stable-derived}
 Following Krause~\cite{Kr}, we define the \emph{stable derived
category} of a Noetherian scheme $X$ as the homotopy category of
acyclic unbounded complexes of injective quasi-coherent sheaves
on~$X$.
 As explained below, this is another (and in some respects better)
``large'' version of the triangulated category of singularities of~$X$;
for this reason, we denote it by $\sD_\Sing^\st(X)$.

 In view of Lemma~\ref{gorenstein-case}(b) (see
also~\cite[Remark~5.4]{Psemi}), one can equivalently define
$\sD_\Sing^\st(X)$ as the quotient category of the homotopy
category of acyclic complexes of quasi-coherent sheaves over $X$
by the thick subcategory of coacyclic complexes, or as the full
subcategory of acyclic complexes in the coderived category
$\sD^\co(X\qcoh)$ of (complexes of) quasi-coherent sheaves over~$X$.
 It is the latter definition that will be used in the sequel.

 Clearly, the category $\sD^\st_\Sing(X)$ has arbitrary infinite
direct sums.
 In~\cite[Corollary~5.4]{Kr}, Krause constructs a fully faithful
functor $\sD^\b_\Sing(X)\rarrow\sD^\st_\Sing(X)$ and proves that
its image is a set of compact generators of the target category.

\begin{thm}
 For any separated Noetherian scheme $Z$ with enough vector bundles,
there is a natural triangulated functor\/ $\sD'_\Sing(Z)\rarrow
\sD^\st_\Sing(Z)$ forming a commutative diagram with the natural
functors from\/ $\sD^\b_\Sing(Z)$ into both these categories.
 The composition\/ $\sD^\b(Z\qcoh)\rarrow\sD'_\Sing(Z)\rarrow
\sD^\st_\Sing(Z)$ preserves those infinite direct sums that exist
in\/ $\sD^\b(Z\qcoh)$.
 When $Z=X_0$ is a divisor in a regular separated Noetherian scheme
of finite Krull dimension, the functor\/ $\sD'_\Sing(X_0)\rarrow
\sD^\st_\Sing(X_0)$ is an equivalence of triangulated categories.
\end{thm}

\begin{proof}
 The construction of the functor $\sD^\b_\Sing(Z)\rarrow
\sD^\st_\Sing(Z)$ in~\cite{Kr} is given in terms of the Verdier
localization functor $Q\:\sD^\co(Z\qcoh)\rarrow\sD(Z\qcoh)$ by
the triangulated subcategory $\sD_\Sing^\st(Z)\subset
\sD^\co(Z\qcoh)$ and its adjoint functors on both sides, which
exist according to~\cite[Corollary~4.3]{Kr}.
 The proof of our Theorem is based on explicit constructions of
the restrictions of these adjoint functors to some subcategories
of bounded complexes in $\sD(Z\qcoh)$.

 It is well known that the Verdier localization functor
$H^0(Z\qcoh)\rarrow\sD(Z\qcoh)$ from the homotopy category of
(complexes of) quasi-coherent sheaves on $Z$ to their derived
category has a right adjoint functor $\sD(Z\qcoh)\rarrow
H^0(Z\qcoh)$.
 The objects in the image of this functor are called
\emph{homotopy injective complexes} of quasi-coherent sheaves
on~$Z$.
 The composition $\sD(Z\qcoh)\rarrow H^0(Z\qcoh)\rarrow
\sD^\co(Z\qcoh)$ provides the functor $Q_\rho\:\sD(Z\qcoh)
\rarrow\sD^\co(Z\qcoh)$ right adjoint to~$Q$.
 In particular, any bounded below complex in $\sD(Z\qcoh)$ has
a bounded below injective resolution and any bounded below
complex of injectives is homotopy injective.
 Furthermore, any bounded below acyclic complex is
coacyclic~\cite[Lemma~2.1]{Psemi}.
 It follows that any bounded below complex from $\sD^+(Z\qcoh)$,
considered as an object of $\sD^\co(Z\qcoh)$, represents its
own image under the functor~$Q_\rho$.

 On the other hand, any bounded above complex from
$\sD(Z\qcoh)$ has a locally free left resolution defined uniquely
up to a quasi-isomorphism of complexes in the exact category of
locally free sheaves, i.~e., there is an equivalence of bounded
above derived categories $\sD^-(Z\qcoh_\lf)\simeq\sD^-(Z\qcoh)$.
 Since the exact category $Z\qcoh_\lf$ has finite homological
dimension, any acyclic complex in it is coacyclic
(and even absolutely acyclic~\cite[Remark~2.1]{Psemi}), so there
are natural functors $\sD^-(Z\qcoh_\lf)\rarrow\sD(Z\qcoh_\lf)
\allowbreak\simeq\sD^\co(Z\qcoh_\lf)\rarrow\sD^\co(Z\qcoh)$.

\begin{lem}
 The composition of the embedding\/ $\sD^-(Z\qcoh)\rarrow
\sD(Z\qcoh)$ with the functor $Q_\lambda\:\sD(Z\qcoh)\rarrow
\sD^\co(Z\qcoh)$ left adjoint to~$Q$ is isomorphic to
the functor\/ $\sD^-(Z\qcoh)\rarrow\sD^\co(Z\qcoh)$ constructed
above.
\end{lem}

\begin{proof}
 We have to show that $\Hom_{\sD^\co(Z\qcoh)}(\L^\bu,\E^\bu)=0$ for
any bounded above complex of locally free sheaves $\L^\bu$ and
any acyclic complex $\E^\bu$ of quasi-coherent sheaves on~$Z$.
 Let us check that any morphism $\L^\bu\rarrow\E^\bu$ in
$H^0(Z\qcoh)$ factorizes through a coacyclic complex of
quasi-coherent sheaves.
 Clearly, we can assume that the complex $\E^\bu$ is bounded
above, too.
 Let $\K^\bu$ be the cocone of a closed morphism of complexes
$\L^\bu\rarrow\E^\bu$; then $\K^\bu$ is bounded above and
the composition $\K^\bu\rarrow\L^\bu\rarrow\E^\bu$ is homotopic
to zero.
 Pick a bounded above complex of locally free sheaves $\F^\bu$
together with a quasi-isomorphism $\F^\bu\rarrow\K^\bu$.
 Then the cone of the composition $\F^\bu\rarrow\K^\bu\rarrow\L^\bu$,
being a bounded above acyclic complex of locally free sheaves,
is coacyclic.
 Since the composition $\F^\bu\rarrow\L^\bu\rarrow\E^\bu$ is
homotopic to zero, the morphism $\L^\bu\rarrow\E^\bu$ factorizes,
up to homotopy, through this cone.
\end{proof}

 Now we can describe the action of the functor $I_\lambda\:
\sD^\co(Z\qcoh)\rarrow\sD^\st_\Sing(Z\qcoh)$ left adjoint to
the embedding $\sD^\st_\Sing(Z\qcoh)\rarrow\sD^\co(Z\qcoh)$
on bounded above complexes in $\sD^\co(Z\qcoh)$.
 If $\K^\bu$ is a bounded above complex of quasi-coherent sheaves
and $\F^\bu$ is its locally free left resolution, then
the cone of the closed morphism $\F^\bu\rarrow\K^\bu$ represents 
the object $I_\lambda(\K^\bu)\in\sD^\st_\Sing(Z\qcoh)$.
 In view of the above Lemma, this cone is functorial and does not
depend on the choice of $\F^\bu$ for the usual
semiorthogonality reasons.

 The embedding of compact generators $\sD^\b_\Sing(Z)\rarrow
\sD^\st_\Sing(Z)$ is constructed in~\cite{Kr} as the functor
induced by the restriction of the composition $I_\lambda\circ Q_\rho\:
\sD(Z\qcoh)\rarrow\sD^\st_\Sing(Z)$ to the full subcategory
$\sD^\b(Z\coh)\subset\sD(Z\qcoh)$.
 Let us explain why this is so.
 By Proposition~\ref{embedding-prop}(d)
(cf.~\cite[Proposition~2.3 and Remark~3.8]{Kr}), the natural functor
$\sD^\b(Z\coh)\rarrow\sD^\co(Z\qcoh)$ is fully faithful and its image
is a set of compact generators in the target.
 This is the image of $\sD^\b(Z\coh)\subset\sD(Z\qcoh)$ under
the functor $Q_\rho$, as constructed above.
 It is clear from the above construction of the functor $Q_\lambda$
that it preserves compactness (and in fact coincides with
the functor $Q_\rho$ on perfect complexes in $\sD(Z\qcoh)$
\cite[Lemma~5.2]{Kr}).
 Since the functors $Q_\lambda$ and $I_\lambda$, being left adjoints,
preserve infinite direct sums, and $I_\lambda$ is a Verdier
localization functor by the image of $Q_\lambda$, it follows that
the image of any set of compact generators of $\sD^\co(Z\qcoh)$
under $I_\lambda$ is a set of compact generators of
$\sD^\st_\Sing(Z)$ \cite[Theorem~2.1(4)]{Neem}.

 In order to define the desired functor $\sD'_\Sing(Z)\rarrow
\sD^\st_\Sing(Z)$, restrict the same composition
$I_\lambda\circ Q_\rho$ to the full subcategory
$\sD^\b(Z\qcoh)\subset\sD(Z\qcoh)$.
 According to the above, this restriction assigns to any bounded
complex of quasi-coherent sheaves $\K^\bu$ on $Z$ the cone of
a morphism $\F^\bu\rarrow\K^\bu$ into it from its locally free
left resolution $\F^\bu$.
 Clearly, the functor $\sD^\b(Z\qcoh)\rarrow\sD^\st_\Sing(Z)$ that
we have constructed preserves those infinite direct sums that
exist in $\sD^\b(Z\qcoh)$ and annihilates the triangulated
subcategory $\sD^\b(Z\qcoh_\lf)\subset\sD^\b(Z\qcoh)$.
 So we have the induced functor $\sD'_\Sing(Z)\rarrow\sD^\st_\Sing(Z)$,
and the first two assertions of Theorem are proven.

 To prove the last assertion, we use the results of
Section~\ref{infinite-matrix}.
 Assume that $Z=X_0$ is the zero locus of a section $w\in\L(X)$
of a line bundle on $X$; as usually, $w\:\O_X\rarrow\L$ has to be
an injective morphism of sheaves.
 Then by Theorem~\ref{infinite-matrix} and Lemma~\ref{infinite-matrix},
the category $\sD'_\Sing(Z)$ admits infinite direct sums and
the image of the fully faithful functor $\sD^\b_\Sing(X_0)\rarrow
\sD'_\Sing(X_0)$ is a set of compact generators in the target.
 Furthermore, it follows from the proof of Theorem~\ref{infinite-matrix}
that any object of $\sD'_\Sing(X_0)$ can be represented by 
a quasi-coherent sheaf on $X_0$ and the direct sum of an infinite
family of such objects is represented by the direct sums of such
sheaves (see Remark~\ref{infinite-matrix}).
 Thus the functor $\sD'_\Sing(Z)\rarrow\sD^\st_\Sing(Z)$, being
an infinite direct sum-preserving triangulated functor identifying
triangulated subcategories of compact generators, is an equivalence
of triangulated categories.
\end{proof}

 We keep the assumptions of Theorem and the notation of the last
paragraph of its proof, i.~e., $X$ is a regular separated Noetherian
scheme of finite Krull dimension with enough vector bundles
and $X_0\subset X$ is the divisor of zeroes of a locally
nonzero-dividing section $w\in\L(X)$.
 The closed embedding $X_0\rarrow X$ is denoted by~$i$.

\begin{cor} \hbadness=1600
 The functor $\Lambda\:\sD^\co((X,\L,w)\qcoh_\lf)\simeq\sD^\co
((X,\L,w)\qcoh)\rarrow\sD^\st_\Sing(X_0)$ assigning to a locally free
(or just $w$\+flat) quasi-coherent matrix factorization $\M$
the acyclic complex of locally free (or quasi-coherent) sheaves
$i^*\M$ on $X_0$ is an equivalence of triangulated categories.
\end{cor}

\begin{proof}
 Given a $w$\+flat matrix factorization $\M$, the complex of
sheaves $i^*\M$ on $X_0$ is acyclic by~\cite[Lemma~1.5]{PV}.
 Clearly, the assignment $\M\longmapsto i^*\M$ defines 
a triangulated functor $\sD^\co((X,\L,w)\qcoh_\wfl)\rarrow
\sD^\st_\Sing(X_0)$.

 To prove that this functor is an equivalence of categories, it
suffices to identify it, up to a shift, with the composition of
the equivalences $\sD^\co((X,\L,w)\qcoh_\lf)\rarrow\sD'_\Sing(X_0)
\rarrow\sD^\st_\Sing(X_0)$.
 Here one simply notices that for any $\M\in\sD^\co((X,\L,w)
\allowbreak\qcoh_\lf)$ the complex $i^*\M$ is isomorphic in
$\sD^\st_\Sing(X_0)$ to its canonical truncation $\tau_{\le1}i^*\M$,
and the latter complex is the cocone of the morphism into
$\Sigma(\M)$ from one of its left locally free resolutions.
 So the functor $\Lambda$ is identified with $\Sigma[-1]$.
\end{proof}

\subsection{Relative stable derived category}
 The goal of this section is to generalize the results of
the previous one to the case of a singular Noetherian scheme~$X$.
 The relative version of stable derived category, defined for
a closed embedding of finite flat dimension $i\:Z\rarrow X$, is
equivalent to the categories $\sD'_\Sing(X_0/X)$ and
$\sD^\co((X,\L,w)\qcoh)$ in the case of the Cartier divisor
$Z=X_0$ corresponding to a locally nonzero-dividing section~$w$
of a line bundle $\L$ on~$X$.

 Let $X$ be a separated Noetherian scheme of finite Krull dimension
and $i\:Z\rarrow X$ be a closed embedding of schemes such that
$i_*\O_Z$ has finite flat dimension as an $\O_X$\+module.
 According to Section~\ref{finite-dim-morphisms}, there is 
a left derived inverse image functor $\boL i^*\:\sD^\co(X\qcoh)
\rarrow\sD^\co(Z\qcoh)$.
 This functor forms a commutative diagram with the similar functor
$\boL i^*\:\sD(X\qcoh)\rarrow\sD(Z\qcoh)$, and consequently, takes
acyclic complexes in $\sD^\co(X\qcoh)$ to acyclic complexes in
$\sD^\co(Z\qcoh)$.

\begin{prop} \label{relative-stable-prop}
 The following four triangulated categories are naturally equivalent:
\par
\textup{(a)} the full subcategory in\/ $\sD^\co(Z\qcoh)$ consisting
of all the objects annihilated by the direct image functor
$i_*\:\sD^\co(Z\qcoh)\rarrow\sD^\co(X\qcoh)$; \par
\textup{(b)} the quotient category of the homotopy category of
complexes over $Z\qcoh$ whose direct images are coacyclic complexes
over $X\qcoh$ by the thick subcategory of coacyclic complexes over
$Z\qcoh$; \par
\textup{(c)} the quotient category of\/ $\sD^\co(Z\qcoh)$ by its
minimal triangulated subcategory, containing the objects in\/
$\boL i^*\sD^\co(X\qcoh)$ and closed under infinite direct sums; \par
\textup{(d)} the quotient category of the full subcategory of
acyclic complexes in\/ $\sD^\co(Z\qcoh)$ by its minimal triangulated
subcategory, containing the left derived inverse images of
acyclic complexes in\/ $\sD^\co(X\qcoh)$ and closed under infinite
direct sums.
\end{prop}

\begin{proof}
 The equivalence of (a) and~(b) is obvious.
 To show that the natural functor from the category~(d) to
the category~(c) is an equivalence, notice that the minimal
triangulated subcategory containing flat quasi-coherent sheaves
and closed under infinite direct sums together with
the triangulated subcategory of acyclic complexes form
a semiorthogonal decomposition of $\sD^\co(X\qcoh)$, and similarly
for~$Z$ \cite[Corollary~A.4.7]{Pcosh}.
 Since flat quasi-coherent sheaves on $Z$ belong to the thick
subcategory in $\sD^\b(Z\qcoh)\subset\sD^\co(Z\qcoh)$ generated
by the inverse images of flat quasi-coherent sheaves from $X$ (see
the proof of Lemma~\ref{infinite-matrix}), the assertion follows.

 Finally, the functor $\boL i^*$ preserves infinite direct sums
and compactness of objects, since its right adjoint functor~$i_*$
preserves infinite direct sums.
 Hence the minimal triangulated subcategory in $\sD^\co(Z\qcoh)$
containing $\boL i^*\sD^\co(X\qcoh)$ and closed under infinite
direct sums is compactly generated by some objects which are 
compact in $\sD^\co(Z\qcoh)$.
 By Brown representability, the quotient category in~(c) is
equivalent to the right orthogonal complement to this
triangulated subcategory, which is the kernel category in~(a).
\end{proof}

 We call any of the equivalent triangulated categories 
in Proposition~\ref{relative-stable-prop} the \emph{relative
stable derived category of\/ $Z$ over\/~$X$} and denote it
by $\sD_\Sing^\st(Z/X)$ (cf.~\cite[Section~2]{Bec}).
 In particular, defining the relative stable derived category by
the construction~(c), we have natural triangulated functors
$\sD^\b(Z\qcoh) \rarrow\sD^\co(Z\qcoh)\rarrow\sD^\st_\Sing(Z/X)$.
 Clearly, the composition $\sD^\b(Z\qcoh)\rarrow
\sD^\st_\Sing(Z/X)$ factorizes through the relative singularity
category $\sD'_\Sing(Z/X)$, providing a natural functor
$\sD'_\Sing(Z/X)\rarrow\sD^\st_\Sing(Z/X)$.

\begin{lem}
 The composition of triangulated functors\/ $\sD^\b_\Sing(Z/X)
\rarrow\sD'_\Sing(Z/X)\rarrow\sD^\st_\Sing(Z/X)$ is fully
faithful and its image forms a set of compact generators for
the triangulated category\/ $\sD^\st_\Sing(Z/X)$.
\end{lem}

\begin{proof}
 By Proposition~\ref{embedding-prop}(d), the full triangulated
subcategory $\sD^\abs(Z\coh)$ compactly generates the triangulated
category $\sD^\co(Z\qcoh)$, and similarly for~$X$.
 In view of the construction~(c) and the argument in the proof of
Proposition, the assertion follows from~\cite[Theorem~2.1]{Neem0}.
\end{proof}

 Now let $\L$ be a line bundle on $X$, let $w\in\L(X)$ be
a locally nonzero-dividing section of~$\L$, and let $i\:X_0
\rarrow X$ be closed embedding of the zero locus of~$w$.
 Defining the category $\sD^\st_\Sing(X_0/X)$ by
the construction~(d), let $\boL\Lambda\:\sD^\co((X,\L,w)\qcoh)
\rarrow\sD^\st_\Sing(X_0/X)$ be the triangulated functor
assigning to a $w$\+flat quasi-coherent matrix factorization
$\M$ the acyclic complex $i^*\M$ over $X_0\qcoh$.

 Since any bounded below acyclic complex over $X_0\qcoh$
is coacyclic, and any bounded above complex belongs
to the minimal triangulated subcategory in $\sD^\co(X_0\qcoh)$
generated by its terms and closed under infinite direct sums,
the following diagram of triangulated functors is commutative
(cf.\ Corollary~\ref{stable-derived}) {\hbadness=1250
$$
\begin{diagram}
\node{\sD^\co((X,\L,w)\qcoh)}\arrow[2]{e,t}{\boL\Xi[-1]}
\arrow{se,b}{\boL\Lambda}
\node[2]{\sD'_\Sing(X_0/X)}\arrow{sw}\\
\node[2]{\sD^\st_\Sing(X_0/X)}
\end{diagram}
$$

\begin{thm}
 For any locally nonzero-dividing section~$w$ of a line bundle $\L$
on a separated Noetherian scheme $X$ of finite Krull dimension,
all the three functors on the above diagram are equivalences of
triangulated categories.
\end{thm}}

\begin{proof}
 The functor $\boL\Xi$ is an equivalence by
Theorem~\ref{infinite-matrix}.
 To show that the functor $\boL\Lambda$ is an equivalence,
let us check that it identifies compact generators.
 By Proposition~\ref{embedding-prop}(d), the category
$\sD^\co((X,\L,w)\qcoh)$ is compactly generated by its full
triangulated subcategory $\sD^\abs((X,\L,w)\coh)$, while
according to Lemma the category $\sD^\st_\Sing(X_0/X)$ is
compactly generated by its full triangulated subcategory
$\sD^\b_\Sing(X_0/X)$.
 The restriction of the functor $\boL\Lambda$ being an equivalence
between these two subcategories (in view of commutativity of
the diagram and) by Theorem~\ref{main-theorem}, it follows that
the functor $\boL\Lambda$ itself is an equivalence, too.
\end{proof}

\begin{rem}
 Another proof of the above Theorem can be obtained using the approach
based on~\cite[Appendix~A]{Ef2}.
 In the notation and assumptions of Remark~\ref{main-theorem}, suppose
that $\sC$ is an abelian category with exact functors of arbitrary
infinite direct sums.
 Then so is the abelian category $\MF(\sC,L,w)$; the full abelian
subcategory $\sC_0\subset\sC$ is closed under infinite direct sums;
and the triangulated functors $i_*$, $\boL i^*$, $\upsilon^n$,
$\boL\xi^n$, $F^n$, $G^{n-}$ act between the coderived categories
$\sD^\co(\sC)$, \ $\sD^\co(\sC_0)$, and $\sD^\co\MF(\sC,L,w)$.

 As in Remark~\ref{main-theorem}, one proves that the functors
$G^{n-}\:\sD^\co(\sC)\rarrow\sD^\co\MF(\sC,L,w)$ and
$\upsilon^n\:\sD^\co(\sC_0)\rarrow\sD^\co\MF(\sC,L,w)$ are fully
faithful and their images form a semiorthogonal decomposition of
the coderived category $\sD^\co\MF(\sC,L,w)$.
 By (the proof of) \cite[Proposition~A.3(3\+4)]{Ef2}, the totalization
functor $\sD^\co\MF(\sC,L,w)\rarrow\sD^\co(\sC,L,w)$ acting between
the coderived category of the abelian category $\MF(\sC,L,w)$ and
the coderived category of matrix factorizations $\sD^\co(\sC,L,w)$
(defined as in Section~\ref{second-kind}) is the Verdier localization
by the minimal triangulated subcategory containing the objects
$G^{n-}(B)$ and $G^{(n+1)-}(B)$ for all $B\in\sC$ and
closed under infinite direct sums.

 It follows that the composition of the functor
$\upsilon^n\:\sD^\co(\sC_0)\rarrow\sD^\co\MF(\sC,L,w)$ with
the totalization functor $\sD^\co\MF(\sC,L,w)\rarrow\sD^\co(\sC,L,w)$
induces an equivalence of triangulated categories
$\sD^\co(\sC_0)/\langle \boL i^*\sD^\co(\sC)\rangle_\oplus\rarrow
\sD^\co(\sC,L,w)$ between the quotient category of the coderived
category $\sD^\co(\sC_0)$ by its minimal triangulated subcategory
containing the image of the functor $\boL i^*\:\sD^\co(\sC)\rarrow
\sD^\co(\sC_0)$ and closed under infinite direct sums, and
the coderived category of matrix factorizations.
 The composition of the functor $\boL\xi^n\:\sD^\co\MF(\sC,L,w)
\rarrow\sD^\co(\sC_0)$ with the Verdier localization functor
$\sD^\co(\sC_0)\rarrow\sD^\co(\sC_0)/\langle\boL i^*\sD^\co(\sC)
\rangle_\oplus$ factorizes through the totalization functor,
providing the inverse equivalence $\sD^\co(\sC,L,w)\rarrow
\sD^\co(\sC_0)/\langle\boL i^*\sD^\co(\sC)\rangle_\oplus$.

 Returning to quasi-coherent matrix factorizations of a global
section $w\in\L(X)$ of a line bundle $\L$ on a separated Noetherian
scheme $X$ with the zero locus $X_0\subset X$, we obtain direct
constructions of two mutually inverse triangulated equivalences
between the coderived category $\sD^\co((X,\L,w)\qcoh)$ and
the relative stable derived category $\sD^\st_\Sing(X_0/X)$
as defined in part~(c) of the above Proposition.
\end{rem}

\Section{Supports, Pull-Backs, and Push-Forwards}

\subsection{Supports}  \label{supports-secn}
 This section paves the ground for the results about preservation
of finite rank or coherence by the push-forwards of matrix
factorizations with proper supports, which will be proven in
Sections~\ref{matrix-push}\+-\ref{finite-dim-matrix-push}.

 Let $X$ be a separated Noetherian scheme and $T\subset X$ be
a Zariski closed subset.
 Denote by $X\coh_T$ the abelian category of coherent sheaves on $X$
with the set-theoretic support in $T$; and similarly for quasi-coherent
sheaves.

 It is a well-known fact (essentially, a reformulation of
the Artin--Rees lemma) that the embedding of abelian categories
$X\qcoh_T\rarrow X\qcoh$ takes injectives to injectives.
 It follows that the functor $\sD^\b(X\coh_T)\rarrow\sD^\b(X\coh)$
is fully faithful.
 Clearly, its image is a thick subcategory and the corresponding
quotient category can be naturally identified with $\sD^\b(U\coh)$,
where $U=X\setminus T$ (cf.\ Section~\ref{cdg-supports}).

 Assume additionally that $X$ has enough vector bundles.
 Let $\Perf_T(X)\subset\Perf(X)$ denote the full subcategory of
perfect complexes with the cohomology sheaves set-theoretically
supported in~$T$.
 By the above result, $\Perf_T(X)$ can be considered as a thick
subcategory in $\sD^\b(X\coh_T)$.
 According to~\cite[Lemma~2.6]{Or2}, the functor
$\sD^\b(X\coh_T)/\Perf_T(X)\rarrow\sD^\b_\Sing(X)$ induced by
the embedding $\sD^\b(X\coh_T)\rarrow\sD^\b(X\coh)$ is fully faithful.
 We denote the source (or the image) category of this functor by
$\sD^\b_\Sing(X,T)$.

 By~\cite[Theorem~1.3]{Ch}, the restriction functor $\sD^\b_\Sing(X)
\rarrow\sD^\b_\Sing(U)$ is the Verdier localization functor by
the triangulated subcategory $\sD^\b_\Sing(X,T)$.
 In particular, the kernel of the restriction functor coincides with
the thick envelope of (i.~e., the minimal thick subcategory containing)
$\sD^\b_\Sing(X,T)$ in $\sD^\b_\Sing(X)$.

 Now we are going to establish the similar results for the triangulated
categories of relative singularities.
 Let $i\:Z\rarrow X$ be a closed subscheme such that
$i_*\O_Z\in\Perf(X)$, and let $\Perf(Z/X)=\sD^\b(\sE_{Z/X})$
(see Remark~\ref{relative-sing}) denote the thick subcategory
in $\sD^\b(Z\coh)$ generated by $\boL i^*\sD^\b(X\coh)$.
 Let $T\subset Z$ be a Zariski closed subset; put $U=X\setminus T$
and $V=Z\setminus T$.
 We denote by $\Perf_T(Z/X)$ the full subcategory of all objects of
$\Perf(Z/X)$ with the cohomology sheaves set-theoretically supported
in~$T$.
 Consider it as a thick subcategory in $\sD^\b(Z\coh_T)$, and denote
by $\sD^\b_\Sing(Z/X,T)$ the quotient category $\sD^\b(Z\coh_T)/
\Perf_T(Z/X)$.

\begin{lem}
 \textup{(a)} The functor\/ $\sD^\b_\Sing(Z/X,T)\rarrow
\sD^\b_\Sing(Z/X)$ induced by the embedding\/ $\sD^\b(Z\coh_T)\rarrow
\sD^\b(Z\coh)$ is fully faithful. \par
\textup{(b)} The restriction functor\/ $\sD^\b_\Sing(Z/X)\rarrow
\sD^\b_\Sing(V/U)$ is the Verdier localization functor by
the triangulated subcategory\/ $\sD^\b_\Sing(Z/X,T)$.
 In particular, the kernel of the restriction functor coincides with
the thick envelope of\/ $\sD^\b_\Sing(Z/X,T)$ in\/ $\sD^\b_\Sing(Z/X)$.
\end{lem}

\begin{proof}
 The proof of part~(a) is similar to that of~\cite[Lemma~2.6]{Or2}.
 One only needs to notice that the tensor product of an object of
$\Perf(Z/X)$ with an object of $\Perf(Z)$ belongs to $\Perf(Z/X)$.
 This follows from the fact that $\Perf(Z)$ as a thick subcategory
in $\sD^\b(Z\coh)$ is generated by the restrictions of vector
bundles from~$X$ (see Section~\ref{relative-sing}).
 Part~(b) is true, since the thick subcategory
$\Perf(V/U)\subset\sD^\b(V\coh)$ is generated by the image of
the restriction functor $\Perf(Z/X)\rarrow\Perf(V/U)$, which is
because any coherent sheaf on $U$ can be extended to a coherent
sheaf on~$X$.
\end{proof}

 Let $\L$ be a line bundle over~$X$ and $w\in\L(X)$ be a section;
set $X_0=\{w=0\}\subset X$.
 The definitions of the set-theoretic and category-theoretic
supports $\Supp\M$ and $\supp\M$ of a coherent matrix factorization
$\M\in(X,\L,w)\coh$ were given (in a greater generality of coherent
CDG\+modules) in Section~\ref{cdg-supports}.

 Given a locally free matrix factorization of finite rank $\M\in
(X,\L,w)\coh_\lf$, define the \emph{(category-theoretic) support}
$\supp\M\subset X$ as the minimal closed subset $T\subset X$
such that the restriction $\M|_U$ of $\M$ to the open subscheme
$U=X\setminus T$ is absolutely acyclic with respect to
$(U,\L|_U,w|_U)\coh_\lf$.
 By Corollary~\ref{exotic-derived-matrix-cor}(i), the definitions
of category-theoretic supports of coherent matrix factorizations and
of locally free matrix factorizations of finite rank agree when
they are both applicable.

 Equivalently, for a locally free matrix factorization $\M$ of
finite rank over~$X$, the open subscheme $X\setminus\supp\M$ is
the union of all affine open subschemes $U\subset X$ such that
the matrix factorization $\M|_U$ is contractible
(see Remark~\ref{second-kind}).
 For any coherent matrix factorization $\M$ one has $\supp\M
\subset X_0$, since any matrix factorization of an invertible
potential is contractible (cf.~\cite[Section~5]{PV}).

 Let $T\subset X$ be a closed subset.
 Denote by $\sD^\abs_T((X,\L,w)\coh_\lf)$ (respectively,
$\sD^\abs_T((X,\L,w)\coh)$)
the quotient category of the homotopy category of locally free
matrix factorizations of finite rank (resp., coherent matrix
factorizations) supported category-theoretically inside $T$ by
the thick subcategory of matrix factorizations absolutely acyclic
with respect to $(X,\L,w)\coh_\lf$ (resp., $(X,\L,w)\coh$).
 Clearly, the functors $\sD^\abs_T((X,\L,w)\coh_\lf)\rarrow
\sD^\abs((X,\L,w)\coh_\lf)$ and $\sD^\abs_T((X,\L,w)\coh)\rarrow
\sD^\abs((X,\L,w)\coh)$ are fully faithful~\cite{PV}.

 By the definition, the thick subcategories $\sD^\abs_T((X,\L,w)
\coh_\lf)\subset\sD^\abs((X,\L,w)\allowbreak\coh_\lf)$ and
$\sD^\abs_T((X,\L,w)\coh)\subset\sD^\abs((X,\L,w)\coh)$ only depend
on the intersection $X_0\cap T$ (rather than the whole of~$T$).
 Equivalently, they can be defined as the full subcategories of
objects annihilated by the restriction functors
$\sD^\abs((X,\L,w)\coh_\lf)\rarrow\sD^\abs((U,\L|_U,w|_U)\coh_\lf)$
and $\sD^\abs((X,\L,w)\coh)\rarrow\sD^\abs((U,\L|_U,w|_U)\coh)$,
where $U=X\setminus T$.

 As in Section~\ref{cdg-supports}, we denote by
$\sD^\abs((X,\L,w)\coh_T)$ the absolute derived category of
coherent matrix factorizations with the set-theoretic support in~$T$.
 The functor $\sD^\abs((X,\L,w)\coh_T)\rarrow\sD^\abs((X,\L,w)\coh)$
is fully faithful by Proposition~\ref{cdg-supports}(d).
 By Corollary~\ref{cdg-supports}(b), the full subcategory
$\sD^\abs_T((X,\L,w)\coh)\subset\sD^\abs((X,\L,w)\coh)$
is the thick envelope of the full subcategory
$\sD^\abs((X,\L,w)\coh_T)$.

 Now assume that $w\:\O_X\rarrow\L$ is an injective morphism
of sheaves.

\begin{prop}
\textup{(a)} The equivalence of categories\/ $\sD^\abs((X,\L,w)\coh)
\simeq\sD^\b_\Sing(X_0/X)$ identifies the triangulated subcategory\/
$\sD^\abs((X,\L,w)\coh_T)$ with the triangulated subcategory\/
$\sD^\b_\Sing(X_0/X\;X_0\cap T)$. 
 In particular, the former triangulated subcategory only depends on
the intersection $X_0\cap T$. \par
\textup{(b)} The full preimage of the thick envelope of
the triangulated subcategory\/ $\sD^\b_\Sing(X_0\;X_0\cap T)\subset
\sD^\b_\Sing(X_0)$ under the fully faithful functor\/
$\Sigma\:\sD^\abs((X,\L,w)\coh_\lf)\allowbreak\rarrow
\sD^\b_\Sing(X_0)$ coincides with the triangulated subcategory\/
$\sD^\abs_T((X,\L,w)\coh_\lf)$.
\end{prop}

\begin{proof}
 Part~(b) follows from the fact that the thick envelope of
$\sD^\b_\Sing(X_0\;X_0\cap T)$ is the kernel of the restriction
functor $\sD^\b_\Sing(X_0)\rarrow\sD^\b_\Sing(X_0\setminus T)$,
the similar fact for $\sD^\abs_T((X,\L,w)\coh_\lf)$, and
the compatibility of the functors~$\Sigma$ with the restrictions
to open subschemes, together with their full-and-faithfulness.

 To prove part~(a), notice first that the functor $\Ups$ obviously
takes $\sD^\b_\Sing(X_0/X\;\allowbreak X_0\cap\nobreak T)$ into
$\sD^\abs((X,\L,w)\coh_T)$.
 Let us check that the functor $\boL\Xi$ takes $\sD^\abs((X,\L,w)
\coh_T)$ into $\sD^\b_\Sing(X_0/X\;X_0\cap T)$.
 Let $\M$ be a coherent matrix factorization supported
set-theoretically in~$T$.
 Present $\M$ as the cokernel of an injective morphism of $w$\+flat
coherent matrix factorizations $\K\rarrow\N$.
 The functor $\boL\Xi$ being triangulated, the object
$\boL\Xi(\M)\in\sD^\b_\Sing(X_0/X)$ is isomorphic to the cone
of the morphism $\Xi(\K)\rarrow\Xi(\N)$ (cf.\
Lemma~\ref{finite-dim-matrix-push}).
 The morphism $\Xi(\K)\rarrow\Xi(\N)$ of coherent sheaves on $X_0$
is an isomorphism outside $T$, so its kernel and cokernel are
supported in $X_0\cap T$.
 Thus the cone is quasi-isomorphic to a two-term complex of coherent
sheaves on $X_0$ with the terms supported set-theoretically in
$X_0\cap T$.
\end{proof}

\subsection{Locality of local freeness}  \label{locality}
 The aim of this section is to show that the property of an object
of $\sD^\abs((X,\L,w)\qcoh_\fl)$ or $\sD^\abs((X,\L,w)\coh)$ to be
a direct summand of an object from $\sD^\abs((X,\L,w)\coh_\lf)$
is local in a separated Noetherian scheme $X$ with a dualizing
complex and enough vector bundles, assuming that the potential
$w\in\L(X)$ is not locally zero-dividing.

 Let $Z$ be a Noetherian scheme of finite Krull dimension with
enough vector bundles.
 Recall that the natural functor $\sD^\b_\Sing(Z)\rarrow
\sD'_\Sing(Z)$ is fully faithful~\cite[Proposition~1.13]{Or1}
(cf.\ Section~\ref{infinite-matrix}).

\begin{prop}
 Let $Z=U\cup V$ be a covering by two open subschemes.
 Then any object of\/ $\sD'_\Sing(Z)$ whose restrictions to $U$ 
and $V$ belong to the full subcategories\/ $\sD^\b_\Sing(U)\subset
\sD'_\Sing(U)$ and\/ $\sD^\b_\Sing(V)\subset\sD'_\Sing(V)$,
respectively, is a direct summand of an object belonging
to the full subcategory\/ $\sD^\b_\Sing(Z)\subset\sD'_\Sing(Z)$.
\end{prop}

\begin{proof}
 Consider the bounded derived category of quasi-coherent sheaves
$\sD^\b(Z\qcoh)$ on $Z$ and two full triangulated subcategories
$\sD^\b(Z\coh)$ and $\sD^\b(Z\qcoh_\fl)$ in it.
 Clearly, the intersection $\sD^\b(Z\coh)\cap\sD^\b(Z\qcoh_\fl)$
coincides with the full subcategory of perfect complexes
$\Perf(Z)=\sD^\b(Z\coh_\lf)\subset\sD^\b(Z\qcoh)$.

\begin{lem}
 Any morphism from an object of the full subcategory\/
$\sD^\b(Z\qcoh_\fl)$ into an object of the full subcategory\/
$\sD^\b(Z\coh)\subset\sD^\b(Z\qcoh)$ factorizes through an object
belonging to\/ $\sD^\b(Z\coh_\lf)$.
\end{lem}

\begin{proof}
 See the proof of~\cite[Proposition~1.13]{Or1}.
\end{proof}

 It follows from Lemma (by the way of the octahedron axiom) that
any object $\K^\bu$ of the full triangulated subcategory
$\sD^\b(Z\qcoh)_\flc$ generated by $\sD^\b(Z\qcoh_\fl)$ and
$\sD^\b(Z\coh)$ in $\sD^\b(Z\qcoh)$ can be included in
a distinguished triangle $\F^\bu\rarrow\K^\bu\rarrow\M^\bu\rarrow
\F^\bu[1]$ with $\F^\bu\in\sD^\b(Z\qcoh_\fl)$ and
$\M^\bu\in\sD^\b(Z\coh)$.
 Besides, the natural functor $\sD^\b(Z\qcoh_\fl)/\sD^\b(Z\coh_\lf)
\rarrow\sD^\b(Z\qcoh)/\sD^\b(Z\coh)$ is fully faithful.

 To prove Proposition, one has to show that any object $\K^\bu\in
\sD^\b(Z\qcoh)$ whose restrictions to $U$ and $V$ belong to
the subcategories $\sD^\b(U\qcoh)_\flc$ and $\sD^\b(V\qcoh_\flc)$,
respectively, is a direct summand of an object from
$\sD^\b(Z\qcoh)_\flc\subset\sD^\b(Z\qcoh)$.
 According to the above, there exist two objects $\F_U^\bu\in
\sD^\b(U\qcoh_\fl)$ and $\F_V^\bu\in\sD^\b(V\qcoh_\fl)$ and
two morphisms $\F_U^\bu\rarrow\K^\bu|_U$ and $\F_V^\bu\rarrow\K^\bu|_V$
whose cones belong to $\sD^\b(U\coh)$ and $\sD^\b(V\coh)$,
respectively.

 Set $W=U\cap V\subset Z$; then the restrictions of $\F_U^\bu$ and
$\F_V^\bu$ to $W$ are isomorphic in $\sD^\b(W\qcoh)/\sD^\b(W\coh)$,
and consequently, in $\sD^\b(W\qcoh_\fl)/\sD^\b(W\coh_\lf)$, too.
 Notice that the category $\Perf(W)=\sD^\b(W\coh_\lf)$ is idempotent
complete, and therefore, a thick subcategory in $\sD^\b(W\qcoh_\fl)$.
 It follows that there exists a finite complex of flat quasi-coherent
sheaves $\F_W^\bu$ on $W$ together with two morphisms
$\F_U^\bu|_W\rarrow\F_W^\bu$ and $\F_V^\bu|_W\rarrow\F_W^\bu$
whose cones are perfect complexes.
 Denote the cocones of these morphisms by $\G_W^\bu$ and $\H_W^\bu$.

 For any object $A$ of a triangulated category $\sD$, let us denote by
${}'\!A$ the object $A\oplus A[1]$.
 For any triangulated subcategory $\sC\subset\sD$, whenever 
an object $A\in\sD$ is a direct summand of an object from $\sC$,
the object ${}'\!A$ belongs to $\sC$, as $A\oplus B\in\sC$ implies
$A\oplus A[1]\in\sC$ in view of the distinguished triangle
$A\oplus B\rarrow A\oplus B\rarrow A\oplus A[1]\rarrow
A[1]\oplus B[1]$ \cite[Theorem~2.1]{Th}.

 By the Thomason--Trobaugh theorem~\cite[Section~5]{TT}, the objects
${}'\G_W^\bu$ and ${}'\H_W^\bu$ can be extended to perfect
complexes on $U$ and $V$, respectively.
 Moreover, these extensions $\G_U^\bu\in\sD^\b(U\coh_\lf)$ and
$\H_V^\bu\in\sD^\b(V\coh_\lf)$ can be chosen in such a way that
the morphisms ${}'\G_W^\bu\rarrow{}'\!\.\F_U^\bu|_W$ and
${}'\H_W^\bu\rarrow{}'\!\.\F_V^\bu|_W$ would be extendable to morphisms
$\G_U^\bu\rarrow{}'\!\.\F_U^\bu$ and $\H_V^\bu\rarrow{}'\!\.\F_V^\bu$
\cite[Theorem~2.1(4\+5)]{Neem}.

 Furthermore, the objects ${}'\G_U^\bu$ and ${}'\H_V^\bu$ can be
extended to perfect complexes $\G^\bu$ and $\H^\bu$ on the whole
scheme $Z$ so that the compositions of morphisms ${}'\G_U^\bu\rarrow
{}''\!\.\F_U^\bu\rarrow{}''\K^\bu|_U$ and ${}'\H_V^\bu\rarrow
{}''\!\.\F_V^\bu\rarrow{}''\K^\bu|_V$ would be extendable to morphisms
$\G^\bu\rarrow{}''\K^\bu$ and $\H^\bu\rarrow{}''\K^\bu$.
 Denote by $\K^\bu_{(1)}$ a cone of the morphism $\G^\bu\oplus\H^\bu
\rarrow{}''\K^\bu$, by $\F^\bu_{U,(1)}$ a cone of the morphism
${}'\G_U^\bu\rarrow{}''\!\.\F_U^\bu$, and by $\F^\bu_{V,(1)}$
a cone of the morphism ${}'\H_V^\bu\rarrow{}''\!\.\F_V^\bu$.
 We have come back to the original situation with an object
$\K_{(1)}^\bu\in\sD^\b(Z\qcoh)$, two objects $\F_{U,(1)}^\bu\in
\sD^\b(U\qcoh_\fl)$ and $\F_{V,(1)}^\bu\in\sD^\b(V\qcoh_\fl)$, and
two morphisms $\F_{U,(1)}^\bu\rarrow\K_{(1)}^\bu|_U$ and
$\F_{V,(1)}^\bu\rarrow\K_{(1)}^\bu|_V$ whose cones belong to
$\sD^\b(U\coh)$ and $\sD^\b(V\coh)$, respectively.
 In addition, the objects $\F_{U,(1)}^\bu|_W$ and $\F_{V,(1)}^\bu|_W$
are now isomorphic in $\sD^\b(W\qcoh_\fl)$.

 The construction does not guarantee commutativity of the diagram
formed by the isomorphism $\F_{U,(1)}^\bu|_W=\F_{W,(1)}^\bu\simeq
\F_{V,(1)}^\bu|_W$ and the restrictions of the morphisms
$\F_{U,(1)}^\bu\rarrow\K_{(1)}^\bu$ and $\F_{V,(1)}^\bu\rarrow
\K_{(1)}^\bu$ to~$W$.
 However, the original choice of the morphisms $\F_U^\bu|_W\rarrow
\F_W^\bu$ and $\F_V^\bu|_W\rarrow\F_W^\bu$ makes this diagram
commute in the quotient category $\sD^b(W\qcoh)/\sD^\b(W\coh)$.
 Hence the difference of two morphisms $\F_{W,(1)}^\bu\rightrightarrows
\K_{(1)}^\bu|_W$ factorizes through a bounded complex of coherent
sheaves on $W$, and consequently (according to Lemma) also
through a perfect complex on~$W$.
 Denote the latter by $\E^\bu\in\sD^\b(W\coh_\lf)$.

 Now let $j\:U\rarrow Z$, \ $k\:V\rarrow Z$, and $h\:W\rarrow Z$
denote the natural open embeddings.
 Consider the square diagram formed by the morphisms
$\boR j_*\F_{U,(1)}^\bu\oplus \boR k_*\F_{V,(1)}^\bu\rarrow
\boR h_*\F_{U,(1)}^\bu|_W$ and $\boR j_*\K_{(1)}^\bu|_U\oplus
\boR k_*\K_{(1)}^\bu|_V\rarrow \boR h_*\K_{(1)}^\bu|_W$.
 According to the above, this diagram is not necessarily commutative;
but it can be made commutative by adding the new direct summand
$\boR h_*\E^\bu$ to the term $\boR j_*\K_{(1)}^\bu|_U\oplus
\boR k_*\K_{(1)}^\bu|_V$ with the morphism $\boR h_*\E^\bu\rarrow
\boR h_*\K_{(1)}^\bu|_W$ induced by the morphism $\E^\bu\rarrow
\K_{(1)}^\bu|_W$ and the morphism $\boR j_*\F_{U,(1)}^\bu\oplus
\boR k_*\F_{V,(1)}^\bu$ equal to zero on the first direct summand
and induced by the morphism $\F_{V,(1)}^\bu|_W\simeq \F_{W,(1)}^\bu
\rarrow\E^\bu$ on the second one.

 Let $\F^\bu$ denote a cocone of the morphism $\boR j_*\F_{U,(1)}^\bu
\oplus \boR k_*\F_{V,(1)}^\bu\rarrow\boR h_*\F_{U,(1)}^\bu|_W$ and
$\L^\bu$ denote a cocone of the morphism $\boR j_*\K_{(1)}^\bu|_U
\oplus\boR k_*\K_{(1)}^\bu|_V\oplus\boR h_*\E^\bu\rarrow
\boR h_*\K_{(1)}^\bu|_W$.
 Then the commutative square can be extended to a morphism of
distinguished triangles, so we obtain a morphism $\F^\bu\rarrow\L^\bu$.
 Since $\K_{(1)}^\bu$ is a cocone of the morphism
$\boR j_*\K_{(1)}^\bu|_U\oplus\boR k_*\K_{(1)}^\bu|_V\rarrow
\boR h_*\K_{(1)}^\bu|_W$, there is also a distinguished triangle
$\K_{(1)}^\bu\rarrow\L^\bu\rarrow\boR h_*\E\rarrow\K_{(1)}^\bu[1]$.

 Notice that the complexes $\F^\bu$ and $\boR h_*\E^\bu$ belong to
$\sD^\b(Z\qcoh_\fl)$ (since the class of bounded complexes of
flat quasi-coherent sheaves is preserved by the derived direct images
with respect to flat morphisms of Noetherian schemes;
cf.\ Proposition~\ref{finite-dim-morphisms}).
 Furthermore, the complex $\boR h_*\E^\bu$ is perfect over~$W$.
 Restricting to $W$ our morphism of distinguished triangles, and
recalling that cones of the morphisms $\F_{U,(1)}^\bu\rarrow
\K_{(1)}^\bu|_U$ and $\F_{V,(1)}^\bu\rarrow\K_{(1)}^\bu|_V$ are
coherent complexes over $U$ and $V$, one easily concludes that
a cone of the morphism $\F^\bu\rarrow\L^\bu$ is a coherent complex
over~$W$.

 Denote this cone temporarily by~$\K_{(2)}^\bu$.
 Clearly, in order to show that the original complex $\K^\bu$ is
a direct summand of an object from $\sD^\b(Z\qcoh)_\flc$ in
$\sD^\b(Z\qcoh)$ (which is our goal) it suffices to check that so
is the complex $\K_{(2)}^\bu$.
 It also follows from the constructions that the restrictions of
the complex $\K_{(2)}^\bu$ to $U$ and $V$ belong to
$\sD^\b(U\qcoh)_\flc$ and $\sD^\b(V\qcoh_\flc)$, respectively.
 Dropping the lower index and redenoting $\K_{(2)}^\bu$ simply
by $\K^\bu$, we are coming back to the situation in the beginning
of the proof with the new knowledge that $\K^\bu$ may be assumed
to be a coherent complex over~$W$.

 The next fragment of our proof is based on the localization theory
for coderived categories of quasi-coherent sheaves on Noetherian
schemes (similar to the Thomason--Trobaugh--Neeman theory for
the conventional derived categories, the difference being that
arbitrary bounded complexes of coherent sheaves play the role of
perfect complexes).
 What we need is a particular case of the theory developed in
Section~\ref{cdg-supports} (corresponding to the choice of 
the quasi-coherent CDG\+algebra $\O_Z$ over~$Z$).

 Specifically, it follows from Proposition~\ref{embedding-prop}(d)
and Theorem~\ref{cdg-supports} together
with~\cite[Theorem~2.1(5)]{Neem} that any morphism from
an object of $\sD^\b(W\coh)$ into a restriction to $W$ of
an object $\K^\bu$ from $\sD^\b(Z\qcoh)$ (or even from
$\sD^\co(Z\qcoh)$) can be extended to a morphism to $\K^\bu$
from an object of $\sD^\b(Z\coh)$.
 Applying this assertion to the identity morphism $\K^\bu|_W\rarrow
\K^\bu|_W$ in the above situation, we obtain a morphism $\M^\bu
\rarrow\K^\bu$ into $\K^\bu$ from a coherent complex $\M^\bu$
over $Z$ that is a quasi-isomorphism over~$W$.
 Passing to a cone of this morphism, we may assume $\K^\bu$ to be
acyclic over~$W$.

 By Corollary~\ref{cdg-supports}, such a complex $\K^\bu$ is
quasi-isomorphic to a (bounded) complex of quasi-coherent sheaves
on $Z$ whose terms are concentrated set-theoretically in
the complement $Z\setminus W$.
 The latter is a disjoint union of two nonintersecting closed
subsets in $Z$, namely, the complements $S=Z\setminus U$ and
$T=Z\setminus V$.
 Now the complex $\K^\bu$ decomposes into a direct sum of two
complexes with the set-theoretic supports inside $S$ and $T$,
respectively.

 One can consider the two direct summands separately.
 We have to show that any bounded complex of quasi-coherent sheaves
$\K^\bu$ on $Z$, which is supported set-theoretically in $T$ and
whose restriction to $U$ belongs to $\sD^\b(U\qcoh)_\flc$, itself
belongs to $\sD^\b(Z\qcoh)_\flc$.
 Arguing as in the beginning of this proof, we have an object
$\G^\bu\in\sD^\b(U\qcoh_\fl)$ together with a morphism
$\G^\bu\rarrow\K^\bu|_U$ whose cone belongs to $\sD^\b(U\coh)$.
 The restriction $\G^\bu|_W$ then belongs to both $\sD^\b(W\qcoh_\fl)$
and $\sD^\b(W\coh)$, and is, therefore, a perfect complex on~$W$.

 Again by the Thomason--Trobaugh theorem, the object ${}'\G^\bu|_W$
can be extended to a perfect complex $\H^\bu$ on~$V$.
 A cocone of the morphism $\boR j_*{}'\G^\bu\oplus\boR k_*\H^\bu
\rarrow\boR h_*{}'\G^\bu|_W$ provides an object
$\F^\bu\in\sD^\b(Z\qcoh_\fl)$ isomorphic to ${}'\G^\bu$ over $U$
and to $\H^\bu$ over~$V$.
 Now the morphism ${}'\G^\bu\rarrow{}'\K^\bu|_U$ over $U$ extends
uniquely to a morphism $\F^\bu\rarrow{}'\K^\bu$ over $Z$, since
the set-theoretic support of ${}'\K^\bu$ is contained in a closed
subset lying inside~$U$.
 A cone of the morphism $\F^\bu\rarrow{}'\K^\bu$ is a coherent
complex on $Z$, since it is so in restrictions to $U$ and~$V$.
 Proposition is proven.
\end{proof}

 Now let $X$ be a separated Noetherian scheme of finite Krull
dimension with enough vector bundles, $\L$ be a line bundle on $X$,
and $w\in\L(X)$ be a locally nonzero-dividing potential.
 Let $X_0\subset X$ be the zero locus of~$w$.

\begin{cor}
 Let $X=U\cap V$ be a covering by two open subschemes.
 Then any object of\/ $\sD^\co((X,\L,w)\qcoh_\fl)$ whose
restrictions to $U$ and $V$ belong to the full triangulated
subcategories\/ $\sD^\abs((U,\L|_U,w|_U)\coh_\lf)\subset
\sD^\co((U,\L|_U,w|_U)\qcoh_\fl)$ and\/
$\sD^\abs((V,\L|_V,w|_V)\coh_\lf)\subset
\sD^\co((V,\L|_V,w|_V)\qcoh_\fl)$, respectively, is a direct summand
of an object from the full triangulated subcategory\/
$\sD^\abs((X,\L,w)\coh_\lf)\subset\sD^\co((X,\L,w)\qcoh_\fl)$.
\end{cor}

\begin{proof}
 By Proposition~\ref{infinite-matrix}, the category
$\sD^\co((X,\L,w)\qcoh_\fl)$ is a full triangulated subcategory
of the triangulated category $\sD'_\Sing(X_0)$.
 The (essential) intersection of the full subcategories
$\sD^\co((X,\L,w)\qcoh_\fl)$ and $\sD^\b_\Sing(X_0)$ in
$\sD'_\Sing(X_0)$ is the triangulated category
$\sD^\abs((X,\L,w)\coh_\lf)$.

 Indeed, an object of $\F\in\sD^\b_\Sing(X_0)$ belongs to
$\sD^\abs((X,\L,w)\coh_\lf)$ if and only if the object
$i_\circ\F$ vanishes in $\sD^\b_\Sing(X)$ (Theorem~\ref{main-theorem});
an object $\F\in\sD'_\Sing(X_0)$ belongs to $\sD^\co((X,\L,w)\qcoh_\fl)$
if and only if the object $i_\circ\F$ vanishes in $\sD'_\Sing(X)$
(Proposition~\ref{infinite-matrix}); and the functor
$\sD^\b_\Sing(X)\rarrow\sD'_\Sing(X)$ is fully faithful.

 Moreover, the (essential) intersection of $\sD^\co((X,\L,w)\qcoh_\fl)$
with the thick envelope of $\sD^\b_\Sing(X_0)$ in $\sD'_\Sing(X_0)$
is the thick envelope of $\sD^\abs((X,\L,w)\coh_\lf)$
in $\sD'_\Sing(X_0)$.
 Indeed, let $\M$ be an object of the intersection; then
$\M\oplus\M[1]$ belongs to both $\sD^\co((X,\L,w)\qcoh_\fl)$ and
$\sD^\b_\Sing(X_0)$, hence also to $\sD^\abs((X,\L,w)\coh_\lf)$,
and consequently $\M$ belongs to the thick envelope of
$\sD^\abs((X,\L,w)\coh_\lf)$.

 Now let $\K$ be our object of $\sD^\co((X,\L,w)\qcoh_\fl)$;
it can be also viewed as an object of $\sD'_\Sing(X_0)$.
 If its restrictions to $U$ and $V$ belong to
$\sD^\abs((U,\L|_U,w|_U)\coh_\lf)$ and $\sD^\abs((V,\L|_V,w|_V)
\coh_\lf)$, they also belong to $\sD^\b_\Sing(U_0)\subset
\sD'_\Sing(U_0)$ and $\sD^\b_\Sing(V_0)\subset\sD'_\Sing(V_0)$
(where we set $U_0=U\cap X_0$ and $V_0=V\cap X_0$).

 Applying Proposition, we can conclude that $\K$ belongs to
the thick envelope of $\sD^\b_\Sing(X_0)$ in $\sD'_\Sing(X_0)$.
 The assertion of Corollary follows from the above.
\end{proof}

 Assume additionally that the scheme $X$ admits a dualizing
complex~$\D_X^\bu$.

\begin{thm}
 Let $X=U\cap V$ be a covering by two open subschemes.
 Then any object of\/ $\sD^\abs((X,\L,w)\coh)$ whose restrictions
to $U$ and $V$ belong to the thick envelopes of the
triangulated subcategories\/ $\sD^\abs((U,\L|_U,w|_U)\coh_\lf)\subset
\sD^\abs((U,\L|_U,w|_U)\coh)$ and\/ $\sD^\abs((V,\L|_V,w|_V)\coh_\lf)
\subset\sD^\abs((V,\L_V,w|_V)\coh)$ itself belongs to the thick
envelope of the triangulated subcategory\/
$\sD^\abs((X,\L,w)\coh_\lf)\subset\sD^\abs((X,\L,w)\coh)$.
\end{thm}

\begin{proof}
 The argument is based on the Serre--Grothendieck duality theory
for matrix factorizations as developed in Section~\ref{serre-duality},
which allows to reduce the question to the result of Corollary.
 Specifically, let $\M$ be our coherent matrix factorization over~$X$.
 Replacing, if necessary, $\M$ with $\M\oplus\M[1]$, we may assume
the restrictions of $\M$ to $U$ and $V$ to be isomorphic to locally
free matrix factorizations of finite rank.

 Let us apply the construction of functor $\Omega\:\sD^\abs((X,\L,w)
\coh)^\sop\rarrow\sD^\co((X,\L,\allowbreak-w)\qcoh_\fl)$ from
Section~\ref{serre-duality} to the matrix factorization~$\M$.
 That is, we pick a left resolution of $\M$ by locally free matrix
factorizations of finite rank, dualize by applying
$\Hom_{X\qc}({-},\O_X)$, and totalize using infinite direct sums.
 By Corollary~\ref{serre-duality}, the functor $\Omega$ is fully
faithful; it also identifies $\sD^\abs((X,\L,w)\coh_\lf)^\sop$ with
$\sD^\abs((X,\L,-w)\coh_\lf)$.
 Hence it suffices to check that the matrix factorization $\Omega(\M)$
belongs to the thick envelope of $\sD^\abs((X,\L,-w)\coh_\lf)$ in
$\sD^\co((X,\L,-w)\qcoh_\lf)$.
 But we know as much from the above Corollary. 
\end{proof}

\subsection{Nonlocalization of local freeness} \label{nonlocalization}
 The lack of a workable notion of the conventional derived category
(as opposed to the coderived category) for quasi-coherent matrix
factorizations stands in the way of a direct extension of
the Thomason--Trobaugh--Neeman localization theory for perfect
complexes~\cite{TT,Neem0,Neem} to locally free matrix factorizations
of finite rank.
 We have seen in Section~\ref{cdg-supports} how the localization theory
can be developed for coherent matrix factorizations.
 In this section we demonstrate a counterexample showing that
the localization theory, in its conventional form, actually does
\emph{not} hold for locally free matrix factorizations.

 In other words, the restriction $\sD^\abs((X,\L,w)\coh_\lf)\rarrow
\sD^\abs((U,\L|_U,w|_U)\coh_\lf)$ for an open subscheme $U\subset X$
is not always a Verdier quotient functor, even up to
the direct summands.
 Moreover, the triangulated category $\sD^\abs((X,\L,w)\coh_\lf)$
may fail to be generated by a single object, unlike in the case of
the categories of perfect complexes on quasi-compact
quasi-separated schemes.

 All the potentials in our example will be simply regular functions,
i.~e., sections of the trivial line bundle~$\O_X$ or $\O_U$, etc.;
so we drop the line bundle $\L$ from our notation in the rest of
the section and write simply $\sD^\abs((X,w)\coh_\lf)$ or
$\sD^\abs((X,w)\coh)$, etc.
 For simplicity, we will work over the basic field of complex
numbers~$\boC$.

 Consider the $3$\+dimensional affine quadratic cone
$$
 X=\{xy=zw\}\.\subset\.\mathbb A^4=\Spec \boC[x,y,z,w].
$$
 Further, let us take the open subset
$$
 U=\{z\ne 0\}\subset X.
$$
 Clearly, we have an isomorphism of pairs (algebraic variety,
regular function on it)
\begin{equation} \label{eq:iso_for_U,w}
 (U,w)\stackrel{\sim}\rarrow
 (\mathbb A^2_{t_1,t_2}\times \mathbb G_m\;t_1t_2),
 \quad (x,y,z,w)\longmapsto ((x,\frac{y}{z}),z),
\end{equation}
where we denote $\mathbb A^2_{t_1,t_2}=\Spec\boC[t_1,t_2]$ and,
as usually, $\mathbb G_m=\mathbb A^1\setminus\{0\}$.

\begin{lem}
\textup{(a)} We have a natural equivalence of triangulated categories
$$
 \sD^\abs((U,w)\coh)\simeq\sD^\abs((\mathbb G_m,0)\coh).
$$
\textup{(b)} The restriction functor
$$
 \sD^\abs((X,w)\coh)\rarrow\sD^\abs((U,w)\coh)
$$ is an equivalence.
\end{lem}

 Here the category of matrix factorizations of the zero potential
$\sD^\abs((Y,0)\coh)$ is, of course, simply the derived
category of $2$\+periodic complexes of coherent sheaves
on a smooth variety~$Y$.

\begin{proof}
Part~(a): by~\eqref{eq:iso_for_U,w}, we have an equivalence
$$
 \sD^\abs((U,w)\coh)\simeq
 \sD^\abs((\mathbb A^2_{t_1,t_2}\times \mathbb G_m,t_1t_2)\coh).
$$
By Kn\"orrer periodicity (cf.~\cite[Theorem~3.1]{Or1a}),
we have an equivalence
$$
 \sD^\abs((\mathbb A^2_{t_1,t_2}\times\mathbb G_m,t_1t_2)\coh)
 \simeq\sD^\abs((\mathbb G_m,0)\coh).
$$

Part~(b): let us put $D=X\setminus U$.
 By Theorem~\ref{cdg-supports}(b) (see also
Section~\ref{supports-secn}), we have a short exact sequence
of triangulated categories
$$
 0\rarrow \sD^\abs_D((X,w)\coh)\rarrow\sD^\abs((X,w)\coh)\rarrow
 \sD^\abs((U,w)\coh)\rarrow 0
$$
 Thus, we need to show that the category $\sD^\abs_D((X,w))\coh)$
is zero.
 It suffices to check that the category $\sD^\abs((D,w)\coh)$ is zero.

 Let us put $S=\{xy=0\}\subset \mathbb A^2$.
 Then we have an isomorphism
$$
 (D,w)\stackrel{\sim}\rarrow(S\times\mathbb A^1_t\;t),
 \quad (x,y,0,w)\longmapsto ((x,y),w).
$$
Since $\sD^\abs((\mathbb A^1_t,t)\coh)=0,$ it follows that
$$
 \sD^\abs((D,w)\coh)\simeq \sD^\abs((S\times\mathbb A^1_t,t)\coh)=0.
$$
\end{proof}

 Since $U$ is smooth, we have an equivalence
$$
 \sD^\abs((U,w)\coh_\lf)\simeq \sD^\abs((U,w)\coh).
$$
 Now we turn to the category $\sD^\abs((X,w)\coh_\lf)$.
 As usually, we put
$$
 X_0=\{w=0\}\.\subset\.X.
$$
 According to Theorem~\ref{main-theorem}, the triangulated category
$\sD^\abs((X,w)\coh_\lf)$ is equivalent to the kernel of the direct
image functor $i_\circ\:\sD^\b_\Sing(X_0)\rarrow\sD^\b_\Sing(X)$ 
acting between the triangulated categories of singularities of
the schemes $X_0$ and~$X$.
 This can be rephrased by saying that $\sD^\abs((X,w)\coh_\lf)$ is
equivalent to the quotient category of the category of bounded
complexes of coherent sheaves on $X_0$ whose direct images are
perfect complexes on $X$ by the category of perfect complexes on~$X_0$.
 Denoting the triangulated category of coherent complexes on $X_0$
whose direct images are perfect on $X$ by $\Perf(X_0,X)\subset
\sD^\b(X_0\coh)$, we have an equivalence of triangulated categories
\begin{equation} \label{eq:equiv_for_D^abs_l.f.(X,w)}
 \sD^\abs((X,w)\coh_\lf)\simeq \Perf(X_0,X)/\Perf(X_0).
\end{equation}
 Note that we have a natural isomorphism
$$
 X_0\simeq S\times \mathbb A^1,\quad (x,y,z,0)\longmapsto ((x,y),z).
$$
 It follows immediately that
\begin{equation} \label{eq:equiv_for_D_sg(X^0)}
 \sD^\b_\Sing(X_0)\simeq \sD^\abs((\mathbb A^1,0)\coh).
\end{equation}

\begin{prop}
\textup{(a)} We have a natural equivalence of triangulated categories
$$
 \Perf(X_0,X)/\Perf(X_0)\simeq\sD^\abs((\mathbb G_m,0)\coh)_\zdim,
$$
where\/ $\sD^\abs((\mathbb G_m,0)\coh)_\zdim\subset
\sD^\abs((\mathbb G_m,0)\coh)$ is the subcategory of complexes with
zero-dimensional support.

 Moreover, we have a commutative diagram of fully faithful
triangulated functors
\dgARROWLENGTH=2em
$$
\begin{diagram}
\node{\Perf(X_0,X)/\Perf(X_0)}\arrow[3]{e,V} \arrow{s,=}
\node[3]{\sD^\b_\Sing(X_0)} \arrow{s,=} \\
\node{\sD^\abs((\mathbb G_m,0)\coh)_\zdim} \arrow[3]{e,t,V}{j_*}
\node[3]{\sD^\abs((\mathbb A^1,0)\coh),}
\end{diagram}
$$
where $j:\mathbb G_m\rarrow\mathbb A^1$ is the open embedding.

\textup{(b)} We have a commutative diagram of fully faithful
triangulated functors and equivalences
\dgARROWLENGTH=2em
$$
\begin{diagram}
\node{\sD^\abs((X,w)\coh_\lf)}\arrow[3]{e,V} \arrow{s,=}
\node[3]{\sD^\abs((U,w)\coh_\lf)} \arrow{s,=} \\
\node{\sD^\abs((\mathbb G_m,0)\coh)_\zdim} \arrow[3]{e,t,V}{\iota}
\node[3]{\sD^\abs((\mathbb G_m,0)\coh),}
\end{diagram}
$$
where $\iota$~is the tautological embedding.
\end{prop}

\begin{proof}
 Part~(a): indeed, from the equivalence~\eqref{eq:equiv_for_D_sg(X^0)}
we have a natural fully faithful triangulated functor
$$
 \Perf(X_0,X)/\Perf(X_0)\rarrow\sD^\abs((\mathbb A^1,0)\coh).
$$
 Let us denote by $\sT\subset\sD^\abs((\mathbb A^1,0)\coh)$
the essential image of this functor.
 For each $z_0\in\boC\setminus\{0\},$ we have a line
$l_{z_0}:=\{y=0\;z=z_0\}\subset X_0.$
 Since $l_{z_0}\subset U$ and $U$ is smooth, the coherent sheaf
$\O_{l_{z_0}}$ is contained in $\Perf(X_0,X).$
 Further, its image in $D^\b_\Sing(X_0)$ corresponds to
the skyscraper $\O_{z_0}\in\sD^\abs((\mathbb A^1,0)\coh)$
under the equivalence~\eqref{eq:equiv_for_D_sg(X^0)}.
 It follows that the triangulated category $\sT$ contains
$j_*(\sD^\abs((\mathbb G_m,0)\coh)_\zdim)$ as a full subcategory.

 Suppose that $\sT$ is strictly bigger than
$j_*(\sD^\abs((\mathbb G_m,0)\coh)_\zdim)$.
Then it contains an object $\F_0=\O_0\oplus\O_0[1]\in
\sD^\abs((\mathbb A^1,0)\coh)$, where $\O_0$ is the structure
sheaf of the origin.
 Denote by $O\in X_0$ the origin $(0,0,0,0)$.
 Then the image of the coherent sheaf $\O_O\in X_0\coh$ in
$\sD^\b_\Sing(X_0)$ corresponds to $\F_0$ under
the equivalence~\eqref{eq:equiv_for_D_sg(X^0)}.
 But the object $\O_O\in \sD^\b(X_0\coh)$ is not relatively perfect
under the inclusion $X_0\rarrow X$ (i.~e., it does not belong to
$\Perf(X_0,X)$), since $O$ is the singular point of $X$.
 We get a contradiction.

 Therefore, we have an equivalence
$\sT\simeq j_*(\sD^\abs((\mathbb G_m,0)\coh)_\zdim)$.
 This proves~(a).

 Part~(b) follows immediately from part~(a) and
the equivalence \eqref{eq:equiv_for_D^abs_l.f.(X,w)}.
\end{proof}

 In particular, we see that the functor 
$\sD^\abs((X,w)\coh_\lf)\rarrow\sD^\abs((U,w)\coh_\lf)$
is not essentially surjective, even up to the direct summands.
 Moreover, the triangulated category $\sD^\abs((X,w)\coh_\lf)$
does not even have a countable set of generators.

\subsection{Pull-backs and push-forwards in singularity categories}
\label{pull-push-singularities}
 Let $f\:Y\rarrow X$ be a morphism of separated Noetherian schemes
with enough vector bundles.
 The morphism $f$~is said to have \emph{finite flat dimension}
if the derived inverse image functor $\boL f^*\:\sD^-(X\qcoh)
\rarrow\sD^-(Y\qcoh)$ takes $\sD^\b(X\qcoh)$ to $\sD^\b(Y\qcoh)$.

 In this case, the functor $\boL f^*$ induces the inverse image
functors on the triangulated categories of singularities
\begin{align*}
 f^\circ\:\sD'_\Sing(X)&\lrarrow\sD'_\Sing(Y) \\
 f^\circ\:\sD^\b_\Sing(X)&\lrarrow\sD^\b_\Sing(Y).
\end{align*}

 Under the same assumption of finite flat dimension, the derived
direct image functor $\boR f_*\:\sD^\b(Y\qcoh)\rarrow\sD^\b(X\qcoh)$
takes $\sD^\b(Y\qcoh_\fl)$ to $\sD^\b(X\qcoh_\fl)$, as one can see by
computing $\boR f_*$ in terms of an affine covering of~$Y$ in
the spirit of the proof of Proposition~\ref{finite-dim-morphisms}.
 When the scheme $X$ has finite Krull dimension, one has
$\sD^\b(X\qcoh_\fl)=\sD^\b(X\qcoh_\lf)$, so the functor $\boR f_*$
induces the direct image functor
$$
 f_\circ\:\sD'_\Sing(Y)\lrarrow\sD'_\Sing(X),
$$
which is right adjoint to~$f^\circ$.

 Whenever the morphism $f$~is proper of finite type and has finite
flat dimension, the functor $\boR f_*$ takes $\sD^\b(Y\coh)$ to
$\sD^\b(X\coh)$ \cite[Th\'eor\`eme~3.2.1]{Groth} and induces
the direct image functor
$$
 f_\circ\:\sD^\b_\Sing(Y)\lrarrow\sD^\b_\Sing(X),
$$
which is right adjoint to~$f^\circ$
\cite[paragraphs before Proposition~1.14]{Or1}.
 More generally, for a morphism $f$ of finite flat dimension and
any closed subset $T\subset Y$ such that (a closed subscheme
structure on) $T$ is proper of finite type over $X$, the functor
$\boR f_*$ takes $\sD^\b(Y\coh_T)$ to $\sD^\b(X\coh)$ and
induces the direct image functor
$$
 f_\circ\:\sD^\b_\Sing(Y,T)\lrarrow\sD^\b_\Sing(X).
$$
 Indeed, the intersection of $\sD^\b(X\qcoh_\fl)$ and $\sD^\b(X\coh)$
in $\sD^\b(X\qcoh)$ is equal to $\sD^\b(X\coh_\lf)$, as any complex
of finite flat dimension with bounded coherent cohomology is
easily seen to be perfect.

 Let $Z\subset X$ and $W\subset Y$ be closed subschemes such that
$\O_Z$ is a perfect $\O_X$\+module, $\O_W$ is a perfect
$\O_Y$\+module, and $f(W)\subset Z$.
 Assume that both morphisms $f\:Y\rarrow X$ and $f|_W\:W\rarrow Z$
have finite flat dimensions.
 Then the derived inverse image functor $\boL f|_W^*\:\sD^\b(Z\qcoh)
\rarrow\sD^\b(W\qcoh)$ induces the inverse image functors on
the triangulated categories of relative singularities
\begin{align*}
 f^\circ\:\sD'_\Sing(Z/X)&\lrarrow\sD'_\Sing(W/Y) \\
 f^\circ\:\sD^\b_\Sing(Z/X)&\lrarrow\sD^\b_\Sing(W/Y).
\end{align*}

 Now let $Z\subset X$ be a closed subscheme; set $W=Z\times_XY$.
 Denote the closed embeddings $Z\rarrow X$ and $W\rarrow Y$ by
$i$ and $i'$, respectively; let also $f'$ denote the morphism
$f|_W\:W\rarrow Z$.
 Assume that $W$ coincides with the derived product of $Z$ and $Y$
over $X$, i.~e., $\boL f^*i_*\O_Z=i'_*\O_W$.
 Assume further that $i_*\O_Z$ is a perfect $\O_X$\+module; then
also $i'_*\O_W$ is a perfect $\O_Y$\+module.

 For any $\M\in\sD^\b(Y\qcoh)$ there is a natural morphism
$\phi_\M\:\boL i^*\boR f_*\M\rarrow\boR f'_*\.\boL i'{}^*\M$
in $\sD^\b(Z\qcoh)$.
 Using the projection formula for tensor products with perfect
complexes, one easily checks that the morphism $i_*\phi_\M$ is
an isomorphism.
 Hence so is the morphism $\phi_\M$, since the functor $i_*$ does not
annihilate any objects of the derived category.
 Hence we obtain the induced functor of direct image
$$
 f_\circ\:\sD'_\Sing(W/Y)\lrarrow\sD'_\Sing(Z/X).
$$
 When the morphism~$f$ is proper of finite type, there is also
the induced functor
$$
 f_\circ\:\sD^\b_\Sing(W/Y)\lrarrow\sD^\b_\Sing(Z/X).
$$

 Assume additionally that the morphism~$f$ has finite flat dimension;
then so does the morphism~$f'$.
 In this case the functor $f_\circ\:\sD'_\Sing(W/Y)\rarrow
\sD'_\Sing(Z/X)$ is right adjoint to the functor 
$f^\circ\:\sD'_\Sing(Z/X)\rarrow\sD'_\Sing(W/Y)$.
 When the morphism~$f$ is proper of finite type, the functor
$f_\circ\:\sD^\b_\Sing(W/Y)\rarrow\sD^\b_\Sing(Z/X)$ is right adjoint
to the functor $f^\circ\:\sD^\b_\Sing(Z/X)\rarrow\sD^\b_\Sing(W/Y)$.

\begin{rem}
 In the case when $Z$ is a Cartier divisor in $X$, we will construct
the functor $f_\circ\:\sD^\b_\Sing(W/Y)\rarrow\sD^\b_\Sing(Z/X)$
under somewhat weaker assumptions below in Section~\ref{matrix-push}.
 Namely, it will suffice that the morphism $f'\:W\rarrow Z$ be proper
of finite type, while the morphism $f\:Y\rarrow Z$ need not be. 
 A generalization to the case of proper support will also be obtained.
\end{rem}

\subsection{Push-forwards of matrix factorizations}
\label{matrix-push}
 Let $f\:Y\rarrow X$ be a morphism of separated Noetherian schemes
with enough vector bundles, $\L$ be a line bundle on $X$, and
$w\in\L(X)$ be a section.

 Set $\B_X=(X,\L,w)$ and $\B_Y=(Y,f^*\L,f^*w)$; then there is
a natural morphism of CDG\+algebras $\B_X\rarrow\B_Y$ compatible
with the morphism of schemes $f\:Y\rarrow X$.
 Therefore, according to Section~\ref{pull-push}, there are
the derived inverse image functors
\begin{alignat*}{2}
 \boL f^*\: &\sD^\co((X,\L,w)\qcoh_\ffd)&&\lrarrow
            \sD^\co((Y,f^*\L,f^*w)\qcoh_\ffd) \\
 \boL f^*\: &\sD^\abs((X,\L,w)\coh_\ffd)&&\lrarrow
            \sD^\abs((Y,f^*\L,f^*w)\coh_\ffd)
\end{alignat*}
and the derived direct image functor
$$
 \boR f_*\:\sD^\co((Y,f^*\L,f^*w)\qcoh)\lrarrow
 \sD^\co((X,\L,w)\qcoh).
$$
 The latter two functors are ``partially adjoint'' to each other.

 Given a triangulated category $\sD$, we denote by
$\overline{\sD\!\.}\.$ its idempotent completion.
 By~\cite[Section~1]{BS}, the category $\overline{\sD\!\.}\.$
has a natural structure of triangulated category.

\begin{lem}
 For any closed subset $T\subset Y$ such that (for a closed subscheme
structure on $T$) the morphism $f|_T\:T\rarrow X$ is proper of finite
type, the functor $\boR f_*$ takes the full subcategory
$\sD^\abs((Y,f^*\L,f^*w)\coh_T)\subset\sD^\co((Y,f^*\L,f^*w)\qcoh)$
into the full subcategory $\sD^\abs((X,\L,w)\coh)\subset
\sD^\co((X,\L,w)\qcoh)$, thus defining a triangulated functor
of direct image
$$
 \boR f_*\:\sD^\abs((Y,f^*\L,f^*w)\coh_T)\lrarrow
 \sD^\abs((X,\L,w)\coh).
$$
Consequently, there is the triangulated functor
$$
 \overline{\boR f_*\!}\,\:\overline{\sD^\abs_T((Y,f^*\L,f^*w)\coh)}
 \lrarrow\overline{\sD^\abs((X,\L,w)\coh)}.
$$
\end{lem}

\begin{proof}
 We will use the construction of the functor $\boR f_*\:
\sD^\co((Y,f^*\L,f^*w)\qcoh)\rarrow\sD^\co((X,\L,w)\qcoh)$
similar to the one in the proof of
Proposition~\ref{finite-dim-morphisms}
(see Remark~\ref{finite-dim-morphisms}).
 According to this construction, given a matrix factorization
$\M\in (Y,f^*\L,f^*w)\qcoh$, the object $\boR f_*\M\in
\sD^\co((X,\L,w)\qcoh)$ is represented by the total matrix
factorization $\boR_{\{U_\alpha\}}f_*\M$ of the finite \v Cech complex
$f_*C_{\{U_\alpha\}}^\bu\M$ of matrix factorizations on~$X$.
 The derived functor of direct image of complexes of
quasi-coherent sheaves $\boR f_*\:\sD^\b(Y\qcoh)\rarrow\sD^\b
(X\qcoh)$ can be constructed in the same way.
{\emergencystretch=1em\par}

 By~\cite[Th\'eor\`eme~3.2.1]{Groth}, the latter functor takes
$\sD^\b(Y\coh_T)$ into $\sD^\b(X\coh)$.
 Hence the cohomology matrix factorizations of the finite complex
of matrix factorizations $f_*C_{\{U_\alpha\}}^\bu \M$ belong to
$(X,\L,w)\coh$ when the matrix factorization $\M$ belongs to
$(Y,f^*\L,f^*w)\coh_T$.
 It follows that the object $\boR f_*\M$ belongs to
$\sD^\abs((X,\L,w)\coh)\allowbreak\subset\sD^\co((X,\L,w)\qcoh)$
in this case.

 To prove the last assertion, it remains to apply
Corollary~\ref{cdg-supports}(b).
\end{proof}

 Now assume that both morphisms of sheaves $w\:\O_X\rarrow\L$ and
$f^*w\:\O_Y\rarrow f^*\L$ are injective.
 Let $X_0\subset X$ and $Y_0\subset Y$ denote the closed subschemes
defined locally by the equations $w=0$ and $f^*w=0$, respectively.
 In this setting, we will compare the constructions of direct
image functors for matrix factorizations and for the triangulated
categories of relative singularities, and prove the assertions
of Lemma in a different way.
 Recall that in Section~\ref{pull-push-singularities} we have
constructed the functor of direct image
$f_\circ\:\sD'_\Sing(Y_0/Y)\rarrow\sD'_\Sing(X_0/X)$.

\begin{prop}
\textup{(a)} Whenever the morphism $f_0=f|_{Y_0}\:Y_0\rarrow X_0$
is proper of finite type, the functor\/ $\boR f_*$ takes the full
subcategory\/ $\sD^\abs((Y,f^*\L,f^*w)\coh)\subset\sD^\co
((Y,f^*\L,f^*w)\qcoh)$ into the full subcategory\/ $\sD^\abs((X,\L,w)
\coh)\subset\sD^\co((X,\L,w)\allowbreak\qcoh)$, thus defining
a triangulated functor
$$
 \boR f_*\:\sD^\abs((Y,f^*\L,f^*w)\coh)\lrarrow
 \sD^\abs((X,\L,w)\coh).
$$ \par
\textup{(b)} For any closed subset $T\subset Y_0$ such that (for a
closed subscheme structure on~$T$) the morphism $f_0|_T\:T\rarrow X_0$
is proper of finite type, the functor $f_\circ$ takes the full
subcategory $\sD^\b_\Sing(Y_0/Y,T)\subset\sD'_\Sing(Y_0/Y)$ into
the full subcategory $\sD^\b_\Sing(X_0/X)\subset\sD'_\Sing(X_0/X)$,
thus defining a triangulated functor
$$
 f_\circ\:\sD^\b_\Sing(Y_0/Y,T)\lrarrow\sD^\b_\Sing(X_0/X).
$$ \par
\textup{(c)} The equivalences of categories\/ $\sD^\abs((Y,f^*\L,f^*w)
\coh_T)\simeq\sD^\b_\Sing(X_0/X,T)$ from
Proposition\/~\textup{\ref{supports-secn}(a)} and\/
$\sD^\abs((X,\L,w)\coh)\simeq\sD^\b_\Sing(X_0/X)$ from
Theorem\/~\textup{\ref{main-theorem}} transform the direct image
functor\/ $\boR f_*\:\sD^\abs((Y,f^*\L,f^*w)\coh_T)\rarrow
\sD^\abs((X,\L,w)\coh)$ from Lemma into the direct image
functor~$f_\circ$ from part~\textup{(b)}. \hbadness=1200
\end{prop}

\begin{proof}
 Part~(a) follows from Lemma and Proposition~\ref{supports-secn}(a),
or alternatively, from part~(b) and the proof of part~(c) below.
 In part~(b), the fact of key importance is that the functor
$\sD^\b_\Sing(X_0/X)\rarrow\sD'_\Sing(X_0/X)$ is fully faithful
(by Theorem~\ref{infinite-matrix}).
 The functor~$f_\circ$ takes $\sD^\b_\Sing(Y_0/Y,T)$ into
$\sD^\b_\Sing(X_0/X)$, because the functor $\boR f_0{}_*\:
\sD^\b(Y_0\qcoh)\rarrow\sD^\b(X_0\qcoh)$ takes
$\sD^\b(Y_0\coh_T)$ into $\sD^\b(X_0\coh)$ \cite{Groth}.
 To prove part~(c), we will check that the equivalences of categories
from Theorem~\ref{infinite-matrix} transform the functor $\boR f_*\:
\sD^\co((Y,f^*\L,f^*w)\qcoh)\rarrow\sD^\co((X,\L,w)\qcoh)$ into
the functor $f_\circ\:\sD'_\Sing(Y_0/Y)\rarrow\sD'_\Sing(X_0/X)$.
 (Together with part~(b) and Proposition~\ref{supports-secn}(a),
this will also provide another proof of Lemma.)

 For this purpose, extend the functor $\Ups_Y\:\allowbreak
\sD^\b(Y_0\qcoh)\rarrow\sD^\co((Y,f^*\L,f^*w)\qcoh)$
to a functor $\widetilde\Ups_Y\:\sD^+(Y_0\qcoh)\rarrow
\sD^\co((Y,\allowbreak f^*\L,f^*w)\qcoh)$ in the obvious way
(taking infinite direct sums of quasi-coherent sheaves in
the construction of the matrix factorization
$\widetilde\Ups_Y(\F^\bu)$).
 The functor $\widetilde\Ups_Y$ is well-defined, since any bounded
below acyclic complex of quasi-coherent sheaves is
coacyclic~\cite[Lemma~2.1]{Psemi}.
 Furthermore, the functor $\widetilde\Ups_Y$ can be presented as
the composition of the ``periodicity summation'' functor
$\sD^+(Y_0\qcoh)\rarrow\sD^\co((Y_0,i'{}^*f^*\L,0)\allowbreak\qcoh)$
taking values in the coderived category of quasi-coherent matrix
factorizations of the zero potential on~$Y_0$, and the functor of
direct image $i'_*\:\sD^\co((Y_0,i'{}^*f^*\L,0)\qcoh)\allowbreak
\rarrow\sD^\co((Y,f^*\L,f^*w)\qcoh)$ with respect to
the closed embedding~$i'$. {\emergencystretch=0em\hfuzz=2pt\par}

 The functors $\boR f_0{}_*\:\sD^+(Y_0\qcoh)\rarrow\sD^+(X_0\qcoh)$
and $\boR f_*\:\sD^\co((Y,f^*\L,f^*w)\allowbreak\qcoh)\rarrow
\sD^\co((X,\L,w)\qcoh)$ form a commutative diagram with the functors
$\widetilde\Ups_X$ and $\widetilde\Ups_Y$.
 Indeed, the ``periodicity summations'' of bounded below complexes
of quasi-coherent sheaves on $Y_0$ and $X_0$, taking injective
resolutions to injective resolutions, obviously commute with
the derived direct images with respect to~$f'$, as the direct image
preserves infinite direct sums.
 Furthermore, the derived direct images of quasi-coherent matrix
factorizations are compatible with the compositions of morphisms
of schemes (see Remark~\ref{pull-push}), hence also commute with
each other.
 It follows that the functors $\boR f_*$ and~$f_\circ$ agree
as they should.
 (Alternatively, one can prove this in the way similar to the proof
of Proposition~\ref{finite-dim-matrix-push} below.)
\end{proof}

\subsection{Push-forwards for morphisms of finite flat dimension}
\label{finite-dim-matrix-push}
 Let $f\:Y\rarrow X$ be a morphism of finite flat dimension between
separated Noetherian schemes with enough vector bundles,
$\L$ be a line bundle on $X$, and $w\in\L(X)$ be a section.
 As in Section~\ref{matrix-push}, we have a natural morphism of
CDG\+algebras $\B_X=(X,\L,w)\rarrow\B_Y=(Y,f^*\L,f^*w)$ compatible
with the morphism of schemes $Y\rarrow X$.

 The quasi-coherent graded algebra $\B_Y$ has finite flat dimension
over $\B_X$.
 Therefore, according to Section~\ref{finite-dim-morphisms},
there are derived inverse image functors
\begin{alignat*}{2}
 \boL f^*\: &\sD^\co((X,\L,w)\qcoh)&&\lrarrow
            \sD^\co((Y,f^*\L,f^*w)\qcoh) \\
 \boL f^*\: &\sD^\abs((X,\L,w)\coh)&&\lrarrow
            \sD^\abs((Y,f^*\L,f^*w)\coh),
\end{alignat*}
the former of which is left adjoint to the functor
$\boR f_*\:\sD^\co((Y,f^*\L,f^*w)\qcoh)\rarrow\sD^\co((X,\L,w)\qcoh)$
from Section~\ref{matrix-push}.

 Furthermore, according to Proposition~\ref{finite-dim-morphisms},
there is a derived direct image functor
\begin{multline*}
 \boR f_*\:\sD^\co((Y,f^*\L,f^*w)\qcoh_\ffd)\.\simeq\.
 \sD^\co((Y,f^*\L,f^*w)\qcoh_\fl) \\
 \lrarrow\sD^\co((X,\L,w)\qcoh_\ffd)\.\simeq\.
 \sD^\co((X,\L,w)\qcoh_\fl)
\end{multline*}
which is right adjoint to the functor
$\boL f^*\:\sD^\co((X,\L,w)\qcoh_\ffd)\rarrow
\sD^\co((Y,f^*\L,\allowbreak f^*w)\qcoh_\ffd)$
from Section~\ref{matrix-push}.

 Now assume that $X$ and $Y$ have finite Krull dimensions.
 Recall that the natural triangulated functors
$\sD^\abs((X,\L,w)\coh_\lf)\rarrow\sD^\co((X,\L,w)\qcoh_\fl)$ and
$\sD^\abs((Y,f^*\L,f^*w)\coh_\lf)\rarrow\sD^\co((Y,f^*\L,f^*w)
\qcoh_\fl)$ are fully faithful by
Corollary~\ref{exotic-derived-matrix-cor}(e) and~(j).

 As in the second half of Section~\ref{matrix-push}, assume that both
morphisms of sheaves $w\:\O_X\rarrow\L$ and $f^*w\:\O_Y\rarrow f^*\L$
are injective, and denote by $f_0\:Y_0\rarrow X_0$ the induced
morphism between the zero loci schemes of $f^*w$ and~$w$.
 Since the morphism $f$ has finite flat dimension, so does
the morphism~$f_0$.

\begin{prop}
\textup{(a)} Whenever the morphism $f_0$~is proper of finite type,
the functor\/ $\boR f_*\:\allowbreak\sD^\co((Y,f^*\L,f^*w)\qcoh_\fl)
\rarrow\sD^\co((X,\L,w)\qcoh_\fl)$ takes the full subcategory\/ 
$\sD^\abs((Y,f^*\L,f^*w)\coh_\lf)\subset\sD^\co((Y,f^*\L,f^*w)
\qcoh_\fl)$ into the full subcategory\/ $\sD^\abs((X,\L,w)\coh_\lf)
\subset\sD^\co((X,\L,w)\allowbreak\qcoh_\fl)$.  
 Besides, the functor $f_0{}_\circ\:\sD^\b_\Sing(Y_0)\allowbreak
\rarrow\sD^\b_\Sing(X_0)$ takes the full subcategory\/
$\sD^\abs((Y,f^*\L,f^*w)\coh_\lf)\subset\sD^\b_\Sing(Y_0)$ into
the full subcategory\/ $\sD^\abs((X,\L,w)\coh_\lf)\subset
\sD^\b_\Sing(X_0)$.
 Both restrictions define the same triangulated functor
$$
 \boR f_*\:\sD^\abs((Y,f^*\L,f^*w)\coh_\lf)\lrarrow
 \sD^\abs((X,\L,w)\coh_\lf).
$$ \par
\textup{(b)} For any closed subset $T\subset Y_0$ such that
(for a closed subscheme structure on~$T$) the morphism
$f_0|_T\:T\rarrow X_0$ is proper of finite type, the functor\/
$\boR f_*\:\sD^\co((Y,f^*\L,f^*w)\qcoh_\fl)\rarrow
\sD^\co((X,\L,w)\qcoh_\fl)$ takes the full subcategory\/ 
$\sD^\abs_T((Y,f^*\L,f^*w)\coh_\lf)\subset\sD^\co((Y,f^*\L,f^*w)
\qcoh_\fl)$ into the thick envelope of the full subcategory\/
$\sD^\abs((X,\L,w)\coh_\lf)\subset\sD^\co((X,\L,w)\qcoh_\fl)$.  
 Besides, the triangulated functor $\overline{f_0{}_\circ}\:
\overline{\sD^\b_\Sing(Y_0,T)}\rarrow\overline{\sD^\b_\Sing(X_0)}$
takes the full subcategory\/
$\sD^\abs_T((Y,f^*\L,f^*w)\coh_\lf)\subset\overline{\sD^\b_\Sing
(Y_0,T)}$ into the thick envelope of the full subcategory\/
$\sD^\abs((X,\L,w)\coh_\lf)\subset\overline{\sD^\b_\Sing(X_0)}$.
 Both restrictions define the same triangulated functor
$$
 \overline{\boR f_*\!}\,\:
 \overline{\sD^\abs_T((Y,f^*\L,f^*w)\coh_\lf)}\lrarrow
 \overline{\sD^\abs((X,\L,w)\coh_\lf)}.
$$
\end{prop}

\begin{proof}
 Both categories $\sD^\co((X,\L,w)\qcoh_\fl)$ and $\sD^\b_\Sing(X_0)$
are full triangulated subcategories of the triangulated category
$\sD'_\Sing(X_0)$ (see Proposition~\ref{infinite-matrix}
and~\cite[Proposition~1.13]{Or1}).
 According to the proof of Corollary~\ref{locality}, the intersection
of $\sD^\co((X,\L,w)\qcoh_\fl)$ with (the thick envelope of)
$\sD^\b_\Sing(X_0)$ in $\sD'_\Sing(X_0)$ (is the thick envelope of)
the subcategory $\sD^\abs((X,\L,w)\coh_\lf)\subset\sD'_\Sing(X_0)$.

 Thus it suffices to show that the direct image functor
$\boR f_*\:\sD^\co((X,\L,w)\qcoh_\fl)\allowbreak\rarrow
\sD^\co((X,\L,w)\qcoh_\fl)$ agrees with the direct image functor 
$f_0{}_\circ\:\sD'_\Sing(Y_0)\rarrow\sD'_\Sing(X_0)$.
 The latter assertion does not depend on any properness
assumptions.

 Recall that the derived functor $\boR f_*$ was constructed in
the proof of Proposition~\ref{finite-dim-morphisms} in terms of
the \v Cech complex whose terms are direct sums of
the CDG\+modules $f|_V{}_*\M|_V$, where $\M\in\sD^\co((Y,f^*\L,f^*w)
\qcoh_\ffd)$ and $V\subset Y$.
 The derived direct image $\boR f_0{}_*\:\sD^\b(Y_0\qcoh)\rarrow
\sD^\b(X_0\qcoh)$ can be constructed in the similar way; moreover,
one can use for this purpose the restriction to $Y_0$ of an affine
open covering $U_\alpha$ of the scheme~$Y$.

 We will make use of the flat dimension analogue of
Corollary~\ref{w-flat-cor}(d).
 Let $\widetilde\Sigma'_X$ and $\widetilde\Sigma'_Y$ denote
the obvious extensions of the functors $\Sigma'$ from
$(X,\L,w)\qcoh_\lf$ to the category of $w$\+flat matrix
factorizations of finite flat dimension $(X,\L,w)
\qcoh_{\wfl\cap\ffd}$ and from $(Y,f^*\L,f^*w)\qcoh_\lf$ to
$(Y,f^*\L,f^*w)\qcoh_{\fwfl\cap\ffd}$
(see the proofs of Proposition~\ref{infinite-matrix} and
Theorem~\ref{main-theorem}).
 Notice that the direct image functors $f|_V{}_*$ take $f^*w$\+flat
sheaves to $w$\+flat sheaves and $(V,f^*\L|_V,f^*w|_V)
\qcoh_{\fwfl\cap\ffd}$ to $(X,\L,w)\qcoh_{\wfl\cap\ffd}$.

 Let $\N$ be a matrix factorization from $(Y,f^*\L,f^*w)\qcoh
_{\fwfl\cap\ffd}$.
 Since the open subschemes $V$ are presumed to be affine, there
are natural isomorphisms $\widetilde\Sigma'_X(f|_V{}_*\N|_V)\simeq
f_0|_{V\cap Y_0}\.{}_*\widetilde\Sigma'_Y(\N)|_{V\cap Y_0}$ of
quasi-coherent sheaves on~$X_0$.
 Now it remains to use the next lemma.
\end{proof}

\begin{lem}
 Let $\M^{-n}\rarrow\dotsb\rarrow\M^N$ be a finite complex of
matrix factorizations from $(X,\L,w)\qcoh_{\wfl\cap\ffd}$ and
$\M$ be its totalization.
 Then the complex\/ $\widetilde\Sigma'(\M^{-n})\rarrow\dotsb\rarrow
\widetilde\Sigma'(\M^N)$ and the quasi-coherent sheaf\/
$\widetilde\Sigma'(\M)$ on $X_0$ represent naturally isomorphic
objects in the triangulated category of singularities\/
$\sD'_\Sing(X_0)$.
 The same applies to a finite complex of matrix factorizations from
$(X,\L,w)\qcoh_\wfl$, the functor\/ $\Xi$, and the triangulated
category of relative singularities\/ $\sD''_\Sing(X_0/X)$.
\end{lem}

\begin{proof}
 For each $-n\le p\le N$, the restriction of the matrix factorization
$\M^p$ to the closed subscheme $X_0\subset X$ is an unbounded
complex of quasi-coherent sheaves $i^*\M^{p,\bu}$.
 By~\cite[Lemma~1.5]{PV}, this complex is acyclic.

 The complex $\widetilde\Sigma'(\M^{-n})\rarrow\dotsb\rarrow
\widetilde\Sigma'(\M^N)$ of quasi-coherent sheaves on $X_0$ is
quasi-isomorphic to the total complex of the bicomplex $\K^{\bu,\bu}$
with the terms $\K^{p,0}=i^*\M^{p,0}$, \ $\K^{p,-1}=i^*\M^{p,-1}$, \
$\K^{p,-2}= \ker(i^*\M^{p,-1}\to i^*\M^{p,0})$, and $\K^{p,q}=0$ for
$q\ne0$, $-1$, $-2$.
 Similarly, the quasi-coherent sheaf $\widetilde\Sigma'(\M)$ on $X_0$
is quasi-isomorphic to the total complex of the bicomplex
$\E^{\bu,\bu}$ with the terms $\E^{p,p}=i^*\M^{p,p}$, \
$\E^{p,p-1}=i^*\M^{p,p-1}$, \ $\E^{p,p-2}=\ker(i^*\M^{p,p-1}\to
i^*\M^{p,p})$, and $\E^{p,q}=0$ for $q-p\ne0$, $-1$, $-2$.

 We can assume that $N$, $n\ge0$.
 Consider the bicomplex $\F^{\bu,\bu}$ with the terms
$\F^{p,q}=i^*\M^{p,q}$ for $-n-1\le q\le N$, \ $\F^{p,-n-2}=
\ker(i^*\M^{p,-n-1}\to i^*\M^{p,-n})$, and $\F^{p,q}=0$ for
$q<-n-2$ or $q>N$.
 Then there are natural surjective morphisms of bicomplexes
$\F^{\bu,\bu}\rarrow\K^{\bu,\bu}$ and $\F^{\bu,\bu}\rarrow
\E^{\bu,\bu}$.
 The kernels of both morphisms are the direct sums of
a finite bicomplex of quasi-coherent sheaves of finite flat
dimension on $X_0$ and a finite bicomplex of quasi-coherent sheaves
on $X_0$ with acyclic columns.
 Thus both morphisms become isomorphisms in $\sD'_\Sing(X_0)$.
\end{proof}

\begin{rem}
 One would like to have a theory of set-theoretic supports for
locally free matrix factorizations of finite rank that would
allow to prove the above Proposition in the way similar to
the proof of Lemma~\ref{matrix-push}.
 However, we do not know how to do this.
 In particular, we do \emph{not} know whether every locally free matrix
factorization of finite rank with the category-theoretic support
in $T$ is isomorphic in the absolute derived category to a direct
summand of an object represented by a coherent matrix factorization
of finite flat dimension with the set-theoretic support in~$T$
(cf.\ Corollary~\ref{cdg-supports} and Section~\ref{nonlocalization}).

 Another alternative approach to proving Proposition would be
to show that the intersection of the full subcategories
$\sD^\abs((X,\L,w)\coh)$ and $\sD^\abs((X,\L,w)\qcoh_\lf)$ in
the absolute derived category $\sD^\abs((X,\L,w)\qcoh)$
coincides with the full subcategory $\sD^\abs((X,\L,w)\coh_\lf)$.
 We do \emph{not} know whether this is true.
\end{rem}

\subsection{Duality and push-forwards} \label{push-duality}
 In the following two sections we discuss the compatibility
properties of the derived direct and inverse image functors for matrix
factorizations with the Serre--Grothendieck duality functors from
Section~\ref{serre-duality}.

 Let $X$ be a separated Noetherian scheme with a dualizing complex
$\D_X^\bu$, and let $f\:Y\rarrow X$ be a separated morphism of
finite type.
 As usually, we set $\D_Y^\bu=f^+\D_X^\bu$, where $f^+$ is the functor
denoted by $f^!$ in~\cite{Har} (right adjoint to $\boR f_*$ for proper
morphisms~$f$ and left adjoint to $\boR f_*$ for open embeddings~$f$;
see~\cite[Example~4.2]{Neem} and~\cite[Remark before
Proposition~V.8.5 and Deligne's Appendix]{Har}).
 This formula defines the dualizing complex $\D_Y^\bu$ up to
a natural quasi-isomorphism only, and we presume this derived
category object (as well as $\D_X^\bu$) to be represented
by a finite complex of injective quasi-coherent sheaves.

\begin{prop}
 Let $T\subset Y_0$ be a closed subset such that (for some
closed subscheme structure on~$T$) the morphism $f|_T\:T\rarrow X_0$
is proper.
 Then the derived direct image functor $\overline{\boR f_*\!}\,\:
\overline{\sD^\abs_T((Y,f^*\L,f^*w)\coh)}\rarrow
\overline{\sD^\abs((X,\L,w)\coh)}$ and the similar functor for
the potential\/~$-w$ form a commutative diagram with the Serre duality
functors\/ $\cHom_{X\qc}({-},\D_X^\bu)\:\sD^\abs((X,\L,-w)\coh)^\sop
\rarrow\sD^\abs((X,\L,w)\coh)$ and\/ $\cHom_{Y\qc}({-},\D_Y^\bu)\:
\sD^\abs_T((Y,f^*\L,-f^*w)\coh)^\sop\rarrow
\sD^\abs_T((Y,f^*\L,f^*w)\coh)$.
\end{prop}

 Two proofs of Proposition are given below.
 One of them is based on the theory of set-theoretic supports of
coherent CDG\+modules developed in Section~\ref{cdg-supports} and
the arguments similar to the proof of Lemma~\ref{matrix-push}.
 It does not depend on the assumption about $w$ and $f^*w$ being
local nonzero-divisors and does not mention the zero loci.
 The other proof is based on the passage to the triangulated
categories of relative singularities and uses
Proposition~\ref{matrix-push}(c).

\begin{proof}[First proof]
 First of all, the duality functor $\cHom_{Y\qc}({-},\D_Y^\bu)\:
\sD^\abs((Y,f^*\L,-f^*w)\allowbreak\qcoh)^\sop\rarrow
\sD^\abs((Y,f^*\L,f^*w)\qcoh)$ obviously takes the full subcategory 
$\sD^\abs((Y,\allowbreak f^*\L,-f^*w)\coh_T)^\sop$ into 
$\sD^\abs((Y,f^*\L,f^*w)\coh_T)$ and vice versa.
 Furthermore, for any quasi-coherent sheaf $\K$ on $Y$ denote
by $\Gamma_T\K\subset\K$ the maximal quasi-coherent subsheaf with
the set-theoretic support in~$T$.
 Then for any matrix factorization
$\M\in\sD^\abs((Y,f^*\L,-f^*w)\coh_T)$ the natural morphism
$\cHom_{Y\qc}(\M\;\Gamma_T\D_Y^\bu)\rarrow
\cHom_{Y\qc}(\M,\D_Y^\bu)$ is an isomorphism in
$\sD^\abs((Y,f^*\L,f^*w)\coh_T)$.

 As in the proof of Lemma~\ref{matrix-push}, we will use
the construction of the functor $\boR f_*\:\sD^\abs((Y,f^*\L,f^*w)
\allowbreak\qcoh)\rarrow\sD^\abs((X,\L,w)\qcoh)$ similar to the one
from the proof of Proposition~\ref{finite-dim-morphisms} (see
Remarks~\ref{pull-push} and~\ref{finite-dim-morphisms}).
 Let $\{U_\alpha\}$ and $\{V_\beta\}$ be two affine open coverings
of the scheme~$Y$.
 For any matrix factorization $\N\in(Y,f^*\L,-f^*w)\qcoh$,
there is a natural morphism of bicomplexes of matrix factorizations
$f_*C_{\{U_\alpha\}}^\bu\cHom_{Y\qc}(\N\;\Gamma_T\D_Y^\bu)\rarrow
\cHom_{X\qc}(f_*\N\;f_*C_{\{U_\alpha\}}^\bu\Gamma_T\D_Y^\bu)$.
 Passing to the total complexes and taking the composition with
the adjunction morphism $f_*C_{\{U_\alpha\}}^\bu\Gamma_T\D_Y^\bu =
\boR f_*(\Gamma_T\D_Y^\bu)\rarrow\D_X^\bu$,
we obtain a natural morphism of complexes of matrix factorizations
$f_*C_{\{U_\alpha\}}^\bu\cHom_{Y\qc}(\N\;\Gamma_T\D_Y^\bu)\rarrow
\cHom_{X\qc}(f_*\N\;\D_X^\bu)$ (cf.~\cite[beginning of
Section~6]{Neem}).

 Substituting $\N=C_{\{V_\beta\}}^\bu\M$ for some
$\M\in(Y,f^*\L,-f^*w)\qcoh$, we get a natural morphism of
bicomplexes of matrix factorizations $f_*C_{\{U_\alpha\}}^\bu
\cHom_{Y\qc}(C_{\{V_\beta\}}^\bu\M\;\Gamma_T\D_Y^\bu)
\allowbreak\rarrow\cHom_{X\qc}(f_*C_{\{V_\beta\}}^\bu\M\;\D_X^\bu)$.
 When $\M$ is a coherent matrix factorization supported
set-theoretically in $T$, the induced morphism of the total complexes
is a quasi-isomorphism of complexes of matrix factorizations by
the conventional Serre--Grothendieck duality theorem for bounded
derived categories of coherent sheaves and proper morphisms of
schemes (see~\cite[Theorem~VII.3.3]{Har} or~\cite[Section~6]{Neem}).
 Hence the induced morphism of the total matrix factorizations is
an isomorphism in $\sD^\abs((X,\L,w)\qcoh)$, and consequently
also in $\sD^\abs((X,\L,w)\coh)$.
\end{proof}

\begin{proof}[Second proof]
 Assume that $w$ and $f^*w$ are locally nonzero-dividing sections
of the respective line bundles.
 Let $i\:X_0\rarrow X$ be the zero locus of~$w$ and 
$i'\:Y_0\rarrow Y$ be the zero locus of~$f^*w$.
 As above, we set $\D_{X_0}^\bu = \boR i^!\D_X^\bu$ and
$\D_{Y_0}^\bu = \boR i'{}^!\D_Y^\bu$
\cite[Proposition~V.2.4]{Har}, and presume all these
dualizing complexes to be finite complexes of injective
quasi-coherent sheaves.

 The duality functor $\cHom_{Y\qc}({-},\D_Y^\bu)\:
\sD^\abs((Y,f^*\L,-f^*w)\coh)^\sop\rarrow\sD^\abs((Y,\allowbreak
f^*\L,f^*w)\coh)$ is compatible with the restrictions to
the open subscheme $Y\setminus T$ and therefore identifies the full
subcategories $\sD^\abs_T((Y,f^*\L,-f^*w)\coh)^\sop$ and
$\sD^\abs_T((Y,f^*\L,f^*w)\coh)$.
 To prove the proposition, we will define the Serre duality functors
on the triangulated categories of relative singularities
$\sD^\b_\Sing(Y_0/Y)$ and $\sD^\b_\Sing(X/X_0)$, then check that
the equivalences of triangulated categories $\boL\Xi=\Ups^{-1}$
commute with the dualities, and finally reduce to the conventional
Serre--Grothendieck duality theorem for bounded complexes of
coherent sheaves. {\hbadness=1700\par}

 The duality functor $\cHom_{X_0\qc}({-},\D_{X_0}^\bu)\:
\sD^\b(X_0\coh)^\sop\rarrow\sD^\b(X_0\coh)$ takes objects of
the form $\boL i^*\K^\bu$, where $\K^\bu\in\sD^\b(X\coh)$,
to similar objects.
 Indeed, one has $\cHom_{X_0\qc}(\boL i^*\K^\bu,\D_{X_0}^\bu)\simeq
\boR i^!\cHom_{X\qc}(\K^\bu,\D_X^\bu)$ \cite[Proposition~V.8.5]{Har}
and $\boR i^!\simeq\L|_{X_0}[-1]\allowbreak\ot_{\O_{X_0}}\boL i^*$
(see the proof of Theorem~\ref{main-theorem}).
 Therefore, we have the induced duality functor $\cHom_{X_0\qc}
({-},\D_{X_0}^\bu)\:\sD^\b_\Sing(X_0/X)^\sop\rarrow\sD^\b_\Sing
(X_0/X)$.
 Similarly, the duality functor $\cHom_{Y_0\qc}({-},\D_{Y_0}^\bu)\:
\sD^\b(Y_0\coh)^\sop\rarrow\sD^\b(Y_0\coh)$ takes the full
subcategory $\sD^\b(Y_0\coh_T)^\sop$ into
$\sD^\b(Y_0\coh_T)$ and $\Perf_T(Y_0/Y)^\sop$ into
$\Perf_T(Y_0/Y)$.
 Hence the induced duality functor $\cHom_{Y_0\qc}
({-},\D_{Y_0}^\bu)\:\sD^\b_\Sing(Y_0/Y,T)^\sop\rarrow
\sD^\b_\Sing(Y_0/Y,T)$.

 Checking that the equivalence of categories
$\sD^\abs((X,\L,w)\coh)\simeq\sD^\b_\Sing(X_0/X)$ commutes with
the dualities is easily done using the functor~$\Ups$.
 It suffices to notice the functorial quasi-isomorphism
$\cHom_{X\qc}(i_*\F^\bu,\D_X^\bu)\simeq i_*
\cHom_{X_0\qc}(\F^\bu,\D_{X_0}^\bu)$ for any complex
$\F^\bu\in\sD^\b(X_0\coh)$ \cite[Theorem~III.6.7]{Har}.
 The same applies to the equivalence of categories
$\overline{\sD^\abs_T((Y,f^*\L,f^*w)\coh)}\simeq
\overline{\sD^\b_\Sing(Y_0/Y,T)}$.
 Furthermore, by Proposition~\ref{matrix-push}(c), the equivalences
of categories $\boL\Xi=\Ups^{-1}$ transform the derived direct image
functor $\overline{\boR f_*\!}\,\:
\overline{\sD^\abs_T((Y,f^*\L,f^*w)\coh)}\rarrow
\overline{\sD^\abs((X,\L,w)\coh)}$ into (the idempotent closure
of) the direct image functor $f_\circ\:\sD^\b_\Sing
(Y_0/Y,T)\rarrow\sD^\b_\Sing(X_0/X)$.

 Finally, the direct image functor $f_\circ\:\sD^\b_\Sing
(Y_0/Y,T)\rarrow\sD^\b_\Sing(X_0/X)$ commutes with the Serre
duality functors, since so do the derived direct image functors
$\boR f|_T{}_*\:\sD^\b(\widetilde T\coh)\rarrow
\sD^\b(X_0\coh)$ for all the closed subscheme structures
$\widetilde T\subset Y_0$ on the closed subset~$T$ and
the similar functors related to the closed embeddings
$\widetilde T'\rarrow\widetilde T''$ of various such subscheme
structures into each other.
 This is the conventional Serre--Grothendieck duality theorem
for proper morphisms of schemes.
\end{proof}

\subsection{Duality and pull-backs}  \label{pull-duality}
 Let $X$ be a separated Noetherian scheme with a dualizing complex
$\D_X^\bu$ and $f\:Y\rarrow X$ be a separated morphism of finite type;
set $\D_Y^\bu = f^+\D_X^\bu$.
 Let $\L$ be a line bundle on $X$ and $w\in\L(X)$ be a section.

 Let us first suppose that the morphism $f$ is smooth of relative
dimension~$n$.
 Then the functor $f^+\:\sD^+(X\qcoh)\rarrow\sD^+(Y\qcoh)$ is
naturally isomorphic to $\omega_{Y/X}[n]\otimes_{\O_Y}\nobreak f^*$,
where $\omega_{Y/X}$ is the line bundle of relative top forms.

 In particular, $\D_Y^\bu\simeq\omega_{Y/X}[n]
\ot_{\O_Y}f^*\D_X^\bu$ (where $f^*\D_X^\bu$ is also presumed to
have been replaced by a complex of injectives).
 Then it is clear that the equivalences of categories
$\D_X^\bu\otimes_{\O_X}{-}\:\sD^\co((X,\L,w)\qcoh_\fl)\rarrow
\sD^\co((X,\L,w)\qcoh)$ and $f^*\D_X^\bu\otimes_{\O_Y}\nobreak{-}\:
\sD^\co((Y,f^*\L,\allowbreak f^*w)\qcoh_\fl)\rarrow
\sD^\co((Y,f^*\L,f^*w)\qcoh)$ from Section~\ref{serre-duality}
transform the inverse image functor for flat matrix factorizations
$f^*\:\sD^\co((X,\L,w)\qcoh_\fl)\rarrow
\sD^\co((Y,f^*\L,f^*w)\qcoh_\fl)$ into the (underived, as
the morphism $f$ is flat) inverse image functor for quasi-coherent
matrix factorizations $f^*\:\sD^\co((X,\L,w)\qcoh)\rarrow
\sD^\co((Y,f^*\L,f^*w)\qcoh)$.

 Furthermore, for any quasi-coherent matrix factorization $\M$ on
$X$ there is a natural morphism of finite complexes of matrix
factorizations $f^*\cHom_{X\qc}(\M,\D_X^\bu)\rarrow
\cHom_{Y\qc}(f^*\M,f^*\D_X^\bu)$ on~$Y$.
 When $\M$ is a coherent matrix factorization, this is
a quasi-isomorphism of complexes of matrix factorizations
(since the similar assertion holds for coherent
sheaves~\cite[Proposition~II.5.8]{Har}), so the related morphism
of total matrix factorizations has an absolutely acyclic cone.
 Thus the anti-equivalences of categories $\cHom_{X\qc}({-},
\D_X^\bu)\:\sD^\abs((X,\L,-w)\coh)^\sop\rarrow\sD^\abs((X,\L,w)\coh)$
and $\cHom_{Y\qc}({-},f^*\D_X^\bu)\:\sD^\abs((Y,f^*\L,-f^*w)\coh)^\sop
\rarrow\sD^\abs((Y,f^*\L,f^*w)\coh)$ form a commutative diagram with
the inverse image functors~$f^*$ for coherent matrix factorizations.

 Now suppose that $f$~is a proper morphism of finite type.
 The following theorem describes the compatibility property of
the covariant Serre--Grothendieck duality with the inverse
images of matrix factorizations (cf.~\cite[Theorem~5.15.3]{Pcosh},
where the similar result is proven for complexes of
quasi-coherent sheaves).

\begin{thm}
 The equivalences of categories
$\D_X^\bu\otimes_{\O_X}{-}\:\sD^\abs((X,\L,w)\qcoh_\fl)\rarrow
\sD^\co((X,\L,w)\qcoh)$ and $\D_Y^\bu\otimes_{\O_Y}\nobreak{-}\:
\sD^\abs((Y,f^*\L,f^*w)\qcoh_\fl)\rarrow\sD^\co((Y,f^*\L,
\allowbreak f^*w)\qcoh)$ transform the inverse image functor
$f^*\:\sD^\abs((X,\L,w)\qcoh_\fl)\rarrow\sD^\abs((Y,f^*\L,f^*w)
\qcoh_\fl)$ into the functor $f^!\:\sD^\co((X,\L,w)\qcoh)\rarrow
\sD^\co((Y,f^*\L,\allowbreak f^*w)\qcoh)$ right adjoint to the direct
image functor $\boR f_*\:\sD^\co((Y,f^*\L,f^*w)\qcoh)\allowbreak
\rarrow\sD^\co((X,\L,w)\qcoh)$ (see the end of
Section~\textup{\ref{pull-push}}).
\end{thm}

\begin{proof}
 For any quasi-coherent matrix factorization $\N$ on $Y$ and any
flat quasi-coherent matrix factorization $\E$ on $X$ we have to
construct an isomorphism
$$
 \Hom_{\sD^\co((X,\L,w)\qcoh)}(\boR f_*\N\;\D_X^\bu\ot_{\O_X}\E)
 \simeq\Hom_{\sD^\co((Y,f^*\L,f^*w)\qcoh)}
 (\N\;\D_Y^\bu\ot_{\O_Y}f^*\E).
$$
 The composition $\Hom_Y(\N\;\D_Y^\bu\ot_{\O_Y}f^*\E)\rarrow
\Hom_X(\boR f_*\N\;\boR f_*(\D_Y^\bu\ot_{\O_Y}f^*\E)) \simeq
\Hom_X(\boR f_*\N\;f_*\D_Y^\bu\ot_{\O_X}\E)\rarrow
\Hom_X(\boR f_*\N\;\D_X^\bu\ot_{\O_X}\E)$ provides a morphism
from the right-hand to the left-hand side.
 Here all the Hom functors are taken in the coderived categories
of quasi-coherent matrix factorizations on $Y$ and $X$; the middle
isomorphism holds since $\D_Y^\bu\ot_{\O_Y}f^*\E$ is an injective
matrix factorization on $Y$ (so the derived direct image can be
computed for it by applying the underived direct image
functor~$f_*$ termwise) and by the projection formula; the last
morphism is induced by the adjunction $f_*\D_Y^\bu\rarrow\D_X^\bu$.

 Furthermore, on both sides of the desired isomorphism we have
injective matrix factorizations in the second arguments of the Hom
functors; hence the Hom can be computed in the homotopy category
of matrix factorizations instead of the coderived category in
both cases.
 Finally, one can assume $\N$ to be an injective matrix
factorization, too, and compute $\boR f_*\N=f_*\N$ termwise
(alternatively, one could use the \v Cech construction).
 Similarly, the tensor products in the second arguments are
totalizations of termwise tensor products.

 Now one can fix the components involved for both matrix
factorizations $\N$ and $\E$, obtaining a morphism of finite
complexes of abelian groups of the same kind as above, but
related to (one-term) complexes of quasi-coherent sheaves rather
than matrix factorizations.
 The latter is an isomorphism by~\cite[Theorem~5.15.3]{Pcosh}.
 It remains to notice that the totalization of an acyclic finite
complex of (unbounded) complexes of abelian groups is acyclic.
\end{proof}

 The next corollary is a matrix factorization version of
the main result of Deligne's appendix to~\cite{Har}
(see also~\cite[Section~5.16]{Pcosh}).

\begin{cor}
 For any morphism of finite type between separated Noetherian
schemes with dualizing complexes $f\:Y\rarrow X$, a line bundle
$\L$ on $X$, and a section $w\in\L(X)$, one can define
a triangulated functor $f^+\:\sD^\co((X,\L,w)\qcoh)\rarrow
\sD^\co((Y,f^*\L,f^*w))\qcoh)$ in such a way that \par
\textup{(i)} for an open embedding~$f$, one has $f^+=f^*$,
and more generally, for a smooth morphism~$f$ of relative
dimension~$n$ one has $f^+=\omega_{Y/X}[n]\ot_{\O_Y}f^*$; \par
\textup{(ii)} for a proper morphism~$f$, the functor $f^+=f^!$
is right adjoint to~$\boR f_*$; \par
\textup{(iii)} the construction is compatible with
the compositions of the morphisms~$f$.
\end{cor}

\begin{proof} \emergencystretch=0em
 It suffices to define $f^+\:\sD^\co((X,\L,w)\qcoh)\rarrow
\sD^\co((Y,f^*\L,f^*w)\qcoh)$ as the functor corresponding to
the inverse image of flat quasi-coherent matrix factorizations
$f^*\:\sD^\abs((X,\L,w)\qcoh_\fl)\rarrow\sD^\abs((Y,f^*\L,
\allowbreak f^*w) \qcoh_\fl)$ under the identifications of
categories $\D_X^\bu\otimes_{\O_X}\nobreak{-}\:\sD^\abs((X,\L,w)
\qcoh_\fl)\allowbreak\rarrow\sD^\co((X,\L,w)\qcoh)$ and
$\D_Y^\bu\otimes_{\O_Y}{-}\:\sD^\abs((Y,f^*\L,f^*w)
\qcoh_\fl)\rarrow\sD^\co((Y,f^*\L,\allowbreak f^*w)\qcoh)$,
where $\D_X^\bu$ is any dualizing complex on $X$ and
$\D_Y^\bu=f^+\D_X^\bu$.
\end{proof}

\appendix \bigskip
\Section{Quasi-Coherent Graded Modules}
\medskip

\subsection{Flat quasi-coherent sheaves}  \label{flat-sheaves}
 I am grateful to A.~Neeman for suggesting that a result of
the following kind can be proven without much difficulty.

\begin{lem}
 On any quasi-compact semi-separated scheme, any quasi-coherent
sheaf is the quotient sheaf of a flat quasi-coherent sheaf. 
\end{lem}

\begin{proof}
 Let $X$ be our scheme.
 Assume that a quasi-coherent sheaf $\M$ over $X$ is flat over
an open subscheme $V\subset X$; given an affine open subscheme
$U\subset X$, we will construct a surjective morphism $\N\rarrow\M$
onto $\M$ from a quasi-coherent sheaf $\N$ over $X$ that is flat
over $U\cup V$.
 Let $j$ denote the embedding $U\rarrow X$.
 There exists a surjective morphism onto $j^*\M$ from a flat
quasi-coherent sheaf $\F$ over~$U$; let $\K$ denote the kernel
of this morphism of sheaves.

 The morphism $j\:U\rarrow X$ being affine and flat, the functor
$j_*$ is exact and preserves flatness.
 Consider the pull-back of the exact triple $j_*\K\rarrow j_*\F
\rarrow j_*j^*\M$ with respect to the morphism $\M\rarrow j_*j^*\M$;
denote the middle term of the resulting exact triple by~$\N$.
 One has $\N|_U=\F|_U$, so $\N$ is flat over~$U$.
 Furthermore, the sheaf $j^*\M$ is flat over~$V\cap U$, hence so
is the sheaf~$\K$.
 The embedding $U\cap V\rarrow V$ is an affine flat morphism, so
the sheaf $j_*\K$ is flat over~$V$.
 From the exact triple $j_*\K\rarrow\N\rarrow\M$ we conclude that
$\N$ is flat over~$V$.
\end{proof}

 It follows immediately that any quasi-coherent graded module over
a quasi-coherent graded algebra $\B$ over $X$ is a quotient module
of a flat quasi-coherent graded module.

\subsection{Locally projective quasi-coherent graded modules}
\label{locally-projective}
 The following result is essentially due to Raynaud and
Gruson~\cite{RG} (for a discussion, see~\cite[Section~2]{Dr});
here we just briefly explain how to deduce the formulation that
interests us from their assertions.

\begin{thm}
 Let $X$ be an affine scheme and $\{U_\alpha\}$ be its finite affine
covering.
 Let $\B$ be a quasi-coherent graded algebra over $X$ and $\P$ be
a quasi-coherent graded module over~$\B$.
 Then the graded $\B(X)$\+module $\P(X)$ is projective if and only
if the graded $\B(U_\alpha)$\+module $\P(U_\alpha)$ is projective
for every~$\alpha$.
\end{thm}

\begin{proof}
 First of all, a graded module $P$ over a graded ring $B$ is
projective if and only if it is projective as an ungraded module.
 Indeed, if $P$ is graded projective, then it is a homogeneous
direct summand of a free graded module, hence $P$ is also ungraded
projective.
 Conversely, pick a homogeneous (of degree~$0$) surjective
homomorphism $F\rarrow P$ onto a given graded module $P$ from
a free graded module~$F$.
 If $P$ is ungraded projective, this homomorphism has a (perhaps
nonhomogeneous) section~$s$, and the homogeneous component
of~$s$ of degree~$0$ provides a homogeneous section.
 Hence it suffices to consider ungraded modules over an ungraded
quasi-coherent algebra~$\B$.

 It is clear that if $\P(X)$ is a projective $\B(X)$\+module,
then $\P(V)$ is a projective $\B(V)$\+module for any affine open
subscheme $V\subset X$.
 Conversely, assume that the $\B(U_\alpha)$\+module $\P(U_\alpha)$
is projective for every~$\alpha$.
 Then by the result of~\cite{Ka} the $\B(U_\alpha)$\+modules
$\P(U_\alpha)$ are direct sums of countably generated modules,
and it follows easily that so is the $\B(X)$\+module $\P(X)$
(essentially, since a connected graph with an at most countable set
of edges at each vertex has a countable number of vertices).
 Hence we can assume the $\B(X)$\+module $\P(X)$ to be countably
generated.

 Besides, the $\B(U_\alpha)$\+modules $\P(U_\alpha)$ are flat,
hence so is the $\B(X)$\+module $\P(X)$.
 By~\cite[Corollaire~II.2.2.2]{RG}, it remains to show that
the $\B(X)$\+module $\P(X)$ satisfies the Mittag-Leffler
condition; this can be easily deduced from the similar property
of the $\B(U_\alpha)$\+modules $\P(U_\alpha)$ using
the formulation of this condition given in Proposition~II.2.1.4(iii)
or Propositions~II.2.1.4(ii) and~II.2.1.1(i) of~\cite{RG}
(cf.\ Sections~II.2.5 and~II.3.1 of the same paper).
\end{proof}

\subsection{Injective quasi-coherent graded modules}
\label{injective-sheaves}
 The following result is a noncommutative generalization of
a theorem of Hartshorne~\cite[Theorem~II.7.18]{Har}
about injective quasi-coherent sheaves on Noetherian schemes.
 Our proof method, based on the Artin--Rees lemma, is different
from the one in~\cite{Har}.

\begin{thm}
 Let $\B$ be a Noetherian quasi-coherent graded algebra over
a Noetherian scheme~$X$.
 Then any injective object in the category of quasi-coherent graded
left modules over $\B$ is also an injective object of the category of
arbitrary sheaves of graded $\B$\+modules over~$X$.

 Consequently, the restriction $\J|_U$ of an injective quasi-coherent
graded module $\J$ over $\B$ to an open subscheme $U\subset X$ is
an injective quasi-coherent graded module over $\B|_U$.
 Conversely, if $U_\alpha$ is an open covering of $X$ and
the quasi-coherent graded $\B|_{U_\alpha}$\+modules $\J|_{U_\alpha}$
are injective, then a quasi-coherent graded $\B$\+module $\J$
is injective.
 Besides, the underlying sheaf of graded abelian groups of any
injective quasi-coherent graded $\B$\+module $\J$ is flabby. 
\end{thm}

\begin{proof}
 First of all, notice that the abelian category $\B\qcoh$ of
quasi-coherent graded modules over $\B$ is a locally Noetherian
Grothendieck category with coherent graded modules forming
the subcategory of Noetherian generators~\cite[Exercise~II.5.15]{Har2};
so in particular $\B\qcoh$ has enough injectives and the assertions
of Theorem are not vacuous.
 The category of sheaves of graded $\B$\+modules $\B\modl$ has
similar properties, with the extensions by zero of the restrictions
of $\B$ to (small) open subschemes of $X$ forming a set
of Noetherian generators~\cite[Theorem~II.7.8]{Har}.

 Secondly, let us check that the main result in the first
paragraph implies the assertions in the second one.
 Indeed, injective sheaves of graded $\B$\+modules have all
the properties we are interested in.
 They remain injective after being restricted to an open subscheme,
since the extension by zero from an open subscheme is an exact
functor.
 They are flabby, since given two open subschemes $U\subset V
\subset X$ and $j_U$, $j_V$ being their identity embeddings
$U$, $V\rarrow X$, the morphism of sheaves of graded $\B$\+modules
$j_U{}_!\B|_U\rarrow j_V{}_!\B|_V$ is injective.
 And their property is local~\cite[Lemma~II.7.16]{Har}, because sheaves
of graded $\B$\+modules supported inside one of the subschemes
$U_\alpha$ form a set of generators of the category $\B\modl$.

 Now let $\J$ be an injective quasi-coherent graded module over~$\B$.
 To prove the main assertion, we have to show that for any open
subscheme $U\subset X$ and a subsheaf of graded $\B$\+modules
$\G\subset j_U{}_!\B|_U$, any homogeneous morphism of sheaves of
graded $\B$\+modules $\G\rarrow\J$ can be extended to a similar
morphism $j_U{}_!\B|_U\rarrow\J$.
 Indeed, $\G$ is a subsheaf of graded $\B$\+modules in
the coherent graded $\B$\+module $\B$, hence according to
the following proposition there exists a quasi-coherent graded
$\B$\+module $\G\subset\F\subset\B$ such that the morphism
$\G\rarrow\J$ can be extended to a homogeneous morphism of
quasi-coherent graded $\B$\+modules $\F\rarrow\J$.

 Since $\J$ is injective in $\B\qcoh$, the latter morphism can in turn
be extended to a similar morphism $\B\rarrow\J$.
 Restricting to $j_U{}_!\B|_U$, we obtain the desired morphism of
sheaves of graded $\B$\+modules $j_U{}_!\B|_U\rarrow\J$.
\end{proof}

\begin{prop}
 In the assumptions of Theorem, let $\E$ be a coherent graded left
$\B$\+module, $\G\subset\E$ be a subsheaf of graded $\B$\+modules,
$\M$ be a quasi-coherent graded $\B$\+module, and
$\phi\:\G\rarrow\M$ be a morphism of sheaves of graded $\B$\+modules.
 Then there exists a coherent graded $\B$\+module
$\G\subset\F\subset\E$ such that the morphism~$\phi$ can be
extended to~$\F$.
\end{prop}

\begin{proof}
 Before proving Proposition, let us reformulate its conclusion as
follows.
 In the same setting, there exists a quasi-coherent graded 
$\B$\+module $\K$ together with an injective morphism $\M\rarrow\K$
and a morphism $\E\rarrow\K$ forming a commutative diagram with
the embedding $\G\rarrow\E$ and the morphism $\phi\:\G\rarrow\M$.
 Indeed, if a coherent $\B$\+module $\F$ exists, one can take $\K$
to be the fibered coproduct of $\E$ and $\M$ over~$\F$;
conversely, if a quasi-coherent $\B$\+module $\K$ exists, one can
take $\F$ to be the full preimage of $\M\subset\K$ under
the morphism $\E\rarrow\K$.
 Notice also that one can always replace $\M$ with its sufficiently
big coherent graded $\B$\+submodule.

 Now let us state the version of Artin--Rees lemma that we will use.

\begin{lem}
 In the assumptions of Theorem, let $\M$ be a coherent graded
$\B$\+module, $\N\subset\M$ a coherent graded $\B$\+submodule,
and $Z\subset X$ a closed subscheme with the sheaf of ideals
$\I_Z\subset\O_X$.
 Then for any $n\ge0$ there exists $m\ge0$ such that the intersection
$\I_Z^m\M\cap\N$ is contained in $\I_Z^n\N$.
\end{lem}

\begin{proof}
 Clearly, the question is local, so it suffices to consider the case
of an affine scheme~$X$.
 Then (the graded version of) the Artin--Rees lemma for
ideals generated by central elements in noncommutative Noetherian
rings~\cite[Theorem~13.3]{GW} applies.
\end{proof}

 Being a Noetherian object, the sheaf of graded $\B$\+modules $\G$
is generated by a finite number of homogeneous sections
$s_n\in\G(U_n)$, where $U_n\subset X$ are some open subschemes.
 If all of these subschemes coincide with $X$, the sheaf $\G$,
being a subsheaf of a coherent sheaf generated by global sections,
is itself coherent, so there is nothing to prove.
 In the general case, we will argue by induction in the number of
open subschemes $U_n$ that are not equal to~$X$.

 Let $U=U_1\subsetneq X$ be one such open subscheme, and
$T=X\setminus U$ be its closed complement.
 We can assume that $\M$ is a coherent graded $\B$\+module.
 Let $\N$ denote its maximal coherent graded $\B$\+submodule
supported set-theoretically in~$T$.
 Applying Lemma to $\N\subset\M$, we conclude that there is a closed
subscheme structure $i\:Z\rarrow X$ on $T$ such that the morphism
$\N\rarrow i_*i^*\M$ is injective.
 Consequently, so is the morphism $\M\rarrow i_*i^*\M\oplus j_*j^*\M$,
where $j$~denotes the open embedding $U\rarrow X$.

 Let us show that there is a thicker closed subscheme structure
$i'\:Z'\rarrow X$ on $T$ such that the kernel of the morphism of
sheaves $i'_*i'{}^*\G\rarrow i'_*i'{}^*\E$ is contained in
the kernel of the morphism of sheaves $i'_*i'{}^*\G\rarrow
i_*i^*\G$.
 Indeed, there exists a finite collection of subsheaves of graded
$\B$\+modules in $\G$, each of them an extension by zero of
a coherent graded $\B|_V$\+module from some open subscheme
$V\subset X$, such that the stalk of $\G$ at each point of $X$
coincides with the stalk of one of these subsheaves.
 So the assertion reduces to the case when $\G$ is a coherent
graded $\B$\+submodule in $\E$, when it is an equivalent
reformulation of Lemma.

 Let $\H\subset i'{}^*\E$ denote the image of the morphism of
sheaves of graded $i'{}^*\B$\+modules $i'{}^*\G\rarrow i'{}^*\E$
over the scheme~$Z'$.
 Let $\iota\:Z\rarrow Z'$ be the natural closed embedding.
 Then, according to the above, the morphism of sheaves of
graded $i'{}^*\B$\+modules $i'{}^*\G\rarrow\iota_*i^*\G$ induces
a morphism $\H\rarrow\iota_*i^*\G$.

 The sheaf of graded $i'{}^*\B$\+modules $\H$ is generated by
the images of the restrictions of the sections $s_n$, \ $n\ge2$,
to the closed subschemes $Z'\cap U_n\subset U_n$.
 Hence the induction assumption is applicable to $\H$, and we can
conclude that there exists a quasi-coherent graded
$i'{}^*\B$\+module $\K$ on the scheme $Z'$ together with
an injective morphism $\iota_*i^*\M\rarrow\K$ and a morphism
$i'{}^*\E\rarrow\K$ forming a commutative diagram with
the embedding $\H\rarrow i'{}^*\E$ and the composition
$\H\rarrow\iota_*i^*\G\rarrow\iota_*i^*\M$.

 Similarly, the sheaf of graded $\B|_U$\+modules $j^*\G$ is generated
by the restrictions of the sections~$s_n$ to the open subschemes
$U_1\cap U_n\subset U_n$, among which the (restriction of)
the section~$s_1$ is a global section over $U=U_1$.
 Hence the induction assumption is applicable to $j^*\G$, and
there exists a quasi-coherent graded $\B|_U$\+module $\L$
together with an injective morphism $j^*\M\rarrow\L$ and
a morphism $j^*\E\rarrow\L$ forming a commutative diagram with
the embedding $j^*\G\rarrow j^*\E$ and the morphism
$j^*\G\rarrow j^*\M$.

 Now the injective morphism $\M\rarrow i'_*\K\oplus j_*\L$
(whose first component is the composition $\M\rarrow
i_*i^*\M\simeq i'_*\iota_*i^*\M\rarrow i'_*\K$) and
the morphism $\E\rarrow i'_*\K\oplus j_*\L$ provide the desired
commutative diagram of morphisms of sheaves of graded
$\B$\+modules over~$X$.
\end{proof}

\Section{Hochschild (Co)Homology of Matrix Factorizations}
\medskip

 This appendix complements the paper~\cite{PP} in two ways.
 Section~\ref{loc-free-hochschild} contains some modifications
and improvements of the main results of~\cite{PP} generally, and
as applied to locally free matrix factorizations of finite rank
in particular.
 The main thrust consists in replacing the finite homological
dimension conditions in~\cite{PP} with the Noetherianness conditions
to the (limited) extent possible.

 Section~\ref{coherent-hochschild}, on the other hand, presents
an elementary approach to computation of Hochschild (co)homology
of coherent matrix factorizations, entirely unrelated to that
in~\cite{PP} and not based on any notion of Hochschild
(co)homology of the second kind, but rather on
the Serre--Grothendieck duality theory.

\subsection{Locally free matrix factorizations of finite rank}
\label{loc-free-hochschild}
 In Sections~\ref{coh-noetherian}\+-\ref{dg-of-cdg},
we start with a bit of categorical nonsense, following the lines
of~\cite[Sections~3.3\+-3.5]{PP}, but with the additional
coherence/Noetherianness conditions imposed from the very beginning.
 We use the notation from~\cite{PP} rather than that of the main body
of this paper.
 Then in Section~\ref{affine-loc-free} we turn to locally free matrix
factorizations of finite rank over certain possibly singular,
affine algebraic varieties.
 The final Section~\ref{smooth-strat} presents an improvement over
the discussion of matrix factorizations over smooth affine varieties
in~\cite[Section~4.8]{PP}.
 An example of application of our techniques to nonaffine
varieties can be found in the preprint~\cite{Ef}.
 
\subsubsection{Coherent and Noetherian CDG-categories}
\label{coh-noetherian}
 Let $(\Gamma,\sigma,\boldsymbol{1})$ be a grading group
data~\cite[Section~1.1]{PP} and $B^\#$ be a small $\Gamma$\+graded
preadditive category~\cite[Section~A.1]{Partin}.
 Both left and right $\Gamma$\+graded $B^\#$\+modules form
abelian categories.

 A $\Gamma$\+graded $B^\#$\+module is said to be \emph{finitely
generated} (respectively, \emph{finitely presented}) if it is
a quotient module of a finitely generated free $\Gamma$\+graded
$B^\#$\+module \cite[Section~1.5]{PP} (respectively, the cokernel of 
a morphism of finitely generated free $\Gamma$\+graded
$B^\#$\+modules).

 A $\Gamma$\+graded preadditive category $B^\#$ is called \emph{left
Noetherian} if any submodule of a finitely generated
$\Gamma$\+graded left $B^\#$\+module is finitely generated, or
equivalently, if the abelian category of $\Gamma$\+graded left
$B^\#$\+modules is locally Noetherian.
 A $\Gamma$\+graded preadditive category $B^\#$ is called \emph{left
coherent} if any submodule of a finitely presented
$\Gamma$\+graded left $B^\#$\+module is finitely presented.

 Let $B$ be a small ($\Gamma$\+graded)
CDG\+category~\cite[Section~1.2]{PP} and $B^\#$ be its underlying
$\Gamma$\+graded preadditive category.
 Following~\cite{PP}, we denote the DG\+categories of left and right
CDG\+modules over $B$ by $B\modlc$ and $\modrc B$.
 The DG\+subcategories of left CDG\+modules whose underlying
$\Gamma$\+graded $B^\#$\+modules are flat or injective are denoted by
$B\modlc_\rfl$ and $B\modlc_\rinj\subset B\modlc$.
 Similarly, the DG\+subcategories of left and right CDG\+modules
over $B$ whose underlying $\Gamma$\+graded $B^\#$\+modules are
projective and finitely generated are denoted by $B\modlc_\rfgp$
and $\modrcfgp B$.

 Assuming that the $\Gamma$\+graded category $B^\#$ is left Noetherian,
the DG\+subcategory of left CDG\+modules whose underlying
$\Gamma$\+graded $B^\#$\+modules are finitely generated is denoted
by $B\modlc_\rfg\subset B\modlc$.
 Assuming that the $\Gamma$\+graded category $B^\#$ is right coherent,
the DG\+subcategory of right CDG\+modules whose underlying
$\Gamma$\+graded $B^\#$\+modules are finitely presented is denoted
by $\modrcfp B$.

 The coderived and contraderived categories of left CDG\+modules 
over $B$ are denoted by $\rD^\rco(B\modlc)$ and $\rD^\rctr(B\modlc)$,
respectively~\cite[Section~3.2]{PP}.
 Assuming that the $\Gamma$\+graded category $B^\#$ is right coherent,
the class of flat $\Gamma$\+graded left
$B$\+modules~\cite[Section~2.2]{PP} is closed under infinite products,
so the contraderived category $\rD^\rctr(B\modlc_\rfl)$ is well-defined.
 The homotopy category of the DG\+category $B\modlc_\rinj$ is
denoted, as usually, by $H^0(B\modlc_\rinj)$.

 In the respective assumptions of left Noetherianness or right
coherence of the $\Gamma$\+graded category $B^\#$, the absolute
derived categories of CDG\+modules with finitely generated or finitely
presented underlying $\Gamma$\+graded $B^\#$\+modules are denoted by
$\rD^\rabs(B\modlc_\rfg)$ and $\rD^\rabs(\modrcfp B)$, respectively.

\subsubsection{Derived functors of the second kind}
\label{coh-noeth-second-kind}
 Let $k$ be a commutative ring and $B$ be a small $k$\+linear
CDG\+category.
 Assume that the $\Gamma$\+graded category $B^\#$ is left Noetherian.
 Let $L$ and $M$ be left CDG\+modules over~$B$; suppose that
the $\Gamma$\+graded left $B^\#$\+module $L^\#$ underlying
the CDG\+module $L$ over $B$ is finitely generated.

 As in~\cite[Sections~2.1\+2]{PP}, we denote by $Z^0(B\modlc)$ and
$Z^0(\modrc B)$ the abelian categories of left and right
CDG\+modules over~$B$.
 Let $Z^0(B\modlc_\rfg)\subset Z^0(B\modlc)$ and
$H^0(B\modlc_\rfg)\subset H^0(B\modlc)$ denote the abelian and
homotopy categories of left CDG\+modules over $B$ with finitely
generated underlying $\Gamma$\+graded $B^\#$\+modules, and
$Z^0(\modrcfp B)\subset Z^0(\modrc B)$ and $H^0(\modrcfp B)\subset
H^0(\modrc B)$ be the similar categories of right CDG\+modules
with finitely presented underlying $\Gamma$\+graded modules.

 Let $J^\bu$ be a right resolution of $M$ in $Z^0(B\modlc)$ such that
the $\Gamma$\+graded left $B^\#$\+modules $J^i{}^\#$ are injective,
and let $J$ be the total CDG\+module of the complex of CDG\+modules
$J^\bu$ constructed by taking infinite direct sums along the diagonals.
 Then the complex $\Tot^\oplus\Hom^B(L,J^\bu)$ computing
$\Ext^{I\!I}_B(L,M)$ \cite[Section~2.2]{PP} is isomorphic to
the complex $\Hom^B(L,J)$ \cite[formula~(6)]{PP}, which computes
the $k$\+modules of morphisms from $L$ into $M[*]$ in the coderived
category $\rD^\rco(B\modlc)$ \cite[Theorems~3.5(a) and~3.7]{Pkoszul}.
 Thus,
$$
 H^*\Ext^{I\!I}_B(L,M)\.\simeq\.\Hom_{\rD^\rco(B\modlc)}(L,M[*]).
$$

 Just as in~\cite[Section~3.3]{PP}, one can lift this isomorphism
from the level of cohomology modules to that of the derived
category $\rD(k\rmodl)$ in the following way.
 Consider the functor
$$
 \Hom^B\:H^0(B\modlc)^\op\times H^0(B\modlc)\lrarrow
 \rD(k\rmodl)
$$
and restrict it to the full subcategory $H^0(B\modlc_\rinj)$ in
the second argument.
 This restriction factorizes through the coderived category
$\rD^\rco(B\modlc)$ in the first argument.
 Taking into account~\cite[Theorem~3.7]{Pkoszul}, we obtain
a right derived functor
$$
 \rD^\rco(B\modlc)^\op\times\rD^\rco(B\modlc)\lrarrow
 \rD(k\rmodl).
$$
 Restricting to the full subcategory $\rD^\rabs(B\modlc_\rfg)^\op
\subset\rD^\rco(B\modlc)^\op$ \cite[Theorem~3.11.1]{Pkoszul} in
the first argument, we have the derived functor
\begin{equation}  \label{ext-second-kind}
 \rD^\rabs(B\modlc_\rfg)^\op\times\rD^\rco(B\modlc)\lrarrow
 \rD(k\rmodl).
\end{equation}
 The composition of this functor with the localization
functors $Z^0(B\modlc_\rfg)\rarrow\rD^\rabs(B\modlc_\rfg)$
and $Z^0(B\modlc)\rarrow\rD^\rco(B\modlc)$ agrees with the derived
functor $\Ext^{I\!I}_B$ where the former is defined.

 Now assume that the $\Gamma$\+graded category $B^\#$ is
right coherent.
 Consider the functor~\cite[formula~(5)]{PP}
$$
 \ot_B\: H^0(\modrc B)\times H^0(B\modlc)\lrarrow\rD(k\rmodl)
$$
and restrict it to the Cartesian product of full subcategories
$H^0(\modrcfp B)\times H^0(B\modlc_\fl)\subset
H^0(\modrc B)\times H^0(B\modlc)$.
 Since the tensor product with a finitely presented $\Gamma$\+graded
right $B^\#$\+module commutes with infinite products of
$\Gamma$\+graded left $B^\#$\+modules, this restriction factorizes
through the contraderived category $\rD^\rctr(B\modlc_\fl)$ in
the second argument.
 Clearly, it also factorizes through the absolute derived category
$\rD^\rabs(\modrcfp B)$ in the first argument.

 By Remark~\ref{embedding-prop} of the main body of this paper
(see also~\cite[Proposition~A.3.1(b)]{Pcosh}), the natural functor
$\rD^\rctr(B\modlc_\fl)\rarrow\rD^\rctr(B\modlc)$ is an equivalence
of triangulated categories.
 Hence we obtain a left derived functor
\begin{equation} \label{tor-second-kind}
 \rD^\rabs(\modrcfp B)\times\rD^\rctr(B\modlc)\lrarrow
 \rD(k\rmodl).
\end{equation}
 Up to composing with the localization functors $Z^0(\modrcfp B)
\rarrow \rD^\rabs(\modrcfp B)$ and $Z^0(B\modlc)\rarrow
\rD^\rctr(B\modlc)$, this functor agrees with the derived functor
$\Tor^{B,I\!I}$ from~\cite[Section~2.2]{PP} where the former
is defined.

 Indeed, let $N$ be an object of $Z^0(\modrcfp B)$.
 Let $P_\bu$ be a left resolution of an object $M\in Z^0(B\modlc)$
by left CDG\+modules over $B$ with flat underlying $\Gamma$\+graded
$B^\#$\+modules, and let $P$ be the total CDG\+module of the complex
$P_\bu$ constructed by taking infinite products along the diagonals.
 Then the complex $\Tot^\sqcap(N\ot_B P_\bu)$ computing 
$\Tor^{B,I\!I}(N,M)$ is isomorphic to the complex $N\ot_B P$
computing the derived functor~\eqref{tor-second-kind} on
the objects $N$ and~$M$.

\subsubsection{Comparison of the two theories}
\label{comparison-subsubsect}
 Let $C$ be a small $k$\+linear ($\Gamma$\+graded) DG\+category.
 The above constructions applicable to CDG\+categories and
CDG\+mod\-ules over them can be also applied to DG\+categories
and DG\+modules as a particular case.
 Following~\cite{PP}, we denote the DG\+categories of left and right
DG\+modules over $C$ by $C\modld$ and $\modrd C$, and generally
use the upper index ``$\mathrm{dg}$'' instead of ``$\mathrm{cdg}$''
in the notation related to DG\+modules.

 As in~\cite[Sections~2.1, 3.1 and~3.4]{PP}, we denote by
$H^0(C\modld)_\rinj$ and $H^0(C\modld)_\rfl$ the homotopy categories
of h\+injective and h\+flat left DG\+modules over~$C$.
 The notation $H^0(C\modld_\rinj)_\rinj$ and $H^0(C\modld_\rfl)_\rfl$
stands for the full triangulated subcategories in $H^0(C\modld)$ formed
by h\+injective DG\+modules over $C$ whose underlying $\Gamma$\+graded
$C^\#$\+modules are injective, or h\+flat DG\+modules whose underlying
$\Gamma$\+graded $C^\#$\+modules are flat, respectively.
 Finally, let $H^0(C\modld_\rfgp)_\rprj\allowbreak\subset H^0(C\modld)$
and $H^0(\modrdfgp C)_\rfl\subset H^0(\modrd C)$ denote the full
triangulated subcategories of h\+projective left and h\+flat right
DG\+modules whose underlying $\Gamma$\+graded $C^\#$\+modules are
projective and finitely generated. {\hbadness=1250\par}

 Assume that the $\Gamma$\+graded category $C^\#$ is left Noetherian.
 Let $L$ be an object of $Z^0(C\modld_\rfg)$.
 Given a left DG\+module $M$ over $C$, pick its injective resolution
$J^\bu$ in the exact category $Z^0(C\modld)$ \cite[Section~2.1]{PP}.
 Let $\Tot^\oplus(J^\bu)\rarrow\Tot^\sqcap(J^\bu)$ be the natural
closed morphism between the total DG\+modules of the complex $J^\bu$
constructed by taking infinite direct sums and infinite products
along the diagonals.
 Then the induced morphism of complexes of $k$\+modules
$$
 \Hom^C(L\;\Tot^\oplus(J^\bu))\lrarrow
 \Hom^C(L\;\Tot^\sqcap(J^\bu))
$$
represents the comparison morphism $\Ext_C^{I\!I}(L,M)\rarrow
\Ext_C(L,M)$ \cite[formula~(10)]{PP} in $\sD(k\rmodl)$ between
the two kinds of $\Ext$ objects for the DG\+modules $L$ and~$M$.

 Similarly, assume that the $\Gamma$\+graded category $C^\#$ is
right coherent.
 Let $N$ be an object of $Z^0(\modrdfp C)$.
 Given a left DG\+module $M$ over $C$, pick its projective
resolution $P_\bu$ in the exact category $Z^0(C\modld)$.
 Let $\Tot^\oplus(P_\bu)\rarrow\Tot^\sqcap(P_\bu)$ be the natural
closed morphism between the total DG\+modules of the complex $P_\bu$
constructed by taking infinite direct sums and infinite products
along the diagonals.
 Then the induced morphism of complexes of $k$\+modules
$$
 N\ot_C \Tot^\oplus(P_\bu)\lrarrow N\ot_C\Tot^\sqcap(P_\bu)
$$
represents the comparison morphism $\Tor^C(N,M)\rarrow
\Tor^{C,I\!I}(N,M)$ \cite[formula~(9)]{PP} in $\rD(k\rmodl)$
between the two kinds of $\Tor$ objects for the DG\+modules
$N$ and~$M$.

\begin{propa}
 Assume that the\/ $\Gamma$\+graded category $C^\#$ is
left Noetherian.
 Let $L$ be a left DG\+module over $C$ whose underlying\/
$\Gamma$\+graded left $C^\#$\+module is finitely generated,
and let $M$ be a left DG\+module over~$C$.
 Then the natural morphism\/ $\Ext_C^{I\!I}(L,M)\rarrow
\Ext_C(L,M)$ is an isomorphism provided that either \par
\textup{(i)} the object $M\in\rD^\rco(C\modld)$ belongs to
the image of the fully faithful functor $H^0(C\modld_\rinj)_\rinj
\rarrow\rD^\rco(C\modld)$; or \par
\textup{(ii)} the object $L\in\rD^\rabs(C\modld)$ belongs to
the image of the fully faithful functor $H^0(C\modld_\rfgp)_\rprj
\rarrow\rD^\rabs(C\modld_\rfg)$.
\end{propa}

\begin{proof}
 Let $J^\bu$ be an injective resolution of the DG\+module $M$ in
the exact category $Z^0(C\modld)$.
 Then the natural morphism $M\rarrow\Tot^\oplus(J^\bu)$ is always
an isomorphism in $\rD^\rco(C\modld)$ \cite[proof of
Theorem~3.7]{Pkoszul}, while the morphism $M\rarrow\Tot^\sqcap(J^\bu)$
is an isomorphism in the conventional derived category
$\rD(C\modld)$ \cite[proofs of Theorems~1.4\+5]{Pkoszul}.
 Furthermore, one has $\Tot^\oplus(J^\bu)\in H^0(C\modld_\rinj)$
and $\Tot^\sqcap(J^\bu)\in H^0(C\modld_\rinj)_\rinj$.

 Part~(i): the functor is fully faithful by~\cite[Theorem~3.5(a)
and Lemma~1.3]{Pkoszul}.
 According to formula~\eqref{ext-second-kind} from
Section~\ref{coh-noeth-second-kind} and~\cite[Section~3.1]{PP},
both kinds of Ext involved are well-defined as functors
of the argument $M\in\rD^\rco(C\modld)$.
 Hence one can assume $M\in H^0(C\modld_\rinj)_\rinj$.
 Then both morphisms $M\rarrow\Tot^\oplus(J^\bu)$ and
$M\rarrow\Tot^\sqcap(J^\bu)$ are homotopy equivalences by
semiorthogonality, hence so is the morphism
$\Tot^\oplus(J^\bu)\rarrow\Tot^\sqcap(J^\bu)$ and
the assertion follows.

 Part~(ii): in view of the first paragraph of this proof, 
a cone $K$ of the morphism $\Tot^\oplus(J^\bu)\rarrow\Tot^\sqcap(J^\bu)$
in $H^0(C\modld)$ is an acyclic DG\+module over $C$ whose underlying
$\Gamma$\+graded $C^\#$\+module is injective.
 Hence the complex of morphisms $\Hom^C({-},K)$ is a well-defined
functor $\rD^\abs(C\modld_\rfg)^\op\rarrow\rD(k\rmodl)$ annihilating
$H^0(C\modld_\rfgp)_\rprj$.
\end{proof}

\begin{propb}
 Assume that the\/ $\Gamma$\+graded category $C^\#$ is right coherent.
 Let $N$ be a right DG\+module over $C$ whose underlying\/
$\Gamma$\+graded right $C^\#$\+module is finitely presented,
and let $M$ be a left DG\+module over~$C$.
 Then the natural morphism\/ $\Tor^C(N,M)\rarrow
\Tor^{C,I\!I}(N,M)$ is an isomorphism provided that either \par
\textup{(i)} there is a closed morphism $P\rarrow M$ into $M$
from a DG\+module $P\in H^0(C\modld_\rfl)_\rfl$ with a cone
contraacyclic with respect to $C\modld$ or completely acyclic
with respect to $C\modld_\rfl$ (see
\textup{\cite[\textit{Sections}~3.2 \textit{and}~4.7]{PP}});
or \par
\textup{(ii)} the object $N\in\rD^\rabs(\modrd C)$ belongs to
the image of the fully faithful functor $H^0(\modrdfgp C)_\rfl
\rarrow\rD^\rabs(\modrdfp C)$.
\end{propb}

\begin{proof}
 Let $P_\bu$ be a projective resolution of the DG\+module $M$ in
the exact category $Z^0(C\modld)$.
 Then the natural morphism $\Tot^\sqcap(P_\bu)\rarrow M$ is always
an isomorphism in $\sD^\ctr(C\modld)$ \cite[proof of
Theorem~3.8]{Pkoszul}, while the morphism $\Tot^\oplus(P_\bu)
\rarrow M$ is an isomorphism in $\rD(C\modld)$ \cite[proof of
Theorem~1.4]{Pkoszul}.
 Furthermore, one has $\Tot^\sqcap(P_\bu)\in H^0(C\modld_\rfl)$
and $\Tot^\oplus(P_\bu)\in H^0(C\modld_\rfl)_\rfl$.

 Part~(i): acyclic DG\+modules in the second argument are annihilated
by the functor $\Tor^C$ by~\cite[Section~3.1]{PP}, while
contraacyclic DG\+modules in the second argument are annihilated
by the functor $\Tor^{C,I\!I}(N,{-})$ according to
the formula~\eqref{tor-second-kind}.
 The latter also applies to DG\+modules completely acyclic with
respect to $C\modld_\rfl$, since the functor of tensor product with
a finitely presented DG\+module preserves infinite direct sums
and products.
 So one can replace $M$ with $P$ and assume that $M\in
H^0(C\modld_\rfl)_\rfl$.

 Then a cone of the morphism $\Tot^\sqcap(P_\bu)\rarrow M$ is
contraacyclic with respect to $C\modld$ with a flat underlying
$\Gamma$\+graded $C^\#$\+module, hence also contraacyclic with
respect to $C\modld_\rfl$.
 On the other hand, a cone of the morphism $\Tot^\oplus(P_\bu)
\rarrow M$ is acyclic and h\+flat.
 It follows that the functor $N\ot_C{-}$ transforms both these
morphisms, and therefore also the morphism $\Tot^\oplus(P_\bu)
\rarrow\Tot^\sqcap(P_\bu)$, into quasi-isomorphisms of
complexes of $k$\+modules.

 Part~(ii): a cone $K$ of the morphism $\Tot^\oplus(P_\bu)\rarrow
\Tot^\sqcap(P_\bu)$ in $H^0(C\modld)$ is an acyclic DG\+module
over $C$ whose underlying $\Gamma$\+graded $C^\#$\+module is flat.
 Hence the tensor product ${-}\ot_C K$ is a well-defined functor
$\rD^\rabs(\modrdfp C)\rarrow\rD(k\rmodl)$ annihilating
$H^0(\modrdfgp C)_\rfl$.
\end{proof}

 In particular, assuming that the category $C^\#$ is left Noetherian,
the natural morphism $\Ext_C^{I\!I}(L,M)\rarrow\Ext_C(L,M)$ is
an isomorphism for all $L\in C\modld_\rfg$ and $M\in C\modld$
provided that the Verdier localization functor $\rD^\rco(C\modld)
\rarrow\rD(C\modld)$ is an equivalence of triangulated categories.
 Assuming that the category $C^\#$ is right coherent, the natural
morphism $\Tor^C(N,M)\rarrow\Tor^{C,I\!I}(N,M)$ is an isomorphism
for all $N\in\modrdfp C$ and $M\in C\modld$ provided that the Verdier
localization functor $\rD^\rctr(C\modld)\rarrow\rD(C\modld)$ is
an equivalence of categories, or alternatively, that any acyclic
DG\+module from $C\modld_\rfl$ is completely acyclic with respect
to $C\modld_\rfl$.

\subsubsection{Comparison for the DG-category of CDG-modules}
\label{dg-of-cdg}
 Let $B$ be a $k$\+linear CDG\+category and $C=\modrcfgp B$ be
the DG\+category of right CDG\+modules over $B$ whose underlying
$\Gamma$\+graded $B^\#$\+modules are projective and finitely
generated.
 The DG\+categories of (left or right) CDG\+modules over $B$ and
DG\+modules over $C$ are naturally equivalent~\cite[Sections~1.5
and~2.6]{PP} (as are the categories of $\Gamma$\+graded modules
over $B^\#$ and $C^\#$).
 Following~\cite[Section~3.5]{PP}, we denote by $M_C$ the DG\+module
over $C$ corresponding to a CDG\+module $M$ over~$B$.

 Let $k\spcheck$ be an injective cogenerator of the abelian
category of $k$\+modules.
 Introduce the notation $B\modlc_\rprj\subset B\modlc$ for
the DG\+category of left CDG\+modules over $B$ with projective
underlying $\Gamma$\+graded $B^\#$\+modules.
 The results below in this section are to be compared with those
from~\cite[Sections~3.5 and~4.7]{PP}.

\begin{propa}
 Assume that the\/ $\Gamma$\+graded category $B^\#$ is left Noetherian.
 Let $L$ be a left CDG\+module over $B$ whose underlying\/
$\Gamma$\+graded left $B^\#$\+module $L^\#$ is finitely generated, and
let $M$ be a left CDG\+module over~$B$.
 Then the natural morphism\/ $\Ext_C^{I\!I}(L_C,M_C)\rarrow
\Ext_C(L_C,M_C)$ is an isomorphism provided that either \par
\textup{(i)} the object $M$ belongs to the minimal triangulated
subcategory of\/ $\rD^\rco(B\modlc)$ containing the objects\/
$\Hom_k(F,k\spcheck)$ for all $F\in H^0(\modrcfgp B)$ and closed under
infinite products; or \par
\textup{(ii)} the object $L$ belongs to the minimal thick
subcategory of\/ $\rD^\rabs(B\modlc_\rfg)$ containing
the image of $H^0(B\modlc_\rfgp)$.
\end{propa}

\begin{proof}
 Part~(i): the equivalence of categories $H^0(C\modld_\rinj)_\rinj
\simeq\rD(C\modld)$ makes the embedding functor
$H^0(C\modld_\rinj)_\rinj\rarrow\rD^\rco(C\modld)$ right adjoint to
the localization functor $\rD^\rco(C\modld)\rarrow\rD(C\modld)$.
 It follows that the functor $H^0(C\allowbreak\modld_\rinj)_\rinj
\rarrow\rD^\rco(C\modld)$ preserves infinite products (also,
all infinite products exist in the coderived category, since
it is compactly generated~\cite[Theorem~3.11.2]{Pkoszul}).
 Since the category $H^0(C\modld_\rinj)_\rinj$ is the minimal
triangulated subcategory of $H^0(C\modld)$ containing the objects
$\Hom_k(F_C,k\spcheck)$ and closed under infinite
products~\cite[Theorem~1.5]{Pkoszul}, the assertion follows
from Proposition~\ref{comparison-subsubsect}.A(i).

 Part~(ii): the equivalence of absolute derived categories
$\rD^\rabs(B\modlc_\rfg)\simeq\rD^\rabs(C\modld_\rfg)$ takes
objects of the full subcategory $H^0(B\modlc_\rfgp)\subset
\rD^\rabs(B\modlc_\rfg)$ to representable (and, consequently,
perfect and h\+projective) DG\+modules in $H^0(C\modld_\rfgp)
\subset\rD^\rabs(C\modld_\rfg)$, so it remains to apply
Proposition~\ref{comparison-subsubsect}.A(ii).
\end{proof}

\begin{propb}
 Assume that the\/ $\Gamma$\+graded category $B^\#$ is right coherent.
 Let $N$ be a right CDG\+module over $B$ whose underlying\/
$\Gamma$\+graded right $B^\#$\+module $N^\#$ is finitely presented,
and let $M$ be a left CDG\+module over~$B$.
 Then the natural morphism\/ $\Tor^C(N_C,M_C)\rarrow
\Tor^{C,I\!I}(N_C,M_C)$ is an isomorphism provided that either \par
\textup{(i)} the object $M$ belongs to the minimal triangulated
subcategory of\/ $H^0(B\modlc_\rprj)\allowbreak\subset
\rD^\rctr(B\modlc)$ containing the image of $H^0(B\modlc_\rfgp)$ and
closed under infinite direct sums; or \par
\textup{(ii)} the object $N$ belongs to the minimal thick
subcategory of\/ $\rD^\rabs(\modrcfp B)$ containing
the image of $H^0(\modrcfgp B)$.
\end{propb}

\begin{proof}
 Similar to that of Proposition~A and based on
Proposition~\ref{comparison-subsubsect}.B.
\end{proof}

 Now assume that the commutative ring~$k$ has finite weak homological
dimension and all the $\Gamma$\+graded $k$\+modules of morphisms in
the category $B^\#$ are flat.
 Clearly, the DG\+categories of left and right CDG\+modules over
the CDG\+category $B\ot_k B^\op$ are naturally equivalent, as are
the DG\+categories of left and right DG\+modules over
the DG\+category $C\ot_k C^\op$.
 The DG\+category of CDG\+modules over $B\ot_k B^\op$ is also
naturally equivalent to the DG\+category of DG\+modules over
$C\ot_k C^\op$ \cite[Section~2.6]{PP}.
 As above, we denote by $M_C$ the DG\+module over $C\ot_k C^\op$
corresponding to a CDG\+module $M$ over $B\ot_k B^\op$.

 To any left CDG\+module $G$ and right CDG\+module $F$ over $B$
one can assign the left CDG\+module $G\ot_k F$ and the right
CDG\+module $F\ot_k G$ (corresponding to each other under
the above equivalence) over the CDG\+category $B\ot_k B^\op$.
 There are also the natural \emph{diagonal} CDG\+module $B$
over $B\ot_k B^\op$ and DG\+module $C$ over $C\ot_k C^\op$
\cite[Section~2.4]{PP}; these also correspond to each other with
respect to the above equivalence of DG\+categories.

 For any DG\+module $M_C$ over $C\ot_k C^\op$, we are interested
in the comparison morphisms between the two kinds of Hochschild
cohomology $HH^{I\!I,\.*}(C,M_C)\rarrow HH^*(C,M_C)$ and homology
$HH_*(C,M_C)\rarrow HH_*^{I\!I}(C,M_C)$ \cite[formula~(23)]{PP}.

\begin{propc}
 Assume that the $\Gamma$\+graded category $B^\#\ot_k B^\#{}^\op$
is Noetherian and the diagonal $\Gamma$\+graded module $B^\#$ over
it is finitely generated.
 Let $M$ be a CDG\+module over $B\ot_k B^\op$.
 Then the natural morphism $HH^{I\!I,\.*}(C,M_C)\rarrow HH^*(C,M_C)$
is an isomorphism provided that either \par
\textup{(i)} the object $M$ belongs to the minimal triangulated
subcategory of\/ $\rD^\rco(B\ot_k B^\op\allowbreak\modlc)$ containing
the CDG\+modules $Hom_k(F\ot_k G\;k\spcheck)$ for all
$F\in H^0(\modrcfgp B)$ and $G\in H^0(B\modlc_\rfgp)$ and closed
under infinite products; or \par
\textup{(ii)} the diagonal CDG\+module $B$ over $B\ot_k B^\op$ 
belongs to the minimal thick subcategory of\/ $\rD^\rabs(B\ot_kB^\op
\modlc_\rfg)$ containing the CDG\+modules $G\ot_k F$ for all
$F\in H^0(\modrcfgp B)$ and $G\in H^0(B\modlc_\rfgp)$.
\end{propc} 

\begin{propd}
 Assume that the $\Gamma$\+graded category $B^\#\ot_k B^\#{}^\op$
is coherent and the diagonal $\Gamma$\+graded module $B^\#$ over
it is finitely presented.
 Let $M$ be a CDG\+module over $B\ot_k B^\op$.
 Then the natural morphism $HH_*(C,M_C)\rarrow HH_*^{I\!I}(C,M_C)$
is an isomorphism provided that either \par
\textup{(i)} the object $M$ belongs to the minimal triangulated
subcategory of $H^0(B\modlc_\rprj)\subset\rD^\rctr(B\ot_k B^\op\modlc)$
containing the CDG\+modules $G\ot_k F$ for all $F\in H^0(\modrcfgp B)$
and $G\in H^0(B\modlc_\rfgp)$ and closed under infinite direct sums;
or \par
\textup{(ii)} the diagonal CDG\+module $B$ over $B\ot_k B^\op$ 
belongs to the minimal thick subcategory of\/ $\rD^\rabs(B\ot_k B^\op
\modlc_\rfg)$ containing the CDG\+modules $G\ot_k F$ for all
$F\in H^0(\modrcfgp B)$ and $G\in H^0(B\modlc_\rfgp)$.
\end{propd} 

\begin{proof}[Proofs of Propositions C\+D]
 Similar to the proofs of Propositions A\+B.
\end{proof}

 In particular, assume that the $\Gamma$\+graded category $B^\#\ot_k
B^\#{}^\op$ is Noetherian and the diagonal $\Gamma$\+graded module
$B^\#$ over it is finitely generated.
 Suppose that the diagonal CDG\+module $B$ over $B\ot_k B^\op$ 
belongs to the minimal thick subcategory of\/ $\rD^\rabs(B\ot_k B^\op
\modlc_\rfg)$ containing the CDG\+modules $G\ot_k F$ for all
$F\in H^0(\modrcfgp B)$ and $G\in H^0(B\modlc_\rfgp)$.
 Then, according to~\cite[formulas~(44\+45) in Section~2.6]{PP}
and parts~(ii) of Propositions~C\+D, there are natural isomorphisms
\begin{gather}
 HH^*(C,M_C)\simeq HH^{I\!I,\.*}(C,M_C) \simeq HH^{I\!I,\.*}(B,M) \\
 HH_*(C,M_C)\simeq HH_*^{I\!I}(C,M_C) \simeq HH_*^{I\!I}(B,M)
\end{gather}
for any CDG\+module $M$ over $B\ot_k B^\op$.
 Specializing to the case of the diagonal CDG\+module $M=B$
and DG\+module $M_C=C$, we obtain
\begin{equation} \label{hochschild-category-iso}
 HH^*(C)\simeq HH^{I\!I,\.*}(C) \simeq HH^{I\!I,\.*}(B)
 \quad\text{and}\quad
 HH_*(C)\simeq HH_*^{I\!I}(C) \simeq HH_*^{I\!I}(B).
\end{equation}

\subsubsection{Locally free matrix factorizations}
\label{affine-loc-free}
 Let $k$ be a regular commutative Noetherian ring of finite Krull
dimension and $X$ be an affine scheme of finite type over $\Spec k$.
 Let $w\in\O(X)$ be a global regular function on~$X$.
 Consider the $\Z/2$\+graded CDG\+algebra $B$ over $k$ with
$B^0=\O(X)$, \ $B^1=0$, \ $d=0$, and $h=-w\in B^0$.
 We will find it convenient to denote the CDG\+algebra $B$
simply by $(X,h)=(X,-w)$ (cf.\ Section~\ref{matrix-subsect} of
the main body of this paper).

 Then $C=\modrcfgp B$ is the $\Z/2$\+graded DG\+category of
locally free matrix factorizations of finite rank of
the potential~$w$ on~$X$.
 Furthermore, one has $B\ot_k B^\op = (X\times_kX\; w_2-w_1)$,
where $w_i=p_i^*w\in\O(X\times_k X)$, \ $i=1$, $2$, and
$p_i\:X\times_k X\rarrow X$ denote the coordinate projections.
 Let $\Delta\:X\rarrow X\times_k X$ be the diagonal embedding
and $\Delta_*\O_X$ be the corresponding coherent sheaf on
$X\times_k X$.

 Consider the coherent matrix factorization of the potential
$w_2-w_1$ on $X\times X$ whose even-degree component is
the sheaf $\Delta_*\O_X$, while the odd-degree component vanishes.
 We will denote this ``diagonal'' matrix factorization simply by
$\Delta_*\O_X\in H^0((X\times_kX\;w_2-w_1)\modlc_\rfg)$.
 Applying the machinery of the previous sections leads to
the following result (cf.~\cite[Sections 4.8\+-4.10]{PP}).

\begin{cor}
 Suppose that the diagonal matrix factorization $\Delta_*\O_X$
belongs to the minimal thick subcategory of
$\rD^\rabs((X\times_kX\;w_2-w_1)\modlc_\rfg)$ containing
the external tensor products of locally free matrix factorizations
of finite rank $p_1^*G\ot_k p_2^*F$ for all $G\in
H^0((X,-w)\modlc_\rfgp)$ and $F\in H^0((X,w)\modlc_\rfgp)$.
 Then the natural isomorphisms~\textup{\eqref{hochschild-category-iso}}
hold for the CDG\+algebra $B=(X,w)$ and the DG\+category
of locally free matrix factorizations $C=\modrcfgp B$. \qed
\end{cor}

 Notice that the condition under which the conclusion of Corollary
has been proven is a rather strong one, particularly when $X$ is
not assumed to be a regular scheme.
 Then it is not even clear when or why the diagonal matrix factorization
$\Delta_*\O_X$ should belong to the thick envelope of the full
triangulated subcategory of locally free matrix factorizations
$H^0((X\times_k X\;w_2-w_1)\modlc_\rfgp)\subset
\rD^\rabs((X\times_kX\;w_2-w_1)\modlc_\rfg)$ on $X\times_kX$, let alone
to the thick subcategory generated by external tensor products of
locally free matrix factorizations from the two copies of~$X$.

\subsubsection{Smooth stratifications}  \label{smooth-strat}
 A scheme $X$ of finite type over a field~$k$ is said to admit
a \emph{smooth stratification}~\cite{Ef2} if it can be presented as
a disjoint union of its locally closed subsets
$X=\bigsqcup_\alpha S_\alpha$ so that each $S_\alpha$, when endowed with
the structure of a reduced locally closed subscheme in $X$, becomes
a smooth scheme over~$k$.
 In particular, every scheme of finite type over a perfect field~$k$
admits a smooth stratification, as any regular scheme of finite
type over a perfect field is smooth over it~\cite[Corollaires~17.15.2
and~17.15.13]{Groth2}.
 Notice that a scheme of finite type over a field admits a smooth
stratification if and only if its maximal reduced closed subscheme does.

 The definition of a \emph{regular stratification} of a Noetherian
scheme is similar, except that the strata $S_\alpha$ are
only required to be regular schemes in their reduced locally
closed subscheme structures.
 Any scheme of finite type over a field admits a regular
stratification~\cite[Scholie~7.8.3(iii\+iv) and \
Proposition~7.8.6(i)]{Groth1}.

 Let $X$ be a smooth affine scheme over a field~$k$ and $w\in\O(X)$ 
be a regular function on $X$.
 Set $X_0=\{w=0\}\subset X$ to be the zero locus of~$w$.
 The following result is a slight generalization
of~\cite[Corollary~4.8.A]{PP} based on the above definitions.

\begin{cor}
 Assume that there exists a closed subscheme $Z\subset X$ such that
$w\:X\setminus\nobreak Z\rarrow \mathbb A^1_k$ is a smooth morphism,
$w|_Z=0$, and the scheme $Z$ admits a smooth stratification over~$k$.
 Then the conditions of Corollary~\textup{\ref{affine-loc-free}} are
satisfied, so its conclusions apply.  \hbadness=1050
\end{cor}

\begin{proof}
 According to the argument in~\cite[Section~4.8]{PP}, it suffices to
show that the bounded derived category of coherent sheaves on
$Z\times Z$ is generated by external tensor products of coherent
sheaves on the two Cartesian factors.
 This is a particular case of the following lemma.
\end{proof}

\begin{lem}
 Let $Z'$ and $Z''$ be schemes of finite type over a field~$k$.
 Assume that the scheme $Z'$ admits a smooth stratification.
 Then the bounded derived category of coherent sheaves
$\sD^\b((Z'\times Z'')\coh)$ on the Cartesian product $Z'\times_k Z''$
coincides with its minimal thick subcategory containing the external
tensor products $\K'\ot_k\K''$ of coherent sheaves on $\K'$ on $Z'$
and $\K''$ on~$Z''$. 
\end{lem}

\begin{proof}
 One proceeds by induction in the total number of strata in a smooth
stratification of $Z'$ and a regular stratification of~$Z''$.
 Clearly, one can replace $Z'$ and $Z''$ with their maximal reduced
closed subschemes.
 Now if $S_{\alpha_0}$ is an open stratum in $Z'$ and $T_{\beta_0}$ is
an open stratum in $Z''$, then $S_{\alpha_0}$ is smooth as an open
subscheme in~$Z'$ and $T_{\beta_0}$ is regular as an open subscheme in
$Z''$, while the induction assumption applies to
$(Z'\setminus S_{\alpha_0})\times_k Z''$ and $Z'\times_k
(Z''\setminus T_{\beta_0})$.
 The scheme $S_{\alpha_0}\times_k T_{\beta_0}$ is regular, since it is
smooth over a regular scheme.
 The rest of the argument is based on~\cite[Proposition~2.7]{Or2}
and follows the lines of~\cite[proof of Theorem~3.7]{PL}.
\end{proof}

\addtocontents{toc}{\smallskip}
\subsection{Coherent matrix factorizations}
\label{coherent-hochschild}
 In this section we return to the notation system typical for
the main body of this paper.
 The notion of a critical value of a regular function on a singular
variety is defined in Section~\ref{noncritical}.
 In Section~\ref{external-tensor} we show that the external tensor
product of coherent matrix factorizations is a fully faithful
functor between the absolute derived categories, and provide
a sufficient condition for the pretriangulated extension of its
DG\+category version to be a quasi-equivalence.
 The Hochschild cohomology of the DG\+category corresponding to
the absolute derived category of coherent matrix factorizations of
a potential having no critical values but zero is computed in
Section~\ref{coh-hochschild-cohomology}.

 The notion of cotensor product of complexes of quasi-coherent
sheaves and quasi-coherent matrix factorizations is discussed
in Sections~\ref{cotensor-complexes}\+-\ref{cotensor-matrix},
and used in order to compute the Hochschild homology of
the (same) DG\+category of coherent matrix factorizations in
Section~\ref{coh-hochschild-homology}.
 The direct sum formula for the Hochschild (co)homology of 
the DG\+categories of coherent matrix factorizations of
a potential with several critical values is established in
Section~\ref{direct-sum-formula}.

 In some sense, the results of this section (as compared to those
of Section~\ref{loc-free-hochschild}) suggest that the DG\+category
corresponding to the absolute derived category of coherent
matrix factorizations on a singular variety may be better behaved
than the similar category of locally free matrix factorizations
of finite rank.
 Other (and in some way related) arguments in support of the same
conclusion are provided by the results of the papers~\cite{Lu,Ef2}
showing that the DG\+category corresponding to the absolute
derived category of coherent matrix factorizations is smooth
(and even homotopically finitely presented), under suitable
conditions on the field~$k$.
 (Cf.\ the counterexample in Section~\ref{nonlocalization}.)

\subsubsection{Noncritical functions} \label{noncritical}
 Let $k$ be a field and $X$ be a scheme of finite type over $\Spec k$.
 Let $f\in\O(X)$ be a global regular function on~$X$.

 Let $Y$ be a scheme of finite type over $\Spec k$ and $g\in\O(Y)$ be
a global regular function.
 Let $p_1\:X\times_k Y\rarrow X$ and $p_2\:X\times_k Y\rarrow Y$ be
the natural projections.
 Consider the regular function $f_1+g_2=p_1^*f+p_2^*g$ on $X\times_k Y$.

 Suppose that $f\:X\rarrow\mathbb A^1_k$ is a flat morphism from $X$ to
the affine line (when $k$ is algebraically closed, this means that
the function $f-c$ is a local nonzero-divisor on $X$ for every
$c\in k$).
 Then the morphism $f_1+g_2\:X\times_k Y\rarrow \mathbb A^1_k$ is also
flat as the composition of two flat morphisms $X\times_k Y\rarrow
\mathbb A^1_k\times_k Y\rarrow \mathbb A^1_k$ (the former morphism being
flat since the morphism $f\:X\rarrow\mathbb A^1_k$ is and the latter one
because the polynomial $x+g$ does not divide zero in $B[x]$ for any
commutative ring $B$ and element $g\in B$).
 In particular, it follows that the function $f_1+g_2$ is a local
nonzero-divisor on the Cartesian product $X\times_k Y$.

 A function $f\in\O(X)$ is said to be \emph{noncritical} (or to
\emph{have no critical values}) if for any regular function $g\in\O(Y)$
on a scheme $Y$ of finite type over $\Spec k$ the absolute derived
category of coherent matrix factorizations
$\sD^\abs((X\times_k Y\;\mathcal O\; f_1+g_2)\coh)$ vanishes
(i.~e., is equivalent to the zero category).
 According to Remark~\ref{second-kind} and
Theorem~\ref{cdg-supports}(b), this condition is local in both
$X$ and~$Y$.

 Therefore, given a scheme $X$ of finite type over $\Spec k$
and a regular function $f\in\O(X)$, there is a unique maximal open
subscheme $X'_f\subset X$ where the function $f$~is noncritical.
 We will see below that the open subscheme $X'_f$ is always dense in
$X$ if the morphism $f\:X\rarrow \mathbb A^1_k$ is flat and
the field~$k$ has zero characteristic.

 Similarly, there is a unique maximal open subscheme $\mathbb A^1_{k,f}
\subset \mathbb A^1_k$ such that the restriction of~$f$ to its full
preimage in $X$ is noncritical.
 The scheme $\mathbb A^1_{k,f}$ is always nonempty if the field~$k$ has
zero characteristic.
 The points in the complement $\mathbb A^1_k\setminus
\mathbb A^1_{k,f}$ are called the \emph{critical values} of~$f$.
 In particular, one says that~$f$ \emph{has no critical values but zero}
if the restriction of~$f$ to $f^{-1}(\mathbb A^1_k\setminus\{0\})
\subset X$ is noncritical.

 Notice that when the schemes $X$ and $Y$ are separated and
the morphism of schemes $f\:X\rarrow\mathbb A^1_k$ is flat,
the category $\sD^\abs((X\times_k Y\;\mathcal O\; f_1+g_2)\coh)$
is equivalent to the triangulated category
$\sD^\b_\Sing(\{f_1+g_2=0\}/X\times_k Y)$ of relative singularities
of the zero locus of the function $f_1+g_2$ on $X\times_k Y$ 
(see Theorem~\ref{main-theorem}).

\begin{rem}
 It would be interesting to have a geometric characterization of
noncriticality of functions on singular schemes.
 For example, how does our definition of noncriticality relate to
the condition that the differential of~$f$ at every closed point
$x\in X$ be a nonzero element of the Zariski cotangent space
$T^*_xX\,$?
 We do \emph{not} know this; cf.\ the smooth stratification
approach below.
\end{rem}

\begin{lem}
 Let $X=\bigsqcup_\alpha S_\alpha$ be a scheme of finite type over\/
$\Spec k$ presented as a disjoint union of its locally closed
subsets, endowed with their reduced locally closed subscheme
structures.
 Let $\L$ be a line bundle on $X$ and $w\in\L(X)$ be its global
section.
 In this setting, if the absolute derived categories\/
$\sD^\abs((S_\alpha,\L|_{S_\alpha},w|_{S_\alpha})\coh)$ vanish for
all\/~$\alpha$, then so does the absolute derived category\/
$\sD^\abs((X,\L,w)\coh)$.
\end{lem}

\begin{proof}
 Proceeding by induction in the number of strata in
the stratification~$S_\alpha$, it suffices to consider the case
when there are only two of them, namely, a closed subset $S\subset X$
and its open complement $X\setminus S$.
 One can also replace $X$ with its maximal reduced closed subscheme.
 Then the desired assertion follows from Theorem~\ref{cdg-supports}(b),
since the triangulated category $\sD^\abs((X,\L,w)\coh_S)$ is generated
by the image of the natural functor $\sD^\abs((S,\L|_S,w|_S)\coh)
\rarrow\sD^\abs((X,\L,w)\coh_S)$.
\end{proof}

\begin{prop}
 Let $X$ be a scheme of finite type over\/ $\Spec k$ and
$f\in\O(X)$ be a regular function on~$X$.
 Let $X=\bigsqcup_\alpha S_\alpha$ be a smooth stratification of
the scheme $X$ over~$k$ (see Section~\textup{\ref{smooth-strat}})
such that the morphisms of schemes $f|_{S_\alpha}\:S_\alpha\rarrow
\mathbb A^1_k$ are smooth for all\/~$\alpha$.
 Then the function~$f$ is noncritical on~$X$.
\end{prop}

\begin{proof}
 Let $Y$ be a scheme of finite type over $\Spec k$ and $g\in\O(Y)$ be
a regular function.
 We have to show that the triangulated category
$\sD^\abs((X\times_k Y\;\O\;f_1+g_2)\allowbreak\coh)$ vanishes.
 Choosing a stratification of $Y$ by regular locally closed subschemes
and applying Lemma, one can assume that $X$ is smooth over~$k$ and
$Y$ is regular.

 Then the scheme $X\times_k Y$ is also regular, the derivative of
the function $f_1+g_2\in\O(X\times_k Y)$, viewed as an element of
the Zariski cotangent space, does not vanish at any points where
the function itself does (and, in a sense, at any other closed points,
too), and it follows that the zero locus of $f_1+g_2$ in $X\times_k Y$
is also a regular scheme.
 It remains to use Theorem~\ref{main-theorem}
(or~\cite[Theorem~3.5]{Or3} and Corollary~\ref{regular-cor}(c)).
\end{proof}

 It follows from Proposition that, for any scheme of finite type $X$
with a smooth stratification $X=\bigsqcup_\alpha S_\alpha$ over
$\Spec k$ and any regular function $f\in\O(X)$, the set of critical
values of the function~$f$ on $X$ is contained in the union of the sets
of critical values of the functions $f|_{S_\alpha}$.
 In particular, if the characteristic of~$k$ is zero, then all
these sets are finite.

\subsubsection{External tensor products} \label{external-tensor}
 Let $X'$ and $X''$ be separated schemes of finite type over
a field~$k$, and let $w'\in\O(X')$ and $w''\in\O(X'')$ be regular
functions.
 Let $X'\times_k X''$ be the Cartesian product, $p_1$ and $p_2$
be its natural projections onto the factors $X'$ and $X''$, and
$w'_1+w''_2=p_1^*w'+p_2^*w''$ be the related regular function
on $X'\times_k X''$.
 Then there is the external tensor product functor
\begin{multline} \label{qcoh-external-tensor}
 \otimes_k\:\sD^\co((X',\O,w')\qcoh)\times\sD^\co((X'',\O,w'')\qcoh)
 \\ \lrarrow\sD^\co((X'\times_k X''\;\O\;w'_1+w''_2)\qcoh),
\end{multline}
which restricts to the similar functor
\begin{multline}  \label{coherent-external-tensor}
 \otimes_k\:\sD^\abs((X',\O,w')\coh)\times\sD^\abs((X'',\O,w'')\coh)
 \\ \lrarrow\sD^\abs((X'\times_k X''\;\O\;w'_1+w''_2)\coh)
\end{multline}
on coherent matrix factorizations.

\begin{prop}
 Let $\K'$ and $\M'$ be coherent matrix factorizations of
the potential~$w'$ on the scheme $X'$, and let $\K''$ and $\M''$
be coherent matrix factorizations of the potential~$w''$ on
the scheme~$X''$.
 Then the natural map of\/ $\Z/2$\+graded $k$\+vector spaces
of morphisms
\begin{multline}
 \Hom_{\sD^\abs((X',\O,w')\coh)}(\K',\M'[*])\ot_k
 \Hom_{\sD^\abs((X'',\O,w'')\coh)}(\K'',\M''[*])\\ \lrarrow
 \Hom_{\sD^\abs((X'\times_k X''\;\O\;w'_1+w''_2)\coh)}
 (\K'\ot_k\K''\;\M'\ot_k\M''[*])
\end{multline}
induced by the additive functor of two
arguments~\textup{\eqref{coherent-external-tensor}}
is an isomorphism.
\end{prop}

\begin{proof}
 By Proposition~\ref{embedding-prop}(d), it suffices to show that
the natural map
\begin{multline} \label{qcoh-external-tensor-morphisms}
 \Hom_{\sD^\co((X',\O,w')\qcoh)}(\K',\M'[*])\ot_k
 \Hom_{\sD^\co((X'',\O,w'')\qcoh)}(\K'',\M''[*])\\ \lrarrow
 \Hom_{\sD^\co((X'\times_k X''\;\O\;w'_1+w''_2)\qcoh)}
 (\K'\ot_k\K''\;\M'\ot_k\M''[*])
\end{multline}
induced by the functor~\eqref{qcoh-external-tensor} is an isomorphism
for any coherent matrix factorizations $\K'$, $\K''$ and quasi-coherent
matrix factorizations $\M'$, $\M''$ of the potentials $w'$ and~$w''$.
 One easily checks that the desired assertion holds for the Hom spaces
in the homotopy categories of matrix factorizations (since it holds
for morphisms between the external tensor products of coherent and
quasi-coherent sheaves).

 Furthermore, one can assume the quasi-coherent matrix factorizations
$\M'$ and $\M''$ to be injective.
 Then the Hom spaces in the left-hand side of
the map~\eqref{qcoh-external-tensor-morphisms} coincide with
the similar Hom spaces computed in the homotopy categories of matrix
factorizations.
 Let $\I^\bu$ be a right resolution of $\M'\ot_k\M''$ in the abelian
category of quasi-coherent matrix factorizations (and closed morphisms
between them) consisting of injective matrix factorizations, and let
$\J$ be the total matrix factorization of the complex $\I^\bu$
constructed by taking infinite direct sums along the diagonals.
 Then the $k$\+vector spaces of morphisms from $\K'\ot_k\K''$ into $\J$
in the homotopy category of matrix factorizations are isomorphic to
the right-hand side of~\eqref{qcoh-external-tensor-morphisms}
\cite[Theorem~3.7]{Pkoszul}.

 It remains to show that the spaces of morphisms from $\K'\ot_k\K''$
to $\M'\ot_k\M''$ in the homotopy category of matrix factorizations
are isomorphic to the similar spaces of morphisms from $\K'\ot_k\K''$
to~$\J$.
 Indeed, taking the termwise Hom from $\K'\ot_k\K''$ preserves
exactness of the sequence $0\rarrow\M'\ot_k\M''\rarrow\I^\bu$, since
the higher Ext spaces from the components of $\K'\ot_k\K''$ into
those of $\M'\ot_k\M''$ in the abelian category of quasi-coherent
sheaves on $X'\ot_k X''$ vanish.
 The latter assertion can be checked for affine schemes $X'$, $X''$
using projective resolutions and then globally for the cohomology of
quasi-coherent sheaves using, e.~g., the \v Cech approach.
\end{proof}

\begin{thm}
 Assume that the morphisms of schemes $w'\:X'\rarrow\mathbb A^1_k$ and
$w''\:X''\rarrow \mathbb A^1_k$ are flat.
 Suppose that there exist closed subschemes $Z'\subset X'$ and
$Z''\subset X''$ such that $w'|_{Z'}=0=w''_{Z''}$, the functions $w'$
and~$w''$ are noncritical on $X'\setminus Z'$ and $X''\setminus Z''$,
and the scheme $Z'$ admits a smooth stratification over~$k$.
 Then the absolute derived category\/
$\sD^\abs((X'\times_k X''\;\O\;w'_1+w''_2)\coh)$ coincides with its
minimal thick subcategory containing the image of
the functor~\textup{\eqref{coherent-external-tensor}}.
\end{thm}

\begin{proof}
 By the definition of noncriticality, one has
$\sD^\abs((((X'\setminus Z')\times_k X'')\;\O\;w'_1+w''_2)\allowbreak
\coh)=0=\sD^\abs(((X'\times_k(X''\setminus Z''))\;\O\;w'_1+w''_2)\coh)$.
 Therefore, any coherent matrix factorization of the potential
$w'_1+w''_2$ on $X'\times_k X''$ has its category-theoretic support
inside $Z'\times_k Z''$, and is consequently isomorphic in
$\sD^\abs((X'\times_k X''\;\O\;w'_1+w''_2)\coh)$ to
a direct summand of an object represented by a coherent matrix
factorization supported set-theoretically inside $Z'\times_k Z''$
(see Corollary~\ref{cdg-supports}(b)).
 It follows that the triangulated category
$\sD^\abs((X'\times_k X''\;\O\;w'_1+w''_2)\coh)$ is generated by
the direct images of coherent matrix factorizations of the zero
potential from the closed embedding $Z'\times_k Z''\rarrow
X'\times_k X''$.

 Furthermore, let $X'_0$, $X''_0$, and $Y_0$ denote the zero loci of
the functions $w'$, $w''$, and $w'_1+w''_2$ on $X'$, $X''$, and
$X'\times_k X''$, respectively.
 Denote the natural closed embeddings by $i'\:X'_0\rarrow X'$, \
$i''\:X''_0\rarrow X''$, \ $\iota\:X'_0\times_k X''_0\rarrow Y_0$,
and $h\:Y_0\rarrow X'\times_k X''$.
 The external tensor product functor (cf.~\cite[Lemma~4.8.B]{PP})
\begin{equation}  \label{sing-external-tensor}
 \ot_k\:\sD^\b_\Sing(X'_0/X')\times\sD^\b_\Sing(X''_0/X'')\lrarrow
 \sD^\b_\Sing(Y_0/(X'\times_kX''))
\end{equation}
is well-defined, since for any bounded complexes of coherent sheaves
$\F^\bu$ on $X'$ and $\K^\bu$ on $X''_0$ one has
$\iota_*(\boL i'{}^*\F^\bu\ot_k\K^\bu)\simeq\boL h^*((\id_{X'}\times
i'')_*(\F^\bu\ot_k\K^\bu))$.
 Indeed, the square diagram of closed embeddings
$$
\begin{diagram}
\node{X_0'\times_k X_0''}\arrow{e}\arrow{s}
\node{X'\times_k X_0''}\arrow{s} \\ \node{Y_0}\arrow{e}
\node{X'\times_k X''}
\end{diagram}
$$
is Cartesian and the higher derived tensor products related to
the construction of this relative Cartesian product of schemes
all vanish.

 The functor $\Ups\:\sD^\b_\Sing(Y_0/(X'\times_kX''))\rarrow
\sD^\abs((X'\times_kX''\;\O\;w'_1+w''_2)\coh)$ (see
Section~\ref{main-theorem}) and the similar functors for
the potentials $w'$ and $w''$ on $X'$ and $X''$ transform
the external product functor~\eqref{coherent-external-tensor} into
the external tensor product functor~\eqref{sing-external-tensor}.
 By the assumption, one has $Z'\subset X'_0$ and $Z''\subset X''_0$.
 It remains to apply Lemma~\ref{smooth-strat} in order to finish
the proof of the theorem.
\end{proof}

\subsubsection{Internal Hom of matrix factorizations}
\label{internal-hom-matrix}
 Let $X$ be a separated Noetherian scheme. 
 Let $\L$ be a line bundle on $X$ and $w'$, $w''\in\L(X)$
be its global sections.
 Then given a matrix factorization $\U^0\rarrow\U^1\ot\L^{\ot1/2}\rarrow
\U^0\ot_{\O_X}\L$ of the potential~$w'$ and a matrix factorization
$\V^0\rarrow\V^1\ot\L^{\ot1/2}\rarrow\V^0\ot_{\O_X}\L$ of
the potential~$w''$ on the scheme $X$ (in the symbolic notation of 
Section~\ref{matrix-subsect}), one can construct the matrix
factorization $\U^0\ot_{\O_X}\V^0\oplus\U^1\ot_{\O_X}\V^1 \rarrow
\U^1\ot\L^{\ot1/2}\ot_{\O_X}\V^0\oplus \U^0\ot_{\O_X}\V^1\ot\L^{\ot1/2}
\rarrow\U^0\ot_{\O_X}\V^0\ot_{\O_X}\L\oplus
\U^1\ot_{\O_X}\V^1\ot_{\O_X}\L$ of the potential $w'+w''$ on~$X$.
 Here the tensor product $\U^1\ot_{\O_X}\V^1$ is defined as
the sheaf $(\U^1\ot\L^{\ot1/2})\ot_{\O_X}(\V^1\ot\L^{\ot1/2})\ot_{\O_X}
\L^{\ot-1}$ on $X$, while the differential on the tensor product of
matrix factorizations is given by the conventional rule
$d(u\ot v) = d(u)\ot v + (-1)^{|u|}u\ot d(v)$. 

 We denote the matrix factorization so obtained by $\U\ot_{\O_X}\V$ and
call it the \emph{tensor product} of two matrix factorizations $\U$ and
$\V$ of two sections $w'$ and~$w''$ of the same line bundle $\L$ on
a scheme~$X$.
 Restricting to the cases when one or both matrix factorizations are
flat, and passing to the coderived categories, one obtains the induced
tensor product functors
\begin{multline}  \label{matrix-tensor-fl-fl}
 \ot_{\O_X}\:\sD^\co((X,\L,w')\qcoh_\fl)\times\sD^\co((X,\L,w'')
 \qcoh_\fl) \\ \lrarrow\sD^\co((X\;\L\;w'+w'')\qcoh_\fl)
\end{multline}
and
\begin{multline}   \label{matrix-tensor-fl-qcoh}
 \ot_{\O_X}\:\sD^\co((X,\L,w')\qcoh_\fl)\times\sD^\co((X,\L,w'')\qcoh)
 \\ \lrarrow\sD^\co((X\;\L\;w'+w'')\qcoh).
\end{multline}
 The functors~(\ref{matrix-tensor-fl-fl}\+-\ref{matrix-tensor-fl-qcoh})
are well-defined, since the tensor product with a flat (quasi-coherent)
matrix factorization takes a short exact sequence of flat matrix
factorizations to a short exact sequence of flat matrix factorizations,
the tensor product with a flat matrix factorization takes a short exact
sequence of quasi-coherent matrix factorizations to a short exact
sequence of quasi-coherent matrix factorizations, and the tensor
product with a quasi-coherent matrix factorization takes a short
exact sequence of flat matrix factorizations to a short exact sequence
of quasi-coherent matrix factorizations.
 Also, the tensor product functor preserves infinite direct sums.

 Given a quasi-coherent matrix factorization $\U^0\rarrow
\U^1\ot\L^{\ot1/2}\rarrow\U^0\ot_{\O_X}\L$ of a potential $w'\in\L(X)$
and a quasi-coherent matrix factorization $\V^0\rarrow\V^1\ot\L^{\ot1/2}
\allowbreak\rarrow\V^1\ot_{\O_X}\L$ of a potential $w''\in\L(X)$ on
the scheme~$X$, one can construct the quasi-coherent matrix
factorization $\cHom_{X\qc}(\U^0,\V^0)\oplus\cHom_{X\qc}(\U^1,\V^1)
\rarrow\cHom_{X\qc}(\U^0\;\V^1\ot\L^{\ot1/2})\oplus
\cHom_{X\qc}(\U^1,\V^0)\ot\L^{\ot1/2}\rarrow\cHom_{X\qc}(\U^0,\V^0)
\ot_{\O_X}\nobreak\L\oplus\cHom_{X\qc}(\U^1,\V^1)\ot_{\O_X}\L$ of
the potential $w''-w'$ on~$X$.
 Here the sheaf $\cHom_{X\qc}(\U^1,\V^0)\ot\L^{\ot1/2}$ is defined
as the tensor product $\cHom_{X\qc}(\U^1\ot\L^{\ot1/2}\;\V^0)\allowbreak
\ot_{\O_X}\L$, while the differential on the internal Hom is given by
the conventional rule $d(g)(u)=d(g(u))-(-1)^{|g|}g(d(u))$.

 We denote the matrix factorization so obtained by
$\cHom_{X\qc}(\U,\V)$ and call it the matrix factorization of
\emph{quasi-coherent internal Hom} between the quasi-coherent
matrix factorizations $\U$ and $\V$ of two sections $w'$ and~$w''$
of the same line bundle $\L$ on a scheme~$X$.
 Restricting to the case when the matrix factorization in
the second argument is injective, one obtains the induced
internal Hom functor
\begin{multline}  \label{matrix-chom-induced}
 \cHom_{X\qc}\:\sD^\abs((X,\L,w')\qcoh)^\sop\times
 H^0((X,\L,w'')\qcoh_\inj) \\
 \lrarrow\sD^\abs((X\;\L\;w''-w')\qcoh),
\end{multline}
which can be also viewed as the right derived internal Hom functor
\begin{multline}  \label{matrix-chom-derived}
 \boR\cHom_{X\qc}\:\sD^\abs((X,\L,w')\qcoh)^\sop\times
 \sD^\co((X,\L,w'')\qcoh) \\
 \lrarrow\sD^\abs((X\;\L\;w''-w')\qcoh).
\end{multline}

\begin{rem}
 Alternatively, one could restrict the quasi-coherent internal Hom
functor to pairs of quasi-coherent matrix factorizations which
are \emph{both} injective, obtaining the triangulated functor
\begin{multline}  \label{matrix-flatvalued-chom-induced}
 \cHom_{X\qc}\:H^0((X,\L,w')\qcoh_\inj)^\sop\times
 H^0((X,\L,w'')\qcoh_\inj) \\
 \lrarrow H^0((X\;\L\;w''-w')\qcoh_\fl),
\end{multline}
which can be also viewed as a derived internal Hom functor
\begin{multline}  \label{matrix-flatvalued-chom-derived}
 \boL\boR\cHom_{X\qc}\:\sD^\co((X,\L,w')\qcoh)^\sop\times
 \sD^\co((X,\L,w'')\qcoh) \\
 \lrarrow\sD^\abs((X\;\L\;w''-w')\qcoh_\fl)
\end{multline}
that is a left derived functor in its first argument and
a right derived functor in the second one.
 Notice that the derived functor so obtained does \emph{not}
agree with the right derived functor defined above, i.~e.,
the composition of the functor~\eqref{matrix-flatvalued-chom-derived}
with the natural fully faithful functor $\sD^\abs((X\;\L\;w''-w')
\qcoh_\fl)\rarrow\sD^\abs((X\;\L\;w''-w')\qcoh)$ and the Verdier
localization functor $\sD^\abs((X,\L,w')\qcoh)\rarrow
\sD^\co((X,\L,w')\qcoh)$ is \emph{not} isomorphic to
the functor~\eqref{matrix-chom-derived}.
\end{rem}

 In particular, when $w'=w''$,
the functors~(\ref{matrix-chom-induced}\+-\ref{matrix-chom-derived})
take values in the absolute derived category of quasi-coherent
matrix factorizations of the zero potential $0\in\L(X)$.
 The objects of this category are simply complexes of quasi-coherent
sheaves $\M^\bu$ on $X$ endowed with a $2$\+periodicity isomorphism
$\M^\bu[2]\simeq\M^\bu\ot_{\O_X}\L$.
 So there is a natural forgetful functor
\begin{equation}  \label{forget-periodicity}
 \sD^\co((X,\L,0)\qcoh)\lrarrow\sD^\co(X\qcoh)
\end{equation}
and the similar functors acting on the homotopy, absolute derived,
etc.\ categories of flat, coherent, locally free, etc.\ matrix
factorizations.

 Furthermore, there is the derived global sections functor
\begin{equation}  \label{derived-global}
 \boR\Gamma(X,{-})\:\sD^\co(X\qcoh)\lrarrow\sD(\Z\modl)
\end{equation}
taking values in the derived category of abelian groups and defined
using either the injective resolutions or the \v Cech construction
(see Sections~\ref{pull-push}\+-\ref{finite-dim-morphisms}).
 In fact, the functor~\eqref{derived-global} factorizes through
the conventional derived category $\sD(X\qcoh)$.

 Composing the forgetful functor with the functor of underived
global sections of complexes of quasi-coherent sheaves, one obtains
a triangulated functor
\begin{equation}  \label{underived-global-matrix}
 \Gamma(X,{-})\:H^0((X,\L,0)\qcoh)\lrarrow\sD(\Z\modl). 
\end{equation}
Alternatively, the functor~\eqref{underived-global-matrix} can be
defined as the functor $\Hom_{(X,\L,0)\qcoh}(\O_X,{-})$, where
the structure sheaf $\O_X$ is viewed as a matrix factorization
$(\U^0,\U^1)$ of the potential $0\in\L(X)$ with the components
$\U^0=\O_X$ and $\U^1\ot\L^{\ot1/2}=0$.

 Similarly, composing the functors~\eqref{forget-periodicity}
and~\eqref{derived-global}, one obtains a triangulated functor
\begin{equation}  \label{derived-global-matrix}
 \boR\Gamma(X,{-})\:\sD^\co((X,\L,0)\qcoh)\lrarrow\sD(\Z\modl),
\end{equation}
which can be alternatively described as the functor
$\Hom_{\sD^\co((X,\L,0)\qcoh)}(\O_X,{-}[*])$.
 In the case when $\L=\O_X$,
the functors~(\ref{underived-global-matrix}\+-%
\ref{derived-global-matrix})
can be viewed as taking values in the derived category of
$\Z/2$\+graded ($2$\+periodic) complexes of abelian groups.

 For any quasi-coherent matrix factorizations $\K$ and $\M$ of
a potential $w\in\L(X)$ on the scheme $X$ there is a natural
isomorphism of complexes of abelian groups
\begin{equation}  \label{matrix-rhom-gamma-chom}
 \Hom_{(X,\L,w)\qcoh}(\K,\M)\simeq\Gamma(X,\cHom_{X\qc}(\K,\M)),
\end{equation}
and more generally, for any quasi-coherent matrix factorizations
$\K$ and $\E$ of potentials $w'$ and $w''\in\L(X)$ and
a quasi-coherent matrix factorization $\M$ of the potential
$w'+w''$ on the scheme $X$ there is a natural isomorphism
of complexes
\begin{equation}
 \Hom_{(X\;\L\;w'+w'')\qcoh}(\K\ot_{\O_X}\E\;\M)\simeq
 \Hom_{(X,\L,w'')\qcoh}(\E,\cHom_{X\qc}(\K,\M)).
\end{equation}

\begin{lem}
 Let $\K$ be a quasi-coherent matrix factorization and
$\M$ be an injective quasi-coherent matrix factorization of
a potential $w\in\L(X)$.
 Let $\cHom_{X\qc}(\K,\M)\rarrow\J$ be a closed morphism
with a coacyclic cone between quasi-coherent matrix factorizations
of the potential\/ $0\in\L(X)$ from the matrix factorization of
quasi-coherent internal Hom into an injective matrix
factorization~$\J$.
 Then the induced morphism
$$
 \Gamma(X,\cHom_{X\qc}(\K,\M))\lrarrow\Gamma(X,\J)
$$
is a quasi-isomorphism of complexes of abelian groups.
\end{lem}

\begin{proof}
 Let $0\rarrow\cHom_{X\qc}(\K,\M)\rarrow\I^\bu$ be a right resolution
of the matrix factorization $\cHom_{X\qc}(\K,\M)$ by injective 
matrix factorization~$\I^i$.
 Then one can take $\J$ to be the total matrix factorization of
the complex $\I^\bu$ constructed by passing to the infinite direct sums
along the diagonals.
 Notice that the functor of global sections of quasi-coherent sheaves
on $X$ commutes with the infinite direct sums.
 It remains to show that the functor $\Gamma(X,{-}) =
\Hom_{(X,\L,w'')\qcoh}(\O_X,{-})$ preserves exactness of the sequence
$0\rarrow\cHom_{X\qc}(\K,\M)\rarrow\I^\bu$
(cf.\ the proof of Proposition~\ref{external-tensor}).
 
 In fact, we claim that the Ext groups from flat quasi-coherent sheaves
to the components of $\cHom_{X\qc}(\K,\M)$ vanish in the abelian
category $X\qcoh$.
 This assertion follows from the results
of~\cite[Lemma~2.5.3(c) and Corollary~4.1.9(b)]{Pcosh}
(the argument is based essentially on the above
Lemma~\ref{flat-sheaves}).
\end{proof}

 We recall the constructions of the (underived and derived) direct
and inverse image functors for matrix factorizations from
Sections~\ref{pull-push}\+-\ref{finite-dim-morphisms}
and~\ref{matrix-push}\+-\ref{finite-dim-matrix-push}.
 In addition to the conventional adjunction of the (underived)
direct and inverse image functors $f_*$ and~$f^*$ (as mentioned
in Section~\ref{pull-push}), there is also the ``internal Hom
adjunction'', formulated as follows.

 Let $f\:Z\rarrow Y$ be a morphism of separated Noetherian
schemes, $\L$ be a line bundle on $Y$, and $w'$, $w''\in\L(Y)$ be
two global sections.
 Let $\K\in H^0((Y,\L,w')\qcoh)$ and $\M\in H^0((Z\;f^*\L\;f^*w'')
\qcoh)$ be quasi-coherent matrix factorizations on $Y$ and~$Z$.
 Then there is a natural isomorphism
\begin{equation}
 f_*\cHom_{Z\qc}(f^*\K,\M)\simeq\cHom_{Y\qc}(\K,f_*\M)
\end{equation}
of quasi-coherent matrix factorizations of the potential
$w''-w'$ on~$Y$.

 Now let $X$ be a separated scheme of finite type over a field~$k$,
and let $w'$, $w''\in\O(X)$ be two global regular functions on~$X$.
 Denote by $p_1$ and~$p_2$ the natural projections $X\times_k X
\rightrightarrows X$, and consider the regular function $w'_1+w''_2=
p_1^*w'+p_2^*w''$ on $X\times_k X$.
% Given a regular function $w\in\O(X)$, we denote by $w_i$, \
%$i=1$,~$2$, the regular functions $p_i^*w$ on $X\times_k X$.
 Let $\Delta\:X\rarrow X\times_k X$ denote the diagonal map.

 Let $\N$ and $\K$ be quasi-coherent matrix factorizations of
the potentials $w'$ and~$w''$ on~$X$.
 Then there is a natural isomorphism $\N\ot_{\O_X}\K\simeq
\Delta^*(\N\ot_k\K)$ of matrix factorizations of the potential
$w'+w''$ on~$X$.
 Therefore, given a quasi-coherent matrix factorization $\M$ of
the potential $w'+w''\in\O(X)$, one has a natural isomorphism of
$\Z/2$\+graded complexes of abelian groups
\begin{multline}  \label{rhom-chom-diagonal}
 \Hom_{(X,\O,w'')\qcoh}(\K,\cHom_{X\qc}(\N,\M)) \\ \.\simeq\.
 \Hom_{(X\times_k X\;\O\;w'_1+w''_2)\qcoh}(\N\ot_k\K\;\Delta_*\M).
\end{multline}

\begin{prop}
\textup{(a)} Assume that the matrix factorization $\N$ is coherent and
the matrix factorization $\M$ is injective.
 Let $\Delta_*\M\rarrow\J$ be a closed morphism with a coacyclic cone
between quasi-coherent matrix factorizations of the potential
$w'_1+w''_2$ on $X\times_k X$ from the direct image $\Delta_*\M$
into an injective matrix factorization~$\J$.
 Then there is a natural closed morphism with a coacyclic cone
$\cHom_{X\qc}(\N,\M)\rarrow p_2{}_*\cHom_{X\times_k X\qc}
(p_1^*\N,\J)$ of quasi-coherent matrix factorizations of
the potential~$w''$ on $X$, and the matrix factorization
$p_2{}_*\cHom_{X\times_k X\qc}(p_1^*\N,\J)$ is injective. \par
\textup{(b)} There is a natural isomorphism of\/ $\Z/2$\+graded
complexes of abelian groups
\begin{multline*}
 \Hom_{(X,\O,w'')\qcoh}(\K\;p_2{}_*\cHom_{X\times_k X\qc}(p_1^*\N,\J))
 \\ \.\simeq\.\Hom_{(X\times_kX\;\O\;w'_1+w''_2)\qcoh}(\N\ot_k\K\;\J).
\end{multline*}
\end{prop}

\begin{proof}
 Part~(a): the desired closed morphism is provided by the composition
\begin{multline*}
\cHom_{X\qc}(\N,\M)\.\simeq\.
p_2{}_*\Delta_*\cHom_{X\qc}(\Delta^*p_1^*\N,\M) \\ \.\simeq\.
p_2{}_*\cHom_{X\times_k X\qc}(p_1^*\N,\Delta_*\M)\lrarrow
p_2{}_*\cHom_{X\times_k X\qc}(p_1^*\N,\J).
\end{multline*}
 To prove that this morphism has a coacyclic cone, pick a right
resolution $\I^\bu$ of the matrix factorization $\Delta_*\M$ on
$X\times_k X$ by injective matrix factorizations, and take $\J$ to be
the totalization of the complex of matrix factorizations $\I^\bu$
constructed by passing to the infinite direct sums along
the diagonals.

 Then the complex of matrix factorizations $0\rarrow
\cHom_{X\times_k X\qc}(p_1^*\N,\Delta_*\M)\rarrow
\cHom_{X\times_k X\qc}(p_1^*\N,\I^\bu)$ is acyclic, since
for any affine open subscheme $U\subset X$ the higher Ext spaces
between the components of the restrictions of $p_1^*\N$ and
$\Delta_*\M$ to $U\times_k U$ vanish.
 The latter assertion follows from the adjunction of derived
functors $\boL\Delta^*$ and $\Delta_*=\boR\Delta_*$ or
$p_1^*=\boL p_1^*$ and $\boR p_1{}_*$ together with the agreement of
the derived direct/inverse images of (complexes of) quasi-coherent
sheaves with the compositions of morphisms of separated
Noetherian schemes.

 It remains to show that our complex will stay acyclic after
applying the direct image functor~$p_2{}_*$.
 According to the argument in the proof of Lemma, the components of
the matrix factorizations $\cHom_{X\times_k X\qc}(p_1^*\N,\I^i)$
are acyclic for the direct image.
 So are the components of the matrix factorization
$\cHom_{X\times_k X\qc}(p_1^*\N,\Delta_*\M)$, in view of the above
local argument and since for any affine open subscheme $U\subset X$
the higher Ext spaces between the components of the restrictions
of $p_1^*\N$ and $\Delta_*\M$ to $X\times_k U$ vanish.
 The latter assertion is checked in the same way as above.

 Finally, the claim that the matrix factorization in question is
injective follows from the computation in part~(b), which shows
that the left-hand side is an exact functor of the argument $\K$,
because the right-hand side is.
 Part~(b) is straightforward:
\begin{multline*}
 \Hom_{(X,\O,w'')\qcoh}(\K\;p_2{}_*\cHom_{X\times_k X\qc}(p_1^*\N,\J))
 \\ \simeq\. \Hom_{(X\times_k X\;\O\;w''_2)\qcoh}(p_2^*\K\;
 \cHom_{X\times_k X\qc}(p_1^*\N,\J)) \\ \simeq\.
 \Hom_{(X\times_k X\;\O\;w'_1+w''_2)\qcoh}
 (p_1^*\N\ot_{\O_{X\times_k X}}p_2^*\K\;\J) \\ \simeq\.
 \Hom_{(X\times_k X\;\O\;w'_1+w''_2)\qcoh}(\N\ot_k\K\;\J). \qed
\end{multline*}  \renewcommand{\qed}{}
\end{proof}

\subsubsection{Hochschild cohomology}  \label{coh-hochschild-cohomology}
 Our goal in the rest of this appendix is to compute the Hochschild
(co)homology of the DG\+category $\sDG^\abs((X,\O,w)\coh)$
corresponding to the triangulated category $\sD^\abs((X,\O,w)\coh)$.
 The word ``corresponding'' here means, first of all, that there is
a natural equivalence of (triangulated) categories
$H^0\sDG^\abs((X,\O,w)\coh)\simeq\sD^\abs((X,\O,w)\coh)$
(see~\cite[Section~1.2]{Pkoszul}).

 As the absolute derived category is constructed from the homotopy
category of matrix factorizations using the Verdier localization
procedure, so the DG\+category $\sDG^\abs((X,\O,w)\coh)$ is obtained
by applying a DG\+version of localization to the DG\+category of
coherent matrix factorizations $(X,\O,w)\coh$ of the potential~$w$
on the scheme~$X$ (see Section~\ref{qcoh-cdg}).
 Several such localization procedures are known, leading to naturally
quasi-equivalent DG\+categories.
 As the Hochschild (co)homology of DG\+categories are preserved by
quasi-equivalences~\cite[Sections~2.1 and~2.4]{PP}, it is not very
important which localization procedure to choose.
 To be specific, let us say that we prefer Drinfeld's
localization~\cite{Dr0}.
 Similarly one localizes the DG\+category $(X,\O,w)\qcoh$ and
obtains a DG\+category $\sDG^\co((X,\O,w)\qcoh)$ ``corresponding''
to the coderived category $\sD^\co((X,\O,w)\qcoh)$.

 Our method will naturally allow to compute the Hochschild cohomology
$HH^*(\sDG^\abs((X,\O,w)\coh))$ together with its structure of
an associative (in fact, supercommutative, but we will neither
prove nor use this fact) $\Z/2$\+graded algebra over~$k$.
 Similarly, the Hochschild homology $HH_*(\sDG^\abs((X,\O,w)\coh))$
will be computed together with its structure of a $\Z/2$\+graded
module over the $\Z/2$\+graded associative algebra
$HH^*(\sDG^\abs((X,\O,w)\coh))$. {\hbadness=2050\par}

 Let $X$ be a separated scheme of finite type over a field~$k$ and
$w\in\O(X)$ be a global regular function.
 Assume that the morphism of schemes $w\:X\rarrow\mathbb A^1_k$ is flat.
 Consider the Cartesian square $X\times_k X$ and endow it with
the potential (global function) $w_2-w_1=p_2^*(w)-p_1^*(w)$.

 Any $\Z/2$\+graded complex of quasi-coherent sheaves $\K^\bu$ on
$X$ can be viewed as a matrix factorization of the potential
$0\in\O(X)$.
 Furthermore, one can take its direct image $\Delta_*\K^\bu$ with
respect to the diagonal embedding $\Delta\:X\rarrow X\times_k X$
and consider it as a quasi-coherent matrix factorization of
the potential $w_2-w_1$ on $X\times_k X$.
 Given a bounded $\Z$\+graded complex of quasi-coherent sheaves
$\K^\bu$, one can associate a $\Z/2$\+complex with it (by taking
direct sums of all terms with the same parity) and then apply
the above constructions.

\begin{thm}
 Assume that there exists a closed subscheme $Z\subset X$ such that
$w|_Z=0$, the function~$w$ is noncritical on $X\setminus Z$, and
the scheme $Z$ admits a smooth stratification over~$k$.
 In particular, if the field~$k$ is perfect, it suffices to require
that the function~$w$ on $X$ have no critical values but zero
(and take $Z=\{w=0\}$).
 Then there is a natural isomorphism between the Hochschild cohomology
algebra $HH^*(\sDG^\abs((X,\O,w)\coh))$ and the Ext algebra\/
$\Hom_{\sD^\co((X\times_kX\;\O\;w_2-w_1)\qcoh)}(\Delta_*\D_X^\bu\;
\Delta_*\D_X^\bu[*])$, where $\D_X^\bu$ denotes a dualizing complex
on~$X$.
\end{thm}

\begin{proof}
 By the definition, the Hochschild cohomology algebra of
a $\Z/2$\+graded DG\+category $\sDG$ is the $\Z/2$\+graded algebra
$\Hom_{\sD(\sDG\ot_k\sDG^\sop\modl)}(\sDG,\sDG[*])$, where
the $\Hom$ is taken in the conventional derived category
$\sD(\sDG\ot_k\sDG^\sop\modl)$ of DG\+bimodules over $\sDG$
(or DG\+modules over $\sDG\ot_k\sDG^\sop$) between two copies of
the diagonal DG\+bimodule $\sDG$ over~$\sDG$
\cite[Sections~2.4 and~3.1]{PP}.

 Specializing to the case of the DG\+category $\sDG_w=
\sDG^\abs((X,\O,w)\coh)$, we notice, first of all,
that the contravariant Serre--Grothendieck duality functor
$\M\longmapsto\cHom_{\O_X}(\M,\D_X^\bu)$ (see
Section~\ref{serre-duality}) provides a quasi-equivalence between
the DG\+categories $\sDG^\abs((X,\O,w)\coh)^\sop$ and
$\sDG^\abs((X,\O,-w)\coh)$.
 Furthermore, the external tensor product is a DG\+functor
\begin{multline}  \label{dg-external}
 \sDG^\abs((X,\O,-w)\coh)\ot_k\sDG^\abs((X,\O,w)\coh)\\
 \lrarrow\sDG^\abs((X\times_kX\;\O\;w_2-w_1)\coh),
\end{multline}
which, according to Proposition~\ref{external-tensor} and
Theorem~\ref{external-tensor}, induces an equivalence between
the derived categories of (left or right) DG\+modules over
the two DG\+categories in the left-hand and right-hand sides.
 Composing the Serre--Grothendieck duality with the external tensor
product, we obtain (perhaps, after replacing our DG\+categories
with naturally quasi-equivalent ones) a DG\+functor
\begin{multline} \label{dg-duality-external}
 \sDG^\abs((X,\O,w)\coh)^\sop\ot_k\sDG^\abs((X,\O,w)\coh)\\
 \lrarrow\sDG^\abs((X\times_kX\;\O\;w_2-w_1)\coh).
\end{multline}
having the same property with respect to the derived categories
of DG\+modules over the left-hand and right-hand sides as
the DG\+functor~\eqref{dg-external}.

 We are interested specifically in the diagonal right DG\+module
over $\sDG_w^\sop\ot_k\sDG_w$, that is the contravariant functor from
$\sDG^\sop_w\ot_k\sDG_w$ to the DG\+category of $\Z/2$\+graded complexes
of $k$\+vector spaces taking an object $(\M^\op,\K)$ to the complex
$\Hom_\sDG(\K,\M)$.
 It is claimed that the diagonal DG\+module is naturally
quasi-isomorphic to the DG\+module obtained by composing
the DG\+functor~\eqref{dg-duality-external} with the right DG\+module
over the right-hand side represented by the object $\Delta_*\D_X^\bu
\in\sDG^\abs((X\times_kX\;\O\;w_2-w_1)\allowbreak\coh)\subset
\sDG^\co((X\times_kX\;\O\;w_2-w_1)\qcoh)$.
{\emergencystretch=0em\hfuzz=3.1pt\par}

 Indeed, for any quasi-coherent matrix factorizations $\K$ and $\N$
of the potentials $w$ and $-w$ on $X$ there is a natural isomorphism
of $\Z/2$\+graded complexes of abelian groups
(see~\eqref{rhom-chom-diagonal})
$$
 \Hom_{(X,\O,w)\qcoh}(\K,\cHom_{\O_X}(\N,\D_X^\bu))\simeq
 \Hom_{(X\times_kX\;\O\;w_2-w_1)\qcoh}(\N\ot_k\K\;\Delta_*\D_X^\bu).
$$
 Proposition~\ref{internal-hom-matrix} shows how one can pass
from this isomorphism to a quasi-isomor\-phism of the similar complexes
of morphisms in the DG\+categories $\sDG^\co((X,\O,w)\allowbreak\qcoh)$
and $\sDG^\co((X\times_kX\;\O\;w_2-w_1)\qcoh)$.

 Now morphisms between representable DG\+modules in the derived
category of DG\+modules over a DG\+category $\sDG$ are computed by
the complex of morphisms in $\sDG$ between the representing objects,
so our proof is finished.
\end{proof}

\begin{rem}
 Given a separated scheme $X$ of finite type over a field~$k$, let
$\D_X^\bu$ be a dualizing complex on $X$ and $\D_{X\times_k X}^\bu$
be a dualizing complex on $X\times_k\nobreak X$ such that $\D_X^\bu
\simeq\boR\Delta^!(\D_{X\times_k X}^\bu)$.
 Then the anti-equivalence of absolute derived categories
$\cHom_{X\times_kX\qc}({-}\;\D_{X\times_kX}^\bu)\:
\sD^\abs((X\times_k X\;\O\;w_1-\nobreak w_2)\coh)^\sop\simeq
\sD^\abs((X\times_k\nobreak X\;\O\;w_2-\nobreak w_1)\coh)$ from
Proposition~\ref{serre-duality} transforms the object $\Delta_*
\O_X\in\sD^\abs((X\times_k\nobreak X\;\O\;w_1-\nobreak w_2)\coh)$
into the object $\Delta_*\D_X^\bu\in\sD^\abs((X\times_k\nobreak X\;
\allowbreak\O\;w_2-\nobreak w_1)\coh)$
(see Proposition~\ref{push-duality}).
 Therefore, in the assumptions of Theorem the Hochschild cohomology
algebra $HH^*(\sD^\abs((X,\O,w)\coh))$ of the DG\+category
$\sDG^\abs((X,\O,w)\coh)$ can be also identified with
the Ext algebra $\Hom_{\sD^\abs((X\times_kX\;\O\;w_1-w_2)\coh)}
(\Delta_*\O_X,\Delta_*\O_X[*])^\op$
(cf.\ Remark~\ref{coh-hochschild-homology} below). \hbadness=1250
\end{rem}

\subsubsection{Cotensor product of complexes of quasi-coherent sheaves}
\label{cotensor-complexes}
 Let $X$ be a separated Noetherian scheme.
 Then there is a tensor product functor on the coderived category
of ($\Z$\+graded complexes of) flat quasi-coherent sheaves on $X$
(cf.~\cite[Chapter~6]{M-th})
\begin{equation}  \label{tensor-fl-fl}
 \ot_{\O_X}\:\sD^\co(X\qcoh_\fl)\times\sD^\co(X\qcoh_\fl)\lrarrow
 \sD^\co(X\qcoh_\fl),
\end{equation}
and a similar functor of tensor product on the coderived categories
of flat and arbitrary quasi-coherent sheaves
(see~\cite[Section~4.12]{Pcosh})
\begin{equation}  \label{tensor-fl-qcoh}
 \ot_{\O_X}\:\sD^\co(X\qcoh_\fl)\times\sD^\co(X\qcoh)\lrarrow
 \sD^\co(X\qcoh).
\end{equation}

 Now let $\D_X^\bu$ be a dualizing complex on~$X$ (viewed, as usually,
as a finite complex of injective quasi-coherent sheaves).
 Then the equivalence of triangulated categories
$\sD^\co(X\qcoh_\fl)\simeq\sD^\co(X\qcoh)$ constructed using
the dualizing complex $\D_X^\bu$ (see~\cite[Chapter~8]{M-th})
transforms the tensor product functor~\eqref{tensor-fl-fl}
into the tensor product functor~\eqref{tensor-fl-qcoh}.
 One can use the same equivalence of categories to define a tensor
triangulated category structure with the unit object $\D_X^\bu$
on the coderived category $\sD^\co(X\qcoh)$.
 We call this operation the \emph{cotensor product} of complexes
of quasi-coherent sheaves on $X$ and denote it by
\begin{equation}  \label{cotensor-product}
 \oc_{\D_X^\bu}\:\sD^\co(X\qcoh)\times\sD^\co(X\qcoh)\lrarrow
 \sD^\co(X\qcoh).
\end{equation} 
 Explicitly, $\N^\bu\oc_{\D_X^\bu}\M^\bu = \D_X^\bu\ot_{\O_X}
\cHom_{X\qc}(\D_X^\bu,\N^\bu)\ot_{\O_X}\cHom_{X\qc}(\D_X^\bu,\M^\bu)$
for any complexes of injective quasi-coherent sheaves $\N^\bu$
and $\M^\bu$ on $X$ (cf.\ Lemma~\ref{gorenstein-case}(b)) and also
$\N^\bu\oc_{\D_X^\bu}\M^\bu = \cHom_{X\qc}(\D_X^\bu,\N^\bu)
\ot_{\O_X}\M^\bu$ for any complex of injective quasi-coherent
sheaves $\N^\bu$ and any complex of quasi-coherent sheaves $\M^\bu$
on~$X$.

 Recall that the full triangulated subcategory of bounded below
complexes in $\sD^\co(X\qcoh)$ is equivalent to $\sD^+(X\qcoh)$
(see~\cite[Lemma~2.1 and Remark~4.1]{Psemi}
or~\cite[Lemma~A.1.2]{Pcosh}).
 Denote by $\sD^+_\ch(X\qcoh)$ the full triangulated subcategory
in $\sD^+(X\qcoh)$ consisting of complexes with coherent
cohomology sheaves; then the category $\sD^+_\ch(X\qcoh)$ can be
also viewed as a full triangulated subcategory in $\sD^\co(X\qcoh)$.

 For any complexes of quasi-coherent sheaves $\N^\bu$ and $\M^\bu$
on $X$ there is a natural morphism of complexes of quasi-coherent
sheaves
\begin{multline}  \label{cotensor-symm-comparison}
 \D_X^\bu\ot_{\O_X}\cHom_{X\qc}(\D_X^\bu,\N^\bu)\ot_{\O_X}
 \cHom_{X\qc}(\D_X^\bu,\M^\bu) \\ \lrarrow
 \cHom_{X\qc}(\cHom_{X\qc}(\N^\bu,\D_X^\bu)\ot_{\O_X}
 \cHom_{X\qc}(\M^\bu,\D_X^\bu)\;\D_X^\bu)
\end{multline}
on $X$ defined in terms of the composition morphisms
$\cHom_{X\qc}(\D_X^\bu,\K^\bu)\ot_{\O_X}\cHom_{X\qc}(\K^\bu,\D_X^\bu)
\rarrow\cHom_{X\qc}(\D_X^\bu,\D_X^\bu)$ for complexes of
quasi-coherent sheaves $\K^\bu$ on $X$ and the natural
quasi-isomorphism $\D_X^\bu\ot_{\O_X}\cHom_{X\qc}(\D_X^\bu,\D_X^\bu)
\ot_{\O_X}\cHom_{X\qc}(\D_X^\bu,\D_X^\bu)\rarrow\D_X^\bu$.

\begin{thm}
 For any bounded below complexes of injective quasi-coherent sheaves
$\N^\bu$ and $\M^\bu$ with coherent cohomology sheaves on a separated
Noetherian scheme $X$ with a dualizing complex $\D_X^\bu$,
the natural morphism~\textup{\eqref{cotensor-symm-comparison}}
is a homotopy equivalence of bounded below complexes of injective
quasi-coherent sheaves on $X$ with coherent cohomology sheaves.
\end{thm}

\begin{proof}
 By Lemma~\ref{serre-duality}(b\+c), both sides
of~\eqref{cotensor-symm-comparison} are bounded below complexes of
injective quasi-coherent sheaves.
 Since the functor $\cHom_{X\qc}({-},\D_X^\bu)\:\sD(X\qcoh)\rarrow
\sD(X\qcoh)$ takes $\sD^+_\ch(X\qcoh)\subset\sD^+(X\qcoh)$ into
$\sD^-(X\coh)\subset\sD^-(X\qcoh)$ and vice versa, while the derived
tensor product functor $\ot^\boL_X\:\sD^-(X\qcoh)\times\sD^-(X\qcoh)
\rarrow\sD^-(X\qcoh)$ takes $\sD^-(X\coh)\times\sD^-(X\coh)$ into
$\sD^-(X\coh)$, the right-hand side has coherent cohomology sheaves.

 It remains to prove the homotopy equivalence claim.
 The homotopy category of bounded below complexes of injectives
being equivalent to $\sD^+(X\qcoh)$, one only has to check that
the map is a quasi-isomorphism. 
 Let us first show that it suffices to do so for complexes of
sheaves on affine open subschemes $U\subset X$.

 Indeed, for any quasi-coherent sheaves $\E$ and $\K$ on $X$
there is a natural morphism of quasi-coherent sheaves
$\cHom_{X\qc}(\E,\K)|_U\rarrow\cHom_{U\qc}(\E|_U,\K|_U)$ on~$U$.
 The morphism of complexes of quasi-coherent sheaves
$\cHom_{X\qc}(\E^\bu,\K^\bu)|_U\rarrow\cHom_{U\qc}(\E^\bu|_U,\K^\bu|_U)$
is a quasi-isomorphism whenever the complex $\E^\bu$ has coherent
cohomology, $\K^\bu$ is a complex of injective quasi-coherent
sheaves, and one of the complexes $\E^\bu$ and $\K^\bu$ is finite.

 Finally, the tensor product in the right-hand side preserves
quasi-isomorphisms of bounded above complexes of flat
quasi-coherent sheaves, while the one in the left-hand side
is well-defined on the coderived category of (complexes of)
flat quasi-coherent sheaves.
 It remains to notice that the functor $\Hom_{X\qc}(\D_X^\bu,{-})$
in the equivalence of categories in Theorem~\ref{serre-duality}
agrees with the restrictions to open subschemes, since so does
its inverse functor $\D_X^\bu\ot_{\O_X}{-}$.

 Now that we are on an affine scheme~$U$, pick bounded above complexes
of vector bundles ${}'\!\N^\bu$ and ${}'\!\M^\bu$ isomorphic to
$\cHom_{U\qc}(\N^\bu,\D_U^\bu)$ and $\cHom_{U\qc}(\M^\bu,\D_U^\bu)$,
respectively, in $\sD^-(U\coh)$.
 Given the isomorphisms $\cHom_{U\qc}(\D_U^\bu,\cHom_{U\qc}({}'\!\N^\bu,
\D_U^\bu))\allowbreak\simeq\cHom_{U\qc}({}'\!\N^\bu,\cHom_{U\qc}
(\D_U^\bu,\D_U^\bu))\simeq\cHom_{U\qc}({}'\!\N^\bu,\O_U)$ and similarly
for ${}'\!\M^\bu$ in $\sD^\co(U\qcoh_\fl)$, the assertion reduces to
the obvious isomorphism of complexes $\D_X^\bu\ot_{\O_X}
\cHom_{U\qc}({}'\!\N^\bu,\O_U)\ot_{\O_U}\cHom_{U\qc}({}'\!\M^\bu,\O_U)
\simeq\Hom_{U\qc}({}'\!\N^\bu\ot_{\O_X}{}'\!\M^\bu\;\D_U^\bu)$.
\end{proof}

 For any complexes of quasi-coherent sheaves $\K^\bu$ and $\M^\bu$
on $X$ we denote by $\cHom_{X\qc}^\oplus(\K^\bu,\M^\bu)$ the complex
of quasi-coherent sheaves on $X$ obtained by totalizing
the bicomplex of quasi-coherent internal Hom sheaves
$\cHom_{X\qc}(\K^i,\M^j)$ by taking infinite direct sums along
the diagonals.
 Assuming that $\M^\bu$ is a complex of injective quasi-coherent
sheaves, the complex $\cHom_{X\qc}^\oplus(\K^\bu,\M^\bu)$ is
absolutely acyclic with respect to $X\qcoh$ whenever
the complex $\K^\bu$~is (see Lemma~\ref{serre-duality}(a)).

 In the same assumption, the complex $\cHom_{X\qc}^\oplus(\K^\bu,\M^\bu)$
is also coacyclic with respect to $X\qcoh$ whenever the complex of
quasi-coherent sheaves $\K^\bu$ is acyclic and bounded from
above~\cite[Lemma~2.1]{Psemi}.
 Therefore, representing the second argument of $\cHom_{X\qc}^\oplus$ by
complexes of injectives, one can construct the right derived functors
\begin{equation}  \label{hom-oplus-abs}
\boR\cHom_{X\qc}^\oplus\:
\sD^\abs(X\qcoh)^\sop\times\sD^\co(X\qcoh)\lrarrow\sD^\abs(X\qcoh)
\end{equation}
and
\begin{equation}  \label{hom-oplus-minus}
\boR\cHom_{X\qc}^\oplus\:
\sD^-(X\qcoh)^\sop\times\sD^\co(X\qcoh)\lrarrow\sD^\co(X\qcoh)
\end{equation}
of the functor $\cHom_{X\qc}^\oplus$.

 For any complexes of quasi-coherent sheaves $\N^\bu$ and $\M^\bu$ on $X$
there is a natural morphism of complexes of quasi-coherent sheaves
\begin{equation}  \label{cotensor-asymm-comparison}
 \N^\bu\ot_{\O_X}\cHom_{X\qc}(\D_X^\bu,\M^\bu) \lrarrow
 \cHom_{X\qc}^\oplus(\cHom_{X\qc}(\N^\bu,\D_X^\bu)\;\M^\bu)
\end{equation}
on $X$ defined in terms of the composition/evaluation morphism
$\N^\bu\ot_{\O_X}\cHom_{X\qc}(\N^\bu,\allowbreak\D_X^\bu)\ot_{\O_X}
\cHom_{X\qc}(\D_X^\bu,\M^\bu)\rarrow\M^\bu$.

\begin{prop}
 For any bounded below complex of injective quasi-coherent sheaves
$\N^\bu$ with coherent cohomology sheaves and any complex of
injective quasi-coherent sheaves $\M^\bu$ on $X$,
the natural morphism~\textup{\eqref{cotensor-asymm-comparison}}
is a homotopy equivalence of complexes of injective quasi-coherent
sheaves on $X$.
\end{prop}

\begin{proof}
 It suffices to check that the morphism~\eqref{cotensor-asymm-comparison}
is an isomorphism in $\sD^\co(X\qcoh)$.
 Both sides of the desired isomorphism being well-defined as functors
of the argument $\N^\bu\in\sD^+(X\qcoh)$ taking values in
$\sD^\co(X\qcoh)$, one can freely replace $\N^\bu$ with any
quasi-isomorphic bounded below complex of quasi-coherent sheaves.
 The same applies to the bounded above complex $\cHom_{X\qc}(\N^\bu,
\D_X^\bu)$ in the right-hand side of~\eqref{cotensor-asymm-comparison}.

 All the functors involved being local in $X$ up to isomorphism in
the relevant triangulated categories, it suffices to consider
complexes of sheaves over affine open subschemes $U\subset X$
(see Remark~\ref{second-kind}).
 Representing the object $\cHom_{U\qc}(\N^\bu,\D_U^\bu)\in\sD^-(U\coh)
\subset\sD^-(U\qcoh)$ by a bounded above complex of vector bundles
${}'\!\N^\bu$, it remains to notice the isomorphism of complexes
$\cHom_{U\qc}({}'\!\N^\bu,\D_U^\bu)\ot_{\O_U}\F^\bu\simeq
\cHom_{U\qc}^\oplus({}'\!\N^\bu\;\D_U^\bu\ot_{\O_U}\F^\bu)$ for any
complex of quasi-coherent sheaves $\F^\bu$ on $U$ and point out that
the functor $\cHom_{U\qc}^\oplus({}'\!\N^\bu,{-})$ takes a homotopy
equivalence $\D_U^\bu\ot_{\O_U}\cHom_{U\qc}(\D_U^\bu,\M^\bu)
\rarrow\M^\bu$ to a homotopy equivalence. 
\end{proof}

 In the particular cases when either $\N^\bu$ is a finite complex
of quasi-coherent sheaves with coherent cohomology sheaves, or
$\N^\bu$ is a bounded below complex of quasi-coherent sheaves with
coherent cohomology sheaves and $\M^\bu$ is a bounded below complex
of quasi-coherent sheaves, the direct sum totalization of
the bicomplex $\cHom_{X\qc}$ in the right-hand side of
the isomorphism~\eqref{cotensor-asymm-comparison} in
the coderived category $\sD^\co(X\qcoh)$ is no different from
the conventional direct product totalization.

 Finally, let $X$ be a separated scheme of finite type over a field~$k$
and $\pi\:X\rarrow\Spec k$ be its structure morphism.
 Then $\D_X^\bu\simeq\pi^+\O_{\Spec k}$ (see Section~\ref{push-duality})
is a natural choice of the dualizing complex on~$X$.
 Let $\pi\times_k\pi\:X\times_k X\rarrow\Spec k$ be the structure
morphism of the Cartesian square of $X$ over~$k$.
 Then the dualizing complex $\D_{X\times_kX}^\bu=(\pi\times_k\pi)^+
\O_{\Spec k}$ on $X\times_kX$ is quasi-isomorphic to the external
tensor product $\D_X^\bu\ot_k\D_X^\bu$, and one has $\D_X^\bu\simeq
\Delta^+(\D_{X\times_kX}^\bu)\simeq\boR\Delta^!(\D_X^\bu\ot_k\D_X^\bu)$,
where $\Delta\:X\rarrow X\times_kX$ denotes the diagonal map.

 The equivalence of categories $\sD^\co(X\times_kX\qcoh_\fl)\simeq
\sD^\co(X\times_kX\qcoh)$ constructed using the dualizing complex
$\D_{X\times_kX}^\bu$ and the similar equivalence $\sD^\co(X\qcoh_\fl)
\simeq\sD^\co(X\qcoh)$ constructed using the dualizing complex
$\D_X^\bu$ transform the external tensor product functor
$$
 \ot_k\:\sD^\co(X\qcoh_\fl)\times\sD^\co(X\qcoh_\fl)\lrarrow
 \sD^\co(X\times_kX\qcoh_\fl)
$$
into the external tensor product functor
$$
 \ot_k\:\sD^\co(X\qcoh)\times\sD^\co(X\qcoh)\lrarrow
 \sD^\co(X\times_kX\qcoh),
$$
since so do the functors $\D_X^\bu\ot_{\O_X}{-}$ and
$\D_{X\times_kX}^\bu\ot_{\O_{X\times_kX}}{-}$.

 Let $\N^\bu$ and $\M^\bu$ be two complexes of injective
quasi-coherent sheaves on~$X$, and let $\J^\bu$ be a complex of
injective quasi-coherent sheaves on $X\times_kX$ isomorphic to
$\N^\bu\ot_k\M^\bu$ in $\sD^\co(X\times_kX\qcoh)$.
 Then in the coderived categories of quasi-coherent sheaves one has
\begin{multline*}
\N^\bu\oc_{\D_X^\bu}\M^\bu = \D_X^\bu\ot_{\O_X}
\Delta^*(\cHom_{\O_X}(\D_X^\bu,\N^\bu)\ot_k
\cHom_{\O_X}(\D_X^\bu,\M^\bu)) \\ \simeq \D_X^\bu\ot_{\O_X}
\Delta^*\cHom_{\O_X}(\D_{X\times_kX}^\bu,\J^\bu)
\simeq\boR\Delta^!(\N^\bu\ot_k\M^\bu)
\end{multline*}
by the result of~\cite[Theorem~5.15.3]{Pcosh} applied to
the proper morphism (actually, closed embedding)~$\Delta$.
 We have obtained the formula
\begin{equation}
 \N^\bu\oc_{\pi^+\O_{\Spec k}}\M^\bu \.\simeq\.
 \boR\Delta^!(\N^\bu\ot_k\M^\bu)
\end{equation}
for the cotensor product of complexes of quasi-coherent sheaves
on the scheme~$X$ (see the end of Section~\ref{pull-push} for
the notation $\boR f^!$ as applied to objects of the coderived
category of quasi-coherent sheaves).

\subsubsection{Cotensor product of matrix factorizations}
\label{cotensor-matrix}
  The equivalences of triangulated categories $\sD^\co((X,\L,w'')
\qcoh_\fl)\simeq\sD^\co((X,\L,w'')\qcoh)$ and $\sD^\co((X\;\L\;
w'+w'')\allowbreak\qcoh_\fl)\simeq\sD^\co((X\;\L\;w'+w'')\qcoh)$
constructed using a dualizing complex $\D_X^\bu$ (see
Section~\ref{serre-duality}) transform the tensor product
functor~\eqref{matrix-tensor-fl-fl} into the tensor product
functor~\eqref{matrix-tensor-fl-qcoh}.
 So one can use the same equivalences of categories together with
the similar equivalence $\sD^\co((X,\L,w')\qcoh_\fl)\simeq
\sD^\co((X,\L,w')\qcoh)$ (constructed using the same dualizing
complex~$\D_X^\bu$) in order to define a triangulated functor
of two arguments
\begin{multline}   \label{matrix-cotensor}
 \oc_{\D_X^\bu}\:\sD^\co((X,\L,w')\qcoh)\times\sD^\co((X,\L,w'')\qcoh)
 \\ \lrarrow\sD^\co((X\;\L\;w'+w'')\qcoh),
\end{multline}
which we call the \emph{cotensor product} of matrix factorizations.

 As in the case of complexes of quasi-coherent sheaves, one
explicitly has $\N\oc_{\D_X^\bu}\M=\D_X^\bu\ot_{\O_X}
\cHom_{X\qc}(\D_X^\bu,\N)\ot_{\O_X}\cHom_{X\qc}(\D_X^\bu,\M)$ for any
injective quasi-coherent matrix factorizations $\N$ and $\M$ on $X$,
and also $\N\oc_{\D_X^\bu}\M = \cHom_{X\qc}(\D_X^\bu,\N)\ot_{\O_X}\M$
for any injective quasi-coherent matrix factorization $\N$ and any
quasi-coherent matrix factorization $\M$ on~$X$.
 As in Section~\ref{internal-hom-matrix}, $\,\N$ and $\M$ must be matrix
factorizations of two sections $w'$ and~$w''$ of the same line bundle
$\L$ on a scheme $X$; then the cotensor product $\N\oc_{\D_X^\bu}\M$ is
a matrix factorization of the section $w'+w''$ of the line bundle~$\L$.
{\hbadness=1150\par}

\begin{rem}
 While a matrix factorization version of
Proposition~\ref{cotensor-complexes} is presented below,
Remark~\ref{internal-hom-matrix} explains the reason why a matrix
factorization version of Theorem~\ref{cotensor-complexes} cannot be
formulated in the way similar to the version for complexes of
quasi-coherent sheaves above.
 Still, let $\N$ and $\M$ be coherent matrix factorizations of
sections $w'$ and~$w''$ of the same line bundle $\L$ on a separated
Noetherian scheme $X$ with enough vector bundles.
 Let $\P$ and $\Q$ be coherent matrix factorizations of
the potentials $-w'$ and $-w''\in\L(X)$ isomorphic to
$\cHom_{X\qc}(\N,\D_X^\bu)$ and $\cHom_{X\qc}(\M,\D_X^\bu)$ in
the respected coderived categories.

 Let $\E_\bu$ and $\F_\bu$ be left resolutions of the matrix
factorizations $\P$ and $\Q$ by locally free matrix
factorizations of finite rank (of the respected potentials).
 Then the totalizations of the bounded below complexes of matrix
factorizations $\cHom_{X\qc}(\E_\bu,\D_X^\bu)$, \
$\cHom_{X\qc}(\F_\bu,\D_X^\bu)$, and $\cHom_{X\qc}
(\E_\bu\ot_{\O_X}\F_\bu\;\D_X^\bu)$ represent objects naturally
isomorphic to $\N$, \ $\M$ and $\N\oc_{\D_X^\bu}\M$ in the coderived
categories of matrix factorizations of the potentials $w'$, $w''$,
and $w'+w''$ (cf.\ Corollary~\ref{serre-duality}).
\end{rem}

 For any quasi-coherent matrix factorizations $\M$ and $\N$ of sections
$w'$ and~$w''$ of the same line bundle $\L$ on the scheme $X$, there is
a natural morphism of quasi-coherent matrix factorizations of
the section $w'+w''$ of the line bundle $\L$ on~$X$
\begin{equation}  \label{cotensor-matrix-asymm-comparison}
 \N\ot_{\O_X}\cHom_{X\qc}(\D_X^\bu,\M) \lrarrow
 \cHom_{X\qc}(\cHom_{X\qc}(\N,\D_X^\bu)\;\M)
\end{equation}
constructed in the same way as it was done for complexes of
quasi-coherent sheaves in Section~\ref{cotensor-complexes}.

\begin{prop}
 For any coherent matrix factorization $\N$ and injective
quasi-coherent matrix factorization $\M$ of sections $w'$ and~$w''$
of the same line bundle $\L$ on a separated Noetherian scheme $X$,
the natural morphism~\textup{\eqref{cotensor-matrix-asymm-comparison}}
is an isomorphism in the coderived category of quasi-coherent matrix
factorizations of the potential $w'+w''\in\L(X)$.
\end{prop}

\begin{proof}
 The argument follows the lines of the proof of
Proposition~\ref{cotensor-complexes}.
 The left-hand side of the desired isomorphism is well-defined as
a functor of the argument $\N\in\sD^\co((X,\L,w'')\qcoh)$ taking
values in $\sD^\co((X\;\L\;w'+w'')\qcoh)$, while the right-hand
side is well-defined as a functor of the argument
$\N\in\sD^\abs((X,\L,w'')\coh)$ taking values, say, in the same
coderived category.
 Besides, the right-hand side, viewed as an object of the coderived
category, only depends on the matrix factorization
$\cHom_{X\qc}(\N,\D_X^\bu)$ viewed as an object of the absolute
derived category.

 Furthermore, the contravariant Serre--Grothendieck duality
$\cHom_{X\qc}({-},\D_X^\bu)$ is well-defined as a functor
$\sD^\abs((X,\L,w'')\qcoh)\rarrow\sD^\abs((X,\L,-w'')\qcoh)$
and takes $\sD^\abs((X,\L,w'')\coh)\subset\sD^\abs((X,\L,w'')\qcoh)$
into $\sD^\abs((X,\L,-w'')\coh)\subset\sD^\abs((X,\L,-w'')\coh)$,
inducing an equivalence between these two subcategories
(see Proposition~\ref{serre-duality}).
 In particular, one can conclude that all the functors involved are
local in $X$, and it suffices to prove the desired assertion for
matrix factorizations over affine open subschemes $U\subset X$.

 Now let $\K$ be a coherent matrix factorization of
the potential $-w''$ isomorphic to $\cHom_{U\qc}(\N,\D_U^\bu)$
in $\sD^\abs((U,\L,-w'')\qcoh)$, and let $\E_\bu$ be its left
resolution by locally free matrix factorizations of
the same potential $-w''\in\L(U)$.
 Then the matrix factorization $\cHom_{X\qc}(\K,\M)$ is isomorphic
in $\sD^\co((X\;\L\;w'+w'')\qcoh)$ to the totalization of
the complex of matrix factorizations $\cHom_{X\qc}(\E_\bu,\M)$
constructed by taking infinite direct sums along the diagonals;
and the matrix factorization $\N\simeq\cHom_{X\qc}(\K,\D_X^\bu)$
can be described similarly (cf.\ the proof of
Corollary~\ref{serre-duality}).

 It remains to notice that the functor of tensoring with
$\cHom_{X\qc}(\E_\bu,\O_X)$ and totalizing by taking infinite
direct sums along the diagonals takes the homotopy equivalence
$\D_X^\bu\ot_{\O_X}\cHom_{X\qc}(\D_X^\bu,\M)\rarrow\M$
to a homotopy equivalence of matrix factorizations.
\end{proof}

 As in Section~\ref{cotensor-complexes}, we finish by discussing
the case of a separated scheme $X$ of finite type over a field~$k$.
 From now on we also assume that $\L=\O_X$.
 So let $w'$, $w''\in\O(X)$ be two global regular functions on~$X$;
as in Section~\ref{external-tensor}, we consider the regular
function $w_1'+w_2''=p_1^*w'+p_2^*w''$ on $X\times_k X$. 
 We use the dualizing complexes $\D_X^\bu=\pi^+\O_{\Spec k}$ and
$\D_{X\times_kX}^\bu = (\pi\times_k\pi)^+\O_{\Spec k}$.

 The equivalence of categories $\sD^\co((X\times_k X\;\O\;w'_1+w''_2)
\qcoh_\fl)\simeq\sD^\co((X\times_k X\;\allowbreak\O\;w'_1+w''_2)\qcoh)$
constructed using the dualizing complex $\D_{X\times_k X}^\bu$
and the similar equivalences of coderived categories of matrix
factorizations of the potentials $w'$ and~$w''$ on $X$ constructed
using the dualizing complex $\D_X^\bu$ transform the external
tensor product functor (cf.~\eqref{qcoh-external-tensor})
\begin{multline*}
 \ot_k\:\sD^\co((X,\O,w')\qcoh_\fl)\times\sD^\co((X,\O,w'')\qcoh_\fl)
 \\ \lrarrow \sD^\co((X\times_k X\;\O\;w'_1+w''_2)\qcoh_\fl)
\end{multline*}
into the external tensor product functor
\begin{multline*}
 \ot_k\:\sD^\co((X,\O,w')\qcoh)\times\sD^\co((X,\O,w'')\qcoh)
 \\ \lrarrow \sD^\co((X\times_k X\;\O\;w'_1+w''_2)\qcoh),
\end{multline*}
since so do the functors $\D_X^\bu\ot_{\O_X}{-}$ and
$\D_{X\times_kX}^\bu\ot_{\O_{X\times_kX}}{-}$.

 Let $\N$ and $\M$ be injective quasi-coherent matrix factorizations
of the potentials $w'$ and~$w''$ on~$X$, and let $\J$ be an injective
quasi-coherent matrix factorization of the potential $w'_1+w''_2$
on $X\times_kX$ isomorphic to $\N\ot_k\M$ in
$\sD^\co((X\times_k X\;\O\;w'_1+w''_2)\qcoh)$.
 Then in the coderived categories of quasi-coherent matrix
factorizations one has
\begin{multline*}
\N\oc_{\D_X^\bu}\M = \D_X^\bu\ot_{\O_X}
\Delta^*(\cHom_{\O_X}(\D_X^\bu,\N)\ot_k\cHom_{\O_X}(\D_X^\bu,\M))
\\ \simeq \D_X^\bu\ot_{\O_X}\Delta^*\cHom_{\O_X}
(\D_{X\times_kX}^\bu,\J) \simeq \boR\Delta^!(\N\ot_k\M)
\end{multline*}
by the result of Theorem~\ref{pull-duality} applied to
the proper morphism~$\Delta$.
 We have obtained the formula
\begin{equation}  \label{cotensor-matrix-diagonal}
 \N\oc_{\pi^+\O_{\Spec k}}\M \.\simeq\.\boR\Delta^!(\N\ot_k\M)
\end{equation}
for the cotensor product of quasi-coherent matrix factorizations
on the scheme~$X$.

\subsubsection{Hochschild homology}  \label{coh-hochschild-homology}
 Let $X$ be a separated scheme of finite type over a field~$k$ and
$\pi\:X\rarrow\Spec k$ be its structure morphism.
 Let $w\in\O(X)$ be a global regular function; as in
Section~\ref{coh-hochschild-cohomology}, we assume that
the morphism of schemes $w\:X\rarrow\mathbb A^1_k$ is flat.
 Consider the scheme $X\times_k X$ and endow it with
the potential $w_2-w_1=p_2^*(w)-p_1^*(w)$.
 Let $\Delta\:X\rarrow X\times_k X$ denote the diagonal morphism.

\begin{thm} \hfuzz=15pt
 In the assumptions of Theorem~\textup{\ref{coh-hochschild-cohomology}},
there is a natural isomorphism between the Hochschild homology module
$HH_*(\sDG^\abs((X,\O,w)\coh))$ over the algebra
$HH^*(\sDG^\abs((X,\O,w)\coh))$ and the Ext module\/
$\Hom_{\sD^\co((X\times_kX\;\O\;w_2-w_1)\qcoh)}\allowbreak
(\Delta_*\O_X\;\Delta_*\D_X^\bu[*])$ over the algebra\/
$\Hom_{\sD^\co((X\times_kX\;\O\;w_2-w_1)\qcoh)}
(\Delta_*\D_X^\bu\;\Delta_*\D_X^\bu[*])$.
 Here $\D_X^\bu$ denotes the dualizing complex $\pi^+\O_{\Spec k}$
on~$X$.
\end{thm}

\begin{proof}
 By the definition, the Hochschild homology of a $\Z/2$\+graded
DG\+category $\sDG$ is the $\Z/2$\+graded vector space
$\Tor^{\sDG\ot_k\sDG^\sop}_*(\sDG,\sDG)$ for the diagonal right
and left DG\+modules $\sDG$ over the DG\+category
$\sDG\ot_k\sDG^\sop$ \cite[Sections~2.4 and~3.1]{PP}.
 This is the conventional derived tensor product (``of the first kind'')
of a left and a right DG\+module over a small DG\+category.
 The Hochschild cohomology algebra $\Hom_{\sD(\sDG\ot_k\sDG^\sop)}
(\sDG,\sDG[*])$ acts on the Hochschild homology space via its action
on, say, the first argument of the Tor.

 As in the proof of Theorem~\ref{coh-hochschild-cohomology}, we set
$\sDG_w=\sDG^\abs((X,\O,w)\coh)$; accordingly,
$\sDG_{-w}=\sDG^\abs((X,\O,-w)\coh)$ and
$\sDG_{w_2-w_1}=\sDG^\abs((X\times_k X\;\O\;w_2-w_1)\coh)$.
 The DG\+functor $\sDG_w^\sop\ot_k\sDG_w\rarrow\sDG_{w_2-w_1}$
\eqref{dg-duality-external} induces a fully faithful functor
between the homotopy categories $H^0(\sDG_w)^\sop\ot_k
H^0(\sDG_w)\rarrow H^0(\sDG_{w_2-w_1})$ such that every object
in the target category can be obtained from objects in the image
using the operations of a cone and the passage to a direct summand.

 Let $\sDG(\modr\sDG_w^\sop\ot_k\sDG_w)$ denote the DG\+category
version of the (conventional) derived category of right DG\+modules
over the DG\+category $\sDG_w^\sop\ot_k\sDG$ (i.~e., contravariant
DG\+functors from $\sDG_w^\sop\ot_k\sDG$ into the DG\+category
$\sDG(k\vect)$ of $\Z/2$\+graded complexes of $k$\+vector spaces).
 Let $\sDG(\modr\sDG_w^\sop\ot_k\sDG_w)^0 \subset
\sDG(\modr\allowbreak\sDG_w^\sop\ot_k\sDG_w)$ denote the full
DG\+subcategory of DG\+modules corresponding to compact objects of
the derived category $\sD(\modr\sDG_w^\sop\ot_k\sDG_w)$
of right DG\+modules.

 The derived tensor product with the left DG\+module $\sDG_w$ over
$\sDG_w^\sop\ot_k\sDG_w$ can be viewed as a covariant DG\+functor
$\sDG(\modr\sDG_w^\sop\ot_k\sDG_w)\rarrow\sDG(k\vect)$.
 We are interested in the restriction of this DG\+functor to
the DG\+subcategory $\sDG(\modr\sDG_w^\sop\ot_k\sDG_w)^0$; let us
denote it by $F\:\sDG(\modr\sDG_w^\sop\ot_k\sDG_w)^0\rarrow
\sD(k\vect)$.

 There is a natural DG\+functor $\sDG_w^\sop\ot_k\sDG_w\rarrow
\sDG(\modr\sDG_w^\sop\ot_k\sDG_w)^0$ assigning to any object of
$\sDG_w^\sop\ot_k\sDG_w$ the contravariant DG\+functor
represented by it.
 Similarly one constructs a DG\+functor $\sDG_{w_2-w_1}\rarrow
\sDG(\modr\sDG_w^\sop\ot_k\sDG_w)^0$ whose composition with
the DG\+functor $\sDG_w^\sop\ot_k\sDG_w\rarrow\sDG_{w_2-w_1}$
is naturally quasi-isomorphic to the DG\+functor
$\sDG_w^\sop\ot_k\sDG_w\rarrow\sDG(\modr\sDG_w^\sop\ot_k\sDG_w)^0$.

 It is claimed that the composition of the DG\+functor
$\sDG_{w_2-w_1}\rarrow\sDG(\modr\sDG_w^\sop\allowbreak\ot_k\sDG_w)^0$
with the DG\+functor $F\:\sDG(\modr\sDG_w^\sop\ot_k\sDG_w)^0\rarrow
\sD(k\vect)$ is naturally quasi-isomorphic to
the DG\+functor $\Hom_{\sDG_{w_2-w_1}}(\Delta_*\O_X\;{-})$.
 The derived categories of left DG\+modules over $\sDG_{w_1-w_2}$
and $\sDG_w^\sop\ot_k\sDG_w$ being equivalent, it suffices to
construct a quasi-isomorphism between the compositions of
the two DG\+functors in question with the DG\+functor
$\sDG_w^\sop\ot_k\sDG_w\rarrow\sDG_{w_2-w_1}$.

 Indeed, let $(\K^\op,\M)$ be an object of $\sDG_w^\sop\ot_k\sDG_w$.
 Then the functor of (derived or underived) tensor product with
the diagonal left DG\+module $\sDG_w$ takes the right DG\+module
over $\sDG_w^\sop\ot_k\sDG_w$ represented by $(\K^\op,\M)$ to
the complex of $k$\+vector spaces $\Hom_{\sDG_w}(\K,\M)$.
 Substituting $\K=\cHom_{X\qc}(\N,\D_X^\bu)$ with $\N\in\sDG_{-w}$
and assuming $\M$ to be represented by an injective matrix 
factorization isomorphic to the given coherent one in
$\sDG^\co((X,\O,w)\qcoh)$, we have to compute the complex
of $k$\+vector spaces $\Hom_{(X,\O,w)\qcoh}
(\cHom_{X\qc}(\N,\D_X^\bu),\M)$.

 Now the formula~\eqref{matrix-rhom-gamma-chom} together with
Lemma~\ref{internal-hom-matrix} allow to interpret this complex
as $\boR\Gamma(X,\cHom_{X\qc}(\cHom_{X\qc}(\N,\D_X^\bu),\M))$.
 According to Proposition~\ref{cotensor-matrix} together with
the formula~\eqref{cotensor-matrix-diagonal}, this is the same as
$\boR\Gamma(X\;\boR\Delta^!(\N\ot_k\M))$, or, in other notation,
$\Hom_{\sDG^\co((X,\O,0)\qcoh)}(\O_X\;\boR\Delta^!(\N\ot_k\M))$.
 Finally, the adjunction of $\Delta_*$ and $\boR\Delta^!$ allows
to rewrite the complex in question as $\Hom_{\sDG^\co((X\times_k X\;
\O\;w_2-w_1)\qcoh)}(\Delta_*\O_X\;\allowbreak\N\ot_k\M)$.
 The desired quasi-isomorphism of DG\+functors is obtained.

 It remains to recall that, according to the proof of
Theorem~\ref{coh-hochschild-cohomology}, the diagonal right
DG\+module $\sDG_w$ over $\sDG_w^\sop\ot_k\sDG_w$ is represented
by the object $\Delta_*\D_X^\bu\in\sDG_{w_2-w_1}$, in order to
finish our proof here.
\end{proof}

\begin{rem} \hfuzz=3.5pt
 According to Remark~\ref{coh-hochschild-cohomology}, the Hochschild
homology module $HH_*((\sDG^\abs\allowbreak(X,\O,w)\coh))$
over the Hochschild cohomology algebra $HH^*((\sDG^\abs(X,\O,w)\coh))$
can be also computed as the Ext module
$\Hom_{\sD^\co((X\times_k X\;\O\;w_1-w_2)\qcoh)}(\Delta_*\O_X\;
\Delta_*\D_X^\bu[*])$ over the Ext algebra
$\Hom_{\sD^\abs((X\times_k X\;\O\;w_1-w_2)\coh)}(\Delta_*\O_X,
\Delta_*\O_X[*])^\op$.
 Moreover, the contravariant Serre duality for matrix factorizations
over $X\times_k X$ can be used in order to obtain an alternative
proof of our Hochschild homology computation.
 Indeed, for any coherent matrix factorizations $\N$ and $\M$ of
the potentials $-w$ and~$w$ on $X$ there are natural quasi-isomorphisms
\begin{multline*}
 \Hom_{\sDG^\abs((X,\O,w)\coh)}(\cHom_{X\qc}(\N,\D_X^\bu),\M) \\
 \simeq \Hom_{\sDG^\abs((X,\O,w)\coh)}(\cHom_{X\qc}(\N,\D_X^\bu)\;
 \cHom_{X\qc}(\cHom_{X\qc}(\M,\D_X^\bu),\D_X^\bu)) \displaybreak[0]\\
 \simeq \Hom_{\sDG^\abs((X\times_kX\;\O\;w_1-w_2)\coh)}
 (\cHom_{X\qc}(\N,\D_X^\bu)\ot_k\cHom_{X\qc}(\M,\D_X^\bu)\;
 \Delta_*\D_X^\bu) \displaybreak[0]\\
 \simeq \Hom_{\sDG^\abs((X\times_kX\;\O\;w_1-w_2)\coh)}
 (\cHom_{X\times_kX\qc}(\N\ot_k\M\;\D_{X\times_k X}^\bu)\;
 \Delta_*\D_X^\bu) \\
 \simeq\Hom_{\sDG^\abs((X\times_kX\;\O\;w_2-w_1)\coh)}
 (\Delta_*\O_X\;\N\ot_k\M)
\end{multline*}
by Proposition~\ref{serre-duality} and the proof of
Theorem~\ref{coh-hochschild-cohomology}.
 In other words, while the right diagonal DG\+module $\sDG_w$ over
$\sDG_w^\sop\ot_k\sDG_w$ is represented by the object $\Delta_*\D_X^\bu
\in\sDG^\abs((X\times_kX\;\O\;w_2-w_1)\coh)$ as a contravariant
DG\+functor on $\sDG_w^\sop\ot_k\sDG_w\subset\sDG_{w_2-w_1}$, the left
diagonal DG\+bimodule $\sDG_w$ over $\sDG_w^\sop\ot_k\sDG_w$ is
represented by the object $\Delta_*\O_X^\bu\in
\sDG^\abs((X\times_kX\;\O\;w_2-w_1)\coh)$ as a covariant DG\+functor
on $\sDG_w^\sop\ot_k\sDG_w\subset\sDG_{w_2-w_1}=
\sDG^\abs((X\times_kX\;\O\;w_2-w_1)\coh)$.
\end{rem}

\subsubsection{Direct sum over the critical values}
\label{direct-sum-formula}
 Let $X$ be a separated scheme of finite type over a field~$k$
and $\pi\:X\rarrow \Spec k$ be its structure morphism.
 As in
Sections~\ref{cotensor-complexes}\+-\ref{coh-hochschild-homology}
(see also Section~\ref{push-duality}), we choose the dualizing complex
$\D_X^\bu\simeq\pi^+\O_{\Spec k}$ on~$X$.
 Let $w\in\O(X)$ be a global regular function on $X$ such that
the morphism of schemes $w\:X\rarrow\mathbb A^1_k$ is flat
(cf.~\cite{Or1,Or3}).

 Let $c_1$,~$\dotsc$, $c_n\in k$ be a finite number of different
elements of the ground field.
 Assume that there exist closed subschemes $Z_i\subset X$ such that
the function~$w$ is noncritical on $X\setminus(Z_1\cup\dotsb\cup Z_n)$,
the restriction of~$w$ to $Z_i$ is equal to the constant~$c_i$,
and the schemes $Z_i$ admit smooth stratifications over~$k$.

 In particular, if the field~$k$ is perfect, it suffices to require
that the function~$w$ has only a finite number of critical values
$c_1$,~$\dotsc$, $c_n\in\mathbb A^1_k$ (i.~e., the open subscheme
$\mathbb A^1_{k,f}\subset\mathbb A^1_k$ is nonempty; see
Section~\ref{noncritical}), and all of these values belong to
the field~$k$ (rather than its algebraic closure).
 When the field~$k$ has zero characteristic, the former condition
holds automatically.
 Then one simply takes $Z_i$ to be the zero locus of the function
$w_i-c_i$ on~$X$.

 Consider the Cartesian square $X\times_k X$ with the global function
$w_2-w_1=p_2^*(w)-p_1^*(w)$ on it.
 Let $\Delta\:X\rarrow X\times_k X$ denote the diagonal morphism.
 The following result is to be compared with~\cite[Corollary~4.10]{PP}.

\begin{cor}
 There are natural isomorphisms of\/ $\Z/2$\+graded $k$\+algebras
\begin{multline} \label{hoch-cohomology-direct-sum}
 \textstyle\bigoplus_{i=1}^n HH^*(\sDG^\abs((X\;\O\;w-c_i)\coh))
 \\ \simeq\.\Hom_{\sD^\co((X\times_kX\;\O\;w_2-w_1)\qcoh)}
 (\Delta_*\D_X^\bu\;\Delta_*\D_X^\bu[*]) \\ \simeq\.
 \Hom_{\sD^\abs((X\times_kX\;\O\;w_1-w_2)\coh)}
 (\Delta_*\O_X,\Delta_*\O_X[*])^\op.
\end{multline}
 There are also natural isomorphisms of\/ $\Z/2$\+graded $k$\+modules
\begin{multline} \label{hoch-homology-direct-sum}
 \textstyle\bigoplus_{i=1}^n HH_*(\sDG^\abs((X\;\O\;w-c_i)\coh))
 \\ \simeq\.\Hom_{\sD^\co((X\times_kX\;\O\;w_2-w_1)\qcoh)}
 (\Delta_*\O_X\;\Delta_*\D_X^\bu[*]) \\ \simeq\.
 \Hom_{\sD^\co((X\times_kX\;\O\;w_1-w_2)\qcoh)}
 (\Delta_*\O_X\;\Delta_*\D_X^\bu[*])
\end{multline}
over the\/ $\Z/2$\+graded
$k$\+algebra~\textup{\eqref{hoch-cohomology-direct-sum}}.
\end{cor}

\begin{proof}
 For each $i=1$,~$\dotsc$, $n$, let $Y_i$ denote the open
subscheme $X\setminus \bigcup_{j\ne i} Z_i\subset X$.
 Let $w_i\in\O(Y_i)$ denote the restriction of the regular function
$w-c_i$ to~$Y_i$.
 The argument is based on the results of
Sections~\ref{coh-hochschild-cohomology}
and~\ref{coh-hochschild-homology} applied to the schemes $Y_i$ (or
their open subschemes) endowed with the potentials~$w_i$.

 The restriction of morphisms (in the coderived
categories) of quasi-coherent matrix factorizations to the open
subschemes $Y_i\subset X$ defines a $\Z/2$\+graded $k$\+algebra
morphism from the (middle or) right-hand side to the left-hand side 
of~\eqref{hoch-cohomology-direct-sum}, and a $\Z/2$\+graded
$k$\+module morphism from the (middle or) right-hand side to
the left-hand side of~\eqref{hoch-homology-direct-sum}.
 It remains to show that these morphisms are isomorphisms.

 For this purpose, one can start with replacing $\Delta_*\D_X^\bu$
or $\Delta_*\O_X$ in the second argument of the Hom spaces in
the middle or right-hand sides of~(\ref{hoch-cohomology-direct-sum}\+-%
\ref{hoch-homology-direct-sum}) with an injective matrix factorization
$\J$ on $X\times_k X$ representing the same object in the coderived
category.
 Then one notices that the restriction from $X\times_k X$ to
its open subscheme $V=\bigcup_{i=1}^n Y_i\times_k Y_i$
does not change the Hom spaces in the right-hand sides, as 
the image of $\Delta$ is contained in~$V$.

 Finally, one writes down the \v Cech resolution of the matrix
factorization $\J|_V$ corresponding to the covering of the scheme $V$
by its open subschemes $Y_i\times_k Y_i$.
 This is a finite acyclic complex of injective matrix factorizations,
so applying the functor $\Hom_{(V\;\O\;(w_2-w_1)|_V)\qcoh}(\K,{-})$
from any quasi-coherent matrix factorization $\K$ preserves its
acyclicity.
 The Hom spaces on any intersection of at least two different open
subschemes in the covering being zero by
Theorems~\ref{coh-hochschild-cohomology}\+-%
\ref{coh-hochschild-homology} (as $w$~is noncritical on $Y_i\cap Y_j$
for any $i\ne j$), the desired isomorphisms follow.
\end{proof}

\begin{rem}
 The Hochschild cohomology algebra and the Hochschild homology module
of the DG\+category version $\sDG^\b(X\coh)$ of the bounded derived
category $\sD^\b(X\coh)$ of (complexes of) coherent sheaves on
a separated scheme $X$ of finite type over a field~$k$ can be computed
in the way similar to (but simpler than) the above.
 The answers are the same as in Theorems~\ref{coh-hochschild-cohomology}
and~\ref{coh-hochschild-homology}:
\begin{multline}
 HH^*(\sDG^\b(X\coh))\.\simeq\.\Hom_{\sD(X\times_kX\qcoh)}
 (\Delta_*\D_X^\bu\;\Delta_*\D_X^\bu[*]) \\ \.\simeq
 \Hom_{\sD^\b(X\times_kX\coh)}(\Delta_*\O_X,\Delta_*\O_X[*])^\op
\end{multline}
and
\begin{equation}
 HH_*(\sDG^\b(X\coh))\.\simeq\.\Hom_{\sD(X\times_kX\qcoh)}
 (\Delta_*\O_X\;\Delta_*\D_X^\bu[*]),
\end{equation}
the only difference being that $\sDG^\b(X\coh)$ is a $\Z$\+graded
DG\+category and the right-hand sides describe the Hochschild
(co)homology as a $\Z$\+graded algebra and module.
 The only assumption is that the scheme $X$ should admit a smooth
stratification over~$k$ (i.~e., it suffices that the field~$k$
be perfect).
\end{rem}

\bigskip\addtocontents{toc}{\medskip}

\hbadness=1100

\end{document}